\documentclass[twoside, 10pt]{article}

\usepackage[hmargin = 0.75in, height =9.5in]{geometry}
\usepackage{amssymb,amsthm, amsfonts,amsmath,mathrsfs,  color, amsxtra, url, hyperref}

\usepackage{fancyhdr}
\renewcommand{\theequation}{\thesection.\arabic{equation}}
 \numberwithin{equation}{section}

\newtheorem {thm}{Theorem}[section]
\newtheorem {prop}{Proposition}[section]
\newtheorem {lemm}{Lemma}[section]
\newtheorem {deff}{Definition}[section]
\newtheorem {cor}{Corollary}[section]
\newtheorem {rem}{Remark}[section]

\def\ba{\begin{array}}
\def\ea{\end{array}}
\def\bea{\begin{eqnarray}}
\def\eea{\end{eqnarray}}
\def\beas{\begin{eqnarray*}}
\def\eeas{\end{eqnarray*}}
\def\bi{\begin{itemize}}
\def\ei{\end{itemize}}
\def\bc{\begin{cases}}
\def\ec{\end{cases}}

\def\bhe{\begin{highlightequation}  }
\def\ehe{\end{highlightequation}  }

\def\ba{\begin{array}}
\def\ea{\end{array}}
\def\bea{\begin{eqnarray}}
\def\eea{\end{eqnarray}}
\def\beas{\begin{eqnarray*}}
\def\eeas{\end{eqnarray*}}
\def\bi{\begin{itemize}}
\def\ei{\end{itemize}}
\def\bc{\begin{cases}}
\def\ec{\end{cases}}

\def\bhe{\begin{highlightequation}  }
\def\ehe{\end{highlightequation}  }

\def\a{\alpha}
\def\ga{\gamma}

\def\d{\delta}
\def\e{\varepsilon}

\def\z{\zeta}
\def\k{\kappa}
\def\l{\lambda}
\def\vr{\varrho}
\def\si{\sigma}

\def\th{\theta}
\def\o{\omega}
\def\vf{\varphi}
\def\vth{\vartheta}

\def\D{\Delta}
\def\Ga{\Gamma}
\def\L{\Lambda}
\def\O{\Omega}

\def\Th{\Theta}
\def\U{\Upsilon}

\def\bF{{\bf F}}
\def\bG{{\bf G}}

\def\bK{{\bf K}}

\def\bQ{{\bf Q}}

\def\bW{{\bf W}}

\def\bz{{\bf 0}}

\def\bw{{\bf w}}
\def\bx{{\bf x}}

\def\bx{{\bf x}}

\def\cA{{\cal A}}
\def\cB{{\cal B}}

\def\cD{{\cal D}}
\def\cE{{\cal E}}
\def\cF{{\cal F}}
\def\cG{{\cal G}}
\def\cH{{\cal H}}
\def\cI{{\cal I}}

\def\cK{{\cal K}}

\def\cM{{\cal M}}
\def\cN{{\cal N}}
\def\cO{{\cal O}}
\def\cP{{\cal P}}
\def\cQ{{\cal Q}}

\def\cS{{\cal S}}

\def\cX{{\cal X}}

\def\hJ{\mathbb{J}}

\def\hN{\mathbb{N}}

\def\hQ{\mathbb{Q}}
\def\hR{\mathbb{R}}

\def\hT{\mathbb{T}}
\def\hU{\mathbb{U}}

\def\hX{\mathbb{X}}

\def\sB{\mathscr{B}}
\def\sC{\mathscr{C}}
\def\sD{\mathscr{D}}

\def\sH{\mathscr{H}}
\def\sI{\mathscr{I}}

\def\sM{\mathscr{M}}
\def\sN{\mathscr{N}}
\def\sO{\mathscr{O}}

\def\sS{\mathscr{S}}
\def\sT{\mathscr{T}}

\def\sV{\mathscr{V}}
\def\sW{\mathscr{W}}
\def\sX{\mathscr{X}}

\def\fA{\mathfrak{A}}
\def\fB{\mathfrak{B}}
\def\fC{\mathfrak{C}}

\def\fE{\mathfrak{E}}
\def\fF{\mathfrak{F}}

\def\fH{\mathfrak{H}}

\def\fN{\mathfrak{N}}

\def\fP{\mathfrak{P}}

\def\fS{\mathfrak{S}}
\def\fT{\mathfrak{T}}

\def\fX{\mathfrak{X}}
\def\fY{\mathfrak{Y}}

\def\fp{\mathfrak{p}}
\def\fq{\mathfrak{q}}
\def\ff{\mathfrak{f}}
\def\fm{\mathfrak{m}}

\def\fy{\mathfrak{y}}
\def\fz{\mathfrak{z}}
\def\fx{\mathfrak{x}}
\def\fl{\mathfrak{l}}

\def\fg{\mathfrak{g}}
\def\fw{\mathfrak{w}}
\def\fra{\mathfrak{a}}
\def\fri{\mathfrak{i}}
\def\fk{\mathfrak{k}}
\def\fn{\mathfrak{n}}
\def\fj{\mathfrak{j}}
\def\ft{\mathfrak{t}}
\def\fs{\mathfrak{s}}

\def\fr{\mathfrak{r}}

\def\ti{\n \times \n}
\def\oti{\n \otimes \n}

\def\df{\n := \n}
\def\ls{\n \le \n}
\def\gs{\n \ge \n}
\def\={\n = \n}
\def\+{\n + \n}
\def\-{\n - \n}
\def\ins{\n \in \n}
\def\ld{\n \land \n}
\def\ve{\n \vee \n}
\def\sb{\n \subset \n}
\def\>{\n > \n}
\def\<{\n < \n}
\def\Cp{\n \cap \n}
\def\cp{\n \cup \n}

\def\nxi{ \rule[-0.5mm]{0.45pt}{2.6mm}\hspace{-0.4pt}\rule[-0.5mm]{1.8mm}{0.45pt}\hspace{-1.4mm}\rule[0.8mm]{1mm}{0.45pt}\hspace{-1.4mm}\rule[2mm]{1.8mm}{0.45pt}\hspace{-0.4pt}\rule[-0.5mm]{0.4pt}{2.6mm}\,}

\def\({\textnormal{(}}
\def\){\textnormal{)}}
\def\[{[\n[}
\def\]{]\n]}
\def\lan{\langle}
\def\ran{\rangle}

\def\no{\noindent}

\def\ss{\smallskip}

\def\q{\quad}
\def\qq{\qquad}
\def\n{\negthinspace}
\def\dn{\n \n}
\def\tn{\n \n \n}

\def\ol{\overline}
\def\ul{\underline}
\def\ua{\mathop{\uparrow}}
\def\da{\mathop{\downarrow}}

\def\Ra {\mathop{\Rightarrow }}

\def\wt{\widetilde}
\def\wh{\widehat}

\def\oMoo{\oM^{\raisebox{3pt}{\scriptsize \hb{$t_\oo$}}}}

\def\loo{{\overset{}{\oo}}}

\def\ogaP{{\oga_{\overset{}{\oP}}}}

\def\hb{\hbox}
\def\dis{\displaystyle}
\def\cd{\cdot}
\def\cds{\cdots}

\def\fa{\,\forall \,}
\def\pa{\partial}
\def\es{\emptyset}

\def\b1{{\bf 1}}
\def\qed{\hfill $\Box$ \medskip}

\def\liminf{\mathop{\ul{\rm lim}}}
\def\limsup{\mathop{\ol{\rm lim}}}
\newcommand{\Sup}[1]{\underset{#1}{\sup}\,}
\newcommand{\Inf}[1]{ \underset{#1}{\inf}\,}

\newcommand{\lsup}[1]{ \underset{#1}{\limsup}}
\newcommand{\linf}[1]{ \underset{#1}{\liminf}}
\newcommand{\lmt}[1]{ \underset{#1}{\lim}}
\newcommand{\lmtu}[1]{ \underset{#1}{\lim} \n \ua \,}
\newcommand{\lmtd}[1]{ \underset{#1}{\lim} \n \da \,}

\newcommand{\ccap}[2]{\underset{#1}{\overset{#2}{\cap}}}
\newcommand{\ccup}[2]{\underset{#1}{\overset{#2}{\cup}}}

\newcommand{\Rho}[1]{\rho_{\overset{}{#1}}}

\def\oP{{\ol{P}}}
\def\oQ{{\ol{Q}}}
\def\oA{{\ol{A}}}
\def\oeta{{\ol{\eta}}}
\def\otau{{\ol{\tau}}}
\def\oga{{\ol{\ga}}}
\def\oz{{\ol{\z}}}

\def\ocP{{\ol{\cP}}}

\def\ocA{{\ol{\cA}}}
\def\wcA{{\wt{\cA}}}
\def\ocF{{\ol{\cF}}}
\def\obF{{\ol{\bF}}}

\def\oo{{\ol{\o}}}
\def\oO{{\ol{\O}}}
\def\oD{{\ol{D}}}

\def\oW{{\ol{W}}}
\def\oK{{\ol{K}}}

\def\oX{{\ol{X}}}
\def\oY{{\ol{Y}}}
\def\oZ{{\ol{Z}}}

\def\oV{{\ol{V}}}
\def\oR{{\ol{R}}}
\def\oT{{\ol{T}}}
\def\oPhi{{\ol{\Phi}}}

\def\oM{{\ol{M}}}
\def\osX{\ol{\sX}}

\def\wA{{\wt{A}}}

\def\mto{\n \mapsto \n}
\def\nto{\n \to \n}

\def\bul{\no $\bullet$ }

\def\oxi{{\ol{\xi}}}
\def\oXi{{\ol{\Xi}}}

\def\ocD{{\ol{\cD}}}

\def\ocN{{\ol{\cN}}}
\def\ofN{{\ol{\fN}}}

\def\ofX{{\ol{\fX}}}

\def\ocO{{\ol{\cO}}}
\def\ocE{{\ol{\cE}}}
\def\ocI{{\ol{\cI}}}

\def\obW{{\ol{\bW}}}
\def\obK{{\ol{\bK}}}
\def\osW{{\ol{\sW}}}

\def\ooV{{\ol{\oV}}}
\def\oocP{{\ol{\ocP}}}

\def\ocp^\sharp{{\ol{\ocP}}}

\def\Wtzo{{\obW_{\oga,\oo}}}
\def\Wtgo{{\obW^t_{\oga,\oo}}}
\def\Ktgo{{\obK^t_{\oga,\oo}}}

\def\gP{{\big[ \n \big[\,\ocP\,\big] \n \big]}}
\def\gcP{{\big\{ \n \big\{ \ocP \big\} \n \big\}}}
\def\nci{\n \circ \n}
\def\OmX{\O_{\overset{}{X}}}
\def\omX{\o_{\overset{}{X}}}

\def\aand{\q \hb{and} \q}

\def\nne{\n \ne \n}

\def\cad{c\`adl\`ag }
\def\btau{\rule[1.45mm]{2.1mm}{0.7pt}\hspace{-3.35pt} \rule[-0.15mm]{0.8pt}{1.8mm}~}
\def\bbtau{\rule[1.3mm]{1.6mm}{0.6pt}\hspace{-2.6pt} \rule[0.1mm]{0.7pt}{1.3mm}\;}

\def\usa{upper semi-analytic}

\begin{document}

  \title{\bf   Optimal Stopping with Expectation Constraints}

\author{
 Erhan Bayraktar\thanks{ \noindent Department of  Mathematics,
  University of Michigan, Ann Arbor, MI 48109; email: {\tt erhan@umich.edu}.}
 \thanks{E. Bayraktar is supported in part  by the National Science Foundation under  DMS-2106556,
 and in part by the Susan M. Smith Professorship.
 Any opinions, findings, and conclusions or recommendations expressed in this material are
 those of the authors and do not necessarily reflect the views of the National Science Foundation.} $\,\,$,
 $~~$Song Yao\thanks{
 \noindent Department of  Mathematics,
  University of Pittsburgh, Pittsburgh, PA 15260; email: {\tt songyao@pitt.edu}. }
  \thanks{S. Yao is supported in part  by the National Science Foundation under DMS-1613208.
} }

\date{}

\maketitle

 \begin{abstract}

\ss

We analyze an optimal stopping problem  with a series of inequality-type and  equality-type expectation constraints in a general non-Markovian framework.
We show that the optimal stopping problem  with expectation constraints (OSEC) in an arbitrary probability setting is equivalent
to the constrained problem in weak formulation (an optimization over joint laws of stopping rules with Brownian motion and state dynamics on an enlarged canonical space) and thus the OSEC value 
is independent of a specific probabilistic setup.
Using a martingale-problem formulation, we make an equivalent characterization of the probability classes  in   weak formulation, which implies that the OSEC value function
is upper semi-analytic.
 Then we exploit a measurable selection argument to establish a dynamic programming principle   in weak formulation for the OSEC value function,
 in which the conditional expected costs  act  as   additional states for constraint levels at the intermediate horizon.

  \ss \no {\bf MSC 2020:}\;  60G40,  49L20, 93E20, 60G44


 \ss  \no   {\bf Keywords:}\; Optimal stopping  with expectation constraints, 
  martingale-problem formulation, enlarged canonical space, Polish space  of stopping times,
  dynamic programming principle,    regular conditional probability distribution,
  measurable selection.

\end{abstract}

  \section{Introduction}

  In this article,  we study a continuous-time   optimal stopping problem with
  a series of inequality-type and  equality-type expectation  constraints   in a general non-Markovian framework.

 Given a historical path $\bx|_{[0,t]}$, let the state  of the game $\cX^{t,\bx}_\cd$ evolve according to some SDE on a probability space
 $(\cQ, \cF, \fp)$   whose   drift and diffusion coefficients depend on the past trajectories   of the solution.
 The player   decides an exercise time $\tau  $ 
 to maximize her expected reward
   while being subject to a series of constraints:  for $i \ins \hN$,
 the expectation  of  some   accumulative cost  $\int_t^\tau \n  g_i ( r,\cX^{t,\bx}_{r \land \cd} ) dr $ 
   should not overpass certain level $y_i$ and the expectation  of  some other  accumulative cost  $\int_t^\tau \n  h_i ( r,\cX^{t,\bx}_{r \land \cd} ) dr $ should exactly hit certain level $z_i$.
 This   optimal stopping problem with expectation constraints (OSEC for short), or  optimization problem over constrained stopping times,
 has many  applications in various economic, engineering and financial areas such as  travel problem with fuel constraint, evaluation of American-type derivatives, quickest detection problem, etc.

 Let $V (t,\bx,y,z)$ denote the OSEC value with $(y,z) \df \big(\{y_i\},\{z_i\}\big)$.
 We aim to study the measurability of this value function   and
   establish an associated  dynamic programming principle (DPP) 
without imposing any   continuity condition on  reward and cost functions in time and state variables. 
  Inspired by \cite{EHJ_1987} and \cite{Elk_Tan_2013b},   we   embed the constrained stopping rule  $\tau$ together with the Brownian and state information into an enlarged canonical space $\oO$ and regard their joint distribution as  a new type of controls.
  Then the optimization of the expected reward  over  constrained stopping times  transforms into
 a maximal  expectation of   reward functional  over a class $ \ocP_{t,\bx}(y,z) $ of probability measures  on $\oO$  under which  three canonical coordinates $(\oW,\oX,\oT)$  serve as  Brownian motion, state process and constrained stopping rules respectively.

  One of our achievements   is to show that  
  the two optimization problems are equivalent:
  the  value $V (t,\bx,y,z)$ of OSEC in   strong formulation (i.e.,  on $\cQ$) is equal to the value $\oV (t,\bx,y,z)$ of OSEC in   weak formulation (i.e.,  over $\oO$).   This result indicates that the OSEC value 
  is actually  a  robust value,    independent of a specific probability model.

A dynamic programming principle   of a stochastic optimization problem allows one to maximize/minimize  the problem stage by stage
in a backward recursive way.  
 It requires the problem value function to be  measurable   so that one can do optimization  at an intermediate horizon first.
To show the measurability of the OSEC value functions, 
 we construct a Polish space of stopping times (which is 
 of independent interest) and exploit the martingale-problem formulation of  \cite{Stroock_Varadhan} to
 describe the   probability class   $\ocP_{t,\bx}(y,z)$ 
  as a series of probabilistic tests on   stochastic behaviors of the canonical coordinates  of $\oO$.
 Under such   countable characterization,
   the set-valued mapping   $(t,\bx,y,z) \mto \ocP_{t,\bx}(y,z)$ has Borel-measurable graph
 and the OSEC value function  $V \= \oV$   is thus  upper semi-analytic in $(t,\bx,y,z)$.

 In the next step  we   establish a DPP for $\oV$ in weak formulation, 
 which takes conditional expectations of the remaining costs as additional states for constraint levels at the intermediate horizon.
 For the subsolution side of this DPP, we use the regular conditional probability distribution
  to indicate that the probability classes $  \ocP_{t,\bx}(y,z) $, $ \fa (t,\bx,y,z)   $ are stable under conditioning.
 For the supersolution side of the DPP, we employ a measurable selection theorem in the analytic-set theory
 to paste  a class of locally $\e-$optimal probability measures.
 By the martingale-problem formulation again,  the canonical coordinates $(\oW,\oX)$ are still Brownian motion and the state process under the pasted probability measure.
 Finally  we make a delicate analysis to show that the third  canonical coordinate $\oT$ serves as  a constrained stopping time   under the pasted probability measure. To wit, the probability classes $  \ocP_{t,\bx}(y,z) $'s are also stable under   pasting (or concatenation).

  \no {\bf Relevant Literature.}

   Since Arrow et al. \cite{ABG_1949} and  Snell \cite{Snell_1952},
  the   theory of (unconstrained) optimal stopping   has been plentifully developed  over decades.
  Expositions of this theory are presented in   monographs \cite{CRS_1971,Shiryayev_1978,El_Karoui_1981,Kara_Shr_MF}.
   For the recent development of the  optimal stopping under  model uncertainty/non-linear expectations and the closely related controller-stopper-games, see  \cite{Karatzas_Sudderth_2001,Kara_Zam_2005,CDK-2006, Delbaen_2006, Kara_Zam_2008, Riedel_2009,OSNE1,OSNE2,OS_CRM,riedel2012,Bayraktar_Huang_2013,ETZ_2014,ROSVU,NZ_2015,RDOSRT,RDG}.

 Kennedy \cite{Kennedy_1982} initiated  the study of optimal stopping problem  with expectation constraint.
 The author used a {\it Lagrange multiplier} method
  to reformulate a discrete-time optimal stopping problem with first-moment constraint as a minimax problem
  and   showed that the optimal value of the dual problem
  is   equal to that of the primal problem.
  Since then, the Lagrangian technique has been prevailing
  in  research  of OSEC 
  (see e.g. \cite{Pontier_Szpirglas_1984,LSMS_1995,Balzer_Jansen_2002,Urusov_2005,Makasu_2009,Tanaka_2019})
  and has been applied to  various economic/financial problems
  such as  Markov decision processes with constrained stopping times \cite{Horiguchi_2001c,Horiguchi_2001b},
  mean-variance optimal control/stopping problem \cite{Pedersen_Peskir_2016,Pedersen_Peskir_2017},  quickest detection problem  \cite{Peskir_2012}, etc.

  Recently,   Ankirchner et al. \cite{AKK_2015} and Miller  \cite{Miller_C_2017a} took   different approaches to optimal stopping problems for diffusion processes with expectation constraints by transforming them to stochastic optimization problems with martingale controls.
  The former   characterizes the value function 
  in terms of a Hamilton-Jacobi-Bellman   equation and obtains a verification theorem, while the latter embeds the optimal stopping problem with first-moment constraint   into  a time-inconsistent (unconstrained) 
  stopping problem. However, the authors only postulate   dynamic programming principles  for  their corresponding  problems.
  In contrast, we rigorously prove in this article   a dynamic programming principle for the  optimal stopping problem   with expectation constraints.

 In their study of a   continuous-time stochastic optimization problem of controlled Markov processes,
 El Karoui, Huu Nguyen and Jeanblanc-Picqu\'e \cite{EHJ_1987}
  regarded   joint laws of    state and control processes
 as {\it control rules} on the product space of canonical state   space and control space.
 Then they  used a measurable selection theorem in the analytic-set theory  to establish a DPP
  without assuming any regularity on the reward functional.
    Nutz \& van Handel   \cite{HN_2012} and Neufeld \& Nutz \cite{Neufeld_Nutz_2013}  came up with a similar idea
  to address a superheging problem  under volatility uncertainty.
  They modeled the ``uncertainty" by path-dependent classes of controlled-diffusion laws  
   and explored the analytic measurability    of these classes. Using the measurable selection techniques,
    the authors obtained   PDD result in a  form of    time-consistency 
     of a sub-linear expectation and they thus established  a duality formula for the robust superhedging of measurable claims.
  The approach of \cite{HN_2012,Neufeld_Nutz_2013} was later extended to derive  DPPs of various non-Markovian   control problems,
  see  \cite{PRT_2013} for a dual formulation of  robust semi-static trading and its application  to   martingale optimal transportation and  see \cite{PTZ_2018} for stochastic control of a class of nonlinear kernels
  and its relation to second-order backward stochastic differential equations.
   Since the class of controlled-diffusion laws is naturally different from the class of stopping-time laws,
  the results of these works are not applicable to our optimal stopping problem with expectation constraints.

   In \cite{Elk_Tan_2013a,Elk_Tan_2013b}, El Karoui and Tan utilized the measurable selection argument to attain the DPP for a general stochastic control/stopping problem by embedding   stopping times with controlled diffusions into an enlarged canonical space  in the spirit of \cite{EHJ_1987}. However, the probability class they considered in weak formulation is not suitable for   optimal stopping with expectation constraints, see our Remark \ref{rem_032322} for  details.
   In this paper, we make a more accurate description  of  probability classes $\ocP_{t,\bx}(y,z)$ in which  the third
    canonical coordinate serves as some constrained stopping time.
  We construct a Polish space of stopping times and use it to show the Borel measurability of the graph $ \gP$.
  We also verify the stability of   probability classes $\ocP_{t,\bx}(y,z)$ under conditioning and concatenation
  so that we can  exploit measurable selection theorem to establish a DPP for the OSEC value function.
  \if{0}
   Instead,  we additionally  require in (D3) of Definition \ref{def_ocP} that under each $\oP$ of $\ocP_{t,\bx}(y,z)$
   the time  canonical coordinate $\oT$ serves as some   stopping time
   (it turns out that such a restriction does not affect the unconstrained optimal stopping problem in weak formulation).
     By  constructing a Polish space of stopping times, we manage to derive the Borel measurability of   graph $ \gP$ and thus obtain the measurability of the OSEC value functions.
    Because of   condition (D3) and   expectation constraints, it is more technically involved to
    verify the stability of our  probability classes $\ocP_{t,\bx}(y,z)$ under conditioning and concatenation
    and thus establish a DPP for the OSEC value function $\oV$.
  \fi

  A closely related topic to our research is optimal stopping   with constraint on the distribution of   stopping times.
 Bayraktar and Miller \cite{Bayraktar_Miller_2016}
 studied the problem of optimally stopping a Brownian motion with the restriction that
   the distribution of the stopping time must equal a given measure with finitely many atoms,
 and obtained a dynamic programming result which relates each of the sequential optimal control problems.
  K\"allblad \cite{Kallblad_2017}  used measure-valued martingales   to transform
  the distribution-constrained optimal stopping problem to a stochastic
 control problem and derived a DPP by measurable selection arguments.
  From the perspective of optimal transport,  Beiglb\"ock et al. \cite{BEES_2016} gave a geometric interpretation  of optimal stopping times  of a Brownian motion with distribution constraint.

 As to the stochastic control problems with expectation constraints,
 Pfeiffer et al. \cite{PTZ_2020} obtained a duality result by   a Lagrange relaxation approach 
 and Yu et al. \cite{CYZ_2020} used the measurable selection argument to derive a DPP result.
 Moreover, for stochastic control problems  with state constraints,   stochastic target problems with controlled losses and   related geometric DPP,  
 see \cite{BEI_2009,BET_2009,Bouchard_Nutz_2012,Soner_Touzi_2002_GDP,Soner_Touzi_2002_SDV,
 Soner_Touzi_2009,Bouchard_Vu_2010,Bouchard_Dang_2013,BMN_2014,BDK_2017}.

 The rest of the paper is organized as follows:
 Section \ref{sec_genprob} introduces the optimal stopping problem with expectation constraints
 in a generic probabilistic setting.
   Section \ref{sec_weak_form}  shows that  the  constrained optimal stopping problem
 can be equivalently embedded into an enlarged canonical space: i.e.,
 the OSEC 
 in   strong formulation has the same value as the OSEC in   weak formulation.
 In Section \ref{sec_Mart_prob}, we construct a Polish space of stopping times and use the martingale-problem formulation to make a countable characterization of the probability class in weak formulation, from which we deduce that  the OSEC  value function 
 is upper semi-analytic. Then in Section \ref{sec_DPP},  we utilize  a measurable selection argument to establish a dynamic programming principle in weak formulation for the OSEC value function. 
 We defer the proofs of our results to  Section \ref{sec_proof} and put some technical lemmata in the appendix.

  We close this section by a description of our notation and a review of the martingale-problem formulation.

\subsection{Notation and Preliminaries}

\label{subsec:preliminary}

  Throughout this paper, let us denote  $a^+   \df a \ve 0$ and $a^-   \df (  -  a) \ve 0$ for any   $ a \ins \hR$.
  We set $\hQ_+ \df \hQ \cap [0,\infty)$, $ \hQ^{2,<}_+ \df \\ \big\{ (s,r) \ins \hQ_+ \ti \hQ_+ \n :   s \< r \big\} $
   and set   $\Re \df   (-\infty,\infty]^\hN $ as the product of countably many copies of $(-\infty,\infty]$.
 On $\hT \df [0,\infty]$  we define a metric   $\Rho{+}(t_1,t_2) \df \big|\arctan(t_1) \- \arctan(t_2)\big|$, $\fa t_1,t_2 \ins \hT$
 and consider the induced topology by $\Rho{+}$.

 For a general topological space $ \big( \hX, \fT(\hX)\big) $, we  denote its  Borel sigma-field by $\sB(\hX)$ and let $\fP(\hX)$ be the set of all probability measures   on $\big(\hX,\sB(\hX)\big)$.  Recall that a topological space $ \hX $ is called a {\it Borel space} if it is homeomorphic to a Borel subset of a complete separable metric space.

  Let $ n \ins \hN$.  For any $x \ins \hR^n$ and $\d \ins (0,\infty)$,   let $O_\d(x)$ denote
 the open ball   centered at $x$ with radius $\d$  and let  $\ol{O}_\d(x)$ be its closure.
  For any $x,\wt{x} \ins \hR^n$
  we denote the usual inner  product by   $x   \cd   \wt{x} \df \sum^n_{i=1} x_i \wt{x}_i$, and for any  $ n \ti n -$real matrices $A,\wt{A}  $
  we denote the Frobenius inner  product by $A \n : \n \wt{A} := trace\big(A \wt{A}^T\big)$, where $\wt{A}^T$ is the transpose of $\wt{A}$.
 Let $ \big\{\cE^n_i\big\}_{i \in \hN}$ be a countable subbase of the Euclidean topology $\fT (\hR^n) $ on $\hR^n$.
 Then $ \sO (\hR^n) \df  \Big\{ \ccap{i=1}{n}  \cE^n_{k_i} \n : \{ k_i  \}^n_{i=1} \sb \hN  \Big\}
 \cp \{\es,\hR^n\}$
 forms a countable base of $\fT (\hR^n) $ and thus $\sB(\hR^n) \= \si\big(\sO (\hR^n)\big)$.
 We also set $\wh{\sO} (\hR^n) \df \ccup{k \in \hN}{} \big( \hQ_+ \ti \sO (\hR^n)   \big)^k $.
  For any $\vf \ins C^2(\hR^n)$,
  let $ D \vf   $ be its gradient, $ D^2 \vf $ be its Hessian matrix and denote $D^0\vf \df \vf$.
 \if{0}
 For any $f \ins C^1(\hR^n)$ and $x_o \ins \hR^n$, let $Df(x_o) \ins \hR^n$ be the gradient of $f$ at $x_o$
 and let $D^2 f(x_o) \ins \hR_{n \ti n}$ be the Hessian matrix of $f$ at $x_o$, i.e.,
 the $i-$the component of $Df(x_o)$ is $D_i f(x_0) \df \frac{\pa f}{\pa x_i} (x_o)$ and
 the element on $i-$th row and $j-$th column of $D^2 f(x_o)$ is  $D_{ij} f(x_0) \df \frac{\pa^2 f}{\pa x_i \pa x_j} (x_o)$
 for $i,y \= 1, \cds, n$.
 \fi
 For $ i  \= 1, \cds, n $,    define $\vf_i(x) \df x_i$,  $\fa x \= (x_1,\cds \n ,x_n) \ins \hR^n $.
 We let  $\fC(\hR^n) $   be the collection of these coordinate  functions   and their products, i.e.,
 $\fC(\hR^n) \df \{\vf_i\}^n_{i=1} \cp \{\vf_i \vf_j\}^n_{i,j=1} $.

  Let  $(\O,\cF,P )$ be a generic probability  space. For   subsets $A_1,A_2$ of $   \O$, we  denote $A_1 \D A_2 \df ( A_1 \Cp A^c_2 ) \cp ( A_2 \Cp A^c_1 ) $.
  For a random variable $\xi$ on $\O$ with values in a measurable space $(\cQ,\cG)$,
  we say   $\xi$ is $\cF/\cG-$measurable if its induced sigma-field $\xi^{-1}(\cG) \df \{\xi^{-1}(\cA) \n : \fa \cA \ins \cG \} $ is
  included in $ \cF $.
 For  a  sub-sigma-field $\fF  $ of $\cF$,
 define $\sN_P(\fF) \df \big\{ \cN \sb \O \n : \cN \sb A $ for some $ A \ins \fF \hb{ with } P (A) \=0  \big\}$,
 which collects all $P-$null sets with respect to   $\fF  $.
 For two sub-sigma-fields $\fF_1, \fF_2  $ of $\cF$,
 we denote   $\fF_1 \ve \fF_2 \df \si(\fF_1 \cp \fF_2)$.
 Let $t \ins [0,\infty)$.
 For a  filtration $\bF \= \{\cF_s\}_{s \in [t,\infty)}$ of $ \cF  $, we  set $\cF_\infty \df \si\Big(\underset{s \in [t,\infty)}{\cup}\cF_s\Big)$
 and refer to   filtration $\bF^P \n \= \big\{\cF^P_s \n \df \si \big( \cF_s \cp \sN_P ( \cF_\infty ) \big) \big\}_{s \in [t,\infty)} $ as the $P-$augmentation of $ \bF $.
 For  a process $X\= \{X_s\}_{s \in [t,\infty)}$  on $ \O  $ with values in a topological space,
 denote its raw filtration by $\bF^X \=  \big\{ \cF^X_s \df \si (X_r ; r \ins [t,s])  \big\}_{s \in [t,\infty)}$
 and denote the $P-$augmentation of $\bF^X $ by $\bF^{X,P} \n \= \big\{\cF^{X,P}_s \n \df \si \big( \cF^X_s \cp \sN_P ( \cF^X_\infty ) \big) \big\}_{s \in [t,\infty)} $.
 We call $X$ a continuous process if its paths are all continuous.
 When the time variable $s$ of $X$ has complicated form, we may write $X(s,\o)$ as $X_s(\o)$   for readability.
 By default, a Brownian motion $\{B_s\}_{s \in [t,\infty)}$ on $(\O,\cF,P )$  is with respect to its raw filtration
 $\bF^B$ unless stated otherwise.

 Fix    $d, l \ins \hN$.  Let $ \O_0 \=  \big\{ \o \ins  C   ( [0,\infty)  ; \hR^d ) \n : \o(0) \= 0 \big\}   $
 be the   space of   all $\hR^d-$valued continuous paths 
 starting  from $\bz$, which  is a Polish space under the topology of locally uniform convergence.
  Let $ P_0  $ be the Wiener measure on $\big(\O_0,  \sB(\O_0) \big)$, under which
   the canonical process $ W \=\{W_s\}_{s \in [0,\infty)} $ of $\O_0$ is a   $d-$dimensional standard Brownian motion.
   For any $t \ins [0,\infty)$,   $W^t_s \df W_s \- W_t$, $  s \ins [t,\infty)$ is also a   Brownian motion on $\big( \O_0,\sB(\O_0),P_0 \big)$.
 Let $\OmX \= C  ( [0,\infty)  ; \hR^l  )$ be the  space of all  $\hR^l-$valued   continuous paths 
 endowed  with the topology of locally uniform convergence.
 The  function $\fl_1 (t,\o_0) \df \o_0(t \ld \cd)  $ is continuous in $ (t,\o_0) \ins  [0,\infty) \ti \O_0  $
 while the  function
 \bea \label{020723_11}
 \fl_2 (t,\omX) \df \omX(t \ld \cd)
 \eea
  is continuous in $ (t,\omX) \ins  [0,\infty) \ti \OmX  $.

  Let   $ b  \n :   (0,\infty) \ti \OmX  \mto  \hR^l  $ and $ \si  \n :   (0,\infty) \ti \OmX  \mto   \hR^{l \times d} $ be two Borel-measurable functions such that   for any  $t \ins (0,\infty)$ and any $\omX, \omX' \ins \OmX$
 \bea
   \big|b(t,\omX)\-b(t,\omX') \big| \+ \big|\si(t,\omX)\-\si(t,\omX') \big| \ls \k(t)   \big\|\omX\-\omX' \big\|_t
    \aand \int_0^t  \big( |b (r,\bz)|^2 \+ |\si (r,\bz)|^2 \big) dr \< \infty,
   \q \label{coeff_cond1}
 \eea
 where $\k \n : (0,\infty) \mto (0,\infty)$ is some non-decreasing   function and $\big\|\omX\-\omX' \big\|_t \df \Sup{s \in [0,t]} \big|\omX(s)\-\omX'  (s) \big| $.
 Under   condition \eqref{coeff_cond1},   SDEs with coefficients $(b,\si)$ are well-posed (see e.g.  Theorem V.12.1 of \cite{Rogers_Williams_2}):

 \begin{prop} \label{prop_122021}
  Let $ (\O,\cF,  P  )$ be  a  probability space.
  Given $t \ins [0,\infty)$, let $\{B_s\}_{s \in  [t,\infty)}$ be a $d-$dimensional Brownian motion
  with respect to a right-continuous complete filtration $ \bF \=  \{\cF_s\}_{s \in  [t,\infty) } $  on $ (\O,\cF,  P  )$.
  For any $ \bx \ins \OmX $,   the SDE
   \bea \label{121621_11}
  X_s  \= \bx(t) + \int_t^s b (r,  X_{r \land \cd}) dr \+ \int_t^s  \si (r,  X_{r \land \cd}) d B_r , \; ~ \; \fa s \in [t,\infty)
  \; \hb{ with   initial condition }  X\big|_{[0,t]} \= \bx|_{[0,t]}
 \eea
  admits a unique strong solution $ X^{t,\bx} \= \{ X^{t,\bx}_s\}_{s \in  [0,\infty)}$ on $ (\O,\cF, \{\cF_s\}_{s \in  [t,\infty) },  P  )$ \big(i.e.,  $X^{t,\bx}$ is an   $\{\cF_{s \vee t}\}_{s \in  [0,\infty) }-$ adapted continuous process satisfying \eqref{121621_11}
  and $P\big\{ X^{t,\bx}_s \= \wt{X}^{t,\bx}_s , \fa s \ins [0,\infty) \big\} \= 1$
  if $ \big\{ \wt{X}^{t,\bx}_s \big\}_{s \in  [0,\infty)}$ is another  $\{\cF_{s \vee t}\}_{s \in  [0,\infty) }-$adapted continuous process satisfying \eqref{121621_11}\big).

 \end{prop}

 Let $ ^o \n X^{t,\bx} \= \{ ^o \n X^{t,\bx}_s\}_{s \in  [0,\infty)}$
 be the unique strong solution of \eqref{121621_11} on $  (\O,\cF,  P   ) \= \big( \O_0,\sB(\O_0),   P_0   \big) $ with $  (B, \bF  ) \= \big( W^t , \bF^{W^t,P_0}  \big) $
 and denote by $\sH_o $ the collection of  all $(-\infty,\infty]-$valued  Borel-measurable functions $\phi$ on $(0,\infty) \ti \OmX$ such that $ E_{P_0} \big[ \int_t^\infty \n  \phi^- (r, \, ^o \n X^{t,\bx}_{r \land \cd} )  dr\big]   \< \infty$ for any $(t,\bx) \ins [0,\infty) \ti \OmX$.

 Moreover, we   take the conventions $ \inf \es \df \infty$, $ \sup \es \df - \infty$   and $(+\infty)\+(-\infty) 
 \= -\infty$.
  In particular, on a measure space  $(\O,\cF,\fm) $,
   one can  define the integral  $\int_\O  \xi \, d \fm
   \df \int_\O  \xi^+ \, d \fm  \- \int_\O  \xi^- \, d \fm  $ for  any  $[-\infty,\infty]-$valued $\cF-$measurable random variable $\xi$  on  $\O$.

\subsection{Review of Martingale-Problem Formulation of SDEs}

 In this subsection, we consider  a general measurable space $(\O, \cF)$.
 Let $  \{B_s\}_{s \in [0,\infty)}$ be an $\hR^d-$valued continuous process on $\O$ with $B_0 \= \bz$
 and let $X  \= \{X_s\}_{s \in [0,\infty)}$ be an $\hR^l-$valued continuous  process on $\O$  such that
 $(B_s,X_s)$ is $\cF-$measurable for each $s \ins [0,\infty)$.

 Let $t \ins [0,\infty)  $. We set $B^t_s \df B_s \- B_t$, $\fa s \ins [t,\infty)$
 and define  filtration $\bF^t  \= \{\cF^t_s\}_{s \in [t,\infty)}$ by $\cF^t_s \df   \cF^{B^t}_s \ve \cF^X_s \=
 \si\big(B^t_r; r \ins [t,s]\big) \ve \si (X_r; r \ins [0,s] )$, $\fa s \ins [t,\infty)$.
 For any $\vf \ins C^2(\hR^{d+l})$,    define
 \beas    
    M^t_s(\vf)   \df   \vf \big(B^t_s  , X_s \big)
    \- \n \int_t^s  \n  \ol{b}  \big( r, X_{r \land \cd} \big) \n \cd \n D \vf \big( B^t_r  , X_r \big) dr
    \-   \frac12 \n \int_t^s  \n  \ol{\si} \, \ol{\si}^T  \big( r,  X_{r \land \cd}  \big) \n : \n D^2 \vf  ( B^t_r, X_r  )   dr  ,
    \q  \fa s \ins [t,\infty),
 \eeas
   where  $ \dis \ol{b}  (r,\omX) \df \binom{0}{  b(r,\omX)} \n \ins \hR^{d+l} $, $ \dis  \ol{\si}  (r,\omX) \df \binom{ I_{d \times d}}{   \si(r,\omX)} \n \ins \hR^{ (d+l)  \times d  } $,   $ \fa  (r,\omX) \ins   (0,\infty) \ti  \OmX $.
      Clearly, $\big\{  M^t_s(\vf) \big\}_{s \in [t,\infty)}$ is an $\bF^t-$adapted continuous process.
   For any $n \ins \hN   $ and $\fra \ins \hR^{d+l}$,     set $\tau^t_n (\fra) \df  \inf\big\{s \ins [t,\infty) \n : |(B^t_s,X_s)\-\fra|   \gs n  \big\} \ld (t\+n) $, which  is an $\bF^t-$stopping time. In particular, we   denote $\tau^t_n(\bz) $ by $\tau^t_n$.


  In virtue of \cite{Stroock_Varadhan},  we have the following  martingale-problem formulation of   SDEs  with coefficients $(b,\si)$ on $\O$.

\begin{prop}  \label{prop_MPF1}

  Let $(t,\bx) \ins [0,\infty) \ti \OmX $ and let $P$ be a probability measure on $(\O,\cF)$ such that
  $P \big\{X_s\=\bx(s), \fa s \ins [0,t] \big\} \= 1$.
  Then   $ \big\{ M^t_{s \land \tau^t_n(\fra) } (\vf) \big\}_{s \in [t,\infty)} $ is a bounded $\bF^t-$adapted continuous process under $P$
  for any $(\vf,n,\fra) \ins  C^2(\hR^{d+l}) \ti \hN \ti \hR^{d+l} $ and
  the following statements are equivalent  on  $(\O, \cF, P)$:

\no \(i\) The process $B^t$ is a   Brownian motion 
    and $P\{ X_s \=   X^{t,\bx}_s,   \fa s \ins [0,\infty)\}=1$, where  $ \big\{ X^{t,\bx}_s \big\}_{s \in [0,\infty)}$
  is the unique $ \big\{ \cF^{B^t,P}_{s \vee t} \big\}_{s \in [0,\infty)}   -$adapted continuous process solving   SDE \eqref{121621_11}.

\no \(ii\)   $ \big\{ M^t_{s \land \tau^t_n (\fra) } (\vf) \big\}_{s \in [t,\infty)} $  is a bounded  $\bF^t-$martingale for any $(\vf,n,\fra) \ins  C^2(\hR^{d+l}) \ti \hN \ti \hR^{d+l}$.

\no \(iii\)   $ \big\{ M^t_{s \land \tau^t_n } (\vf) \big\}_{s \in [t,\infty)} $  is a bounded $\bF^t-$martingale for any $(\vf,n) \ins \fC(\hR^{d+l}) \ti \hN $.

 Under either of these situations, one clearly has    $P \big\{ X^{t,\bx}_s \= ^o \n X^{t,\bx}_s ( B)     , \, \fa s \ins [0,\infty) \big\} \= 1$ and $ E_P  \big[ \int_t^\infty \n  \phi^- (r, X^{t,\bx}_{r \land \cd} )   dr\big]  \= E_{P_0} \big[ \int_t^\infty \n \phi^-(r,^o \n X^{t,\bx}_{r \land \cd} ) dr\big]  \< \infty$ for any $\phi \ins \sH_o$.


\end{prop}

\section{Optimal Stopping with Expectation Constraints}

\label{sec_genprob}

  Let $(\cQ, \cF, \fp)$ be   a    probability space
  equipped with a $d-$dimensional standard  Brownian motion $  \{\cB_s\}_{s \in [0,\infty)}$. 

  Let $t \ins [0,\infty)   $.  We set $\cB^t_s \df \cB_s \- \cB_t$, $\fa s \ins [t,\infty)$,
  which is also  a   Brownian motion on $(\cQ, \cF, \fp)$. 
  For any $\bx \ins \OmX $,    Proposition \ref{prop_122021} shows that   the SDE
  \bea   \label{121621_12}
  \cX_s  \= \bx(t) + \int_t^s b (r,  \cX_{r \land \cd}) dr \+ \int_t^s  \si (r,  \cX_{r \land \cd}) d \cB_r , \; ~ \; \fa s \in [t,\infty)
  \; \hb{ with   initial condition }  \cX\big|_{[0,t]} \= \bx|_{[0,t]}
  \eea
   admits   a unique strong  solution $  \cX^{t,\bx} \= \big\{ \cX^{t,\bx}_s\big\}_{s \in [0,\infty)} $
  on $\big(\cQ, \cF,\bF^{\cB^t,\fp}, \fp \big)$ \big(i.e.,  $\cX^{t,\bx}$ is the unique $ \big\{ \cF^{\cB^t,\fp}_{s \vee t} \big\}_{s \in [0,\infty)}   -$adapted continuous process solving   SDE \eqref{121621_12}\big).
  Let $\cS_t $ collect all $[t,\infty]-$valued $ \bF^{\cB^t,\fp}  -$stopping times.

 Let $f \ins \sH_o $, $ \{ g_i, h_i \}_{i \in \hN} \sb \sH_o$
 and let $\pi \n : [0,\infty) \ti \OmX \mto (-\infty,\infty]$ be a Borel-measurable function bounded from below
 by some $c_\pi \ins (-\infty,0)$.

 Given a historical path $\bx|_{[0,t]}$, 
 the state of the  game  then evolves   along process  $\{\cX^{t,\bx}_s\}_{s \in [t,\infty)}$.
 The player of the game need to select  an exercise time $\tau \ins \cS_t$ to cease  the game,
 at which she will receive an accumulative reward $\int_t^\tau   \n  f  \big(r, \cX^{t,\bx}_{r \land \cd}  \big) \, dr $ plus
 a terminal reward $   \pi \big(  \tau  , \cX^{t,\bx}_{ \tau   \land \cd}  \big) $ (both random rewards can take negative values).
 The player intends to  maximize the expectation of her    total wealth,
 but her choice $\tau$   is subject to a series of  expectation constraints
 \bea  \label{111920_11}
   E_\fp \Big[ \int_t^\tau  g_i ( r,\cX^{t,\bx}_{r \land \cd} ) dr  \Big] \ls y_i,   \q
  E_\fp \Big[ \int_t^\tau  h_i ( r,\cX^{t,\bx}_{r \land \cd} ) dr  \Big] \= z_i,   \q \fa i \ins \hN
 \eea
 for some  $(y,z) \= \big(\{y_i\}_{i \in \hN}, \{z_i\}_{i \in \hN}\big)
 \ins \Re \ti \Re$. One can regard   $\int_t^\tau  g_i ( r,\cX^{t,\bx}_{r \land \cd} ) dr$ or $\int_t^\tau h_i ( r,\cX^{t,\bx}_{r \land \cd} ) dr$ as   certain accumulative costs.
 So   the  value  of this   optimal stopping problem   with expectation constraints (OSEC for short)  is
 \bea    \label{081820_11}
 V (t,\bx,y,z)   \df    \Sup{\tau \in \cS_{t,\bx}(y,z) } E_\fp \Big[ \int_t^\tau  f \big( r, \cX^{t,\bx}_{r \land \cd}  \big) dr
 \+ \b1_{\{\tau < \infty\}} \pi \big( \tau  , \cX^{t,\bx}_{ \tau  \land \cd}  \big)  \Big] ,
\eea
  where $ \cS_{t,\bx}(y,z) \df \big\{ \tau  \ins   \cS_t \n : E_\fp \big[ \int_t^\tau  g_i ( r,\cX^{t,\bx}_{r \land \cd} ) dr  \big] \ls y_i,   \, E_\fp \big[ \int_t^\tau  h_i ( r,  \cX^{t,\bx}_{r \land \cd} ) dr  \big] \= z_i, \, \fa i \ins \hN \big\} $. 

 \if{0}
 For any $(t,\bx,y,z) \ins [0,\infty) \ti \OmX \ti \Re \ti \Re$   and for any $y' \= \{y'_i\}_{i \in \hN} \ins \Re$ such that
 $y_i \ls y'_i$ for any $i \ins \hN$, it is clear that
 \beas
 V (t,\bx,y,z) \ls V (t,\bx,y',z) .
 \eeas
 \fi

\begin{rem} \label{rem_112220}
 Let $(t,\bx) \ins [0,\infty) \ti \OmX$.

  \no 1\) \(finitely many constraints\) For   $i \ins \hN$,  the constraint $E_\fp \big[ \int_t^\tau  g_i ( r,\cX^{t,\bx}_{r \land \cd} ) dr  \big] \ls y_i$ holds for any $\tau  \ins   \cS_t$ 
  if  $ y_i \= \infty $, and the constraint $E_\fp \big[ \int_t^\tau  h_i ( r,\cX^{t,\bx}_{r \land \cd} ) dr  \big] \= z_i$
  holds for any $\tau  \ins   \cS_t$ 
  if $\big(h_i(\cd,\cd),z_i\big) \= (0,0)$.

  \no 1a\) If we take  $\big(y_i,h_i(\cd,\cd),z_i\big) \= (\infty,0,0)$, $\fa i \ins \hN$, there is no expectation constraint at all.

  \no 1b\) If one takes $y_i \= \infty$, $\fa i \gs 2$ and
   $\big(h_i(\cd,\cd),z_i\big) \= (0,0)$, $\fa i \ins \hN$,   \eqref{111920_11} reduces  to
  a single constraint $ E_\fp \big[ \int_t^\tau  g_1 ( r,\cX^{t,\bx}_{r \land \cd} ) dr  \big] \\ \ls y_1 $.
  In addition, if   $y_1 \gs 0$,
  then $t \ins  \cS_{t,\bx}(y,\bz) $.    

   \no 1c\) If one takes $y_i \= \infty$, $\fa i \ins \hN $ and
   $\big(h_i(\cd,\cd),z_i\big) \= (0,0)$, $\fa i \gs 2$,   \eqref{111920_11} degenerates to  
  $ E_\fp \big[ \int_t^\tau  h_1 ( r,    \cX^{t,\bx}_{r \land \cd} ) dr  \big]    \= z_1 $.

    \no 1d\) If we take  $\big(y_i,h_i(\cd,\cd),z_i\big) \= (\infty,0,0)$, $\fa i \gs 2$, \eqref{111920_11} becomes a couple of constraints
  $ E_\fp \big[ \int_t^\tau  g_1 ( r,\cX^{t,\bx}_{r \land \cd} ) dr  \big] \ls y_1 $ and $ E_\fp \big[ \int_t^\tau  h_1 ( r,\cX^{t,\bx}_{r \land \cd} ) dr  \big] \= z_1$.

  \no 1e\) If we take  $g_2 \= -g_1$, $ y_2 \gs - y_1$; $y_i \= \infty$, $\fa i \gs 3$ and
   $\big(h_i(\cd,\cd),z_i\big) \= (0,0)$, $\fa i \ins \hN$,   \eqref{111920_11} becomes a range constraint  $ - y_2 \ls E_\fp \big[ \int_t^\tau  g_1 ( r,   \cX^{t,\bx}_{r \land \cd} ) dr  \big]   \ls y_1 $.

\no 2\) \(moment constraints\) Let $i \ins \hN$, $ a  \ins (0,\infty)$ and $q \ins [1,\infty) $.
If   $g_i(s,\bx)  \= a q s^{q-1}   $, $\fa (s,\bx) \ins (0,\infty) \ti \OmX  $
\big(resp.  $h_i(s,\bx)  \= a q s^{q-1}   $, $\fa (s,\bx) \ins (0,\infty) \ti \OmX  $\big),
then the expectation constraint  $E_\fp \big[ \int_t^\tau  g_i ( r,\cX^{t,\bx}_{r \land \cd} ) dr  \big] \ls y_i$
\big(resp.  $E_\fp \big[ \int_t^\tau  h_i ( r,\cX^{t,\bx}_{r \land \cd} ) dr  \big] \= z_i$\big)
specifies  as a moment constraint   $E_\fp \big[ a  ( \tau^q \- t^q )    \big] \ls y_i$
\big(resp.  $E_\fp \big[ a  (\tau^q \- t^q )   \big] \= z_i$\big).

\if{0}
If   $g_i(s,\bx)  \= a q s^{q-1}   $, $\fa (s,\bx) \ins (0,\infty) \ti \OmX  $,
then the expectation constraint  $E_\fp \big[ \int_t^\tau  g_i ( r,\cX^{t,\bx}_{r \land \cd} ) dr  \big] \ls y_i$
specify as a moment constraint   $E_\fp \big[ a  ( \tau^q \- t^q )    \big] \ls y_i$;
If    $h_i(s,\bx)  \= a q s^{q-1}   $, $\fa (s,\bx) \ins (0,\infty) \ti \OmX  $,
then the expectation constraint   $E_\fp \big[ \int_t^\tau  h_i ( r,\cX^{t,\bx}_{r \land \cd} ) dr  \big] \= z_i$
specify as a moment constraint    $E_\fp \big[ a  (\tau^q \- t^q )   \big] \= z_i$.
\fi

\end{rem}

 We would like to study the measurability of   value function $V$  and  derive a  dynamic programming principle for $V$ without imposing any   continuity condition on 
 functions $f,\pi,g_i$'s and $h_i$'s in time and state variables.
 Inspired by \cite{EHJ_1987}, we will use   mapping $   \o \mto \big(\cB_\cd(\o),\cX^{t,\bx}_\cd(\o),\tau(\o)\big)  $
 to  transfer the OSEC 
 onto an enlarged canonical space
   and regard   joint laws of $(\cB_\cd,\cX^{t,\bx}_\cd,\tau)$  as  a new type of controls.

\section{Weak Formulation}
\label{sec_weak_form}

 In this section, we study   the  optimal stopping problem with expectation constraints
 in a weak formulation or over an  enlarged canonical space
 \beas
 \oO \df \O_0   \ti \OmX  \ti \hT .
 \eeas
 Clearly, $\oO$ is a Borel space  under the product topology.
 Let $\fP\big(\oO\big)$ be  the space  of all probability measures on $\big(\oO, \sB(\oO) \big)$
 equipped with the topology of weak convergence, which is also a Borel space (see e.g. Corollary 7.25.1 of \cite{Bertsekas_Shreve_1978}).
  For any $\oP \ins \fP\big(\oO\big)$,   set $\sB_\oP(\oO) \df \si \big( \sB(\oO) \cp \sN_\oP  (\sB(\oO) ) \big)$.
 We define the canonical coordinates on $\oO$ by
\beas
\oW_{\n s} (\oo) \df \o_0(s) ,   \q \oX_s(\oo) \df \omX(s) ,  \q    s \ins [0,\infty)
 \q \hb{and} \q \oT(\oo) \df   \ft , \q \fa \oo \= \big(\o_0,\omX,\ft\big) \ins \oO ,
\eeas
 in which one can regard $\oW$ as a canonical coordinate for Brownian motion, $\oX$ as a canonical coordinate for the state process, and $\oT$ as a canonical coordinate for stopping rules.
 Given $t \ins [0,\infty)$, we also  set $\oW^t_{\n s} \df \oW_{\n s} \- \oW_{\n t}$, $\fa s \ins [t,\infty)$.

 The weak formulation of the OSEC 
 relies on the following probability classes of $\fP\big(\oO\big)$.

\begin{deff} \label{def_ocP}
 For any $(t,\bx) \ins [0,\infty) \ti \OmX$,
 let $\ocP_{t,\bx}$ be the collection of  all probability measures $ \oP \ins \fP\big(\oO\big) $ satisfying:

\no \(D1\) The process $ \oW^t  $ is a  $d-$dimensional Brownian motion on $\big(\oO , \sB(\oO) , \oP\big)$.

\no \(D2\)  $  \oP\big\{ \oX_s \= \osX^{t,\bx}_{\n s}, ~ \fa s \ins [0,\infty) \big\} \= 1$,
  where  $ \big\{\osX^{t,\bx}_{\n s}\big\}_{s \in [0,\infty)}$ is an $\big\{ \cF^{\oW^t,\oP}_{s \vee t } \big\}_{s \in [0,\infty)} -$adapted continuous process   that
  uniquely  solves   the following SDE  on $\big(\oO , \sB\big(\oO\big) , \oP\big) \n : $
 \bea    \label{Ju01_01}
 \osX_{\n s} = \bx(t) + \int_t^s b \big( r, \osX_{\n r \land \cd}  \big)dr \+ \int_t^s \si \big( r, \osX_{\n r \land \cd}  \big) d \oW_{\n r}, \q \fa s \ins [t,\infty) \; \hb{ with initial condition } \osX\big|_{[0,t]} \= \bx\big|_{[0,t]} .
 \eea

  \no \(D3\)  There exists   a $[t,\infty]-$valued   $ \bF^{W^t,P_0} -$stopping time $\wh{\tau}$   on $\O_0$ such that
 $ \oP \big\{  \oT  \=  \wh{\tau} (\oW  )   \big\} \= 1$.


 \end{deff}

 Let $ t  \ins [0,\infty)   $. For any $s \in [t,\infty)$, define
 $\ocF^t_s \df \cF^{\oW^t}_s \ve \cF^\oX_s \=  \si \big(\oW^t_{\n r}; r \ins [t,s]\big) \ve \si \big(\oX_r;r \ins [0,s]\big)$, which is countably generated by
 $ \big\{\oX^{-1}_r (\cO) \n: r \ins \hQ \Cp [0,t], \cO \ins \sO(\hR^l) \big\} \cp \big\{ (\oW^t_r,\oX_r)^{^{-1}} (\cO' ) \n: r \ins \hQ \Cp (t,s], \cO'  \ins \sO(\hR^{d+l}) \big\}   $.   We   denote the filtration  $  \big\{\ocF^t_s   \big\}_{s \in [t,\infty)}$ by $\obF^t$.
 For any  $ (\vf,n,\fra )  \ins   C^2 (\hR^{d+l}) \ti \hN \ti \hR^{d+l} $,
 \beas  
  \oM^t_{\n s}(\vf)   \df  \vf \big(\oW^t_{\n s}   , \oX_s \big)
    \- \n \int_t^s  \n  \ol{b}  \big( r, \oX_{r \land \cd} \big) \n \cd \n D \vf \big( \oW^t_{\n r}  , \oX_r \big) dr
    \-   \frac12 \n \int_t^s  \n  \ol{\si} \, \ol{\si}^T  \big( r,  \oX_{r \land \cd}  \big) \n : \n D^2 \vf  ( \oW^t_{\n r}  , \oX_r  )   dr  ,
    \q \fa s \ins [t,\infty)
 \eeas
 is an $\obF^t-$adapted continuous process  and   $\otau^t_n (\fra) \df  \inf\big\{s \ins [t,\infty) \n : \big|(\oW^t_{\n s} ,   \oX_s) \- \fra \big|    \gs n   \big\} \ld (t\+n) $ is an $\obF^t-$stopping time.
 We will simply denote $ \otau^t_n (\bz) $ by $\otau^t_n  $.

 Let us also define a shifted canonical  process on $\oO$ by $\osW^t_{\n \fs} (\oo) \df \oW_{t+\fs} (\oo) \- \oW_t (\oo) \= \oW^t_{t+\fs} (\oo)$, $\fa (\fs,\oo) \ins [0,\infty) \ti \oO$.  \big(Note: the subscript $\fs \ins [0,\infty)$ of $\osW^t$  is the relative time after $t$ while the subscript $s \ins [t,\infty)$ of $\oW^t $ is the real time.\big)

   \if{0}

   \begin{rem} \label{rem_061222}
Given $t \ins [0,\infty)$, as Lemma \ref{lem_010922} shows that $\sW^t$   is a continuous mapping  from $\O_0$ to $\O_0$,  it follows that
   $\osW^t \= \sW^t(\oW)$ is a continuous mapping  from $ \oO $ to $\O_0$.

\end{rem}

   \fi

 According to the martingale-problem formulation of   SDEs (Proposition \ref{prop_MPF1}), we have an alternative description  of
 the probability class  $\ocP_{t,\bx}$:

 \begin{rem} \label{rem_ocP}

 Let $(t,\bx) \ins [0,\infty) \ti \OmX$.
 In   Definition  \ref{def_ocP}  of $\ocP_{t,\bx}$,   \(D1\) and \(D2\)   is equivalent to

   \no   \(D1\,$'$\)  $\oP\{\oX_s \= \bx(s), \fa s \ins [0,t]\} \= 1$ and
  $\big\{ \oM^t_{s \land \otau^t_n } (\vf) \big\}_{s \in  [t,\infty) } $  is a bounded $ \big( \obF^t , \oP \big)   -$martingale for any $(\vf,n) \ins \fC(\hR^{d+l}) \ti \hN $,

  \no  while 
 \(D3\)   is equivalent to

   \no   \(D3\,$'$\)  There exists   a $[0,\infty]-$valued   $ \bF^{W,P_0} -$stopping time $\ddot{\tau}$   on $\O_0$ such that
 $ \oP \big\{  \oT  \=  t \+ \ddot{\tau} \big(\osW^t  \big)   \big\} \= 1$.

   \end{rem}

 \begin{rem} \label{rem_ocP2}

 Let $(t,\bx) \ins [0,\infty) \ti \OmX$
 and let $\oP \ins \fP(\oO)$ satisfy  \(D1\) and \(D2\) of Definition  \ref{def_ocP}. 

\no \(1\) For any $\phi \ins \sH_o$, Proposition \ref{prop_MPF1} shows that  $ E_\oP  \big[ \int_t^\infty \n  \phi^- (r, \oX_{r \land \cd} )   dr\big] \= E_\oP  \big[ \int_t^\infty \n  \phi^- \big(r, \osX^{t,\bx}_{r \land \cd} \big)   dr\big]   \< \infty$.

\no \(2\) Let $ (\vf,n,\fra )  \ins   C^2 (\hR^{d+l}) \ti \hN \ti \hR^{d+l} $. As $\big\{ \oM^t_{s \land \otau^t_n (\fra)} (\vf) \big\}_{s \in  [t,\infty) } $  is a bounded $ \big( \obF^t , \oP \big)   -$martingale,
the optional sampling theorem implies that
  for any   two $[t,\infty]-$valued $\obF^t-$stopping times $ \oz_1, \oz_2 $ with $\oz_1 \ls \oz_2$, $\oP-$a.s.,
 \bea   \label{020522_17}
  E_\oP  \Big[  \big( \, \oM^t_{\oz_2 \land \otau^t_n(\fra)} (\vf ) \-  \oM^t_{\oz_1 \land \otau^t_n(\fra)} (\vf ) \big) \b1_\oA  \Big]
  \=  E_\oP  \Big[  E_\oP  \Big[ \oM^t_{\oz_2 \land \otau^t_n(\fra)} (\vf ) \-  \oM^t_{\oz_1 \land \otau^t_n(\fra)} (\vf ) \Big| \ocF^t_{\oz_1 }\Big] \b1_\oA  \Big]
  \= 0 , ~ \;  \fa \oA \ins \ocF^t_{\oz_1} . \q
  \eea

  \end{rem}

  Let $(t,\bx) \ins [0,\infty) \ti \OmX$, $(y,z) \= \big(\{y_i\}_{i \in \hN}, \{z_i\}_{i \in \hN}\big) \ins \Re \ti \Re$
  and set   $\oR  (t) \df \int_{\oT \land t}^\oT \n  f   (r, \oX_{r \land \cd}  ) dr  \+ \b1_{\{\oT < \infty\}} \pi   \big(\oT, \oX_{\oT \land \cd}\big) $.
  Given a historical state  path  $\bx|_{[0,t]}$,
  the   value  of the optimal stopping problem with    expectation constraints
 \bea \label{111920_14}
 E_\oP \bigg[ \int_t^\oT   g_i ( r,\oX_{r \land \cd} ) dr  \bigg] \ls y_i,   \q  E_\oP \bigg[ \int_t^\oT   h_i ( r,\oX_{r \land \cd} ) dr  \bigg] \= z_i, \q \fa i \ins \hN
 \eea
   in weak formulation is
   \beas
    \oV (t,\bx,y,z) \df \underset{\oP \in \ocP_{t,\bx}(y,z)}{\sup} E_\oP \big[ \, \oR  (t) \big]
    \= \underset{\oP \in \ocP_{t,\bx}(y,z)}{\sup} E_\oP \Big[ \int_t^\oT \n  f   (r, \oX_{r \land \cd}  ) dr  \+ \b1_{\{\oT < \infty\}} \pi   \big(\oT, \oX_{\oT \land \cd}\big) \Big]  ,
   \eeas
  where $   \ocP_{t,\bx}(y,z) \df \Big\{\oP \ins \ocP_{t,\bx} \n : E_\oP \big[ \int_t^\oT   g_i ( r,\oX_{r \land \cd} ) dr  \big] \ls y_i,   \,
   E_\oP \big[ \int_t^\oT   h_i ( r,\oX_{r \land \cd} ) dr  \big] \= z_i, \, \fa i \ins \hN  \Big\} $.
   We will simply call $ \oV (t,\bx,y,z) $ the {\it weak} value of the   optimal stopping problem with expectation constraints.
   In case $ \ocP_{t,\bx}(y,z) \=  \es$, $\oV (t,\bx,y,z)   \= -\infty $ by the convention   $ \sup \es \df -\infty$.

 We can consider another weak value function of the  OSEC: 
 Let  $ \bw  \ins   \O_0    $ and  define
   $ \ocP_{t,\bw,\bx} \df \big\{\oP \ins \ocP_{t,\bx} \n :   \oP  \big\{  \oW_{\n s}  \=   \bw(s)   ,
    \fa s \ins [0,t]   \big\} \= 1 \big\}$
     as the subclass of $\ocP_{t,\bx}$   given the historical Brownian path $ \bw  |_{[0,t]}$.
   The weak value  of the optimal stopping problem with    expectation constraints \eqref{111920_14}
  given $(\bw , \bx) \big|_{[0,t]}$ 
  is
  $ 
  \oV (t,\bw,\bx,y,z) \df \underset{\oP \in \ocP_{t,\bw,\bx}(y,z)}{\sup} E_\oP \big[  \, \oR  (t)  \big] $,
 where  $ \ocP_{t,\bw,\bx}(y,z) \df \big\{\oP \ins \ocP_{t,\bx}(y,z)   \n : \oP  \big\{  \oW_{\n s}  \=   \bw(s)   ,
    \fa s \ins [0,t]   \big\} \= 1 \big\} $.

One of our main results  (Theorem \ref{thm_V=oV} below) exposes that
the value function  $V(t,\bx,y,z)$ in \eqref{081820_11} coincides with  the weak value function $\oV(t,\bx,y,z)$,
and is even equal to   $ \oV (t,\bw,\bx,y,z) $.

 \begin{thm} \label{thm_V=oV}
 Let $(t,\bw,\bx,y,z) \ins [0,\infty) \ti \O_0   \ti \OmX \ti \Re \ti \Re $.
 Then $   V (t,\bx,y,z ) \= \oV (t,\bx,y,z ) = \oV (t, \bw,\bx, y,z ) $, and
   $ \cS_{t,\bx}(y,z)  \nne \es \Leftrightarrow \ocP_{t,\bx}(y,z) \nne \es \Leftrightarrow \ocP_{t,\bw,\bx}(y,z) \nne \es $.

 \end{thm}

 Theorem \ref{thm_V=oV} demonstrates that   the value of the OSEC 
is independent of a specific probabilistic setup and is also indifferent to the Brownian  history.
This result even allows us to deal with the robust case:

 \begin{rem} \label{rem_robust}

 Let $  \big\{(\cQ_\a, \cF_\a, \fp_\a)\big\}_{\a \in \fA} $ be a family of probability spaces,
 where $\fA$ is a countable or uncountable  index set \big(e.g. one can consider   a non-dominated class $\{\fp_\a\}_{\a \in \fA}$ of probability measures on a measurable space $(\cQ,\cF)$\big).

  Given  $\a \ins \fA$, let $\cB^\a \= \{\cB^\a_s\}_{s \in [0,\infty)}$ be a $d-$dimensional standard  Brownian motion on $(\cQ_\a, \cF_\a, \fp_\a)$.
  For any $ (t,\bx)  \ins [0,\infty) \ti \OmX  $, set $\cB^{\a,t}_s \df \cB^\a_s \- \cB^\a_t$, $s \ins [t,\infty) $ and let $\cX^{\a,t,\bx} \= \big\{ \cX^{\a,t,\bx}_s  \big\}_{s \in [0,\infty)}$ be the unique $\big\{\cF^{\cB^{\a,t},\fp_\a}_{s \vee t}\big\}_{s \in [0,\infty)}-$adapted continuous process solving the SDE
 \beas
 \cX_s \= \bx(t)  + \int_t^s b ( r, \cX_{r \land \cd}) dr \+ \int_t^s  \si ( r, \cX_{r \land \cd}) d \cB^\a_r , ~ \; \fa s \in [t,\infty)
 \hb{ with initial condition } \cX \big|_{[0,t]} \= \bx \big|_{[0,t]}
 \eeas
  on $\big(\cQ_\a, \cF_\a, \bF^{\cB^{\a,t},\fp_\a},  \fp_\a\big)$.

 Then we know from Theorem \ref{thm_V=oV} that  for any $(t,\bx) \ins [0,\infty)    \ti \OmX  $ and
    $(y,z) \= \big(\{y_i\}_{i \in \hN}, \{z_i\}_{i \in \hN}\big) \ins \Re \ti \Re $
 \beas
 \;\; \oV (t,\bx,y,z) \= \Sup{\a \in \fA} \, \Sup{ \tau_\a  \in \cS^\a_{t,\bx}(y,z) } E_{\fp_\a} \Big[ \int_t^{ \tau_\a} f \big( r,\cX^{\a,t,\bx}_{r \land \cd}  \big) dr
 \+ \b1_{\{\tau_\a < \infty\}} \pi \big(  \tau_\a, \cX^{\a,t,\bx}_{  \tau_\a  \land \cd}  \big) \Big] ,
 \eeas
  where $ \cS^\a_{t,\bx}(y,z) $ collects all $[t,\infty]-$valued  $\bF^{B^{\a,t},\fp_\a}-$stopping times $\tau_\a$ satisfying
  $ E_{\fp_\a} \big[ \int_t^{ \tau_\a} g_i \big( r,\cX^{\a,t,\bx}_{r \land \cd} \big) dr  \big] \ls y_i$ and $ E_{\fp_\a} \big[ \int_t^{ \tau_\a} h_i \big( r,\cX^{\a,t,\bx}_{r \land \cd} \big) dr  \big] \= z_i$ for any $ i \ins \hN $.
  To wit, the weak value  $ \oV (t,\bx,y,z)$ is also equal to the robust value of the OSEC
  under model uncertainty.

\end{rem}

  The equivalence between   strong and weak formulation of an  (unconstrained)  optimal stopping problem was discussed in
\cite{Elk_Tan_2013b}. However, their argument may not be applicable to optimal stopping with expectation constraints:

\begin{rem} \label{rem_032322}

 When $\big(y_i,h_i(\cd,\cd),z_i\big) \= (\infty,0,0)$, $\fa i \ins \hN$, 
the unconstrained version of Theorem \ref{thm_V=oV} states that for any $(t,\bx ) \ins [0,\infty) \ti  \OmX   $,
$ V (t,\bx ) \df  \Sup{\tau \in \cS_t } E_\fp \big[ \int_t^\tau \n f \big( r, \cX^{t,\bx}_{r \land \cd}  \big) dr
 \+ \b1_{\{\tau < \infty\}} \pi \big( \tau  , \cX^{t,\bx}_{ \tau  \land \cd}  \big)  \big] $ is equal to $ \oV (t,\bx  ) \df \underset{\oP \in \ocP_{t,\bx} }{\sup} E_\oP \big[\, \oR(t)\big]$.
On the other hand,  \cite{Elk_Tan_2013b} showed  that for any $(t,\bx ) \ins [0,\infty) \ti  \OmX   $,
$ V (t,\bx )  $   equals   $  \ooV (t,\bx) \df \underset{\oP \in \oocP_{t,\bx}}{\sup} E_\oP \big[\, \oR(t)\big] $,
  where $\oocP_{t,\bx}$ collects all $ \oP \ins \fP(\oO)$ satisfying \(D1\), \(D2\)  of Definition \ref{def_ocP}
  and  ``\,$\oP \big\{ \oT \gs t \big\} \= 1$"
  \(We summarize \cite{Elk_Tan_2013b}'s result in our terms 
  for an easy comparison with our work\).
As $  \ocP_{t,\bx} \sb \oocP_{t,\bx}$,
 the equality $ V (t,\bx ) \= \underset{\oP \in \ocP_{t,\bx} }{\sup} E_\oP \big[\, \oR(t)\big]
 \= \underset{\oP \in \oocP_{t,\bx}}{\sup} E_\oP \big[\, \oR(t)\big] $
 indicates that the probability classes $\ocP_{t,\bx}$'s are more accurate than 
 $\oocP_{t,\bx}$'s  to describe the \(unconstrained\)  optimal stopping problem in weak formulation.

 The condition \(D3\) of   Definition \ref{def_ocP} 
 is necessary for  the expectation-constraint case.
 Without it, the weak value $\ooV(t,\bx,y,z)   \df \underset{\oP \in \oocP_{t,\bx}(y,z) }{\sup} E_\oP \big[\, \oR(t)\big]$
 \Big(with $ \oocP_{t,\bx}(y,z)   \df \Big\{\oP \ins \oocP_{t,\bx} \n : E_\oP \big[ \int_t^\oT   \n  g_i ( r,\oX_{r \land \cd} ) dr  \big] \ls y_i,   \,
   E_\oP \big[ \int_t^\oT   \n  h_i ( r,\oX_{r \land \cd} ) dr  \big]   \= z_i, \, \fa i \ins \hN  \Big\} $\Big)
 may not be equal to $V (t,\bx,y,z )  $ for the following reason:

 In Proposition 4.3 of \cite{Elk_Tan_2013b}, the key  to show $\ooV (t,\bx) \ls  V (t,\bx) $, or $E_\oP \big[\, \oR(t)\big] \ls  V (t,\bx) $ for a given $\oP \ins \oocP_{t,\bx} $,
 relies on   transforming  the hitting times of process $\big\{E_\oP   \big[   \b1_{\{\oT \in [t,s]\}} \big| \cF^{\oW^t,\oP}_\infty \big]\big\}_{s \in  [t,\infty)}$
 to a member of $\cS_t$. More precisely,
  the so-called {\it Property \(K\)} assures
   an $\bF^{W^t,P_0}-$adapted \cad   process $ \wh{\vth}_\cd $  such that
$\wh{\vth}_s  ( \oW  ) \=  E_\oP   \big[   \b1_{\{\oT \in [t,s]\}} \big| \cF^{\oW^t}_s \big] \=  E_\oP   \big[   \b1_{\{\oT \in [t,s]\}} \big| \cF^{\oW^t,\oP}_\infty \big]  $, $\oP-$a.s. for any $s \ins [t,\infty)$.
 It follows  that
 $   E_\oP \big[ \b1_{\{\oT \in [t,s]\}}  \n  \b1_{\{\osX^{t,\bx}   \in A\}}   \big| \cF^{\oW^t,\oP}_\infty \big]
   \= \b1_{\{\osX^{t,\bx}_\cd \in A\}}  \wh{\vth}_s  ( \oW )
    \=   \int_t^s   \b1_{\{\osX^{t,\bx}  \in A\}}   \wh{\vth}  (dr ,\oW ) $, $\oP-$a.s.  for any $(s,A) \ins [t,\infty)  \ti  \sB(\OmX)$,
    where $ \osX^{t,\bx} \= \big\{\osX^{t,\bx}_s\big\}_{s \in [0,\infty)} $ is the unique solution of SDE \eqref{Ju01_01}.
 Let $\Phi $ be a nonnegative Borel-measurable  function  on $ [0,\infty) \ti \OmX$.
     Then   
     a standard approximation argument and the ``change-of-variable" formula 
     yield that
    $ E_\oP   \big[  \Phi (\oT  ,  \oX )   \big| \cF^{\oW^t,\oP}_\infty \big]
  \= \int_t^\infty \n \Phi  (  r,  \oX )    \wh{\vth}  (dr ,\oW )
  \= \int_0^1   \Phi (   \vr(\oW ,\l) , \oX )   d \l $, $\oP-$a.s.,
  where $\vr (\o_0,\l) \df \inf\big\{s \ins [t,\infty) \n : \wh{\vth}_s (  \o_0) \>  \l \big\}  $, $ \fa  (\o_0,\l) \ins \O_0 \ti (0,1) $.
  Since the  joint $\oP-$distribution of $ \big(\oW ,\osX^{t,\bx} \big)$ is equal to the joint $\fp-$distribution of $(\cB,\cX^{t,\bx})$,
   \bea \label{031922_11}
    E_\oP   \big[  \Phi (\oT  ,  \oX )    \big]
   \= \int_0^1  E_\oP  \big[ \Phi (   \vr(\oW ,\l) , \osX^{t,\bx} ) \big]  d \l
   \=   \int_0^1  E_\fp  \big[ \Phi (   \vr(\cB ,\l) , \cX^{t,\bx} ) \big]  d \l  .
   \eea
  As   $\tau_\l \df \vr(\cB ,\l) \ins \cS_t $ for each $\l \ins (0,1)$, taking $\Phi$ to  be the total reward  function
  implies  that
  \bea \label{031922_14}
  E_\oP \big[\, \oR(t)\big] \= \int_0^1 E_\fp \Big[ \int_t^{\tau_\l}  \n  f \big( r, \cX^{t,\bx}_{r \land \cd}  \big) dr
 \+ \b1_{\{\tau_\l < \infty\}} \pi \big( \tau_\l  , \cX^{t,\bx}_{ \tau_\l  \land \cd}  \big)  \Big]   d \l
 \ls \int_0^1 V(t,\bx) d \l \= V(t,\bx) .
 \eea

 However, this argument is not applicable to the expectation-constraint case: Given a $\oP \ins \oocP_{t,\bx} (y,z)$,
 since   $ \tau_\l $ may not belong to $\cS_{t,\bx} (y,z) $ for a.e. $\l \ins (0,1)$,   one can not get $E_\oP \big[\, \oR(t)\big]  \ls V(t,\bx,y,z)$ like \eqref{031922_14}.
  Actually, for each $\l \ins (0,1)$, $ \tau_\l $ is only of $\cS_{t,\bx} (y_\l,z_\l) $ with
 $(y_\l,z_\l) \= \big(\{y^i_\l\}_{i \in \hN} , \{z^i_\l\}_{i \in \hN} \big)$ and
 $(y^i_\l,z^i_\l) \df \Big(  E_\fp \big[ \int_t^{\tau_\l} \n  g_i ( r,\cX^{t,\bx}_{r \land \cd} ) dr  \big] ,   E_\fp \big[ \int_t^{\tau_\l} \n  h_i ( r, \\  \cX^{t,\bx}_{r \land \cd} ) dr  \big]  \Big) $.
 For $i \ins \hN$, choosing   accumulative cost functions for $\Phi$ in \eqref{031922_11}
 renders that
 $ \int_0^1   E_\fp \big[ \int_t^{\tau_\l} \n  g_i ( r,\cX^{t,\bx}_{r \land \cd} ) dr  \big]   d \l \= E_\oP \big[ \int_t^\oT   g_i ( r,\oX_{r \land \cd} ) dr  \big]   \ls y_i $ and similarly
   $ \int_0^1   E_\fp \big[ \int_t^{\tau_\l}    h_i ( r,\cX^{t,\bx}_{r \land \cd} ) dr  \big]  d \l
    \= z_i $, so $ V \big(t,\bx,\{ \int_0^1 y_\l d \l \}_{i \in \hN} ,    \{ \int_0^1 z_\l d \l \}_{i \in \hN}\big)  \ls V(t,\bx,y,z)$.
  Then  the attempt to show   $ E_\oP \big[\, \oR(t)\big] \ls V(t,\bx,y,z) $ reduces to deriving  a Jensen-type inequality:
 \beas
  \int_0^1 V(t,\bx,y_\l,z_\l)  d \l \ls V \Big(t,\bx, \Big\{ \int_0^1 y_\l d \l \Big\}_{i \in \hN} , \Big\{ \int_0^1 z_\l d \l \Big\}_{i \in \hN}\Big)   .
  \eeas
 But this does not hold  since the value function $V $ is not concave in level $z$ of equality-type expectation constraints.

\end{rem}

\section{The Measurability of OSEC Values}
\label{sec_Mart_prob}

   In this section, using the martingale-problem formulation of SDEs,
   we   characterize the probability class $ \ocP_{t,\bx} $ by countably many stochastic behaviors of the
   canonical coordinates $(\oW,\oX,\oT)$ of $\oO$.
   This will enable us to analyze   the measurability  of   value functions 
   of the  optimal stopping problem   with expectation constraints.

 Let $ \fS $ be the equivalence classes of  all   $[0,\infty]-$valued
 $ \bF^{W,P_0} -$stopping times on $\O_0$
 in the sense that $\tau_1, \tau_2 \ins \fS $ are equivalent if
 $ P_0 \{   \tau_1  \= \tau_2  \} \= 1 $.
 We  endow $\fS$   with the   metric
 \beas
  \Rho{\fS} (\tau_1,\tau_2) \df E_{P_0} \big[ \Rho{+}  ( \tau_1 , \tau_2 )  \big]   ,  \q \fa \tau_1,   \tau_2 \ins \fS  .
  \eeas

 \begin{lemm} \label{lem_082020_15}
   $\big(\fS, \Rho{\fS}\big)$ is a complete separable metric  space, i.e.,   a Polish space.

\end{lemm}

 For any $\tau \ins \fS$, we   define its joint distribution with $W$ under $P_0$ by $\Ga  (   \tau   ) \df P_0 \nci  ( W,  \tau    )^{-1} \ins \fP \big( \O_0 \ti \hT \big) $.

\begin{lemm} \label{lem_082020_17}

The mapping $\Ga \n :   \fS \mto \fP \big(\O_0 \ti \hT\big)$ is a continuous injection
    from $  \fS$  into $ \fP \big(\O_0 \ti \hT\big) $.

\end{lemm}

   \if{0}

We know from   Lemma \ref{lem040221} and Remark \ref{rem_061222} that for $t \ins [0,\infty)$,
the mapping $ \big( \osW^t,  \oT \- t\big)  \n : \oO \mto \O_0   \ti \ddot{\hR}$
is $\sB(\oO) \big/ \sB(\O_0)   \oti \sB\big(\ddot{\hR}\big) -$measurable.
So $(\osW^t, \oT\-t)^{-1} ( D ) \ins \sB(\oO)$ for any $  D \ins \sB ( \O_0 \ti \hT)$.
The probability measure $\oQ_{t,\oP} \df
\oP \nci \big( \osW^t, \oT \- t\big)^{-1} $ is well-defined.

 \fi

 We can use   Remark \ref{rem_ocP} and Lemma \ref{lem_082020_17}   to decompose
 the probability class $\ocP_{t,\bx}$ as the intersection of   countable many 
 action sets of processes $(\oW,\oX,\oT)$:

\begin{prop} \label{prop_Ptx_char}
For any $(t,\bx) \ins [0,\infty) \ti \OmX $,
the probability class $\ocP_{t,\bx}$ is the intersection of the following three subsets of  $\fP\big(\oO\big)$:

\no i\)    $ \ocP^1_{t,\bx} \df \big\{ \oP \ins \fP\big(\oO\big) \n : \oP  \{  \oX_s \= \bx(s), \;
\fa s \ins [0,t]   \} \= 1 \big\}$.

\no  ii\)
 $ \ocP^2_t \df \Big\{ \oP \ins \fP\big(\oO\big) \n :
  E_\oP \Big[  \Big( \oM^t_{\otau^t_n \land (t+\fr)} (\vf )  \- \oM^t_{\otau^t_n \land (t+\fs)} (\vf ) \Big) \underset{i=1}{\overset{k}{\prod}}   \b1_{    \{(\oW^t_{ t+s_i  },\oX_{ t+s_i   }) \in \cO_i  \}    }    \Big]  \n \= 0  ,  ~ \fa (\vf,n) \ins \fC(\hR^{d+l}) \ti  \hN , \,  \fa (\fs,\fr)  \\ \ins   \hQ^{2,<}_+  , \,\fa
  \{(s_i,\cO_i )\}^k_{i=1} \sb \big(\hQ \Cp [0,\fs]\big) \ti \sO (\hR^{d+l})   \Big\}$.

\no iii\)   $ \ocP^3_t \df \big\{ \oP \ins \fP\big(\oO\big) \n :  \oP  \nci (\osW^t, \oT\-t)^{-1} \ins \Ga (\fS) \big\}$.

 \if{0}
Given $t \ins [0,\infty)$,  Lemma \ref{lem_oWUX} (1) implies that $ \big( \osW^t, \oT \- t\big)  \n : \oO \mto \O_0 \ti \hT$
is $\sB(\oO) \big/ \sB(\O_0)  \oti \sB(\hT) -$measurable.  So it holds  for any $\oP \ins \fP\big(\oO\big)$
that $\oP \nci \big( \osW^t, \oT \- t\big)^{-1} \ins \fP \big( \O_0 \ti \hT \big)$.
 \fi

\end{prop}

 Based on the countable decomposition of the probability class $\ocP_{t,\bx}$ by Proposition \ref{prop_Ptx_char},
 the next proposition shows that   the graph of   probability classes $ \big\{\ocP_{t,\bx}\big\}_{(t,\bx) \in [0,\infty) \times \OmX} $
 is a Borel subset of   $ [0,\infty) \ti \OmX \ti \fP\big(\oO\big) $,  which is crucial for the measurability of the value functions $V \= \oV$.

\begin{lemm} \label{lem_082020_19}
The mapping $ \ol{\Ga} (t,\oP) \df \oP \nci \big(\osW^t, \oT\-t\big)^{-1}
\ins \fP \big( \O_0 \ti \hT \big) $, $\fa (t,\oP) \ins [0,\infty) \ti \fP\big(\oO\big)$ is continuous.
\end{lemm}

\begin{prop} \label{prop_graph_ocP}
 The graph  $\big\lan\n\big\lan  \ocP  \big\ran\n\big\ran   \df \big\{\big(t,\bx,\oP \big)  \ins [0,\infty) \ti \OmX \ti \fP\big(\oO\big) \n :  \oP \ins \ocP_{t,\bx} \big\}$ is a Borel subset of $[0,\infty) \ti \OmX \ti \fP\big(\oO\big)$.
\end{prop}

  Set $ \oD \df \big\{(t,\bx,y,z) \ins [0,\infty) \ti \OmX \ti \Re \ti \Re \n :   \ocP_{t,\bx}(y,z) \nne \es \big\} $
 and   $ \ocD \df \big\{(t,\bw,\bx,y,z) \ins [0,\infty) \ti \O_0 \ti \OmX \ti \Re \ti \Re \n :   \ocP_{t,\bw,\bx}(y,z) \nne \es \big\} $.

 \if{0}
    According to Remark \ref{rem_112220} (1b) and Theorem \ref{thm_V=oV}, if   $ h_i  \= 0$ for all $  i \ins \hN$,
     $\hb{Proj}  \big(\oD\big) \df \{(t,\bx) \ins [0,\infty) \ti \OmX \n : (t,\bx,y,z) \ins \oD \hb{ for some } (y,z) \ins \Re \ti \Re  \} \= [0,\infty) \ti \OmX$ and $\hb{Proj}  \big(\ocD\big) \df \{(t,\bw,\bx) \ins [0,\infty) \ti \O_0 \ti \OmX \n : (t,\bw,\bx,y,z) \ins \ocD \hb{ for some } (y,z) \ins \Re \ti \Re \} \= [0,\infty) \ti \O_0 \ti \OmX$.
 \fi

 \begin{cor} \label{cor_graph_ocP}

   The graph $\gP \df \big\{ \big(t,\bx,y,z, \oP\big) \ins \oD \ti \fP\big(\oO\big)  \n :
 \oP \ins \ocP_{t,\bx}(y,z)   \big\}$   is a Borel subset of $ \oD \ti \fP\big(\oO\big)$
  and the graph $\gcP \df \big\{ \big(t,\bw,\bx,y,z, \oP\big) \ins \ocD \ti \fP\big(\oO\big)  \n :
 \oP \ins \ocP_{t,\bw,\bx}(y,z)   \big\}$   is a Borel subset of $\ocD \ti \fP\big(\oO\big)$.


 \end{cor}

By Corollary \ref{cor_graph_ocP}, the value function $\oV$ is   \usa ~ and is thus universally measurable.

\begin{thm} \label{thm_V_usa}

  The value function $\oV (t,\bx,y,z)$ is \usa ~ on $   \oD$
  and the value function $\oV (t,\bw,\bx,y,z) $ is \usa ~ on $   \ocD $.

\end{thm}

\section{Dynamic Programming Principle for $\oV$ }
\label{sec_DPP}

In this section, we 
explore a   dynamic programming principle (DPP)    for the value function  $\oV$ in weak formulation, which    takes the conditional expected integrals of constraint functions   as   additional states.

  Given $t \ins [0,\infty)$, let  $\oga$ be  a  $[t,\infty)-$valued $\bF^{\oW^t}-$stopping time
  and let $\oP \ins \fP\big(\oO\big)$.
  According to Lemma 1.3.3 and  Theorem 1.1.8 of \cite{Stroock_Varadhan},
    $ \cF^{\oW^t}_\oga  $ is  countably generated and  there is thus a family $\big\{ \oP^t_{\oga,\oo} \big\}_{\oo \in \oO}$
 of probability measures in $\fP\big(\oO\big)$, called  the {\it regular conditional probability distribution} (r.c.p.d.) of $\oP$ with respect to $\cF^{\oW^t}_\oga$,  such that
\bea
\hspace{-0.8cm} &&\hb{\no (R1) for any $\oA \ins \sB(\oO)$, the mapping $\oo \mto \oP^t_{\oga,\oo} \big(\oA\big)$ is $\cF^{\oW^t}_\oga-$measurable;} \nonumber\\
\hspace{-0.8cm}  &&\hb{(R2) for any $(-\infty,\infty]-$valued, $\sB_\oP(\oO)-$measurable   random variable $ \oxi $ 
 that is bounded from below under $\oP$,
it holds} \nonumber\\
\hspace{-5cm} && \hb{for all $\oo \ins \oO$ except on a  $ \ocN_\oxi   \ins \sN_\oP\big(\cF^{\oW^t}_\oga\big)   $ that
 $ \oxi$ is $ \sB_{\oP^t_{\oga,\oo}}(\oO)   -$measurable and $ E_{\oP^t_{\oga,\oo}} \big[ \, \oxi \, \big] \= E_\oP \big[ \, \oxi \, \big| \cF^{\oW^t}_\oga \big] (\oo) $;} \label{R2}\\
 \hspace{-6cm} && \hb{(R3) for some   $ \ocN_0  \ins \sN_\oP \big(\cF^{\oW^t}_\oga\big) $,
   $ \oP^t_{\oga,\oo} \big(\oA\big) \= \b1_{\{\oo \in \oA\}} $,  $\fa \big(\oo,\oA\big) \ins \ocN^c_0  \ti \cF^{\oW^t}_\oga $.} \label{R3}
\eea

    Let  $\oo \ins \oO  $ and set $   \Wtgo   \df \big\{\oo' \ins \oO \n : \oW^t_{\n r} (\oo') \= \oW^t_{\n r} (\oo), ~ \fa r \ins [t,\oga(\oo)] \big\} $.
    We know from  Galmarino's test  that
   \bea \label{090520_11}
   \oga(\oo') \= \oga(\oo) , \q \fa  \oo' \ins \Wtgo,
  \eea
  and   $\Wtgo $ is thus  $ \cF^{\oW^t}_\oga -$measurable. 
  Since $\oo \ins \Wtgo$ for any $\oo \ins \oO $, (R3) in \eqref{R3} shows that
  \bea \label{Jan11_03}
   \oP^t_{\oga,\oo} \big(\Wtgo  \big) \= \b1_{\big\{\oo \in \Wtgo\big\}} \= 1  , \q \fa  \oo \ins \ocN^c_0.
  \eea

 For any $i \ins \hN$, define
  $   \oY^i_{\n \oP} (\oga)   \df   E_\oP \Big[   \int_{\oT \land  \oga }^\oT    g_i (r,\oX_{r \land \cd} ) dr \Big| \cF^{\oW^t}_\oga \Big]
  $ and
    $  \oZ^i_\oP (\oga)  \df   E_\oP \Big[   \int_{\oT \land  \oga }^\oT    h_i (r,\oX_{r \land \cd} ) dr \Big| \cF^{\oW^t}_\oga \Big] $.
   So  $ \big(\oY_{\n \oP} (\oga), \oZ_\oP (\oga) \big)  \df \Big( \big\{\oY^i_{\n \oP} (\oga)\big\}_{i \in \hN}, \big\{\oZ^i_\oP (\oga)\big\}_{i \in \hN} \Big)  $   is   an $ \Re \ti \Re -$valued $ \cF^{\oW^t}_\oga -$measurable  random variable. 

 In terms of the r.c.p.d. $\big\{ \oP^t_{\oga,\oo} \big\}_{\oo \in \oO}$,
 the probability class $ \big\{ \ocP_{t,\bx}(y,z) \n : (t,\bx,y,z) \ins \oD \big\}$
 is stable under conditioning as follows.
 It will play an important role in deriving the sub-solution side  of the DPP for $\oV$.

 \begin{prop} \label{prop_flow}

 Given $ (t,\bx ) \ins [0,\infty) \ti \OmX   $,
 let  $\oga $ be a  $[t,\infty)-$valued $\bF^{\oW^t}-$stopping time and let $\oP \ins \ocP_{t,\bx}$.
 There exists a $\oP-$null set $\ocN    $ such that
  \bea   \label{062821_11}
  \oP^t_{\oga,\oo} \ins \ocP_{ \oga(\oo),\oX_{\oga \land \cd}  (\oo)} \Big(  \big(\oY_{\n \oP}  (\oga )\big)    (\oo), \big(\oZ_\oP  (\oga )\big)    (\oo) \Big), \q \fa \oo \ins  \big\{\oT \gs   \oga  \big\}  \Cp  \ocN^c  .
  \eea

 \end{prop}

 Now, we are ready to present a dynamic programming principle in weak formulation for
the value function $\oV$, 
in which  $ \big(\oY_{\n \oP} (\oga),\oZ_\oP (\oga)\big) $  act  as   additional states for constraint levels at the intermediate horizon $\oga$.

\begin{thm} \label{thm_DPP1}

 Given $ (t,\bx,y,z) \ins \oD  $,
 let  $\big\{\ogaP   \big\}_{\oP \in \ocP_{t,\bx}(y,z)}$ be a   family of $[t,\infty)-$valued  $\bF^{\oW^t}  - $stopping times. Then
 \bea
 &&\hspace{-1.5cm}
 \oV(t,\bx,y,z)   \=   \Sup{\oP \in \ocP_{t,\bx}(y,z)} \n E_\oP \bigg[
 \b1_{\{\oT <  \ogaP \}} \Big( \n \int_t^\oT \n f(r,\oX_{r \land \cd}) dr  \+ \pi \big(\oT, \oX_{\oT \land \cd}\big) \Big) \nonumber \\
 &&\hspace{2.3cm} + \b1_{\{\oT \ge  \ogaP \}} \bigg( \n \int_t^\ogaP \n f(r,\oX_{r \land \cd}) dr  \+ \oV \Big(  \ogaP ,\oX_{ \ogaP  \land \cd} ,    \oY_{\n \oP}  \big( \ogaP \big) , \oZ_\oP  \big( \ogaP \big)   \Big) \bigg) \bigg] .      \label{020422_14}
  \eea

 \end{thm}

\section{Proofs}
\label{sec_proof}

\no {\bf Proof of Proposition \ref{prop_MPF1}: 1)} Set $\cN \df \big\{\o \ins \O \n: X_s (\o) \nne \bx(s) \hb{ for some } s \ins [0,t] \big\} \ins \sN_P(\cF^X_t) $ and let $(\vf,n,\fra) \ins  C^2(\hR^{d+l}) \ti \hN \ti \hR^{d+l}$.
 We denote $ c^\vf_n(\fra) \df \underset{|( w, x )| \le  n + \fra }{\sup} \big( \sum^2_{i=0} | D^i \vf ( w, x ) |   \big)
 \+ \big|\vf (0,\bx(t) )\big|   \< \infty$ and $c^n_{t,\bx}(\fra) \df \big[ d / 2 \+  \k( t\+n)  ( \|\bx\|_t  \+  n \+ \fra    )  \+     \k^2( t\+n)  ( \|\bx\|_t  \+  n \+ \fra     )^2      \big] n \+     \int_t^{t+n}   \big(  |b( r,\bz)| \+ |\si( r,\bz)|^2 \big)  dr \< \infty $.
   Given $\o \ins \cN^c$, since $\big\|X_{r \land \cd} (\o)\big\|_r \= \Sup{r' \in [0,r]} \big|X_{r'} (\o) \big|
    \ls \|\bx\|_t  \ve ( n \+ \fra)  $, $\fa r \ins \big[t,  (\tau^t_n (\fra)) (\o)\big]$,
   we can deduce from   \eqref{coeff_cond1}  and Cauchy-Schwarz inequality   that
   \bea
   && \hspace{-1.2cm} \Sup{s \in [t,(\tau^t_n (\fra)) (\o)]} \big| \big(M^t_s(\vf)\big)(\o) \big|
    \ls   \Sup{s \in [t,(\tau^t_n (\fra)) (\o)]} \big| \vf  \big( B^t_s (\o) , X_s (\o) \big) \big| \+ c^\vf_n(\fra) \n \int_t^{(\tau^t_n (\fra)) (\o)} \n \Big( \big|   b  (r, X_{r \land \cd}(\o) ) \big| \+  \frac12   \big(d\+\big|    \si   (r, X_{r \land \cd}(\o) ) \big|^2\big) \Big)  dr \nonumber  \\
   & & \hspace{-0.7cm}  \ls c^\vf_n(\fra) \+ c^\vf_n(\fra)  \n \int_t^{(\tau^t_n (\fra)) (\o)} \Big(     \k( r)   \big\|X_{r \land \cd}(\o)\big\|_r
    \+   |b( r,\bz )|    \+        d/2   \+    \k^2( r)   \big\|X_{r \land \cd}(\o)\big\|^2_r
    \+   |\si( r,\bz )|^2 \Big) dr \ls c^\vf_n(\fra) (1\+c^n_{t,\bx}(\fra))    . \q   \label{122721_11}
    \eea
     So $ \big\{ M^t_{s \land \tau^t_n (\fra) } (\vf)   \big\}_{s \in [t,\infty)} $ is a bounded $\bF^t-$adapted continuous process under $P$.

 \no {\bf 2)} We next show that  (i) implies  (ii):  Suppose that (i) holds and let   $ (\vf,n,\fra)  \ins  C^2(\hR^{d+l}) \ti \hN \ti \hR^{d+l}  $.
 We simply denote  $\Xi^{t,\bx}_s \df (B^t_s ,  X^{t,\bx}_s) $, $ \fa s \ins [t,\infty)$
 and set $\tau^{t,\bx}_n (\fra)  \df  \inf\big\{s \ins [t,\infty) \n :  |\Xi^{t,\bx}_s \- \fra|    \gs n   \big\} \ld (t\+n) $,
 which is an $\bF^{B^t,P} -$stopping time.
 Applying  It\^o's formula yields   that   $P-$a.s.
\beas
 M^{t,\bx}_s(\vf)  & \tn \df & \tn   \vf  ( \Xi^{t,\bx}_s  )
 \- \n \int_t^s \ol{b}  \big( r, X^{t,\bx}_{r \land \cd} \big)  \n \cd \n D \vf  ( \Xi^{t,\bx}_r  ) dr
    \-   \frac12 \int_t^s   \ol{\si} \, \ol{\si}^T  \big( r, X^{t,\bx}_{r \land \cd} \big) \n : \n D^2 \vf ( \Xi^{t,\bx}_r  )   dr \nonumber \\
 & \tn \= & \tn  \vf \big(0,\bx(t)\big)   \+ \int_t^s    D  \vf ( \Xi^{t,\bx}_r  )  \n \cd \n
 \ol{\si}  \big(r, X^{t,\bx}_{r \land \cd} \big) d B_r , \q s \ins [t,\infty) . 
\eeas

 For any $\o \ins \O$, an analogy to \eqref{122721_11} shows that
    $\Sup{s \in [t,(\tau^{t,\bx}_n(\fra))(\o)]} \big| \big(M^{t,\bx}_s(\vf)\big)(\o) \big| \ls c^\vf_n(\fra) (1\+c^n_{t,\bx}(\fra))  $ and   $  \int_t^{(\tau^{t,\bx}_n(\fra))(\o)} \\ \big|  D  \vf (   \Xi^{t,\bx}_r (\o)  )   \cd
 \ol{\si}  (r, X^{t,\bx}_{r \land \cd} (\o) ) \big|^2 d r
  \ls    (c^\vf_n(\fra))^2   \big[ d   \+ 2    \k^2( t\+n)  ( \|\bx\|_t \+ n \+ \fra    )^2      \big] n \+ 2   (c^\vf_n(\fra))^2 \n   \int_t^{t+n}   |\si( r,\bz)|^2  dr    \< \infty $.
 So
 \bea \label{122921_14}
  \hb{$\big\{ M^{t,\bx}_{s \land \tau^{t,\bx}_n (\fra)   }(\vf) \big\}_{s \in [t,\infty)}$ is   a    bounded  $  \bF^{B^t,P}  -$martingale.}
  \eea
 Set $ \cN_{t,\bx} \df \{\o \ins \O \n: X_s (\o) \nne   X^{t,\bx}_s (\o)  \hb{ for some } s \ins [0,\infty)\} \ins \sN_P\big(\cF^{B^t,P}_\infty \ve \cF^X_\infty\big) $.
 For any $(s,\o) \ins [0,\infty) \ti \cN^c_{t,\bx} $,
 \bea \label{122021_17}
 X^{t,\bx}_s (\o) \= X_s (\o)   ,\q \big( M^{t,\bx}_{s \vee t}(\vf) \big) (\o) \=   \big( M^t_{s \vee t}(\vf) \big) (\o)
 \q \hb{and thus} \q (\tau^{t,\bx}_n (\fra)) (\o) \= (\tau^t_n (\fra)) (\o) .
 \eea

  Fix $t_1,t_2 \ins [t,\infty)$ with $t_1 \< t_2$.
  Let $ \big\{(s_i,\cE_i)\big\}^m_{i=1} \sb [t,t_1] \ti \sB(\hR^d)$ and $ \big\{(r_j,A_j)\big\}^k_{j=1} \sb [0,t_1] \ti \sB(\hR^l)$.
 We can derive from \eqref{122921_14} and \eqref{122021_17}   that
 $
    E_P \Big[ \b1_{\cN^c_{t,\bx}} \big( M^t_{t_2 \land \tau^t_n(\fra)}(\vf) \-  M^t_{t_1 \land \tau^t_n(\fra)}(\vf) \big)  \prod^m_{i=1}   \b1_{(B^t_{s_i})^{-1}(\cE_i)} \prod^k_{j=1} \b1_{ X_{r_j}^{-1}(A_j)}   \Big]
     \=  E_P \Big[ \b1_{\cN^c_{t,\bx}} \big( M^{t,\bx}_{t_2 \land \tau^{t,\bx}_n(\fra)}(\vf) \-  M^{t,\bx}_{t_1 \land \tau^{t,\bx}_n(\fra)}(\vf) \big)  \prod^m_{i=1} \b1_{(B^t_{s_i})^{-1}(\cE_i)} \prod^k_{j=1} \b1_{ (X^{t,\bx}_{r_j})^{-1}(A_j)}  \Big]
 \= 0 $.
   So   the Lambda-system $ \L \df \big\{ A \ins \cF \n : E_P \big[   \big( M^t_{t_2 \land \tau^t_n(\fra) } (\vf) \- M^t_{t_1 \land \tau^t_n(\fra) } (\vf)  \big)   \b1_A \big] \= 0  \big\}$   contains the  Pi-system
 $ \Big\{  \Big( \underset{i=1}{\overset{m}{\cap}}  (B^t_{s_i})^{-1}(\cE_i) \Big) \Cp \Big( \underset{j=1}{\overset{k}{\cap}} X_{r_j}^{-1}(A_j) \Big) \n :  \big\{(s_i,\cE_i)\big\}^m_{i=1} \sb [t,t_1] \ti \sB(\hR^d) , \, \big\{(r_j,A_j)\big\}^k_{j=1} \sb [0,t_1] \ti \sB(\hR^l)   \Big\} $,
  which generates $\cF^t_{t_1}$.    Dynkin's Pi-Lambda Theorem (see e.g Theorem 3.2 of \cite{Billingsley_PM}) renders   $ \cF^t_{t_1} \sb \L $, i.e.,
  $E_P \big[   \big( M^t_{t_2 \land \tau^t_n(\fra) } (\vf) \- M^t_{t_1 \land \tau^t_n(\fra) } (\vf)   \big)   \b1_A \big] \= 0 $, $\fa A \ins \cF^t_{t_1}$.
  Hence, $\big\{ M^t_{s \land \tau^t_n(\fra) } (\vf) \big\}_{s \in [t,\infty)}$ is a bounded $ \bF^t-$martingale.

 \no {\bf 3)} As $\fC(\hR^{d+l}) \sb C^2(\hR^{d+l}) $, (ii) \n $\Ra$ \n (iii)   is straightforward. It remains to show that (iii) gives rise to (i).

 \no {\bf 3a)} Let  $\bF^{t,P} \= \big\{ \cF^{t,P}_s \big\}_{s \in [t,\infty)}$ be the $P-$augmentation of $ \bF^t  $
   \Big(i.e.,  $\cF^{t,P}_s \df \si(\cF^t_s \cp \sN_P(\cF^t_\infty))$
   with $\cF^t_\infty \df \si \Big(\underset{s \in [t,\infty)}{\cup}\cF^t_s \Big) $\Big)
   We define $  \cF^{t,P}_{s+} \df \ccap{\e>0}{} \cF^{t,P}_{s+\e} $, $\fa s \ins   [t,\infty)  $
   and set $\bG^{t,P} \= \big\{ \cG^{t,P}_s \df \cF^{t,P}_{s+}\big\}_{s \in   [t,\infty)}$.

  Let $i,j \ins \{ 1,\cds \n ,d \}$. We  set $  \phi_i(w,x) \df w_i $ and $\phi_{ij}(w,x) \df w_i w_j $
 for any $w \= (w_1,\cds \n ,w_d) \ins \hR^d$ and $  x \ins \hR^l$. Clearly,  $ \phi_i , \phi_{ij} \ins \fC(\hR^{d+l})$.
 One can calculate that   $ M^t_s(\phi_i) \= B^{t,i}_s   $,   $ M^t_s(\phi_{ij}) \=  B^{t,i}_s B^{t,j}_s \- \d_{ij} (s\-t) $, $\fa s \ins [t,\infty)$,
 where  $B^t_s \= \big( B^{t,1}_s ,\cds \n , B^{t,d}_s  \big)$ and $\d_{ij}$ is the $(i,j)-$element of the identity matrix $I_{d \times d}$.

 Let $n \ins \hN $.  By (iii),   $ \big\{ M^t_{s \land \tau^t_n } (\phi_i) \big\}_{s \in [t,\infty)} $ and $ \big\{ M^t_{s \land \tau^t_n } (\phi_{ij}) \big\}_{s \in [t,\infty)} $ are bounded $   \bF^t-$martingales and are thus
  bounded  $   \bF^{t,P}  -$martingales.  
  \if{0}

  For any $s_1,s_2 \ins [t,\infty)$ with $s_1 \< s_2$, the {\it optional sampling} theorem (e.g. Theorem 1.3.22 of \cite{Kara_Shr_BMSC}) shows that
  $ E_P \Big[ M^t_{s_2 \land \tau^t_n } (\phi_i) \big| \cF^{t,P}_{s_1+}\Big] \= M^t_{s_1 \land \tau^t_n } (\phi_i) $ and
  $ E_P \Big[ M^t_{s_2 \land \tau^t_n } (\phi_{ij}) \big| \cF^{t,P}_{s_1+}\Big] \= M^t_{s_1 \land \tau^t_n } (\phi_{ij})$,  $P-$a.s.
  So   $ \big\{ M^t_{s \land \tau^t_n } (\phi_i) \big\}_{s \in [t,\infty)} $ and $ \big\{ M^t_{s \land \tau^t_n } (\phi_{ij}) \big\}_{s \in [t,\infty)} $
  are further  $   \bG^{t,P}-$martingales.

  \fi
  The {\it optional sampling} theorem (e.g. Theorem 1.3.22 of \cite{Kara_Shr_BMSC}) implies that
  they are further  $   \bG^{t,P}-$martingales.
  Since  $\lmtu{n \to \infty} \tau^t_n \= \infty$,
  we see that $\big\{M^t_s(\phi_i) \= B^{t,i}_s   \big\}_{s \in [t,\infty)}  $
  and $\big\{M^t_s(\phi_{ij}) \= B^{t,i}_s B^{t,j}_s \- \d_{ij} (s\-t) \big\}_{s \in [t,\infty)}$   are $   \bG^{t,P}-$local martingales.
    L\'evy's characterization theorem then yields that $ B^t   $ is a    Brownian motion with respect to  filtration $   \bG^{t,P} $
   and is thus a Brownian motion with respect to   filtration $ \bF^{B^t} $. 

 \no {\bf 3b)}
 We simply denote   $\Xi_s \df (B^t_s ,  X_s) $, $\beta_s \df \ol{b}  \big(s, X_{s \land \cd} \big)$ and
 $ \a_s \df \ol{\si} \, \ol{\si}^T  \big( s,  X_{s \land \cd} \big) $, $ \fa  s \ins [t,\infty)$.
 Let  $i,j  \ins \{ 1,\cds \n ,d\+l\}$. We set $  \psi_i(\nxi) \df \nxi_i  $
 and $\psi_{ij}(\nxi) \df \nxi_i  \nxi_j  $
 for any  $  \nxi \= \big( \nxi_1   ,\cds \n , \nxi_{d + l} \big)  \ins \hR^{d+l}$.
 Similar to $ M^t_\cd (\phi_i) $ and $ M^t_\cd   (\phi_{ij}) $,
  the processes  $  M^t_s(\psi_i) \=    \Xi^{(i)}_s   \-  \int_t^s \beta^{(i)}_r  dr$  and
$ M^t_s(\psi_{ij}) \= \Xi^{(i)}_s  \Xi^{(j)}_s
 \- \int_t^s \beta^{(i)}_r \Xi^{(j)}_r dr
 \- \int_t^s  \beta^{(j)}_r   \Xi^{(i)}_r dr
 \- \int_t^s  ( \a_r  )_{ij} dr $, $ s \ins [t,\infty)$
 are  $ \bG^{t,P}-$local martingales.
 Using the {\it integration by parts} formula, we obtain   that  $P-$a.s.
 \beas
 && \hspace{-1.5cm}
 \Xi^{(i)}_s  \Xi^{(j)}_s  \-  M^t_s(\psi_i) M^t_s(\psi_j)
   \=      M^t_s(\psi_i) \int_t^s  \beta^{(j)}_r    dr \+M^t_s(\psi_j) \int_t^s \beta^{(i)}_r  dr \+ \int_t^s \beta^{(i)}_r  dr \cd  \int_t^s  \beta^{(j)}_r    dr \\
  &&  \hspace{-0.5cm} \=    \int_t^s M^t_r(\psi_i)  \beta^{(j)}_r    dr
  \+  \int_t^s   \Big(  \int_t^r \beta^{(j)}_{r'}  dr' \Big) d M^t_r(\psi_i)
  \+ \int_t^s M^t_r(\psi_j) \beta^{(i)}_r  dr
  \+  \int_t^s   \Big(  \int_t^r \beta^{(i)}_{r'}  dr' \Big) d M^t_r(\psi_j) \\
  && \hspace{-0.5cm} \q
   +  \int_t^s \Big(  \int_t^r \beta^{(i)}_{r'}  dr' \Big)   \beta^{(j)}_r    dr
  \+ \int_t^s  \Big(  \int_t^r \beta^{(j)}_{r'}  dr' \Big) \beta^{(i)}_r  dr \\
  && \hspace{-0.5cm} \=       \int_t^s \Big[ \Xi^{(i)}_r  \beta^{(j)}_r
  \+ \Xi^{(j)}_r \beta^{(i)}_r \Big] dr
  \+  \int_t^s  \n \Big(  \int_t^r \beta^{(i)}_{r'}  dr' \Big) d M^t_r(\psi_i)
  \+  \int_t^s  \n \Big(  \int_t^r \beta^{(j)}_{r'}  dr' \Big) d M^t_r(\psi_j)  , \q  s \ins [t,\infty) .
 \eeas
   \if{0}
 Note
 \beas
 \int_t^s \beta^{(i)}_r  dr \int_t^s  \beta^{(j)}_r
 \= \int_t^s \Big(  \int_0^r \beta^{(i)}_{r'}  dr' \Big)   \beta^{(j)}_r    dr
  \+ \int_t^s  \Big(  \int_0^r \beta^{(j)}_{r'}  dr' \Big) \beta^{(i)}_r  dr ,
   \q s \ins [t,\infty) .
 \eeas
   \fi
 So $ M^t_s(\psi_i) M^t_s(\psi_j) \- \int_t^s  ( \a_r  )_{ij} dr
 \= M^t_s(\psi_{ij}) \-  \int_t^s   \big(  \int_t^r \beta^{(i)}_{r'}  dr' \big) d M^t_r(\psi_i)
\- \int_t^s   \big(  \int_t^r \beta^{(j)}_{r'}  dr' \big) d M^t_r(\psi_j) $,
  $s \ins [t,\infty)$  is also an $ \bG^{t,P}-$local martingale, which implies that
 the quadratic variation of  the $ \bG^{t,P}-$local martingale
  $  M^t_s  \df \big( M^t_s(\psi_1), \cds \n , \\  M^t_s(\psi_{d+l})\big)
  \= \Xi_s  \-  \int_t^s \beta_r   dr $, $ s \ins [t,\infty)$   is
  $ \big\lan  M^t ,  M^t    \big\ran_s  \= \int_t^s   \a_r   dr $, $ s \ins [t,\infty) $.

  Let $ n \ins \hN $, $a \ins \hR^l$ and set $ \dis \cH^a_s \df   \binom{-   \si^T \big( s, X_{s \land \cd} \big) a }{a}$, $ \fa s \ins (t,\infty)$.
   The stochastic exponential   of the  $ \bG^{t,P}-$ martingale
  $\big\{\int_t^{\tau^t_n \land s} \cH^a_r \n \cd \n d M^t_r\big\}_{s \in [t,\infty)}$ is
  \beas
 && \hspace{-1.2cm} \exp\Big\{\int_t^{\tau^t_n \land s} \cH^a_r \n \cd \n d M^t_r \- \frac12 \int_t^{\tau^t_n \land s}   (\cH^a_r)^T  \a_r  \cH^a_r  dr \Big\}
  \= \exp\Big\{\int_t^{\tau^t_n \land s} \cH^a_r \n \cd \n d \Xi_r  \-  \int_t^{\tau^t_n \land s} \cH^a_r \n \cd \n  \beta_r   dr  \Big\} \\
 && \= \exp\bigg\{ a \n \cd \n \Big( \int_t^{\tau^t_n \land s}        d X_r  \- \int_t^{\tau^t_n \land s}    \si \big( r, X_{r \land \cd} \big)   d B_r  \-  \int_t^{\tau^t_n \land s}      b \big( r, X_{r \land \cd} \big)   dr \Big)   \bigg\} , \q s \ins [t,\infty) .
  \eeas
  \if{0}
  We can deduce from \eqref{coeff_cond1} that  $\big\{\int_0^{\tau^t_n \land s} \cH^a_r \n \cd \n d M^t_r\big\}_{s \in [t,\infty)}$ is a $  \bF^t-$BMO martingale
  and thus its stochastic exponential is a  $   \bF^t-$uniformly integrable martingale.
  \fi
  Letting $a $ vary over $\hR^l$ yields that $P-$a.s.,
  $ X_{ \tau^t_n \land s} \=   \bx(t)  \+  \int_t^{\tau^t_n \land s}     b \big( r, X_{r \land \cd} \big)   dr \+ \int_t^{\tau^t_n \land s}    \si \big( r, X_{r \land \cd} \big)   d B_r   $,  $ \fa s \ins [t,\infty) $.
   Sending $n \nto \infty$ then renders  that      $P-$a.s.,
  $ X_s \=   \bx(t)  \+  \int_t^s    b \big( r, X_{r \land \cd} \big)   dr \+ \int_t^s    \si \big( r, X_{r \land \cd} \big)   d B_r   $,  $ \fa s \ins [t,\infty) $.
  Viewing  SDE \eqref{121621_11}  on   $\big(\O,\cF, \bG^{t,P},P\big)$,
  we know from Proposition \ref{prop_122021} that  there is a unique   $ \big\{\cG^{t,P}_{s \vee t}   \big\}_{s \in [0,\infty)}  -$adapted continuous process satisfying   \eqref{121621_11}. Hence, $P\{X_s \= X^{t,\bx}_s , \;  \fa s \ins [0,\infty)\} \= 1$.  \qed
  \if{0}

  {\bf 4)} Let  $B^t$ be  a   Brownian motion on  $(\O,\cF,P)$ and let $(\vf,n) \ins \fC(\hR^{d+l}) \ti \hN $.  On $\O_0$,
  \beas    
  \q ^o \n M^t_s(\vf)   \df   \vf \big(W^t_s  , \, ^o \n X^{t,\bx}_s \big)
    \- \n \int_t^s  \n  \ol{b}  \big( r, \, ^o \n X^{t,\bx}_{r \land \cd} \big) \n \cd \n D \vf \big( W^t_r  , \, ^o \n X^{t,\bx}_r \big) dr
    \-   \frac12 \n \int_t^s  \n  \ol{\si} \, \ol{\si}^T  \big( r, \, ^o \n X^{t,\bx}_{r \land \cd}  \big) \n : \n D^2 \vf  ( W^t_r, \, ^o \n X^{t,\bx}_r  )   dr  ,     ~ \fa s \ins [t,\infty)
 \eeas
 is an $\bF^{W^t,P_0}-$adapted continuous process
 and  $^o   \tau^t_n   \df  \inf\big\{s \ins [t,\infty) \n : \big|(W^t_s  , \, ^o \n X^{t,\bx}_s ) \big|   \gs n  \big\} \ld (t\+n) $  is an $\bF^{W^t,P_0}-$stopping time.
 As   $ ^o \n X^{t,\bx}  $ is
  the unique strong solution of \eqref{121621_11} on $   \big( \O_0,\sB(\O_0),   P_0   \big) $ against Brownian  motion  $   W^t  $,
  Part \(iii\) shows that     $ \big\{ ^o \n M^t_{s \land ^o\tau^t_n } (\vf) \big\}_{s \in [t,\infty)} $  is a bounded $\bF^{W^t,P_0}-$martingale.

   Applying Lemma \ref{lem_122921_11} with $t_0 \= t$,  $(\O_1, \cF_1, P_1,B^1)   \= \big(\O, \cF, P , B \big) $, $(\O_2, \cF_2, P_2,B^2) \= \big(\O_0, \sB(\O_0) , P_0 , W\big) $ and $\Phi \= B$ implies that
     $ \breve{X}^{t,\bx}_s \df ^o \n X^{t,\bx}_s  ( B)$, $  \breve{M}^t_s(\vf) \df  \big( ^o \n M^t_s(\vf) \big) (B)$, $s \ins [0,\infty)$
     are $\bF^{B^t,P}-$adapted continuous processes and  $\breve{\tau}^t_n \df ^o   \tau^t_n  (B)   $  is an $\bF^{B^t,P }-$stopping time.
  Define   filtration $\breve{\bF}^{t,\bx}  \= \{\breve{\cF}^{t,\bx}_s\}_{s \in [t,\infty)}$ by $\breve{\cF}^{t,\bx}_s \df 
 \si\big(B^t_r; r \ins [t,s]\big) \ve \si (\breve{X}^{t,\bx}_r; r \ins [0,s] ) \sb \cF^{B^t,P}_s$, $\fa s \ins [t,\infty)$.

      Let $0 \ls s \< r \< \infty$ and $ \big\{(s_i,\cE_i)\big\}^k_{i=1} \sb [0,s] \ti \sB(\hR^d)  $.
     We can also deduce from  Lemma \ref{lem_122921_11}   that
    \beas
      0 & \tn \= & \tn   E_{P_0} \Big[ \Big( {^o \n M^t}_{\tn r \land ^o   \tau^t_n   } (\vf) \- ^o \n M^t_{ s \land  ^o   \tau^t_n   } (\vf) \Big) \b1_{\ccap{i=1}{k} (W^t_{s_i})^{-1} (\cE_i)} \Big] \\
      & \tn  \= & \tn  \int_{\o_0 \in \O_0}  \Big( \big( {^o \n M^t} (\vf) \big) \big(r \ld ^o   \tau^t_n (\o_0),\o_0 \big) \- \big( {^o \n M^t} (\vf) \big) \big(s \ld ^o   \tau^t_n (\o_0),\o_0 \big) \Big)   \b1_{\ccap{i=1}{k} \{  W^t_{s_i} (\o_0) \in \cE_i \}} (P \nci B^{-1}) (d \o_0) \\
      & \tn  \= & \tn  \int_{\o  \in \O }  \Big( \big( {^o \n M^t} (\vf) \big) \big(r \ld ^o   \tau^t_n (B(\o)),B(\o) \big) \- \big( {^o \n M^t} (\vf) \big) \big(s \ld ^o   \tau^t_n (B(\o)),B(\o) \big) \Big)  \b1_{\ccap{i=1}{k} \{  W^t_{s_i} (B(\o)) \in \cE_i \}}  P    (d \o ) \\
       & \tn  \= & \tn  \int_{\o  \in \O }  \Big( \big( {^o \n M^t} (\vf) \big) \big(r \ld  \breve{\tau}^t_n  (\o) ,B(\o) \big) \- \big( {^o \n M^t} (\vf) \big) \big(s \ld \breve{\tau}^t_n  (\o),B(\o) \big) \Big)  \b1_{\ccap{i=1}{k} \{  B^t_{s_i}  (\o)  \in \cE_i \}}  P    (d \o ) \\
      & \tn  \= & \tn  E_P \Big[ \Big( \breve{M}^t_{  r \land \breve{\tau}^t_n   } (\vf) \- \breve{M}^t_{ s \land \breve{\tau}^t_n   } (\vf) \Big) \b1_{\ccap{i=1}{k} (B^t_{s_i})^{-1} (\cE_i)} \Big] .
      \eeas
      So the Lambda-system $\breve{\L}_{s,r} \df \big\{ A \ins \cF^{B^t,P}_\infty \n : E_P \big[ \big( \breve{M}^t_{  r \land \breve{\tau}^t_n   } (\vf) \- \breve{M}^t_{ s \land \breve{\tau}^t_n   } (\vf) \big) \b1_A \big] \= 0 \big\}$ includes
      the Pi-system $\Big\{ \ccap{i=1}{k} (B^t_{s_i})^{-1} (\cE_i) \n : \big\{(s_i,\cE_i)\big\}^k_{i=1} \sb [0,s] \ti \sB(\hR^d) \Big\} \Cp \sN_P(\cF^{B^t}_\infty)$.
      An application of Dynkin's Pi-Lambda Theorem shows that
      $ \cF^{B^t,P}_s \sb \breve{\L}_{s,r} $, i.e.,  $  E_P \big[ \big( \breve{M}^t_{  r \land \breve{\tau}^t_n   } (\vf) \- \breve{M}^t_{ s \land \breve{\tau}^t_n   } (\vf) \big) \b1_A \big] \= 0  $ for any $A \ins \cF^{B^t,P}_s$.
      To wit, the $\breve{\bF}^{t,\bx}-$adapted continuous process
  $\breve{M}^t_s(\vf) \=    \vf \big( B^t_s  , \breve{X}^{t,\bx}_s \big)
    \- \n \int_t^s  \n  \ol{b}  \big( r, \breve{X}^{t,\bx}_{r \land \cd} \big) \n \cd \n D \vf \big( B^t_r  , \breve{X}^{t,\bx}_r \big) dr
    \-   \frac12 \n \int_t^s  \n  \ol{\si} \, \ol{\si}^T  \big( r, \breve{X}^{t,\bx}_{r \land \cd}  \big) \n : \n D^2 \vf  ( B^t_r, \breve{X}^{t,\bx}_r  )   dr$, $\fa s \ins [t,\infty)$ stopped by the $\breve{\bF}^{t,\bx}-$stopping time
    $ \breve{\tau}^t_n   \=  \inf\big\{s \ins [t,\infty) \n : \big|(B^t_s  , \breve{X}^{t,\bx}_s ) \big|   \gs n  \big\} \ld (t\+n)$
    is a bounded $\breve{\bF}^{t,\bx}-$martingale.  Using Part (i) yields that $P\{   \breve{X}^{t,\bx}_s \=   X^{t,\bx}_s,   \fa s \ins [0,\infty)\}=1$.

   Moreover, for any $ \phi \ins \sH_o $, one has
   $    E_P  \big[ \int_t^\infty \n \phi^-(r,  X^{t,\bx}_{r \land \cd}  ) dr\big] \= E_P  \big[ \int_t^\infty \n \phi^-(r, \breve{X}^{t,\bx}_{r \land \cd}  ) dr\big]
  \= E_P  \big[ \int_t^\infty \n \phi^- \big(r, ^o \n X^{t,\bx}_{r \land \cd} (B) \big) dr\big] \= E_{P_0} \big[ \int_t^\infty \n \phi^-(r,^o \n X^{t,\bx}_{r \land \cd} ) dr\big]  \< \infty $.

 \fi

   \if{0}

 \no {\bf Proof of Remark \ref{rem_ocP}:} As the equivalence between \(D1\)+\(D2\) and  \(D1\,$'$\)   is straightforward by the martingale-problem formulation of SDEs on $\oO$,
 we only verify the equivalence between (D3) and (D3') under (D1).
 Define two processes on $\O_0$ by:
\beas
    \beta_s (\o_0) \df W_{s \vee t} (\o_0) \- W_t (\o_0) \aand
  \ddot{\beta}_s(\o_0) \df W_{(s-t)^+}(\o_0)   ,  \q \fa (s,\o_0) \ins [0,\infty) \ti \O_0 .
\eeas

\no {\bf (1)} We first show   (D3) \n $\Ra$ \n (D3') under (D1):   Let $\wh{\tau}$ be   a $[t,\infty]-$valued $\bF^{W^t,P_0} -$stopping time   on $\O_0$. We claim that
\beas
  \wh{\tau}\= \wh{\tau}(\beta)  , \q \hb{$P_0- $a.s.}
  \eeas
  Since it holds for any $(r,\cE) \ins [t,\infty) \ti \sB(\hR^d)$ that    $\beta^{-1} \big(\{W^t_r \ins \cE\}\big) \= \big\{W^t_r(\beta) \ins \cE \big\} \=  \{ W_r(\beta) \- W_t(\beta) \ins \cE \big\}
 \=  \{  \beta_r \-  \beta_t \ins \cE \big\}
\= \{W_r  \- W_t  \ins \cE\} \= \{W^t_r   \ins \cE\} $,
  the sigma-field $\{A_0 \sb \O_0 : \beta^{-1}(A_0) \= A_0 \} $ contains all generating sets of $\cF^{W^t}_\infty$ and thus includes
$\cF^{W^t}_\infty$.
 For any $\cN_0 \ins \sN_{P_0} \big(\cF^{W^t}_\infty\big)$, there exists $A_0 \ins \cF^{W^t}_\infty$ such that $\cN_0  \sb A_0  $ and  $ P_0(A_0) \= 0$.
As $ \beta^{-1}(\cN_0) \sb \beta^{-1}(A_0) \= A_0 $, we see that $ \beta^{-1}(\cN_0) \ins \sN_{P_0}\big(\cF^{W^t}_\infty\big)$. Hence,
 \bea \label{010722_11}
 \beta^{-1}(A_0) \= A_0 , \q \fa A_0 \ins \cF^{W^t}_\infty \aand
 \beta^{-1}(\cN_0) \sb \sN_{P_0} \big(\cF^{W^t}_\infty\big) , \q \fa  \cN_0 \ins \sN_{P_0} \big(\cF^{W^t}_\infty \big) .
 \eea

 Let $n \ins \hN$ and set $s^n_i \df t \+ i2^{-n}$, $\fa i \ins \hN \cp \{0\}   $.
  We denote $A^n_i  \df  \{ s^n_{i-1}   \ls  \wh{\tau}  \< s^n_i  \}  \ins \cF^{W^t,P_0}_{ s^n_i }  $, $\fa i \ins \hN   $
  and $A^n_\infty  \df  \{    \wh{\tau}  \= \infty  \}  \ins \cF^{W^t,P_0}_\infty  $.
 Define $\wh{\tau}_n \df \sum_{i \in \hN} s^n_i \b1_{A^n_i}   \+ \infty \b1_{A^n_\infty}$.
 Clearly,
    $\wh{\tau}_n(\beta) \= \sum_{i \in \hN} s^n_i \b1_{\beta^{-1}(A^n_i)}  \+ \infty \b1_{\beta^{-1}(A^n_\infty)} $
    is equal to $ \wh{\tau}_n $  on $\cA_n \df \ccup{i \in \hN \cup \{\infty\}}{} \big(  A^n_i  \Cp \beta^{-1}(A^n_i)   \big) $.
  For any $i \ins \hN \cp \{\infty\}$, since there exists  some    $\wA^n_i   \ins  \cF^{W^t}_{ s^n_i }$
   such that   $ \cN^n_i \df  A^n_i   \D  \wA^n_i  \ins  \sN_{P_0}\big(\cF^{W^t}_\infty\big)  $ (see e.g.   Problem 2.7.3 of \cite{Kara_Shr_BMSC}),
one can deduce from \eqref{010722_11} that
 $ (A^n_i)^c \Cp \beta^{-1}(A^n_i)   \sb 
\big( (A^n_i)^c \Cp \beta^{-1} \big(\wA^n_i  \big)   \big) \cp \big( (A^n_i)^c \Cp \beta^{-1}(  \cN^n_i)   \big)
\= \big(  (A^n_i)^c \Cp \wA^n_i    \big) \cp \big( (A^n_i)^c \Cp \beta^{-1}(  \cN^n_i)   \big)
\sb \cN^n_i \cp \beta^{-1}(  \cN^n_i) \ins \sN_{P_0} \big(\cF^{W^t}_\infty \big)  $.
  So $\cA^c_n \= \ccup{i \in \hN \cup \{\infty\}}{} \big( ( A^n_i )^c \Cp \beta^{-1}(A^n_i)   \big) \ins \sN_{P_0} \big(\cF^{W^t}_\infty\big)$.

 Set $\cN_* \df \ccup{n \in  \hN}{} \cA^c_n   \ins \sN_{P_0} \big(\cF^{W^t}_\infty\big)$.
 Given $\o_0 \ins \cN^c_* \= \ccap{n \in  \hN}{} \cA_n $, as $\wh{\tau}   \= \lmtd{n \to \infty} \wh{\tau}_n   $, we have
 $ \wh{\tau} \big(\beta(\o_0)\big) \= \lmtd{n \to \infty} \wh{\tau}_n  \big(\beta(\o_0)\big) \= \lmtd{n \to \infty} \wh{\tau}_n  ( \o_0 )
 \= \wh{\tau} (\o_0) $, proving the claim.

 Since $\oW^t$ is a Brownian motion under $\oP$ by (D1), applying Lemma \ref{lem_122921_11} with $t_0 \= t$, $(\O_1, \cF_1, P_1,B^1)   \= \big(\oO ,  \sB(\oO ),  \oP, \oW \big) $, $(\O_2, \cF_2, P_2,B^2) \= \big(\O_0,  \sB(\O_0),  P_0,  W \big) $  and $\Phi \= \oW$ yields that
 $  \ocN_* \df \oW^{-1} (\cN_*) \ins \sN_\oP \big(\cF^{\oW^t}_\infty\big) $. So
 \bea \label{010722_14}
 \wh{\tau} \big(\beta(\oW(\oo))\big)  \= \wh{\tau} \big(\oW(\oo)\big) , \q \fa \oo \ins \ocN^c_* \=  \oW^{-1} (\cN^c_*) .
 \eea

  As a shifted canonical process on $\O_0$,   $\sW^t_\fs(\o_0) \df W_{t+\fs}(\o_0) \- W_t (\o_0) \= W^t_{t+\fs}(\o_0) $, $  (\fs,\o_0) \ins [0,\infty) \ti \O_0$  is also a Brownian motion under $P_0$.  
   Since it holds for any $(s,\o_0) \ins [0,\infty) \ti \O_0$ that
   $ \sW^t_\fs \big( \ddot{\beta} (\o_0)\big) \= W_{t+\fs} \big( \ddot{\beta} (\o_0)\big) \- W_t \big( \ddot{\beta} (\o_0)\big)
    \= \ddot{\beta}_{t+\fs} (\o_0)  \-   \ddot{\beta}_t (\o_0) 
    \=  W_\fs(\o_0) $,
 using Lemma \ref{lem_122921_11} again with $t_0 \= 0$, $(\O_1, \cF_1, P_1,B^1)   \= \big(\O_0,  \sB(\O_0),  P_0,W\big) $, $(\O_2, \cF_2, P_2,B^2) \= \big(\O_0,  \sB(\O_0),  P_0, \sW^t\big) $ and $\Phi \= \ddot{\beta}$   shows that
  \bea  \label{010622_11}
  \ddot{\beta}^{-1} \big(\cF^{\sW^t,P_0}_\fs\big) \sb \cF^{W,P_0}_\fs  , \q \fa \fs \ins [0,\infty] .
  \eea
 As  $\cF^{\sW^t}_\fs \= \si\big(\sW^t_\fr; \fr \ins [0,\fs] \Cp \hR \big) 
 \= \si\big( W^t_r; r \ins [t,t\+\fs]  \Cp \hR \big) \= \cF^{W^t}_{t+\fs} $ for any $\fs \ins [0,\infty]$, we also obtain that
 $  
  \cF^{\sW^t,P_0}_\fs \= \si\big( \cF^{\sW^t}_\fs \cp \sN_{P_0} (\cF^{\sW^t}_\infty) \big)
 \= \si\big( \cF^{W^t}_{t+\fs} \cp \sN_{P_0} \big(\cF^{W^t}_\infty\big) \big) \= \cF^{W^t,P_0}_{t+\fs} $, $ \fa \fs \ins [0,\infty) $.

 Define $\ddot{\tau}(\o_0) \= \wh{\tau} \big( \ddot{\beta} (\o_0)\big) \- t$, $\fa \o_0 \ins \O_0$.
 For any  $\fs \ins [0,\infty)$, since $ \big\{\wh{\tau} \ins [t,t\+\fs]\big\} \ins \cF^{W^t,P_0}_{t+\fs} \= \cF^{\sW^t,P_0}_\fs $,
   \eqref{010622_11} renders that $\big\{ \ddot{\tau} \ins [0,\fs]\big\} \= \big\{\wh{\tau} \big(\ddot{\beta}) \ins [t,t\+\fs] \big\} \= \ddot{\beta}^{-1}\big(\{\wh{\tau} \ins [t,t\+\fs]\}\big)  \ins \cF^{W,P_0}_\fs$.
 So $ \ddot{\tau} $ is a  $[0,\infty]-$valued $\bF^{W,P_0}-$stopping time on $\O_0$.
 For any $ s \ins [0,\infty) $,
 $ \ddot{\beta}_s (\osW^t)   \= \osW^t_{(s-t)^+} \= \oW_{t + (s-t)^+} \- \oW_t  \= \oW_{s \vee t} \- \oW_t \=    \beta_s (\oW)   $.
Then \eqref{010722_14} implies that
  $ \ddot{\tau} \big(\osW^t(\oo)\big) \= \wh{\tau} \big( \ddot{\beta} \big(\osW^t(\oo)\big)\big) \- t
   \= \wh{\tau} \big(\beta (\oW(\oo))\big) \- t \= \wh{\tau}  \big( \oW (\oo) \big) \- t $, $ \fa \oo \ins \ocN^c_* $.
 If   $   \oP \big\{ \oT \=   \wh{\tau}(\oW ) \big\} \= 1$, we also have
 $  \oP \big\{ \oT \= t \+ \ddot{\tau}  \big(\osW^t \big)   \big\} \= 1 $. Hence,  (D3) gives rise to   (D3') under (D1).

\no {\bf (2)} Next, we   show  (D3') $\Ra$  (D3): Let $ \ddot{\tau} $ be a  $[0,\infty]-$valued $\bF^{W,P_0}-$stopping time on $\O_0$.

 Since   ${\ddot{\beta}}^t_s \df \ddot{\beta}_s \- \ddot{\beta}_t \= W_{s-t}$, $s \ins [t,\infty)$ is   a Brownian motion on $\big(\O_0,  \sB(\O_0),  P_0\big)$
 and since $  {\ddot{\beta}}^t_s\big(\sW^t(\o_0)\big) \= \sW^t_{s-t}(\o_0) \= W_s(\o_0) \- W_t(\o_0) \= W^t_s (\o_0)$ for any $(s,\o_0) \ins [t,\infty) \ti \O_0 $,
 applying Lemma \ref{lem_122921_11} with $t_0 \= t$,  $(\O_1, \cF_1, P_1,B^1)   \= \big(\O_0,  \sB(\O_0),  P_0,W\big) $, $(\O_2, \cF_2, P_2,B^2) \= \big(\O_0,  \sB(\O_0),  P_0, \ddot{\beta}\big) $ and $\Phi \= \sW^t$   shows that
  \bea \label{010722_17}
   (\sW^t)^{-1} ( \cF^{{\ddot{\beta}}^t,P_0}_s) \sb \cF^{W^t,P_0}_s , \q \fa s \ins [t,\infty) .
   \eea
       Since $\cF^{{\ddot{\beta}}^t}_s \= \si\big({\ddot{\beta}}^t_r; r \ins [t,s] \Cp \hR \big)   \= \si\big( W_{r-t}; r \ins [t,s]  \Cp \hR \big)
 \= \si\big( W_\fr; \fr \ins [0,s\-t]  \Cp \hR  \big) \= \cF^W_{s-t} $ for any $s \ins [t,\infty]$,  we see that
 $ \cF^{{\ddot{\beta}}^t,P_0}_s \= \si\big( \cF^{{\ddot{\beta}}^t}_s \cp \sN_{P_0} (\cF^{{\ddot{\beta}}^t}_\infty) \big)
 \= \si\big( \cF^W_{s-t} \cp \sN_{P_0} (\cF^W_\infty) \big) \= \cF^{W,P_0}_{s-t} $, $ \fa s \ins [t,\infty) $.

 Define $\wh{\tau} (\o_0) \df  t \+ \ddot{\tau} \big(\sW^t (\o_0)\big)   $, $\fa \o_0 \ins \O_0$. For any $s \ins [t,\infty)$,
 as $\big\{ \ddot{\tau}  \ins [0,s\-t]\big\} \ins \cF^{W,P_0}_{s-t} \= \cF^{{\ddot{\beta}}^t,P_0}_s $,  \eqref{010722_17} yields that
 $ \big\{\wh{\tau} \ins [t,s]\big\} \= \{ \ddot{\tau} (\sW^t) \ins [0,s\-t]\}
\= (\sW^t)^{-1} \big(\{ \ddot{\tau}  \ins [0,s\-t]\}\big) \ins 
\cF^{W^t,P_0}_s  $.
So $\wh{\tau}  $ is a  $[t,\infty]-$valued $\bF^{W^t,P_0}-$stopping time on $\O_0$.
Clearly, $\wh{\tau} \big(\oW(\oo)\big)   \= t \+ \ddot{\tau} \big(\sW^t (\oW(\oo))\big)
\=  t \+ \ddot{\tau} \big(\osW^t (\oo)\big) $ for any $\oo \ins \oO$. It means that (D3') directly implies (D3). \qed

   \fi

 \no {\bf Proof of Theorem \ref{thm_V=oV}:}
  Fix  $(t,\bw,\bx) \ins [0,\infty) \ti \O_0   \ti \OmX  $ and $(y,z) \= \big(\{y_i\}_{i \in \hN}, \{z_i\}_{i \in \hN}\big)
 \ins \Re \ti \Re$.

  \no {\bf 1)} We first show   that $V (t,\bx,y,z )  \ls \oV (t, \bw,\bx, y,z ) $:
 If $\cS_{t,\bx}(y,z) \= \es$, then   $V (t,\bx,y,z) \= -\infty \ls \oV (t, \bw,\bx, y,z )$. 

 So we assume   $\cS_{t,\bx}(y,z)  \nne \es$ and let $\tau \ins \cS_{t,\bx}(y,z) $.
 Define a  process
 $ \cB^{t,\bw}_s  (\o) \df   \bw(s \ld t)  \+ \cB^t_{s \vee t} (\o)     $, $ \fa (s,\o) \ins [0,\infty) \ti \cQ  $
  and define  a mapping $\Psi  \n : \cQ \mto \oO$ by
  $ \Psi (\o) \df \big( \cB^{t,\bw}  (\o),  \cX^{t,\bx}  (\o) ,  \tau(\o) \big) \ins \oO $,  $ \fa \o \ins \cQ $.  
 It holds for any  $(s,\o) \ins [t,\infty) \ti \cQ$ that
  $\oW^t_s(\Psi (\o)) \= \oW_s(\Psi (\o)) \-  \oW_t(\Psi (\o)) \=  \cB^{t,\bw}_s(\o) \- \cB^{t,\bw}_t(\o) \= \cB^t_s(\o)     $.
  Since $   \cX^{t,\bx} \=
  \big\{ \cX^{t,\bx}_s\big\}_{s \in [0,\infty)} $ is an $ \big\{ \cF^{\cB^t,\fp}_{s \vee t} \big\}_{s \in [0,\infty)}   -$adapted continuous process and since $\tau$ is an $ \bF^{\cB^t,\fp}  -$stopping time,
 we can deduce  that    the mapping  $ \Psi $ 
 is $\cF^{\cB^t,\fp}_\infty \big/ \sB(\oO)-$measurable
 and is $\cF^{\cB^t,\fp}_s \big/ \ocF^t_s-$measurable for any $s \ins [t,\infty)$.
 \if{0}


 Let $W^X \= \{W^X_s\}_{s \in [0,\infty)}$ be the canonical process of $\OmX$.
 For any $s \ins [0,\infty)$, $\cE_1 \ins \sB(\hR^d)$ and $\cE_2 \ins \sB(\hR^l)$, we can deduce that
  $ \big\{\o \ins \cQ \n :  \cB^{t,\bw} (\o) \ins W^{-1}_s (\cE_1) \big\}
  \= \big\{\o \ins \cQ \n : W_s \big(\cB^{t,\bw} (\o)\big) \ins \cE_1 \big\}
  \= \big\{\o \ins \cQ \n : \cB^t_{s \vee t} (\o) \ins \cE^{t,\bw}_1 \big\} \ins \cF^{\cB^t}_{s \vee t} \sb \cF^{\cB^t}_\infty$
   with $\cE^{t,\bw}_1 \df \{a \-  \bw(s \ld t)   \n : a \ins \cE_1\} \ins  \sB(\hR^d)$,
  and that  $ \big\{\o \ins \cQ \n :  \cX^{t,\bx} (\o) \ins (W^X_s)^{-1} (\cE_2) \big\}
  \= \big\{\o \ins \cQ \n : W^X_s \big(\cX^{t,\bx} (\o)\big) \ins \cE_2 \big\}
  \= \big\{\o \ins \cQ \n : \cX^{t,\bx}_s (\o) \ins \cE^{t,\bw}_2 \big\} \ins \cF^{\cB^t,\fp}_{s \vee t}  \sb \cF^{\cB^t,\fp}_\infty$.
  So the sigma-field $\big\{A \sb \cQ \n : (\cB^{t,\bw})^{-1}(A) \ins \cF^{\cB^t }_\infty  \big\}$ of $\cQ$ contains all generating sets of
  $ \cF^W_\infty \= \sB(\O_0) $
  and   the sigma-field $\big\{A \sb \cQ \n : (\cX^{t,\bx})^{-1}(A) \ins \cF^{\cB^t,\fp}_\infty  \big\}$ of $\cQ$ contains all generating sets of
  $ \cF^{W^X}_\infty \= \sB(\OmX) $. It follows that
  the mapping $\cB^{t,\bw} \n : \cQ \mto \O_0 $ is $\cF^{\cB^t}_\infty \big/ \sB(\O_0)-$measurable
  and  the mapping $\cX^{t,\bx} \n : \cQ \mto \OmX $ is $\cF^{\cB^t,\fp}_\infty \big/ \sB(\OmX)-$measurable.

     Since the  $\bF^{\cB^t,\fp}-$stopping time  $\tau$ satisfies that   $    \big\{\o  \ins \cQ \n : \tau(\o) \ins [t,s] \big\} \ins \cF^{\cB^t,\fp}_s \sb \cF^{\cB^t,\fp}_\infty$ for any $s \ins [t,\infty]$,
      the sigma field of $[t,\infty]$,
   $\L_\tau \df \big\{\cE \sb [t,\infty] \n :   \tau^{-1}(\cE) \ins \cF^{\cB^t,\fp}_\infty \big\}$ contains all closed intervals $[t,s]$, $s \ins [t,\infty]$, which generates $\sB[t,\infty]$. It follows that
   $\sB[t,\infty] \sb \L_\tau $ or $  \tau^{-1}(\cE) \ins \cF^{\cB^t,\fp}_\infty $, $ \fa \cE \ins \sB[t,\infty] $.
   So $\tau \n: \cQ \mto [t,\infty]$ is $\cF^{\cB^t,\fp}_\infty   -$measurable,
   the mapping $\Psi \n : \cQ \mto \oO$ is eventually $\cF^{\cB^t,\fp}_\infty \big/ \sB(\oO)-$measurable.

  Let $s \ins [t,\infty)$. For any $(r_1,\cE_1) \ins [t,s] \ti \sB(\hR^d)$ and $(r_2,\cE_2) \ins [0,s] \ti \sB(\hR^l)$,
  one has  $ \big\{\o \ins \cQ \n :  \Psi (\o) \ins (\oW^t_{r_1})^{-1} (\cE_1) \big\}
  \= \big\{\o \ins \cQ \n :  \oW^t_{r_1}\big(\Psi (\o)\big) \ins  \cE_1  \big\}
  \= \big\{\o \ins \cQ \n :  \cB^t_{r_1} (\o)  \ins  \cE_1  \big\} \ins \cF^{\cB^t}_s $
  and  $ \big\{\o \ins \cQ \n :  \Psi (\o) \ins  \oX_{r_2}^{-1} (\cE_2) \big\}
  \= \big\{\o \ins \cQ \n :  \oX_{r_2}\big(\Psi (\o)\big) \ins  \cE_2  \big\}
  \= \big\{\o \ins \cQ \n :  \cX^{t,\bx}_{r_2} (\o)  \ins  \cE_2  \big\} \ins \cF^{\cB^t,\fp}_s $.
   So the sigma-field $\big\{A \sb \cQ \n : \Psi^{-1}(A) \ins \cF^{\cB^t,\fp}_s  \big\}$ of $\cQ$ contains all generating sets of
   $\ocF^t_s$ and thus includes $\ocF^t_s$.
   Then   the mapping $\Psi  $ is also $\cF^{\cB^t,\fp}_s \big/ \ocF^t_s -$measurable.

  \fi
 Let $\oP_\Psi \ins \fP\big(\oO\big)  $ be the probability measure induced by $\Psi$, i.e.,
 $  \oP_\Psi \big(\oA\big) \df \fp \big( \Psi^{-1}\big(\oA\big) \big) $, $ \fa \oA \ins \sB(\oO)  $.

  Fix  $ (\vf,n)  \ins   \fC(\hR^{d+l}) \ti \hN $. We define an $\bF^{\cB^t,\fp}-$adapted continuous process
  $   
    \cM^{t,\bx}_s(\vf)   \df   \vf \big(\cB^t_s  , \cX^{t,\bx}_s \big)
    \- \n \int_t^s  \n  \ol{b}  \big( r, \cX^{t,\bx}_{r \land \cd} \big) \n \cd \n D \vf \big( \cB^t_r  , \cX^{t,\bx}_r \big) dr
    \-   \frac12 \n \int_t^s  \n  \ol{\si} \, \ol{\si}^T  \big( r,  \cX^{t,\bx}_{r \land \cd}  \big) \n : \n D^2 \vf  ( \cB^t_r, \cX^{t,\bx}_r  )   dr  $, $  \fa s \ins [t,\infty)$
   and define an  $\bF^{\cB^t,\fp}-$stopping time $\btau^{t,\bx}_n  \df  \inf\big\{s \ins [t,\infty) \n : \big| (\cB^t_s ,\cX^{t,\bx}_s) \big|    \gs n   \big\} \ld (t\+n) $.
   Applying Proposition \ref{prop_MPF1} with $(\O,\cF,P,B,X) \= (\cQ,\cF,\fp,\cB,\cX^{t,\bx}) $ yields that
   $ \big\{\cM^{t,\bx}_{s \land \bbtau^{t,\bx}_n}(\vf)\big\}_{s \in [t,\infty)} $
   is a bounded $ \big( \bF^{\cB^t,\fp} , \fp \big) -$martingale.  
 \if{0}

 Actually, $ \big\{\cM^{t,\bx}_{s \land \bbtau^{t,\bx}_n}(\vf)\big\}_{s \in [t,\infty)} $
   is a  martingale with respect to the filtration $\bF^{\cB^t,\cX^{t,\bx}} \=
   \Big\{ \cF^{\cB^t,\cX^{t,\bx}}_s \df \si \big(\cB^t_r;r \ins [t,s]\big) \ve \si\big( \cX^{t,\bx}_r;r \ins [0,s]\big)\Big\}_{s \in [t,\infty)}$ and is thus a martingale with respect to the filtration $\bF^{\cB^t,\cX^{t,\bx},\fp} \=
   \Big\{ \cF^{\cB^t,\cX^{t,\bx},\fp}_s \df \si \big(\cF^{\cB^t,\cX^{t,\bx}}_s \cp \sN_\fp(\cF^{\cB^t,\cX^{t,\bx}}_\infty)\big)\Big\}_{s \in [t,\infty)}$.  Since $\cF^{\cB^t,\fp}_s \=  \si \big(\cF^{\cB^t}_s   \cp \sN_\fp(\cF^{\cB^t}_\infty)\big) \sb \cF^{\cB^t,\cX^{t,\bx},\fp}_s $ for any $s \ins [t,\infty)$, we see that $ \big\{\cM^{t,\bx}_{s \land \bbtau^{t,\bx}_n}(\vf)\big\}_{s \in [t,\infty)} $
   is also an $ \bF^{\cB^t,\fp} -$martingale.

 \fi

   Since $\oP_\Psi \big\{   \oX_s   \= \bx(s) , \fa s \ins [0,t]  \big\} \= \fp \big\{  \oX_s (\Psi ) \= \bx(s), \fa s \ins [0,t] \big\}
   \= \fp \big\{  \cX^{t,\bx}_s    \= \bx(s), \fa s \ins [0,t] \big\} \=1$, using Proposition \ref{prop_MPF1}
   with $(\O,\cF,P,B,X) \= \big(\oO,\sB(\oO),\oP_\Psi,\oW,\oX \big) $
    shows that  $\big\{\oM^t_{s \land \otau^t_n } (\vf) \big\}_{s \in  [t,\infty) }  $ is a bounded $\obF^t-$adapted continuous process under $ \oP_\Psi $.
   Given $\o \ins \cQ$, since $\oW^t_s \big(\Psi (\o)\big) \=  \cB^t_s(\o)     $, $\fa s \ins [t,\infty)$, we see that
  $ \big(\oM^t_s(\vf)\big)\big(\Psi(\o)\big) \=  \big(\cM^{t,\bx}_s(\vf) \big) (\o)$, $ \fa s \ins [t,\infty) $ and
  $\otau^t_n \big(\Psi(\o)\big) \= \btau^{t,\bx}_n (\o)$.  Then
  \bea
   \big(\oM^t_{s \land \otau^t_n }(\vf) \big) \big(\Psi(\o)\big)
  & \tn \= & \tn  \big(\oM^t (\vf) \big) \big( s \ld \otau^t_n (\Psi(\o)),  \Psi(\o)  \big)
  \= \big(\oM^t (\vf) \big) \big( s \ld \btau^{t,\bx}_n  (\o) ,  \Psi(\o)  \big) \nonumber \\
  & \tn  \= & \tn  \big(\cM^{t,\bx} (\vf) \big) \big( s \ld \btau^{t,\bx}_n  (\o) ,  \o   \big)
  \= \big(\cM^{t,\bx}_{s \land \bbtau^{t,\bx}_n} (\vf) \big)   (\o) , \q   \fa (s,\o) \ins [t,\infty) \ti \cQ . \label{122621_17}
  \eea

  Let $t_1,t_2 \ins [t,\infty)$ with $t_1 \< t_2$ and let $\oA \ins \ocF^t_{t_1} $.
    As  $\Psi^{-1}(\oA) \ins \cF^{\cB^t,\fp}_{t_1}$, the $ \big(\bF^{\cB^t,\fp}, \fp\big) -$martingality of
    $ \big\{\cM^{t,\bx}_{s \land \bbtau^{t,\bx}_n}(\vf)\big\}_{s \in [t,\infty)} $ and \eqref{122621_17}    imply that
 \beas
   E_{\oP_\Psi} \Big[ \big(\oM^t_{t_2 \land \otau^t_n }(\vf) \-  \oM^t_{t_1 \land \otau^t_n }(\vf)  \big) \b1_\oA    \Big]
   & \tn \=    & \tn    E_\fp \Big[ \Big( \big(\oM^t_{t_2 \land \otau^t_n }(\vf) \big) (\Psi ) \-  \big(\oM^t_{t_1 \land \otau^t_n }(\vf) \big) (\Psi ) \Big) \b1_{\Psi^{-1}(\oA)}  \Big] \\
      & \tn \=    & \tn   E_\fp \Big[ \big(  \cM^{t,\bx}_{t_2 \land \bbtau^t_n }(\vf)      \-  \cM^{t,\bx}_{t_1 \land \bbtau^t_n }(\vf)   \big)  \b1_{\Psi^{-1}(\oA)} \Big]   \= 0  .
      \eeas
  So  $ \big\{ \oM^t_{s \land \otau^t_n } (\vf) \big\}_{s \in [t,\infty)} $  is a bounded $ \big(\obF^t,\oP_\Psi\big)   -$martingale.
 By Remark \ref{rem_ocP}, $\oP_\Psi$ satisfies (D1) and (D2) of Definition  \ref{def_ocP}. 

 Since   $   W^t_s( \cB^{t,\bw} (\o))   \n \= \n   \cB^{t,\bw}_s(\o) \n \- \n  \cB^{t,\bw}_t(\o)    \n \= \n    \cB^t_s(\o)    $ for any $(s,\o) \ins [t,\infty) \ti \O $,
 taking $(\O,\cF,P,B,\Phi) \= (\cQ,\cF,\fp,\cB,\cB^{t,\bw}) $ in Lemma \ref{lem_M31_01} (2) shows that
 $\fp \big\{ \tau \= \wh{\tau}(\cB^{t,\bw}) \big\} \= 1$ for some  $[t,\infty]-$valued $ \bF^{W^t,P_0}-$stopping time $\wh{\tau}$ on $\O_0$,
 it follows that $ \oP_\Psi \big\{\oT \= \wh{\tau}(\oW) \big\} \= \fp \big\{\oT(\Psi) \= \wh{\tau}(\oW(\Psi))\big\}
 \= \fp \big\{\tau \= \wh{\tau}(\cB^{t,\bw})\big\} \= 1 $. 
 As  $ \oW_{\n s} \big(\Psi (\o)\big) \= \cB^{t,\bw}_s(\o) \= \bw(s) $, $ \fa (s,\o) \ins [0,t] \ti \cQ $,
 it is clear that
 $ \oP_\Psi \big\{\oW_s \= \bw(s) , \fa s \ins [0,t] \big\} \= \fp \big\{ \oW_s (\Psi) \= \bw(s) , \fa s \ins [0,t]  \big\}
 \= 1 $. Thus  $\oP_\Psi \ins \ocP_{t,\bw,\bx}$.
  For any $i \ins \hN$,
  \bea
 E_{\oP_\Psi}   \bigg[ \int_t^\oT g_i \big(r,\oX_{r \land \cd} \big) dr    \bigg]
   \=    E_\fp \bigg[    \int_t^{\oT (\Psi   )  } g_i \big(r,\oX_{r \land \cd} (\Psi   ) \big) dr  \bigg]
  \=  E_\fp \Big[ \int_t^\tau g_i ( r, \cX^{t,\bx}_{r \land \cd} ) dr  \Big] \ls y_i    \label{Feb12_01}
 \eea
 and 
 $ E_{\oP_\Psi}   \big[ \int_t^\oT h_i \big(r,\oX_{r \land \cd} \big) dr    \big]
   \=  E_\fp \big[ \int_t^\tau h_i ( r, \cX^{t,\bx}_{r \land \cd} ) dr  \big] \= z_i $,
 which means that  $ \oP_\Psi \ins \ocP_{t,\bw,\bx}(y,z)   $.
 Then an analogy to  \eqref{Feb12_01}   renders that
  $  E_\fp \big[ \int_t^\tau f ( r, \cX^{t,\bx}_{r \land \cd} ) dr \+ \b1_{\{\tau < \infty\}} \pi\big(\tau, \cX^{t,\bx}_{\tau \land \cd} \big) \big]
 \=  E_{\oP_\Psi}  \Big[ \int_t^\oT f \big(r,\oX_{r \land \cd} \big) dr \+ \b1_{\{\oT < \infty\}} \pi \big( \oT, \oX_{\oT \land \cd} \big) \Big] \ls \oV (t,\bw,\bx,y,z) $.
 Taking supremum over $  \tau \ins \cS_{t,\bx}(y,z)   $  yields  that $V(t,\bx,y,z) \ls \oV (t,\bw,\bx,y,z)  $.

  \no {\bf 2)} As $\ocP_{t,\bw,\bx}(y,z) \sb \ocP_{t,\bx}(y,z)$, we automatically have  $\oV (t, \bw,\bx, y,z ) \ls \oV (t,\bx,y,z )  $.
  It remains to demonstrate  that $\oV (t,\bx,y,z ) \ls V (t,\bx,y,z )   $.
  If $\ocP_{t,\bx}(y,z) \= \es$, then $\oV (t,\bx,y,z ) \= -\infty \ls V (t,\bx,y,z )   $.

   Assume   $\ocP_{t,\bx}(y,z) \nne \es$ and let $\oP \ins \ocP_{t,\bx}(y,z) $.
   Given $ (\vf,n)  \ins  \fC(\hR^{d+l}) \ti \hN   $,
$ \oM^{t,\bx}_s(\vf)    \df    \vf  ( \oW^t_{\n s} ,  \osX^{t,\bx}_{\n s} )
 \- \n \int_t^s \ol{b}  \big( r, \osX^{t,\bx}_{\n r \land \cd} \big)  \n \cd \n D \vf  ( \oW^t_{\n r} ,  \osX^{t,\bx}_{\n r}  ) dr
    \-   \frac12 \int_t^s   \ol{\si} \, \ol{\si}^T  \big( r, \osX^{t,\bx}_{\n r \land \cd} \big) \n : \n D^2 \vf ( \oW^t_{\n r} ,  \osX^{t,\bx}_{\n r} )   dr  $, $s \ins [t,\infty)$   is  an $\bF^{\oW^t,\oP}-$adapted continuous process
 and   $\otau^{t,\bx}_n  \df  \inf\big\{s \ins [t,\infty) \n :  \big|(\oW^t_s,\osX^{t,\bx}_s ) \big|   \gs n  \big\} \ld (t\+n) $
   is an $\bF^{\oW^t,\oP} -$stopping time. Since $\oW^t$ is a Brownian motion under $\oP$ by (D1) of Definition \ref{def_ocP},
  applying Proposition \ref{prop_MPF1}    with $(\O,\cF,P,B,X) \= \big(\oO,\sB(\oO),\oP,\oW,\osX^{t,\bx} \big) $ shows that
 \bea \label{123021_11}
 \big\{\oM^{t,\bx}_{s \land \otau^{t,\bx}_n} (\vf) \big\} \hb{ is a bounded $\bF^{\oW^t,\oP}-$martingale.}
 \eea

 Let $(\fx_o,\ft_o)$ be an arbitrary pair in $\OmX \ti [t,\infty]$ and define  a mapping $\Psi_o \n : \cQ \mto \oO$ by
 $ \Psi_o (\o) \df \big( \cB (\o),  \fx_o, \ft_o \big)  \ins \oO $, $ \fa \o \ins \cQ $.
 (Actually, we are indifferent to the second and third components  of $\Psi_o(\o)$.)
 Since   $   \oW^t_s  ( \Psi_o (\o) )    \=    \oW_s  (\Psi_o (\o)) \- \oW_t  (\Psi_o (\o)) \=   \cB^t_s (\o)    $ for any $(s,\o) \ins [t,\infty) \ti \cQ$ and since $\oW^t$ is a Brownian motion under $\oP$ by (D1) of Definition \ref{def_ocP},
    applying Lemma \ref{lem_122921_11} with $t_0 \= t$,  $(\O_1, \cF_1, P_1,B^1)   \= \big(\cQ, \cF, \fp , \cB \big) $, $(\O_2, \cF_2, P_2,B^2) \= \big(\oO , \sB(\oO ) , \oP  , \oW\big) $ and $\Phi \= \Psi_o$ yields that
   \bea
   \Psi^{-1}_o \big( \cF^{\oW^t }_s \big) \= \cF^{\cB^t }_s , ~
    \Psi^{-1}_o \big( \cF^{\oW^t,\oP}_s \big) \sb  \cF^{\cB^t,\fp}_s  , ~ \fa s \ins [t,\infty] \aand
     \big(\fp \nci \Psi_o^{-1}\big) (\oA) \= \oP(\oA)  , ~ \; \fa \oA \ins \cF^{\oW^t, \oP}_\infty   .  \label{Oct01_07b}
 \eea
 Then $ \sX^{t,\bx}_s (\o) \df  \osX^{t,\bx}_{\n s} ( \Psi_o(\oo))$, $s \ins [0,\infty)$ defines  an $\big\{\cF^{\cB^t,\fp}_{s \vee t}\big\}_{s \in [0,\infty)}-$adapted continuous process.

    Let $ (\vf,n)  \ins  \fC(\hR^{d+l}) \ti \hN   $. We define an $\bF^{\cB^t,\fp}-$adapted continuous process
$ \sM^{t,\bx}_s(\vf)    \df    \vf  ( \cB^t_s ,  \sX^{t,\bx}_{\n s} )
 \- \n \int_t^s \ol{b}  \big( r, \sX^{t,\bx}_{\n r \land \cd} \big)  \n \cd \n D \vf  ( \cB^t_r ,  \sX^{t,\bx}_{\n r}  ) dr
    \-   \frac12 \int_t^s   \ol{\si} \, \ol{\si}^T  \big( r, \sX^{t,\bx}_{\n r \land \cd} \big) \n : \n D^2 \vf ( \cB^t_r ,  \sX^{t,\bx}_{\n r} )   dr  $, $ \fa s \ins [t,\infty)$  and define an  $\bF^{\cB^t,\fp}-$stopping time $\z^{t,\bx}_n  \df  \inf\big\{s \ins [t,\infty) \n : \big|(\cB^t_s , \sX^{t,\bx}_s) \big|   \gs n   \big\} \ld (t\+n) $.

   Applying Proposition \ref{prop_MPF1}    with $(\O,\cF,P,B,X) \= \big(\cQ,\cF,\fp,\cB,\sX^{t,\bx} \big) $ and using an analogy to \eqref{122621_17} renders that
   $ \big\{ \sM^{t,\bx}_{s \land \z^{t,\bx}_n} (\vf) \big\}_{s \in [t,\infty)}  $ is a bounded $\bF^{\cB^t,\fp}-$adapted continuous process under $ \fp  $ satisfying
   \if{0}
 It holds  for any $\o \ins \cQ$ that
  \beas
  \big(\oM^{t,\bx}_s(\vf)\big) \big(\Psi_o(\o)\big) \= \big( \sM^{t,\bx}_s(\vf) \big) (\o) , \q \fa s \ins [t,\infty) \aand
  \otau^{t,\bx}_n \big(\Psi_o(\o)\big) \= \z^{t,\bx}_n(\o) .
  \eeas
   \fi
  \bea \label{123021_14}
 \big( \oM^{t,\bx}_{s \land \otau^{t,\bx}_n} (\vf) \big) \big(\Psi_o(\o)\big) \= \big( \sM^{t,\bx}_{s \land \z^{t,\bx}_n} (\vf) \big) (\o) ,
 \q \fa (s,\o) \ins [t,\infty) \ti \cQ.
 \eea

 Let $t_1,t_2 \ins [t,\infty)$ with $t_1 \< t_2$ and let $A \ins \cF^{\cB^t }_{t_1}$.
 Since $\Psi_o^{-1}\big(\oA\big) \= A   $ for some  $\oA \ins \cF^{\oW^t }_{t_1}$ by \eqref{Oct01_07b},
 we can derive  from \eqref{123021_11}, \eqref{Oct01_07b} and \eqref{123021_14} that
 \beas
  0  & \tn \=  & \tn    E_\oP  \Big[ \big(\oM^{t,\bx}_{t_2 \land \otau^{t,\bx}_n }(\vf) \-  \oM^{t,\bx}_{t_1 \land \otau^{t,\bx}_n }(\vf)  \big) \b1_\oA    \Big]
   \= E_\fp  \Big[ \Big( \big(\oM^{t,\bx}_{t_2 \land \otau^{t,\bx}_n }(\vf) \big) (\Psi_o) \- \big( \oM^{t,\bx}_{t_1 \land \otau^{t,\bx}_n }(\vf)  \big)(\Psi_o) \Big) \b1_{\Psi_o^{-1}(\oA)}   \Big]  \\
    & \tn  \=   & \tn   E_\fp  \Big[ \big(\sM^{t,\bx}_{t_2 \land \z^{t,\bx}_n }(\vf) \-  \sM^{t,\bx}_{t_1 \land \z^{t,\bx}_n }(\vf)  \big) \b1_A    \Big]  ,
    \eeas
 which implies that $ \big\{\sM^{t,\bx}_{s \land \z^{t,\bx}_n }(\vf)\big\}_{s \in [t,\infty)}$ is a bounded $\bF^{\cB^t,\fp}-$martingale.
 Then an application of  Proposition \ref{prop_MPF1} with $(\O,\cF,P,B,X) \= \big(\cQ,\cF,\fp,\cB,\sX^{t,\bx}\big) $ shows that
   \if{0}

 Given $s \ins [t,\infty)$, we set   $\fF^{t,\bx}_s \df \cF^{\cB^t}_s \ve \cF^{\sX^{t,\bx}}_s \= \si\big(\cB^t_r; r \ins [t,s]\big) \ve \si (\sX^{t,\bx}_r; r \ins [0,s] )$.  Clearly, $\fF^{t,\bx}_s \sb \cF^{\cB^t,\fp}_s$.
 Actually,  
 $ \big\{ \sM^{t,\bx}_{s \land \z^{t,\bx}_n} (\vf) \big\}_{s \in [t,\infty)}  $ is a bounded $\{\fF^{t,\bx}_s\}_{s \in [t,\infty)}-$adapted continuous process under $ \fp  $ and is thus a bounded $\bF^{\cB^t,\fp}-$adapted continuous process under $ \fp  $.

  We see from this equality that the Lambda-system $\L \df \Big\{A \ins \cF^{\cB^t, \fp }_\infty  \sb \cQ \n :  E_\fp  \big[ \big(\sM^{t,\bx}_{t_2 \land \z^{t,\bx}_n }(\vf) \-  \sM^{t,\bx}_{t_1 \land \z^{t,\bx}_n }(\vf)  \big) \b1_A    \big] \= 0 \Big\} $ contains the Pi-system
 $ \cF^{\oW^t }_{t_1} \Cp \sN_\fp \big(\cF^{\oW^t }_\infty \big)$. Then Dynkin's Pi-Lambda Theorem implies that
 $ \cF^{\oW^t }_{t_1,\fp} \= \si \Big( \cF^{\oW^t }_{t_1} \Cp \sN_\fp \big(\cF^{\oW^t }_\infty \big) \Big) \sb \L $, i.e.,
 $  E_\fp  \Big[ \big(\sM^{t,\bx}_{t_2 \land \z^{t,\bx}_n }(\vf) \-  \sM^{t,\bx}_{t_1 \land \z^{t,\bx}_n }(\vf)  \big) \b1_A    \Big] \= 0 $ for any $A \ins \cF^{\oW^t }_{t_1,\fp}$.
 Hence, $ \big\{\sM^{t,\bx}_{s \land \z^{t,\bx}_n }(\vf)\big\}_{s \in [t,\infty)}$ is a bounded $\bF^{\cB^t,\fp}-$martingale.
 In particular, $ \big\{\sM^{t,\bx}_{s \land \z^{t,\bx}_n }(\vf)\big\}_{s \in [t,\infty)}$ is a bounded $\{\fF^{t,\bx}_s\}_{s \in [t,\infty)}-$martingale.

   \fi
 \bea \label{012521_23}
 \fp\{ \sX^{t,\bx}_s \= \cX^{t,\bx}_s, ~ \fa s \ins [0,\infty) \} \= 1 .
 \eea

    By (D3) of Definition \ref{def_ocP},
 there exists   a $[t,\infty]-$valued $\bF^{W^t,P_0}-$stopping time $\wh{\ga}$  on $\O_0$
 such that $   \oP \big\{  \oT  \=    \wh{\ga} (\oW  )   \big\} \= 1 $.
 Lemma \ref{lem_M31_01} (1) renders that  $ \ga \df \wh{\ga}(\cB ) $ is an $\bF^{\cB^t,\fp}-$stopping time on $\cQ$
 while $   \wh{\ga}(\oW ) $ is an $\bF^{\oW^t,\oP}-$stopping time on $\oO$.
 For any $ i \ins \hN $,   we can deduce from (D2) of Definition \ref{def_ocP},  \eqref{Oct01_07b}  and  \eqref{012521_23}     that
  \bea
    y_i & \tn \gs  & \tn  E_\oP \bigg[ \int_t^\oT g_i \big( r,\oX_{r \land \cd} \big) dr\bigg]
   \=  E_\oP \bigg[ \int_t^{ \wh{\ga}(\oW )} g_i \big( r,\osX^{t,\bx}_{\n r \land \cd} \big) dr\bigg]
     \=  E_\fp \bigg[ \int_t^{ \wh{\ga}(\oW (\Psi_o))} g_i \big( r,\osX^{t,\bx}_{\n r \land \cd}(\Psi_o) \big) dr\bigg] \nonumber \\
        & \tn  \=  & \tn  E_\fp \bigg[ \int_t^{ \wh{\ga}(\cB )} g_i ( r, \sX^{t,\bx}_{r \land \cd}  ) dr\bigg]
         \=  E_\fp \bigg[ \int_t^{ \ga } g_i ( r, \cX^{t,\bx}_{r \land \cd}  ) dr\bigg]  , \label{012521_15}
 \eea
 and similarly that  $  E_\fp \big[ \int_t^{  \ga } h_i ( r, \cX^{t,\bx}_{r \land \cd}  ) dr\big]
 \=  E_\oP \big[ \int_t^\oT h_i ( r,\oX_{r \land \cd} ) dr\big] \= z_i$.
 So    $ \ga   \ins \cS_{t,\bx} (y,z)$. Analogous  to \eqref{012521_15},
 \beas
  E_\oP   \bigg[ \int_t^\oT f \big(r,\oX_{r \land \cd} \big) dr \+ \b1_{\{\oT < \infty\}}
 \pi \big( \oT,\oX_{\oT \land \cd} \big)  \bigg]
 \=  E_\fp \Big[ \int_t^{ \ga} f \big( r,\cX^{t,\bx}_{r \land \cd}  \big) dr
 \+ \b1_{\{\ga < \infty\}} \pi \big(  \ga,\cX^{t,\bx}_{ \ga  \land \cd}  \big)  \Big] \ls V(t,\bx,y,z).
 \eeas
 Taking supremum over $ \oP \ins  \ocP_{t,\bx} (y,z) $   yields that   $ \oV(t,\bx,y,z) \ls V(t,\bx,y,z)$. \qed

  \no {\bf Proof of Lemma \ref{lem_082020_15}: 1)}
  \if{0}

   Let $\tau_1,\tau_2 \ins \fS$ with $ \Rho{\fS} (\tau_1,\tau_2) \= 0 $,
  then $ \big|\arctan(\tau_1(\o_0)) \- \arctan(\tau_2(\o_0))\big| \= \Rho{+}  \big( \tau_1(\o_0) , \tau_2(\o_0) \big) \= 0$ and thus $ \tau_1(\o_0) \= \tau_2(\o_0) $ for $P_0-$a.s. $\o_0 \ins \O_0$, i.e.,  $\tau_1 \= \tau_2 $ in $\fS$.
 Also, it holds for any   $ \tau_1,\tau_2,\tau_3 \ins \fS$ that
  \beas
  \Rho{\fS} (\tau_1,\tau_3) & \tn \=  & \tn  E_{P_0} \big[ \Rho{+}  ( \tau_1 , \tau_3 )  \big]
  \= E_{P_0} \big[ \big|\arctan(\tau_1 ) \- \arctan(\tau_3 )\big| \big]
  \ls  E_{P_0} \big[ \big|\arctan(\tau_1 ) \- \arctan(\tau_2 )\big|
  \+    \big|\arctan(\tau_2 ) \- \arctan(\tau_3 )\big| \big] \\
  & \tn \=  & \tn   E_{P_0} \big[ \Rho{+}  ( \tau_1 , \tau_2 ) \big]  \+  E_{P_0} \big[ \Rho{+}  ( \tau_2 , \tau_3 )  \big]
  \= \Rho{\fS} (\tau_1,\tau_2) \+ \Rho{\fS} (\tau_2,\tau_3) .
  \eeas
 So $\Rho{\fS}$ is a metric on $\fS$.

  \fi
 We first show that the   metric  space $\big(\fS, \Rho{\fS}\big)$ is   complete.
  Let $\{\tau_n\}_{n \in \hN}$ be a Cauchy sequence in $\big(\fS, \Rho{\fS}\big)$
      such that
  $ \Sup{k \in \hN} \, \Rho{\fS}  ( \tau_n, \tau_{n+k}  ) \< 2^{-n} $ for any $  n \ins \hN$.

 For any $n  \ins \hN$, 
 the monotone convergence theorem implies that
  \beas
     E_{P_0} \Big[ \, \underset{k \in \hN}{\sup} \, \big|\arctan(\tau_n)     \-  \arctan(\tau_{n+k})    \big| \Big]
    & \tn \ls  & \tn     E_{P_0} \Big[ \, \sum_{k \in \hN}    \big|\arctan(\tau_{n+k-1})     \-  \arctan(\tau_{n+k})    \big| \Big]  \\
   & \tn  \=  & \tn    \sum_{k \in \hN} \Rho{\fS}  \big( \tau_{n+k-1}, \tau_{n+k} \big)
      \ls \sum_{k \in \hN} 2^{1-n-k} \= 2^{1-n }   .
      \eeas
  So $\lmt{n \to \infty}   E_{P_0} \Big[ \, \underset{k \in \hN}{\sup} \, \big|\arctan(\tau_n)     \-  \arctan(\tau_{n+k})    \big| \Big] \= 0$.
 Then one can extract  a subsequence $\big\{ \tau_{n_j} \big\}_{j \in \hN}$ of $ \{ \tau_n  \}_{n \in \hN}$  such that
 $ \lmt{j \to \infty}   \Big( \underset{k \in \hN}{\sup} \, \big|\arctan(\tau_{n_j}(\o_0))     \-  \arctan(\tau_{ k + n_j }(\o_0))    \big|   \Big)   \= 0 $
   for all $ \o_0  \ins \O_0  $ except on a $  P_0-$null set $\cN  $.
   Given $ \o_0  \ins \cN^c  $, we see that
   $   \lmt{j \to \infty}    \Big( \underset{\ell \in \hN}{\sup} \,  \big| \arctan(\tau_{n_j}(\o_0) ) \-  \arctan( \tau_{n_{\overset{}{j+\ell}}}(\o_0) ) \big| \Big)  \= 0$, i.e.,
     $ \big\{ \arctan \big( \tau_{n_j}(\o_0) \big) \big\}_{j \in \hN} $ is a Cauchy sequence
   in $[0,\pi/2]$. Let $ \xi_* (\o_0) $ be the limit of $ \big\{ \arctan \big( \tau_{n_j}(\o_0) \big) \big\}_{j \in \hN} $ in $[0,\pi/2]$.

 As $\bF^{W,P_0}$ is a right-continuous complete filtration,  Lemma 1.2.11 of \cite{Kara_Shr_BMSC} implies that
 $ \tau_* \df \linf{j \to \infty} \tau_{n_j}$  is an  $\bF^{W,P_0}-$stopping time on $\O_0$ satisfying
 \beas
 \arctan \big(\tau_*(\o_0)\big) & \tn \=  & \tn 
 \arctan  \Big( \Sup{j \in \hN} \, \Inf{\ell \ge j} \tau_{n_\ell} (\o_0)\Big)
 \= \Sup{j \in \hN} \arctan  \Big(    \Inf{\ell \ge j} \tau_{n_\ell} (\o_0)\Big)
 \= \Sup{j \in \hN} \,   \Inf{\ell \ge j} \arctan  \Big(   \tau_{n_\ell} (\o_0)\Big) \\
  & \tn \= & \tn  \linf{j \to \infty} \arctan  \big(  \tau_{n_j} (\o_0) \big)
 \= \lmt{j \to \infty} \arctan  \big(  \tau_{n_j} (\o_0) \big)
 \= \xi_* (\o_0) , \q \fa \o_0 \ins   \cN^c .
 \eeas
   Applying the bounded convergence theorem renders that
  $ \lmt{j \to \infty} \Rho{\fS} ( \tau_{n_j} , \tau_* ) \=  \lmt{j \to \infty}  E_{P_0} \big[ \big|\arctan( \tau_{n_j} ) \- \arctan(\tau_* )\big| \big] \= 0 $.

   We next  let   $\{\tau_n\}_{n \in \hN}$ be a general Cauchy sequence in $\big(\fS, \Rho{\fS}\big)$.
   For any $j \ins \hN$, there exists $n_j \ins \hN$ such that
     $ \Sup{k \in \hN} \, \Rho{\fS}  ( \tau_{n_j}, \tau_{k+n_j}  ) \< 2^{-j} $.
     In particular, the subsequence $\big\{\wt{\tau}_j \df \tau_{n_j}\big\}_{j \in  \hN}$ of $\{\tau_n\}_{n \in \hN}$ satisfies that
     $ \Sup{\ell \in \hN} \, \Rho{\fS}  \big( \wt{\tau}_j, \wt{\tau}_{j+\ell}  \big) \< 2^{-j} $ for any $  j \ins \hN$
     and thus has a limit $\wt{\tau}_*$ in $\big(\fS, \Rho{\fS}\big)$ by the above argument.
     Let $\e \ins (0,1)$. There exists a $\fk \ins \hN$ with $\fk \gs 1 \- \log_2 \e$
     such that $ \Rho{\fS} \big( \wt{\tau}_\fk , \wt{\tau}_* \big) \ls \e/2 $.
     Then it holds for any $j \gs n_\fk$ that
     $ \Rho{\fS} \big(  \tau_j , \wt{\tau}_* \big)
     \ls  \Rho{\fS} \big(  \tau_j , \wt{\tau}_\fk  \big) \+  \Rho{\fS} \big( \wt{\tau}_\fk  , \wt{\tau}_* \big)
     \ls \Sup{\ell \in \hN} \, \Rho{\fS}  \big(  \tau_{n_\fk},  \tau_{n_\fk+\ell}  \big)  \+ \e/2 \<
      2^{-\fk} \+ \e/2 \ls  \e $. So $ \lmt{j \to \infty} \Rho{\fS} \big( \tau_j , \wt{\tau}_* \big) \= 0 $,
     which shows the completeness of  $\big(\fS, \Rho{\fS}\big)$.

  \no {\bf 2)}  We need some technical preparation for constructing   a countable dense  subset   of $\fS$.

   Fix $s \ins [0,\infty)$.      Given $ \d \ins \hQ_+  $,   set $ O^s_\d (\o_0) \df \Big\{\o'_0  \ins  \O_0 \n :  \Sup{r \in [0,s]} \big|\o'_0 (r)\- \o_0 (r) \big|  \<   \d  \Big\} $. Since $\O_0$ is a continuous-path space, we can deduce that
  \beas 
     O^s_\d (\o_0) \=  \underset{n \in \hN }{\cup} \underset{r \in (0,s) \cap \hQ  }{\cap}
     \big\{ \o'_0 \ins \O_0 :  | \o'_0 (r) \- \o_0 (r) | \ls \d \- \d/n \big\}
   \= \underset{n \in \hN }{\cup} \underset{r \in (0,s) \cap \hQ }{\cap}
   \big\{ \o'_0 \ins \O_0 : W_r(\o'_0) \ins \ol{O}_{\d-\d/n}  \big( \o_0(r) \big)  \big\} \ins \cF^W_s .
    \eeas

   Let $\fT_s(\O_0)$  collect the empty set $\es$ and all subsets $\cO$ of $\O_0$ such that for any
  $\o_0 \ins \cO $ there exists some $ \d \ins (0,1)$ satisfying $ O^s_\d (\o_0) \sb \cO $.
  Obviously, $\fT_s (\O_0)$ forms a topology on $\O_0$.

    We claim that for any   $A \ins \cF^W_s $ and $\e \ins (0,1)$,
  \bea \label{010122_11}
  \hb{there exist
     $\cO_1, \cO_2 \ins \fT_s(\O_0) $  such that $  \cO^c_1 \sb A \sb \cO_2$ and
     $P_0(A \Cp \cO_1  ) \ve P_0(A^c \Cp \cO_2 ) \< \e $.}
  \eea
  To see this, we define $\L_s \df \{A \ins \sB(\O_0) \n : $
   for any $\e \ins (0,1)$  there exist
     $\cO_1, \cO_2$ in $ \fT_s(\O_0) $  such that $  \cO^c_1 \sb A \sb \cO_2$ and
     $P_0(A \Cp \cO_1  ) \ve P_0(A^c \Cp \cO_2 ) \< \e $\}.
     Clearly, $\es, \O_0 \ins \L_s$ as they both belong to $ \fT_s(\O_0) $.   It is also easy to see that
   $A^c \in \L_s$ if $A \in \L_s$.

    Let $\{A_n\}_{n \in \hN} \sb \L_s $ and $\e \ins (0,1)$. For any $n \in \hN$, there
   exist $\cO^1_n, \cO^2_n$ in $ \fT_s(\O_0) $     such that $(\cO^1_n)^c \sb A_n \sb \cO^2_n$
   and  $P_0 (A_n \Cp \cO^1_n) \ve P_0(A^c_n \Cp \cO^2_n  ) \< \e 2^{- 1-n }  $.
   The set  $\cO_2 \df \underset{n \in \hN}{\cup} \cO^2_n \ins \fT_s(\O_0) $ contains $A \df \underset{n \in \hN}{\cup} A_n $
   and satisfies $P_0 (A^c \Cp \cO_2)
   \ls \underset{n \in \hN}{\sum} P_0  (A^c \Cp \cO^2_n  )  \ls \underset{n \in \hN}{\sum} P_0  (A^c_n \Cp \cO^2_n)
   \< \e/2$.  Similarly, it holds for $\cE \df \underset{n \in \hN}{\cap} \cO^1_n $ that
   $  P_0 (A \cap \cE) 
    \ls \underset{n \in \hN}{\sum}  P_0  (A_n \cap \cO^1_n) \< \e/2 $.
   We can find an $N \ins \hN$ such that
   $ P_0 \Big(\underset{n =1 }{\overset{N}{\cap}}  \cO^1_n  \Big)
     \<    P_0  (  \cE   ) \+ \e/2$.
   Then $\cO_1 \df  \ccap{n =1 }{N} \cO^1_n $ is a set of $ \fT_s(\O_0)$ satisfying that
     $\cO^c_1 \=  \ccup{n =1 }{N} (\cO^1_n)^c \sb \ccup{n \in \hN }{ } A_n \= A$ and
    $    P_0 (A \Cp \cO_1 )  \= P_0 (A \Cp \cE ) \+ P_0 \big(A \Cp(\cO_1 \backslash \cE)\big)
      \ls  P_0 (A \Cp \cE ) \+ P_0(\cO_1 \backslash \cE) \< \e   $,
    which shows   $   \underset{n \in \hN}{\cup} A_n \= A \ins \L_s  $. Hence $ \L_s $ is a sigma-field of $\O_0$.

 Let $r \ins [0,s]$ and let  $ O $ be a nonempty  open subset of $\hR^d$.
 Given $\o_0 \ins W^{-1}_r(O)$, there exists $\d   \ins (0,1)$ such that  $ O_\d \big(W_r(\o_0)\big) \sb O$.
  As $ W_r \big( O^s_\d (\o_0) \big) \sb O_\d \big(W_r(\o_0)\big) $,
  we obtain that   $O^s_\d (\o_0) \sb  W^{-1}_r(O) $ and thus $ W^{-1}_r(O) \ins  \fT_s(\O_0) $.
 Let $\e \ins (0,1)$ and define   closed sets $ D_n \df \{x \ins \hR^d \n: \hb{dist}(x,O^c) \gs 1/n \}    $,
 $\fa n \ins \hN$.
 Since $ \ccap{n \in \hN}{} W^{-1}_r(O \backslash D_n)
  \= W^{-1}_r \Big( \ccap{n \in \hN}{} (O \backslash D_n) \Big) \= \es $,
    there exists $N$ such that
   $  P_0 \big( W^{-1}_r(O \backslash D_N) \big) \< \e  $.
  Similar to the inclusion $  W^{-1}_r(O) \ins \fT_s(\O_0) $,
  one has    $ \cO_1 \df W^{-1}_r(D^c_N) \ins \fT_s(\O_0)$.
  Since  $ \cO^c_1 \= W^{-1}_r(D_N) \sb W^{-1}_r(O) $
  and   $ P_0( W^{-1}_r(O)\Cp \cO_1) \= P_0 \Big( W^{-1}_r(O \Cp D^c_N) \Big) \< \e$,
  we see that $ W^{-1}_r(O) \ins \L_s $.
  It follows that $\cF^W_s \= \si \big( W^{-1}_r(O); r \ins [0,s], \, \hb{open subset $O$ of }\hR^d \big) \sb \L_s$.
  So \eqref{010122_11} holds.

   Let $ \big\{ \o^i_0 \big\}_{i \in \hN}$ be a countable dense subset of   $\O_0$  
  and let $s \ins [0,\infty)$.
  We   set   $\Th_s \df \big\{O^s_\d (\o^i_0) \n : \, \d \ins \hQ_+, \, i \ins \hN \big\} \sb \cF^W_s $.
  Let $A \ins  \cF^W_s$ and $\e \ins (0,1) $. By \eqref{010122_11},
  there exists    $\cO_2 \ins \fT_s(\O_0) $  such that $   A \sb \cO_2 $ and $P_0( \cO_2 ) \- P_0(A)   \< \e $.
  As usual,    $\cO_2$ is   the union of some sequence $\big\{ O_i \big\}_{i \in \hN}$
     in  $\Th_s$. So $A$ satisfies that
  \if{0}

     We demonstrate  that any   $ \cO $ of $ \fT_s(\O_0) $  is a union of some sequence $\big\{ O_i \big\}_{i \in \hN}$   in  $\Th_s$:
     For any $\o_0,\o'_0 \ins \O_0$, set $\{\n|\o'_0\- \o_0 |\n\}_s \df \Sup{r \in [0,s]}  |\o'_0(r)\- \o_0(r)|$.

     Let $i \in \hN$.  If $\o^i_0 \notin \cO$, we set $O_i \df \es$. Suppose next that $\o^i_0 \ins \cO$.
     There exists   $\d_i \ins (0,1)$ such that $ O^s_{\d_i} (\o^i_0) \sb \cO $.
     For any $\o'_0 \ins \cO^c \sb \big( O^s_{\d_i} (\o^i_0) \big)^c $, we have
     $ \{\n|\o'_0\- \o^i_0 |\n\}_s \gs \d_i$. So $  \l_i \df
        \underset{\o'_0 \in \cO^c}{\inf} \{\n|\o'_0 \- \o^i_0 |\n\}_s   \gs \d_i \> 0 $.
     We  choose a $q_i \in \hQ  \Cp (  \l_i/2, \l_i)$ and set      $ O_i \df O^s_{q_i} (\o^i_0 )   $.
     For any $\o'_0 \ins O^s_{q_i} (\o^i_0 ) $, since $ \{\n|\o'_0 \-\o^i_0 |\n\}_s \< q_i \< \l_i $,
     we see that $\o'_0$ can not be in $\cO^c$ and thus $O^s_{q_i} (\o^i_0 ) \sb \cO $.
     \if{0}
     We  claim that $O^s_{q_i} (\o^i_0 ) \sb \cO$: If not, let $\o'_0 \ins O^s_{q_i} (\o^i_0 ) \Cp \cO^c$.
     Then $\l_i \ls \{\n|\o'_0\- \o^i_0 |\n\}_s \ls q_i \< \l_i $. A contradiction appears.
     So $O^s_{q_i} (\o^i_0 ) \sb \cO$
     \fi

    Set $\fs \df \lceil s \rceil$ and let $\o_0 \ins \cO$. Similar to $\l_i$'s,
    we set $\l \n = \n \l(\o_0)   \df \underset{\o'_0 \in \cO^c}{\inf} \{\n|\o'_0 -\o_0 |\n\}_s \> 0   $.
      As $ \big\{ \o^i_0 \big\}_{i \in \hN}$ is a countable dense subset of   $\O_0$,
    there exists a  $\fj  \n = \n  \fj(\o_0) \ins  \hN$ such that
     $ \Rho{\O_0} (\o^\fj_0, \o_0) \ls 2^{-\fs}(\l/3)  $.
    It follows that $ 2^{-\fs} \ld \{\n|\o^\fj_0 \-\o_0 |\n\}_\fs \< 2^{-\fs}(\l/3)$ and thus
     \bea \label{010122_14}
     \{\n|\o^\fj_0 \-\o_0 |\n\}_s  \ls  \{\n|\o^\fj_0 \-\o_0 |\n\}_\fs \< 2^{-\fs}(\l/3) \ls \l/3  .
     \eea
 This shows
   $\o^\fj_0 \in   O^s_{\l/3} (\o_0) \sb \cO $ (Analogous to $O^s_{q_i} (\o^i_0 ) \sb \cO $).

     Since $\l_\fj \= \underset{\o'_0 \in \cO^c}{\inf} \{\n|\o'_0 \-\o^\fj_0 |\n\}_s  \gs \underset{\o'_0 \in \cO^c}{\inf} \{\n|\o'_0 \-\o_0 |\n\}_s
     \-  \{\n| \o^\fj_0 -\o_0  |\n\}_s > \frac23 \l    $, one has  $q_{\overset{}{\fj}} > \l_\fj /2 >  \l /  3 $.
     Then we can deduce from \eqref{010122_14} that
       $\o_0 \ins  O^s_{\l/3} (\o^\fj_0 ) \sb O^s_{q_{\overset{}{\fj}}} (\o^\fj_0 ) = \cO_\fj   $.
      It follows that    $ \cO = \underset{i \in \hN}{\cup} O_i $.

 \fi
     \bea \label{011222_11}
       A \sb   \ccup{i \in \hN}{ } O_i \aand P_0(A) \>   P_0 \Big( \underset{i \in \hN}{\cup} O_i   \Big) \- \e  .
      \eea

  \no {\bf 3)}
Now we are ready to demonstrate the separability of $\big(\fS, \Rho{\fS}\big)$.

 Given $q \ins \hQ_+$,
    let us simply denote  by  $  \{O^q_j\}_{j \in \hN}$ the  countable sub-collection $\Th_q \= \big\{O^q_\d (\o^i_0) \n : \, \d \ins \hQ_+, \, i \ins \hN \big\}$ of $\cF^W_q$  and
   define  $\U^q_{k,\a}   \df      \Big\{  q  \b1_{\underset{j \in I }{\cup } O^q_j}
       +      \a \b1_{\underset{j \in I }{\cap }   ( O^q_j )^c}     \n  :
     I \sb \{1,\cds     , k  \} \Big\}      \sb  \fS $, $ \fa k       \in      \hN$, $\fa \a \ins \hN \Cp [q,\infty)$.
       For any $ k, n        \in      \hN $,    we set
  $  \wh{\U}_{k,n }      \df  \ccup{\a \in \hN}{}     \Big\{ \underset{i = 1 }{\overset{2^n \a}{\land}}  \tau_i     \n  :
  \tau_i      \ins      \U^{i2^{-n}}_{k,\a} ,\,   i \= 1,\cds \n , 2^n \a \Big\}   $, which is a countable subset  of   $      \fS  $.
 Then   $\wh{\U} \df   \underset{ k,n   \in \hN}{\cup} \, \wh{\U}_{k,n }  $ is also a countable subset of $  \fS  $.
   To show   $\wh{\U}$ is dense in $\big(\fS, \Rho{\fS}\big)$,
 we let $\tau \ins \fS$, $\e \ins (0,1)$ and try  to pick   $\ga \ins \wh{\U}$ such that
 \bea \label{031221_17}
 \Rho{\fS} (\tau , \ga) \< \e .
 \eea

  Since 
   $ \lmt{\a \to \infty} \Rho{\fS}   ( \tau,\tau \ld \a  )
   \= \lmt{\a \to \infty} E_{P_0}  \big[   \big| \arctan(\tau)  \-  \arctan(\tau \ld \a)   \big|  \big] \= 0$,
    one can find $\wh{\a} \ins \hN$ such that $\Rho{\fS}   ( \tau,\tau \ld \wh{\a}  ) \< \e/4$.

   Let $ n \ins \hN $ and set $s^n_0 \df 0$. Given $i \ins \{ 0, 1, \cds \n , 2^n \wh{\a} \}   $, we set $s^n_i \df i2^{-n}$
  and    $A^n_i  \df  \{ s^n_{i-1}   \ls  \tau  \< s^n_i  \}  \ins \cF^{W,P_0}_{ s^n_i }  $.
  By e.g.   Problem 2.7.3 of \cite{Kara_Shr_BMSC},
  there exists   $\cA^n_i   \ins  \cF^W_{ s^n_i }$   such that
  $ \cN^n_i \df  A^n_i   \D   \cA^n_i   \ins  \sN_{P_0}(\cF^W_\infty)  $.
     Define    $\wcA^n_i  \df  \cA^n_i \Big\backslash \Big(\underset{j < i}{\cup}
     \cA^n_j \Big)  \ins  \cF^W_{ s^n_i }  $
     and $ \wcA_n  \df  \ccup{i = 1}{2^n \wh{\a}}  \wcA^n_i
      \= \ccup{i = 1}{2^n \wh{\a}}  \cA^n_i  \ins  \cF^W_{\wh{\a}} $.
     The $\bF^{W,P_0}-$stopping time $ \dis  \tau_n  \df  \sum^{2^n \wh{\a}}_{i = 1} s^n_i \b1_{A^n_i}   \+ \wh{\a} \b1_{\{\tau \ge \wh{\a}\}}$ coincides with the $\bF^W-$stopping time
     $ \dis  \wt{\tau}_n  \df  \sum^{2^n \wh{\a}}_{i = 1} s^n_i \b1_{\wcA^n_i} \+   \wh{\a} \b1_{ \wcA^c_n }  $   over
     $\O_n \df \Big( \ccup{i = 1}{2^n \wh{\a}} \big(  A^n_i \Cp \wcA^n_i \big) \Big) \cp \Big(\{\tau \gs \wh{\a}\} \Cp \wcA^c_n \Big) $.
     We  can deduce that
     \bea
     A^n_i   \backslash    \wcA^n_i  & \tn \= &  \tn  A^n_i \Cp
      \Big[ \big( \cA^n_i \big)^c \cp \Big( \underset{j < i}{\cup} \cA^n_{j} \Big)  \Big]
      \= \big( A^n_i \Cp ( \cA^n_i  )^c \big) \cp \Big( \underset{j < i}{\cup} \big(  \cA^n_{j} \Cp A^n_i \big)  \Big)
      \nonumber \\
     & \tn  \sb & \tn  \big( A^n_i \D \cA^n_i \big) \cp \Big( \underset{j < i}{\cup} \big( \cA^n_{j} \Cp (A^n_{j})^c  \big)  \Big)
      \sb \underset{j \le i}{\cup}   \cN^n_j   \ins  \sN_{P_0}(\cF^W_\infty) ,
      \q  \hb{ for } i \= 1, \cds, 2^n \wh{\a}   , \qq  \label{eq:cc131} \\
     \hb{and} \q  \{\tau \gs \wh{\a}\} \Cp \wcA_n  & \tn \= & \tn  \ccup{i = 1}{2^n \wh{\a}} \big( \{\tau \gs \wh{\a}\} \Cp \cA^n_i \big)
      \sb \ccup{i = 1}{2^n \wh{\a}} \big( (A^n_i)^c \Cp \cA^n_i \big)
      \sb \ccup{i = 1}{2^n \wh{\a}} \cN^n_i  \ins  \sN_{P_0}(\cF^W_\infty) .   \nonumber
     \eea
     Putting them together shows that
     $\O^c_n \=  \Big( \ccup{i = 1}{2^n \wh{\a}} \big(  A^n_i \backslash \wcA^n_i \big) \Big) \cp
     \Big(  \{\tau \gs \wh{\a}\} \Cp \wcA_n  \Big)  $   belongs to   $ \sN_{P_0}(\cF^W_\infty) $.
   To wit,
   \bea \label{031221_11}
    \tau_n \= \wt{\tau}_n   , \q  P_0 \n - \n a.s.
    \eea
    Since  $ \lmtd{n \to \infty}  \tau_n  \= \tau \ld \wh{\a}   $,
    one has   $ \lmt{n \to \infty} \Rho{\fS}   ( \tau \ld \wh{\a},\tau_n  )
   \= \lmt{n \to \infty} E_{P_0}  \big[   \big| \arctan(\tau \ld \wh{\a})  \-  \arctan(\tau_n)   \big|  \big] \= 0$.
   So there exists $\fn \ins \hN$ such that
   $ \Rho{\fS}   \big( \tau \ld \wh{\a},\tau_\fn  \big)   \< \e/4$.

  Given $i \ins \{ 1, \cds \n , 2^\fn \wh{\a} \} $, we know from \eqref{011222_11} that
  for some sequence $  \{ O^i_j  \}_{j \in \hN}$
  in $\Th_{ s^\fn_i } = \big\{ O^{s^\fn_i}_j \big\}_{j \in \hN}  $
   \bea   \label{eq:bb431}
      \wcA^\fn_i \sb \underset{ j \in \hN }{\cup} O^i_j
  \q \hb{and} \q   P_0 \big( \wcA^\fn_i \big) \> P_0 \Big( \underset{ j \in \hN }{\cup} O^i_j  \Big)  \-   \frac{\e}{2^{2+2\fn} \pi \wh{\a}^2  } .
  \eea
 And we can find   $  \ell_i  \ins  \hN$
  such that $ \cO_i  \df  \underset{ j = 1 }{\overset{ \ell_i  }{\cup}} O^i_j  \ins \cF^W_{s^\fn_i} $ satisfies
  \bea    \label{eq:bb433}
    P_0  (  \cO_i   )  \>  P_0 \Big( \underset{ j \in \hN }{\cup} O^i_j  \Big)
     \-  \frac{\e}{2^{2+2\fn} \pi \wh{\a}^2   } .
     \eea
      Clearly,     $ \ga_i  \df  s^\fn_i   \b1_{ \cO_i }   \n + \n  \wh{\a} \b1_{  \cO^c_i  }
   \ins  \U^{s^\fn_i}_{k_i,\wh{\a}}$ for some $k_i \ins \hN$.

 Let $i \ins \{ 1, \cds \n , 2^\fn \wh{\a} \}  $ and set   $\wt{\cO}_i \df \cO_i \backslash  \ccup{\fri < i}{ }  \cO_\fri   \in \cF^W_{s^\fn_i}  $.
 Analogous to  \eqref{eq:cc131},
 $ \wt{A}^\fn_i \backslash \wt{\cO}_i  \= \wt{A}^\fn_i \cap \Big[  \cO_i^c \cp \big( \ccup{\fri < i}{ }  \cO_\fri \big)  \Big]
 \sb \Big( \big( \underset{ j \in \hN }{\cup} O^i_j \big) \Cp \cO^c_i  \Big) \cp \Big( \ccup{\fri<i}{} \big(    \cO_\fri \Cp (\wt{A}^\fn_\fri)^c \big) \Big) $.  So   \eqref{eq:bb431} and \eqref{eq:bb433} yield  that
  \bea \label{eq:aa115}
   P_0 (\wt{A}^\fn_i \backslash \wt{\cO}_i)\le P_0 \Big( \big( \underset{ j \in \hN }{\cup} O^i_j \big) \Cp \cO^c_i  \Big) +  \sum_{\fri < i} \, P_0   \Big(  \big( \underset{ j \in \hN }{\cup} O^\fri_j \big)   \Cp (\wt{A}^\fn_{\fri})^c \Big)
   \< \sum_{\fri \le i} \frac{\e}{2^{2+2\fn} \pi \wh{\a}^2  } \= \frac{i\e}{2^{2+2\fn} \pi \wh{\a}^2 } \ls \frac{ \e}{ 2^{2+\fn} \pi \wh{\a}   }   .
    \eea

  Define $ \wt{\cO}  \df  \ccup{i=1}{2^\fn \wh{\a} } \, \wt{\cO}_i
   \n = \n  \ccup{i=1}{2^\fn \wh{\a} }  \cO_i \ins \cF^W_{\wh{\a}} $
  and $\fk   \df  \max\{   k_i  \n : i \n = \n 1, \cds \n , 2^\fn \wh{\a}   \}$.  Then
         $  \ga \df  \underset{i=1}{\overset{ 2^\fn \wh{\a}  }{\land}}  \ga_i $
     is a stopping time of      $ \wh{\U}_{\fk,\fn } $  \big(and is thus of $ \wh{\U} $\big).
      In particular, $\ga  \=  \underset{i=1}{\overset{ 2^\fn \wh{\a}  }{\sum}} s^\fn_i \b1_{\wt{\cO}_i }
      \n + \n \wh{\a} \b1_{\wt{\cO}^c}  $ is equal  to
  $   \wt{\tau}_\fn \= \underset{i=1}{\overset{ 2^\fn \wh{\a}  }{\sum}} s^\fn_i \b1_{\wcA^\fn_i} \+   \wh{\a} \b1_{ \wcA^c_\fn }   $   over  $\wh{\cA} \df \Big(\underset{i = 1}{\overset{ 2^\fn \wh{\a}}{\cup}} \big( \wt{\cO}_i \Cp \wt{A}^\fn_i  \big)\Big) \cp \big( \wt{\cO}^c \Cp \wcA^c_\fn   \big) \\ \ins \cF^W_{\wh{\a}} $. Since \eqref{eq:bb431} implies that
 \beas
 P_0 \big(   \wt{\cO} \Cp \wcA^c_\fn \big)
 \ls \underset{i = 1}{\overset{  2^\fn \wh{\a}}{\sum}} P_0 \big( \wt{\cO}_i \Cp \wcA^c_\fn    \big)
 \ls \underset{i = 1}{\overset{  2^\fn \wh{\a}}{\sum}} P_0 \big( \wt{\cO}_i \Cp (\wcA^\fn_i)^c   \big)
\ls \underset{i = 1}{\overset{  2^\fn \wh{\a}}{\sum}} P_0 \big(   \cO_i \backslash \wcA^\fn_i \big)
\ls \underset{i = 1}{\overset{  2^\fn \wh{\a}}{\sum}} P_0 \Big\{ \Big( \ccup{ j \in \hN }{}  O^i_j \Big) \Big\backslash \wcA^\fn_i    \Big\}
\< \underset{i = 1}{\overset{  2^\fn \wh{\a}}{\sum}}  \frac{\e}{2^{2+2\fn} \pi \wh{\a}^2  }  \= \frac{\e}{2^{2+ \fn} \pi \wh{\a}  } ,
\eeas
    \eqref{eq:aa115} renders  that
      $ \dis  P_0 \big( \wh{\cA}^c \big)  \=
        P_0 \Big\{ \Big( \underset{i = 1}{\overset{ 2^\fn \wh{\a}}{\cup}}
      \big( \wt{A}^\fn_i \backslash \wt{\cO}_i \big) \Big) \cp
        \big( \wt{\cO} \Cp \wcA^c_\fn   \big) \Big\}
      \=  \underset{i = 1}{\overset{  2^\fn \wh{\a}}{\sum}}
     P_0 \big( \wt{A}^\fn_i \backslash \wt{\cO}_i \big)
  \+ P_0 \big( \wt{\cO} \Cp \wcA^c_\fn   \big)
   \< ( 2^\fn \wh{\a} \+ 1)\frac{\e}{2^{2+ \fn} \pi \wh{\a}  } \ls \frac{\e}{2 \pi  } $
   and thus
   $ \Rho{\fS} \big(\wt{\tau}_\fn,\ga \big) 
    =  E_{P_0}  \big[ \b1_{\wh{\cA}^c}   \big| \arctan(\wt{\tau}_\fn )   \-  \arctan(\ga )    \big| \big]
      \le \pi  P_0 \big( \wh{\cA}^c \big)   \< \e/2 $.
      By \eqref{031221_11},   it follows that
   $ \Rho{\fS} (\tau ,\ga)   \ls
    \Rho{\fS} (\tau ,\tau \ld \wh{\a}) \+     \Rho{\fS}   \big( \tau \ld \wh{\a},\tau_\fn  \big) \+  \Rho{\fS} \big(\wt{\tau}_\fn , \ga\big)
   \< \e $, proving \eqref{031221_17}.
 Therefore, $\big(\fS, \Rho{\fS}\big)$ is   a complete separable space.    \qed

\no {\bf Proof of Lemma \ref{lem_082020_17}: 1)} We first show that   $\Ga$ is injective:
Set $\hQ_\pi \df \big(\hQ  \Cp [0,\pi/2) \big) \cp \{\pi/2\} $ and let $ \tau_1 , \tau_2  \ins   \fS$ such that $\Ga( \tau_1) \= \Ga( \tau_2)$.

  Given $q   \ins \hQ_\pi   $ and $n \ins \hN$, we define $\cE^q_n \df (q\-1/n,q\+1/n) \Cp [0,\pi/2]$
 and  $ A^i_{n,q} \df    \big\{ \arctan(\tau_i) \ins \cE^q_n \big\}   \ins \cF^{W,P_0}_\infty$ for $i \= 1,2$.
 Then   $A_{n,q} \df  A^1_{n,q} \Cp (A^2_{n,q})^c$  satisfies
 \bea
 \hspace{-0.8cm}
  P_0 (A_{n,q} ) & \tn \= & \tn    P_0 \big\{ \o_0 \ins (A^2_{n,q})^c :    \tau_1 (\o_0)  \ins \tan(\cE^q_n) \big\}
  \= \big( \Ga ( \tau_1) \big) \big( (A^2_{n,q})^c   \ti \tan(\cE^q_n) \big)
     \=   \big( \Ga ( \tau_2) \big) \big( (A^2_{n,q})^c   \ti \tan(\cE^q_n) \big)   \nonumber  \\
   & \tn   \=  & \tn    P_0 \big\{ \o_0 \ins (A^2_{n,q})^c :   \tau_2 (\o_0) \ins \tan(\cE^q_n) \big\}
     \= P_0(\es) \= 0   .   \label{031821_a01}
 \eea

  Clearly, $ \underset{n  \in \hN}{\cup} \, \ccup{q \in \hQ_\pi}{} A_{n,q} \sb \{ \tau_1   \nne \tau_2 \} $.
  To see the reverse inclusion, we let $\o_0 \ins \Big(  \underset{n  \in \hN}{\cup} \, \ccup{q \in \hQ_\pi}{}  A_{n,q} \Big)^c $
  and let $n   \ins \hN$.
  There exists   $\fq \= \fq(n) \ins \hQ_\pi$ such that    $\arctan\big(\tau_1(\o_0)\big) \ins  \cE^\fq_n$,
  or  $ \o_0 \ins   \big\{\arctan(\tau_1) \ins  \cE^\fq_n \big\}   \= A^1_{n,\fq}  $.
  As $ \o_0 \ins  A^c_{n,\fq}  $, we see that $ \o_0 \ins A^2_{n,\fq} $, i.e.,
  $    \arctan(\tau_2 (\o_0))   $ also belongs to $    \cE^\fq_n  $. It follows that
    $\Rho{+} \big(\tau_1(\o_0), \tau_2(\o_0) \big) \= \big|\arctan(\tau_1(\o_0)) \- \arctan(\tau_2(\o_0)) \big| \< 2/n$.
  Letting $n \nto \infty$ yields that   $ \tau_1(\o_0) \= \tau_2(\o_0) $.
  So $ \underset{n  \in \hN}{\cup} \, \ccup{q \in \hQ_\pi}{} A_{n,q} \= \{ \tau_1   \nne \tau_2 \} $.
  It follows from   \eqref{031821_a01} that   $  P_0 \big\{  \tau_1   \nne \tau_2   \big\}
  \= 0$,   which means that  $\tau_1 \= \tau_2$ in $\fS$.
  Hence,  the mapping $\Ga \n :   \fS \mto \fP \big(\O_0 \ti \hT\big)$   is injective.

\no  {\bf 2)} We next discuss the continuity of $\Ga$:
    Let $\{\tau_n \}_{n \in \hN}  $ be a sequence of $\fS$ that converges to a $\tau \ins \fS$ under $\Rho{\fS}$.
 We need to show   that $ P^n \df \Ga( \tau_n)$ converges to $ P \df \Ga( \tau)$ under
 the weak topology of $\fP \big(\O_0 \ti \hT\big)$,
 i.e.
 \bea \label{Au18_01}
  \lmt{n \to \infty} \int_{(\o_0,\ft) \in \O_0  \times \hT} \phi (\o_0,\ft) P^n \big(d(\o_0,\ft)\big)
  \= \int_{(\o_0,\ft) \in \O_0  \times \hT} \phi (\o_0,\ft) P \big(d(\o_0,\ft)\big)
 \eea
 for any bounded continuous function $\phi \n :  \O_0 \ti \hT \mto \hR$.

    Let $\phi$ be   a bounded continuous function on $ \O_0 \ti \hT $.
  For \eqref{Au18_01}, it suffices to show that for any subsequence $ \{ \tau_{n_k}  \}_{k \in \hN}$ of
  $ \{ \tau_n  \}_{n \in \hN}  $, we can find a subsequence $\big\{ \tau_{n'_k} \big\}_{k \in \hN}$   of   $ \{ \tau_{n_k}  \}_{k \in \hN}$ satisfying \eqref{Au18_01}.

   Let $ \{ \tau_{n_k}  \}_{k \in \hN}$ be an arbitrary subsequence of  $ \{ \tau_n  \}_{n \in \hN}  $.
  As $ \dis 0 \=  \lmt{k \to \infty} \Rho{\fS} \big(\tau_{n_k},\tau \big) \= \lmt{k \to \infty}   E_{P_0} \big[ \Rho{+} (\tau_{n_k},\tau )   \big]   $, one can extract a subsequence $\big\{  {n'_k} \big\}_{k \in \hN} $ from $ \phantom{\Big|} \big\{  {n_k} \big\}_{k \in \hN} $ such that $ \lmt{k \to \infty}   \Rho{+} \big(\tau_{n'_k}(\o_0),\tau(\o_0) \big)   \= 0 $ for all $\o_0 \ins \O_0$  except on a $P_0-$null set $\cN $.
 Given $\o_0 \ins  \cN^c   $,
 since    $ \lmt{k \to \infty}   \Rho{+} \big(\tau_{n'_k}(\o_0),\tau(\o_0) \big)   \= 0 $,
   the continuity of $\phi$ renders that $ \lmt{k \to \infty} \phi \big( \o_0,   \tau_{n'_k} (\o_0) \big)
  \= \phi \big( \o_0,  \tau (\o_0) \big) $.   Applying the bounded convergence theorem   yields  that
 $ \lmt{k \to \infty} \int_{ (\o_0,\ft) \in \O_0  \times \hT} \phi (\o_0,\ft) P^{n'_k} \big(d(\o_0,\ft)\big)
    \=    \lmt{k \to \infty} \int_{\O_0 } \phi \big( \o_0,\tau_{n'_k}(\o_0)\big) P_0 (d \o_0)
   \= \int_{\O_0 } \phi \big( \o_0,\tau(\o_0)\big) P_0 (d \o_0)
     \=    \int_{ (\o_0,\ft) \in \O_0  \times \hT} \phi (\o_0,\ft) P \big(d(\o_0,\ft)\big) $.  \qed

\no {\bf Proof of Proposition \ref{prop_Ptx_char}:}
Fix $(t,\bx) \ins [0,\infty) \ti \OmX$.

  \no {\bf 1)} Let $\oP \ins \ocP_{t,\bx}$. It is   clear that   $\oP \ins \ocP^1_{t,\bx}$.
 Let $(\vf,n) \ins \fC(\hR^{d+l}) \ti \hN$.  By (D1$'$) of Remark \ref{rem_ocP}, $\big\{\oM^t_{s  \land \otau^t_n  } (\vf)\big\}_{s  \in [t,\infty)}$ is a bounded $(\obF^t,\oP)-$martingale.
 For any  $(\fs,\fr) \ins \hQ^{2,<}_+  $  and $  \{(s_i,\cO_i )\}^k_{i=1} \sb \big(\hQ \cap [0,\fs]\big) \ti \sO (\hR^{d+l})  $,
 as $ \big\{ (\oW^t_{ t+s_i  },\oX_{ t+s_i   }) \ins \cO_i \big\} \ins 
 \ocF^t_{t+\fs}$ for $i \= 1, \cds \n , k$,
  one directly   has  $ E_\oP \Big[ \big(\oM^t_{  \otau^t_n \land (t+\fr) } (\vf) \- \oM^t_{  \otau^t_n \land (t+\fs) } (\vf) \big) \underset{i=1}{\overset{k}{\prod}}   \b1_{    \{(\oW^t_{ t+s_i  },\oX_{ t+s_i   }) \in \cO_i  \}    }  \Big] \= 0 $.
  So $\oP \ins \ocP^2_t$.

   By (D3$'$) of  Remark \ref{rem_ocP},
 there exists    a $[0,\infty]-$valued $\bF^{W,P_0} -$stopping time $\ddot{\tau}$  on $\O_0$  such that
 $   \oP \big\{ \oT \=  t \+ \ddot{\tau} \big(\osW^t \big) \big\} \= 1$.
 Since $\osW^t_{\n \fs} \= \oW_{t+\fs} \- \oW_t$, $\fs \ins [0,\infty)$ is a Brownian motion under $\oP$ by (D1) of Definition  \ref{def_ocP},
 applying Lemma \ref{lem_122921_11}   with $t_0 \= 0$, $(\O_1, \cF_1, P_1,B^1)   \= \big(\oO ,  \sB(\oO ),  \oP , \osW^t\big) $, $(\O_2, \cF_2, P_2,B^2) \= \big(\O_0,  \sB(\O_0),  P_0, W\big) $ and $\Phi \= \osW^t$   shows that
   \if{0}

 \bea \label{010722_21}
 \oP \nci \big(\osW^t\big)^{-1}  ( A_0 ) \= P_0  ( A_0 )  , \q \fa A_0 \ins \cF^{W,P_0}_\infty  .
 \eea
 For any   $\cA_0 \ins \sB(\O_0) \= \cF^W_\infty$ and   $\cE \ins \sB(\hT)$,
 since $ \ddot{\tau}^{-1}(\cE) \ins \cF^{W,P_0}_\infty $, applying \eqref{010722_21} with $A_0 \= \cA_0 \Cp \ddot{\tau}^{-1}(\cE) $ yields that
 \beas
  && \hspace{-1cm} \oP \nci \big(\osW^t, \oT \- t\big)^{-1} (\cA_0  \ti \cE)
  \=   \oP \big\{     \big(\osW^t  , \oT \- t  \big) \ins \cA_0  \ti \cE \big\}
  \= \oP \big\{    \big(\osW^t  ,  \ddot{\tau} (\osW^t   )   \big) \ins \cA_0  \ti \cE \big\}
     \=      \oP \nci (\osW^t)^{-1}  \big\{     (W,  \ddot{\tau}  ) \ins \cA_0  \ti \cE\big\} \\
  && \=      \oP \nci (\osW^t)^{-1}  \big(     \cA_0  \Cp \ddot{\tau}^{-1}(\cE)\big)
   \= P_0  \big(     \cA_0  \Cp \ddot{\tau}^{-1}(\cE)\big)
   \=  P_0 \big\{     (W,  \ddot{\tau}  ) \ins \cA_0  \ti \cE \big\}
    \=       P_0 \nci (W, \ddot{\tau})^{-1}    (\cA_0  \ti \cE) .
 \eeas
 Then the Lambda-system $\big\{ D \ins \sB(\O_0   \ti \hT) \n : \oP \nci \big(\osW^t,\oT \- t\big)^{-1}(D)
 \=   P_0 \nci (W, \ddot{\tau})^{-1}   (D) \big\} $
 contains all measurable rectangles of $   \sB(\O_0  )   \oti \sB (\hT) $
 and is thus equal  to  $ \sB(\O_0  )   \oti \sB (\hT) \= \sB\big(\O_0   \ti \hT \big)$ thanks to Dynkin's Pi-Lambda Theorem.
 So  $ \oP \nci \big(\osW^t,\oT \- t \big)^{-1}   \= P_0 \nci (W, \ddot{\tau})^{-1}    \= \Ga ( \ddot{\tau}) \ins \Ga(  \fS)  $, i.e.,  $\oP \ins \ocP^3_t$.

   \fi
 \bea \label{010922_14}
 \oP \nci \big(\osW^t\big)^{-1} (A_0)  \= P_0 (A_0)  , \q \fa   A_0  \ins  \cF^{W,P_0}_\infty   .
 \eea
 For any   $ \cA_0 \ins \sB(\O_0) \= \cF^W_\infty$ and   $\cE \ins \sB(\hT)$, since $\ddot{\tau}^{-1}(\cE) \ins \cF^{W,P_0}_\infty$,
 we can derive  that
 \beas
 && \hspace{-1cm} \oP \nci \big(\osW^t, \oT \- t\big)^{-1} (\cA_0  \ti \cE)
   \=      \oP \big\{     \big(\osW^t  , \oT \- t  \big) \ins \cA_0  \ti \cE \big\}
  \= \oP \big\{    \big(\osW^t  ,  \ddot{\tau} (\osW^t   )   \big) \ins \cA_0  \ti \cE \big\}
     \=      \oP \nci (\osW^t)^{-1}  \big\{     (W,  \ddot{\tau}  ) \ins \cA_0  \ti \cE\big\} \\
   & &  \=       \oP \nci (\osW^t)^{-1}  \big(     \cA_0  \Cp \ddot{\tau}^{-1}(\cE)\big)
   \= P_0  \big(     \cA_0  \Cp \ddot{\tau}^{-1}(\cE)\big)
   \=  P_0 \big\{     (W,  \ddot{\tau}  ) \ins \cA_0  \ti \cE \big\}
     \=       P_0 \nci (W, \ddot{\tau})^{-1}    (\cA_0  \ti \cE) .
 \eeas
 Then Dynkin's Pi-Lambda Theorem implies that
 $ \oP \nci \big(\osW^t, \oT \- t\big)^{-1} \=  P_0 \nci (W, \ddot{\tau})^{-1}    $ on  $ \sB\big(\O_0   \ti \hT \big) $.
  i.e.,  $   \oP \nci \big(\osW^t,\oT \- t \big)^{-1}  
  \n  \= \Ga ( \ddot{\tau})   \ins \Ga(  \fS)  $.
 So  $\oP $ also belongs to $ \ocP^3_t$, which shows $ \ocP_{t,\bx} \sb \ocP^1_{t,\bx} \Cp \ocP^2_t \Cp \ocP^3_t   $.

  \no {\bf 2a)} Let $\oP \ins \ocP^1_{t,\bx} \Cp \ocP^2_t   $.
 To see that $\oP$ satisfies (D1$'$) of  Remark \ref{rem_ocP}, we take  $(\vf,n) \ins \fC(\hR^{d+l}) \ti \hN$.
 As $\oP  \big\{   \oX_s   \= \bx(s) , \fa s \ins [0,t]  \big\}   \=1$,
applying Proposition \ref{prop_MPF1} with $\big(\O,\cF,P,B,X  \big)\=\big(\oO,\sB(\oO),\oP ,\oW,\oX  \big) $ implies that $\big\{\oM^t_{s \land \otau^t_n } (\vf) \big\}_{s \in [t,\infty)} $ is a bounded $\obF^t-$adapted continuous process under $ \oP  $.

  Let $ (\fs,\fr) \ins \hQ^{2,<}_+  $,  $ \big\{(t_i,\cO_i )\big\}^k_{i=1} \n \sb  \big( \hQ \cap [0,t] \big)  \ti \sO (\hR^l)  $ and    $   \big\{(s_j,\cO'_j )\big\}^m_{j=1} \n \sb \big(\hQ \cap (0,\fs]\big) \ti \sO (\hR^{d+l})   $.
  If $\bx(t_\fri) \notin \cO_\fri$ for some $\fri \ins \{1,  \cds \n , k\}$, then   $\oP\{\oX_{t_\fri} \ins \cO_\fri\} \= 0$ and thus
  $ E_\oP \Big[  \big( \oM^t_{\otau^t_n \land (t+ \fr)} (\vf ) \- \oM^t_{\otau^t_n \land (t+ \fs)} (\vf ) \big)   \underset{i=1}{\overset{k}{\prod}} \b1_{    \{  \oX_{t_i}  \in  \cO_i\}    } \underset{j=1}{\overset{m}{\prod}} \b1_{    \{ (\oW^t_{t+s_j},\oX_{t+s_j}) \in  \cO'_j\}    }   \Big]  \= 0 $.
  On the other hand, if $\bx(t_i) \ins \cO_i$ for each $i \ins \{1,  \cds \n , k\}$, then
  $  E_\oP \Big[  \big( \oM^t_{\otau^t_n \land (t+ \fr)} (\vf ) \- \oM^t_{\otau^t_n \land (t+ \fs)} (\vf ) \big)     \underset{i=1}{\overset{k}{\prod}} \b1_{    \{  \oX_{t_i}  \in  \cO_i\}    }    \underset{j=1}{\overset{m}{\prod}} \\  \b1_{    \{ (\oW^t_{t+s_j},\oX_{t+s_j}) \in  \cO'_j\}    }  \Big]
    \=  E_\oP \Big[  \big( \oM^t_{\otau^t_n \land (t+ \fr)} (\vf ) \- \oM^t_{\otau^t_n \land (t+ \fs)} (\vf ) \big)   \underset{j=1}{\overset{m}{\prod}} \b1_{    \{ (\oW^t_{t+s_j},\oX_{t+s_j}) \in  \cO'_j\}    }     \Big] \= 0  $.
 So the Lambda-system 
 $ \ol{\L}^{t,n}_{\fs,\fr}  \df \Big\{\oA \ins \sB\big(\oO\big) \n : E_\oP \Big[  \big( \oM^t_{\otau^t_n \land (t+ \fr)} (\vf ) \- \oM^t_{\otau^t_n \land (t+ \fs)} (\vf ) \big)  \b1_{\oA} \Big]   \= 0 \Big\} $
  includes the   Pi-system 
  $   \Big\{ \Big( \underset{i=1}{\overset{k}{\cap}}    \oX_{t_i}^{\,-1}(\cO_i)   \Big) \Cp \Big( \underset{j=1}{\overset{m}{\cap}}   (\oW^t_{t+s_j},\oX_{t+s_j})^{-1}(\cO'_j)   \Big)
       \n :    \big\{(t_i,\cO_i )\big\}^k_{i=1} \n \sb  \big( \hQ \Cp [0,t]  \big) \ti \sO (\hR^l) , \; \big\{(s_j,\cO'_j )\big\}^m_{j=1} \n \sb \big(\hQ \Cp (0,\fs]\big) \ti \sO (\hR^{d+l})  \Big\} $,
  which generates $\ocF^t_{t+\fs} $.
 \if{0}

 As $\sO(\hR^{d+l})$ is closed under intersection, $ \ol{\sC}^t_s $ is a Pi-system  of $\oO$.
 It is clear that $ \si \big( \ol{\sC}^t_s \big) \sb \ocF^t_{t+s} $.
 The   continuity of process $ \{\oW^t_r\}_{r \in [t,\infty)} $ and process $ \{\oX_r\}_{r \in [0,\infty)} $ implies that
  \bea
   \ocF^t_{t+s}   \=    \si \big( \big(\oW^t_r,\oX_r\big) ; r \ins   (t,t\+s]\big) \ve \si \big(\oX_r; r \ins [0,t]\big)
   \=  \si \big( \big(\oW^t_{t+r},\oX_{t+r}\big) ; r \ins  \hQ  \Cp  (0,s]   \big) \ve \si \big(\oX_r; r \ins [0,t]\big) .   \label{Jan13_21}
  \eea

   Let $r \ins  \hQ  \Cp  (0,s]  $. Since
 $ \big(\oW^t_r,\oX_r\big)^{-1}(\cO)   \ins \ol{\sC}^t_s $, $ \fa \cO \ins \sO(\hR^{d+l}) $,
 the sigma-field 
 $ \L^t_s \df \big\{ \cE \sb \hR^{d+l} \n : \big(\oW^t_r,\oX_r\big)^{-1} (\cE) \ins \si \big(\ol{\sC}^t_s\big)\big\}$
 contains $ \sO(\hR^{d+l}) $. Then
 $ \sB(\hR^{d+l}) \= \si \big( \sO(\hR^{d+l}) \big) \sb \L^t_s $ or
 $ \big(\oW^t_r,\oX_r\big)^{-1} (\cE) \ins \si \big(\ol{\sC}^t_s\big) $ for any $   \cE \ins \sB(\hR^{d+l}) $.
 It follows from \eqref{Jan13_21} that
 $  \ocF^t_{t+s}  \=  \si \big( \big(\oW^t_{t+r},\oX_{t+r}\big) ; r \ins  \hQ  \Cp  (0,s]   \big) \ve \si \big(\oX_r; r \ins [0,t]\big)
 \sb \si \big(\ol{\sC}^t_s\big) \sb \ocF^t_{t+s} $.

 \fi
     Dynkin's Pi-Lambda Theorem renders that $ \ocF^t_{t+\fs} \sb \ol{\L}^{t,n}_{\fs,\fr}  $, i.e.,
  \bea \label{010922_11}
   E_\oP \Big[  \big( \oM^t_{\otau^t_n \land (t+ \fr)} (\vf ) \- \oM^t_{\otau^t_n \land (t+ \fs)} (\vf ) \big)  \b1_{\oA} \Big]
   \= 0 ,  \q \fa \oA \ins \ocF^t_{t+\fs}  .
  \eea

   Let $t \ls s \< r \< \infty$ and $\oA \ins \ocF^t_s $.
   \if{0}

   Given $ k \ins \hN $  with $k \> -\log_2(r\-s) $,
   set $s_k \df \frac{\lceil s 2^k \rceil}{2^k} \ins \hQ_+ $ and $r_k \df \frac{\lceil r 2^k \rceil}{2^k} \ins \hQ_+ $.
   As $\oA \ins \ocF^t_{t+s_k}$,  taking $(s,r) \= (s_k,r_k)$ in \eqref{010922_11} yields that
   $  E_\oP \Big[  \big( \oM^t_{\otau^t_n \land (t+ r_k)} (\vf ) \- \oM^t_{\otau^t_n \land (t+ s_k)} (\vf ) \big)   \b1_{\oA} \Big]   \= 0 $.
   Letting $k \nto \infty$, we can deduce from the continuity of bounded process $\big\{\oM^t_{s \land \otau^t_n } (\vf) \big\}_{s \in [t,\infty)} $     and the bounded convergence theorem that
   $  E_\oP \Big[ \big( \oM^t_{\otau^t_n \land (t+ r )} (\vf ) \- \oM^t_{\otau^t_n \land (t+ s )} (\vf ) \big) \b1_{\oA} \Big]
   \= \lmt{k \to \infty} E_\oP \Big[  \big( \oM^t_{\otau^t_n \land (t+ r_k)} (\vf ) \- \oM^t_{\otau^t_n \land (t+ s_k)} (\vf ) \big) \b1_{\oA} \Big]   \= 0 $.

   \fi
   Taking $  (\fs,\fr) \= \Big(\frac{\lceil (s -t) 2^k \rceil}{2^k}, \frac{1+\lceil (r -t) 2^k \rceil}{2^k}\Big)$, $ k \ins \hN $  in \eqref{010922_11}
    and sending   $k \nto \infty$,  we can deduce from
      the continuity of bounded process $\big\{\oM^t_{s \land \otau^t_n } (\vf) \big\}_{s \in [t,\infty)} $
   that   $  E_\oP \Big[ \big( \oM^t_{\otau^t_n \land r} (\vf ) \- \oM^t_{\otau^t_n \land s} (\vf ) \big) \b1_{\oA} \Big]
   \n \= 0 $.
   So $\big\{\oM^t_{s \land \otau^t_n  } (\vf )  \big\}_{s \in [t,\infty)}$   is an $\big(\obF^t,\oP\big)-$martingale.
   By Remark \ref{rem_ocP}, $\oP$ satisfies (D1) and (D2) of Definition  \ref{def_ocP}. 

     \no {\bf 2b)} Let $\oP \ins \ocP^1_{t,\bx} \Cp \ocP^2_t \Cp    \ocP^3_t $.
 There exists   a $[0,\infty]-$valued $ \bF^{W,P_0} -$stopping time $\ddot{\tau}$ on $\O_0$ such that
 $ \oP \nci \big(\osW^t, \oT \- t\big)^{-1} \= \Ga  (\ddot{\tau} ) \= P_0 \nci  ( W, \ddot{\tau} )^{-1} $.
 We still have  \eqref{010922_14} since $\osW^t$ 
 is a Brownian motion under $\oP$ by (D1) of Definition  \ref{def_ocP}.
 Given $D \ins \sB(\O_0   \ti \hT)$, 
 taking $A_0 \= \big( W, \ddot{\tau}\big)^{-1} (D) \ins \cF^{W,P_0}_\infty $ in \eqref{010922_14}  yields that
\beas
 \hspace{2mm}
 \oP \big\{  \big(\osW^t  ,\oT\-t \big) \ins D \big\}
 \=  \oP \nci \big(\osW^t, \oT \- t\big)^{-1} (D)  \= P_0 \nci  ( W, \ddot{\tau}  )^{-1} (D)
 \= \oP \nci \big(\osW^t\big)^{-1} \big( \big( W, \ddot{\tau}\big)^{-1} (D) \big)
   \=   \oP \big\{  \big(  \osW^t   ,   \ddot{\tau} (\osW^t  ) \big) \ins D \big\} .
\eeas
 So the  joint distribution of $ (\osW^t,\oT\-t)$ is the same as that of $ \big(\osW^t,   \ddot{\tau} (\osW^t  )\big)$ under $\oP$.
 In particular,   the $\oP-$law of $\oT  $ is equal  to the $\oP-$law of $ t \+ \ddot{\tau}  (\osW^t) $ and therefore
 $\oP$ satisfies  (D3') of  Remark \ref{rem_ocP} or equivalently (D3)  of Definition  \ref{def_ocP}. 
  \if{0}
 So the  joint $\oP-$distribution of $ (\osW^t,\oT )$ is equal to the joint $\oP-$distribution of $ \big(\osW^t, t \+ \ddot{\tau} (\osW^t  )\big)$.

 Set $\otau \df t\+ \ddot{\tau} (\osW^t)$, which is a   $[t,\infty]-$valued  $  \bF^{\oW^t,\oP}-$stopping time.
 We claim that
  \bea \label{071622_14}
 \oP \big( \oA \Cp  \big\{\oT \ins \cE \big\}\big) \= \oP\big( \oA \Cp \big\{ \otau \ins \cE \big\} \big) , \q \fa \oA \ins \cF^{\oW^t,\oP}_\infty ,~  \fa \cE \ins \sB(\hT) .
\eea
 To see this, we fix $\cE \ins \sB(\hT)$ and define $\L \df \big\{\oA \in \sB_\oP(\oO) \n : \oP \big( \oA \Cp  \big\{\oT \ins \cE \big\} \big)
\= \oP \big( \oA \Cp \big\{ \otau \ins \cE \big\} \big) \big\}$.
 As $\oP      \big\{\oT \ins \cE \big\} \= \oP \{(\osW^t,\oT) \ins \O_0 \ti \cE\}
 \= \oP \{(\osW^t,\otau) \ins \O_0 \ti \cE\}
\= \oP    \big\{ \otau \ins \cE \big\}  $, we see that $\oO \ins \L$ and $\L$ is thus a  Lambda-system.

 For any $(\fs,A) \ins [0,\infty) \ti \sB(\hR^d)$, since $W_\fs \n : \O_0 \mto \hR^d$ is a continuous function,
 one has $ W^{-1}_\fs (A) \ins \sB(\O_0) $.
 Then it holds for any $\{(\fs_i,A_i)\}^N_{i=1} \sb  [0,\infty) \ti \sB(\hR^d)$ that
 \beas
 \oP \Big( \Big(\ccap{i=1}{N} (\osW^t_{\fs_i})^{-1} (A_i)\Big) \Cp \{ \oT \ins \cE   \}\Big)
& \tn \=& \tn  \oP \Big( \Big\{  \osW^t  \ins  \ccap{i=1}{N}   W^{-1}_{\fs_i} (A_i) \Big\} \Cp \{ \oT \ins \cE   \}\Big)
\= \oP \Big\{ \big(\osW^t,\oT\big) \ins  \ccap{i=1}{N}   W^{-1}_{\fs_i} (A_i) \ti \cE  \Big\} \\
& \tn \= & \tn   \oP \Big\{ \big(\osW^t,\otau\big) \ins \ccap{i=1}{N}   W^{-1}_{\fs_i} (A_i) \ti \cE \Big\}
 \=   \oP \Big( \Big(\ccap{i=1}{N} (\osW^t_{\fs_i})^{-1} (A_i)\Big) \Cp \{ \otau  \ins \cE   \}\Big) .
\eeas
So $\L$ contains the Pi-system $\Big\{ \ccap{i=1}{N} (\osW^t_{\fs_i})^{-1} (A_i) \n : (\fs_i,A_i) \ins [0,\infty) \ti \sB(\hR^d), i \= 1,\cds,N \Big\}$, which generates $\cF^{\osW^t}_\infty \= \si\big(\osW^t_\fs;   \fs \ins [0,\infty)\big)
\=  \si\big(\oW^t_{t+\fs};   \fs \ins [0,\infty)\big)  \=   \si\big(\oW^t_s;   s \ins [t,\infty)\big) \= \cF^{\oW^t}_\infty $.
 Dynkin's Pi-Lambda Theorem implies that   $ \cF^{\oW^t,\oP}_\infty \sb \L $,
\if{0}

 Dynkin's Pi-Lambda Theorem implies that  $ \cF^{\oW^t}_\infty   \sb \L $.

 For any $\fs \ins [0,\infty)$, $\cF^{\oW^t}_\fs \cp \sN_\oP\big(\cF^{\oW^t}_\infty\big)$ is another Pi-system included in $\L$,
 Dynkin's Pi-Lambda Theorem  renders that $\cF^{\oW^t,\oP}_\fs \sb \L$.
 As $\ccup{\fs \in [0,\infty)}{} \cF^{\oW^t,\oP}_\fs $ is also a Pi-system,
 applying Dynkin's Pi-Lambda Theorem again yields   $ \cF^{\oW^t,\oP}_\infty \sb \L $.

 \fi
 proving the claim \eqref{071622_14}.

Let $c \ins (t,\infty)$. Applying \eqref{071622_14} with $(\cE,\oA) \= \big( [t,c),  \{c  \< \otau\}  \big)$ yields
$   \oP  \{\oT \< c \<  \otau\}   \=  \oP  \{\otau \< c \<  \otau\}   \= 0 $,
and using \eqref{071622_14} with $(\cE,\oA) \= \big( (c,\infty],  \{  \otau \< c\}  \big)$ renders
$   \oP  \{\oT \> c \>  \otau\}   \=  \oP  \{\otau \> c \>  \otau\}  \= 0 $.
It follows that $\oP\{\oT \< \otau\} \= \oP \Big( \ccup{c \in  (t,\infty) \cap \hQ}{} \{\oT \< c \<  \otau\} \Big) \= 0$
and $\oP\{\oT \> \otau\} \= \oP \Big( \ccup{c \in  (t,\infty) \cap \hQ}{} \{\oT \> c \>  \otau\} \Big) \= 0$.
Hence, $\oP\{\oT   \= t\+ \ddot{\tau} (\osW^t) \} \= 1$. Namely, $\oP$ satisfies  (D3$'$) of  Remark \ref{rem_ocP} or equivalently (D3)  of Definition  \ref{def_ocP}. 
\fi
  \qed

\no {\bf Proof of Lemma \ref{lem_082020_19}:} Let $\{ t_n \}_{n \in \hN} \sb [0,\infty) $  converge to $t \ins [0,\infty) $
 and let $\{ \oP_n \}_{n \in \hN} \sb \fP\big(\oO\big)$ converge to $\oP \ins \fP\big(\oO\big)$ under
 the weak topology of $\fP\big(\oO\big)$
    \big(i.e.,  $\lmt{n \to \infty} \int_{ \oo \in \oO } \phi (\oo) \oP_n (d\oo)
    \= \int_{ \oo \in \oO } \phi (\oo) \oP (d\oo) $ for any bounded continuous function $\phi \n : \oO \mto \hR$\big).
  To see that
$\big\{\ol{\Ga} (t_n,\oP_n) \= \oP_n \nci \big(\osW^{t_n},\oT\-t_n \big)^{-1} \big\}_{n \in \hN}$ converges to
$ \ol{\Ga} (t,\oP) \= \oP \nci \big(\osW^t,\oT\-t \big)^{-1}$
under the weak topology of $\fP\big(\O_0 \ti \hT\big) $,
we let $\psi \n : \O_0 \ti \hT \mto \hR$ be a bounded continuous function and   show that
$ \lmt{n \to \infty} \int_{(\o_0,\l) \in  \O_0  \times \hT} \psi (\o_0,\l ) \big(\oP_n \nci (\osW^{t_n},\oT\-t_n)^{-1}\big) \big(d (\o_0,\l)\big) \= \int_{(\o_0,\l) \in  \O_0  \times \hT} \psi (\o_0,\l ) \big(\oP  \nci (\osW^t,\oT\-t)^{-1}\big) \big(d (\o_0,\l)\big) $.

Set $\|\psi\|_\infty \df \Sup{(\o_0,\l) \in \O_0  \times \hT} \big| \psi (\o_0,\l) \big|$ and let $\e \ins (0,1)$. 
Since the weakly convergent sequence $\{\oP_n\}_{n \in \hN}$ 
is relatively compact in $ \fP\big(\oO\big) $,
Prohorov's Theorem yields that $ \{\oP_n\}_{n \in \hN} $ is tight, i.e.,
$\dis \Sup{n \in \hN}\oP_n \big(\ol{\cK}^c_\e\big) \ls \frac{\e}{ 4 \|\psi\|_\infty}  $ for some compact subset $\ol{\cK}_\e$ of $\oO$.

 \if{0}
Let $t,t' \ins [0,\infty)$, $\oo \= \big(\o_0,\omX,\ft\big) \ins \oO$ and $\oo' \= \big(\o'_0,\omX',\ft'\big) \ins \oO$.
When $(t',\oo') \nto (t,\oo) $,  Lemma \ref{lem_010922} shows that
\beas
\Rho{\O_0} \big( \osW^{t'} (\oo'), \osW^t (\oo) \big) \= \Rho{\O_0} \big( \sW^{t'} \n \big(\oW(\oo')\big), \sW^t \big(\oW(\oo)\big) \big)
\=   \Rho{\O_0} \big( \sW^{t'} \n  (\o'_0) , \sW^t  (\o_0)   \big) \nto 0 .
\eeas
 \fi

 The topology of locally uniform convergence on $\O_0$ implies that   $(s,\oo ) \mto \osW^s (\oo )$ is a continuous mapping from $[0,\infty) \ti \oO  $ to $\O_0$ and  $ \ol{\Phi}(s,\oo) \df \big(\osW^s (\oo),\oT(\oo)\-s\big) $ 
 is thus  a continuous mapping from $[0,\infty) \ti \oO$ to $\O_0 \ti \hT$. 
 There exists   $\d 
 \ins (0,1)$ such that
 $ 
 \big|\psi \nci \ol{\Phi} (s,\oo) \- \psi \nci \ol{\Phi} (s',\oo') \big| \< \e / 4  $ for any $ (s,\oo) , (s',\oo') \ins [0,t\+1] \ti \ol{\cK}_\e $   with $ |s\-s'| \vee \Rho{\oO} (\oo,\oo') \< \d  $.
 And one  can find   $N  
 \ins \hN$ such that
 $ 
 \big| \int_{ \oo \in \oO } \psi \nci \ol{\Phi} (t,\oo) \oP_n (d\oo)
\- \int_{ \oo \in \oO } \psi \nci \ol{\Phi} (t,\oo) \oP (d\oo) \big| \< \frac{\e}{4} $  
 and $|t_n \- t| \< \d$ for any $n \gs N $.

  For any $n \gs N$, we can deduce that  
 \beas
 \hspace{1.3cm}
 && \hspace{-1.2cm} \Big| \int_{(\o_0,\l) \in  \O_0  \times \hT} \psi (\o_0,\l ) \big(\oP_n \nci (\osW^{t_n},\oT\-t_n)^{-1}\big) \big(d (\o_0,\l)\big)   - \int_{(\o_0,\l) \in  \O_0  \times \hT} \psi (\o_0,\l ) \big(\oP  \nci (\osW^t,\oT\-t)^{-1}\big) \big(d (\o_0,\l)\big) \Big| \\
 &&  \hspace{-0.7cm}  \ls \Big| \int_{\oo \in  \oO} \Big( \psi \big( \ol{\Phi} (t_n,\oo) \big)
 \-   \psi \big( \ol{\Phi} (t,\oo) \big) \Big) \oP_n  (d\oo) \Big|
 \+ \Big| \int_{\oo \in  \oO} \psi \big( \ol{\Phi} (t,\oo) \big) \oP_n (d\oo)
 \- \int_{\oo \in  \oO} \psi \big( \ol{\Phi} (t,\oo) \big) \oP  (d\oo) \Big| \\
 &&  \hspace{-0.7cm}   \<   \int_{\oo \in \ol{\cK}_\e} \Big| \psi \big( \ol{\Phi} (t_n,\oo) \big)
 \-   \psi \big( \ol{\Phi} (t,\oo) \big) \Big| \oP_n  (d\oo)
 \+   \int_{\oo \in  \ol{\cK}^c_\e} \Big| \psi \big( \ol{\Phi} (t_n,\oo) \big) \Big| \oP_n (d\oo)
 \+ \int_{\oo \in  \ol{\cK}^c_\e} \Big| \psi \big( \ol{\Phi} (t,\oo) \big) \Big| \oP_n  (d\oo) \+ \e / 4  \\
 &&  \hspace{-0.7cm}  \ls \frac{\e}{4} \oP_n (\ol{\cK}_\e) \+ 2 \|\psi\|_\infty \oP_n \big(\ol{\cK}^c_\e\big) \+ \e/4 \ls \e.
 \hspace{10.5cm} \hb{\qed}
 \eeas

\no {\bf Proof of Proposition \ref{prop_graph_ocP}:}
   According to Proposition \ref{prop_Ptx_char}, $\big\lan\n\big\lan  \ocP  \big\ran\n\big\ran $
 is the intersection of   $ \big\lan\n\big\lan  \ocP  \big\ran\n\big\ran_1 \df \big\{\big(t,\bx,\oP \big)  \ins [0,\infty) \ti \OmX \ti \fP\big(\oO\big) \n :  \oP \ins \ocP^1_{t,\bx} \big\}$   and
 $ \big\lan\n\big\lan  \ocP  \big\ran\n\big\ran_i \df \big\{\big(t,\bx,\oP \big)  \ins [0,\infty) \ti \OmX \ti \fP\big(\oO\big) \n :  \oP \ins \ocP^i_t \big\} $ for $i \= 2,3$.

  \no {\bf 1)} Since the  function $\fl_2 (t,\omX) \df \omX(t \ld \cd)  $ is continuous in $ (t,\omX) \ins  [0,\infty) \ti \OmX  $,
   the mapping  $  \psi_{\overset{}{X}}  (t, \bx,\oo) \df     \b1_{\{ \fl_2(t,\oX(\oo) )  - \fl_2(t,\bx )   = 0 \}}
   $, $   (t,\bx, \oo) \ins [0,\infty) \ti \OmX \ti  \oO  $  is $\sB[0,\infty) \oti   \sB(\OmX) \oti \sB(\oO)  -$measurable.
      Lemma \ref{lem_A1} implies that 
 $ \Psi_{\n X}  (t,\bx,\oP) \df \int_{\oo \in \oO} \psi_{\overset{}{X}}  (t, \bx,\oo)  \oP(d \, \oo)
   \= \oP  \big\{  \oX_s \= \bx(s),   \fa s \ins [0,t]   \big\}  $, $   (t,\bx, \oP) \ins [0,\infty) \ti \OmX  \ti   \fP\big(\oO\big) $
  is $ \sB[0,\infty) \oti \sB(\OmX)   \oti \sB\big(\fP\big(\oO\big)\big) 
  - $measurable. So    $  \big\lan\n\big\lan  \ocP  \big\ran\n\big\ran_1 \=
  \big\{ (t,\bx,\oP) \ins [0,\infty) \ti \OmX \ti \fP\big(\oO\big) \n : 
  \Psi_{\n X}  (t,\bx,\oP) \= 1  \big\} \ins  \sB[0,\infty) \oti \sB(\OmX)   \oti \sB\big(\fP\big(\oO\big)\big)$.

  \no {\bf 2)}
  Since 
  $W(s,\o_0) \df \o_0(s)$ is continuous in $ (s,\o_0) \ins  [0,\infty) \ti \O_0 $
  and 
    $W^X(s,\omX) \df \omX(s)$ is continuous in $ (s,\omX) \ins  [0,\infty) \ti \OmX $,
  the function $\Xi(t,\fs,\o_0,\omX) \df \big( W(t \+ \fs,\o_0) \- W(t,\o_0), W^X(t \+ \fs,\omX)  \big) $   
 is continuous in $   (t,\fs,\o_0,\omX) \ins [0,\infty) \ti [0,\infty) \ti \O_0 \ti \OmX $. 

 Let $ (\vf,n)  \ins   \fC(\hR^{d+l}) \ti  \hN$.
  The measurability of   functions $b,\si, \fl_2,\Xi$ imply that the mapping
 \beas
   \cH_\vf (t,\fs,r,\o_0,\omX) \df \b1_{\{t \le r \le t + \fs \}} \Big\{ \ol{b} \big( r, \fl_2( r,\omX)\big) \n \cd \n D \vf \big( \Xi(t,r\-t,\o_0,\omX) \big) \+
 \frac12 \ol{\si} \ol{\si}^T \big( r, \fl_2( r,\omX)\big) \n : \n D^2 \vf \big( \Xi(t,r\-t,\o_0,\omX) \big) \Big\} ,
 \eeas
  \if{0}

 \beas
   \cH_\vf (t,\fs,r,\o_0,\omX) \df \b1_{\{t \le r \le t + \fs \}} \Big\{ \ol{b} \big( r, \fl_2( r,\omX)\big) \n \cd \n D \vf \big( \Xi(t,(r\-t)^+,\o_0,\omX) \big) \+
 \frac12 \ol{\si} \ol{\si}^T \big( r, \fl_2( r,\omX)\big) \n : \n D^2 \vf \big( \Xi(t,(r\-t)^+,\o_0,\omX) \big) \Big\}
 \eeas

  \fi
 $  (t,\fs,r,\o_0,\omX) \ins [0,\infty) \ti [0,\infty) \ti (0,\infty) \ti \O_0 \ti \OmX$  is Borel-measurable.
 \if{0}

  Given $(t,s,\o_0,\omX) \ins [0,\infty) \ti [0,\infty) \ti \O_0 \ti \OmX$, 
   the continuity of
   $ \Xi(t,r,\o_0,\omX)$  in $r  $ shows that
  \beas
   c_{t,\fs} \df \Sup{r \in [t,t +\fs]} \Big( \big| D \vf \big(\Xi(t,r,\o_0,\omX) \big)  \big| \+ \big| D^2 \vf \big(\Xi(t,r,\o_0,\omX) \big)  \big|  \Big) \< \infty .
  \eeas
  An analogy to \eqref{122721_11} renders that
   \beas
      \hspace{-0.3cm} \int_0^\infty \big|\cH_\vf (t,s,r,\o_0,\omX)\big| dr
         \ls    \int_t^{t+\fs}   \big| \ol{b} \big( r,\fl_2( r,\omX) \big) \big|\big| D \vf \big( \Xi(t,r,\o_0,\omX) \big)\big| dr    \+   \frac12 \int_t^{t+\fs} \big| \ol{\si}   \big( r, \fl_2( r,\omX) \big)\big|^2    \big| D^2 \vf \big( \Xi(t,r,\o_0,\omX) \big) \big|  dr \q \nonumber \\
      \eeas

      \vspace{-1.18cm}
      \bea
  & &  \ls    c_{t,s} \bigg(  \int_t^{t+\fs}  \Big(   \k( r)    \|\fl_2( r,\omX) \|_r  \+    |b( r,\bz)|  \Big) dr
  \+ \frac12 \int_t^{t+\fs}  \Big(d   \+ 2 \big( \k( r)  \|\fl_2( r,\omX) \|_r \big)^2   \+ 2  |\si( r,\bz)|^2 \Big) dr \bigg)  \nonumber  \\
  & &  \ls     c_{t,s} \Big(  \big( \k( t\+\fs) \| \omX  \|_{t+\fs}   \big) \fs \+ \int_t^{t+\fs}    |b( r,\bz)|   dr
  \+    \big[ d /2   \+   \big( \k( t\+\fs) \| \omX  \|_{t+\fs} \big)^2 \big] \fs \+  \int_t^{t+\fs}    |\si( r,\bz)|^2   dr  \Big)    \< \infty .
  \label{080421_11}
 \eea
    So $\cI_\vf (t,s,\o_0,\omX) \df \int_0^\infty \cH_\vf (t,s,r,\o_0,\omX) dr$ is a well-defined  Borel-measurable function
    on $[0,\infty) \ti [0,\infty) \ti \O_0 \ti \OmX$.

 As the mapping $\cH_\vf (t,s,r,\o_0,\omX)$ is   Borel-measurable in $(t,s,r,\o_0,\omX) \ins [0,\infty) \ti [0,\infty) \ti (0,\infty) \ti \O_0 \ti \OmX$, its integration with respect to $r \ins (0,\infty)$,
$\cI_\vf (t,s,\o_0,\omX) \= \int_0^\infty \cH_\vf (t,s,r,\o_0,\omX) dr $, is then Borel-measurable in  $   (t,s,\o_0,\omX) \ins [0,\infty) \ti [0,\infty) \ti \O_0 \ti \OmX$.

  and 
   $\cI_\vf (t,\fs,\o_0,\omX) \df \int_0^\infty \cH_\vf (t,\fs,r,\o_0,\omX) dr $  is thus  Borel-measurable
     in  $   (t,\fs,\o_0,\omX) \ins
     [0,\infty) \ti [0,\infty) \ti \O_0 \ti \OmX$ (By an analogy to \eqref{122721_11},   $ |\cI_\vf (t,\fs,\o_0,\omX)| \ls \int_0^\infty \big|\cH_\vf (t,s,r,\o_0,\omX)\big| dr
     \< \infty$).

 \fi
  For each $   (t,\fs,\o_0,\omX) \ins  [0,\infty) \ti [0,\infty) \ti \O_0 \ti \OmX$,
  an analogy to \eqref{122721_11} renders that   $  \int_0^\infty \big|\cH_\vf (t,\fs,r,\o_0,\omX)\big| dr
     \< \infty$. So  $\cI_\vf (t,\fs,\o_0,\omX) \df \int_0^\infty \cH_\vf (t,\fs,r,\o_0,\omX) dr $
     is a real-valued,  Borel-measurable mapping  on $[0,\infty) \ti [0,\infty) \ti \O_0 \ti \OmX$.
 Define $\sT_n (t,\o_0,\omX) \df  \inf\big\{ \fs \ins [0,\infty) \n : |\Xi(t, \fs,\o_0,\omX)|    \gs n  \big\}$,
 $   (t, \o_0,\omX) \ins [0,\infty)   \ti \O_0 \ti \OmX $, which is also Borel-measurable since for any $a \ins [0,\infty)$,
  \beas
   && \hspace{-1.2cm}   \big\{(t, \o_0,\omX) \ins [0,\infty)   \ti \O_0 \ti \OmX \n : \sT_n (t,\o_0,\omX) \> a \big\}
    \= \Big\{(t, \o_0,\omX) \ins [0,\infty)   \ti \O_0 \ti \OmX  \n : \Sup{ a' \in  [0,a]} |\Xi(t, a',\o_0,\omX)|   \< n   \Big\}   \\
   && \hspace{-0.5cm} 
    \= \Big(\ccup{k \in \hN}{} \, \ccap{q \in \hQ \cap [0,a] }{} \big\{(t, \o_0,\omX) \ins [0,\infty)   \ti \O_0 \ti \OmX  \n : |\Xi(t, q,\o_0,\omX)| \ls n  \- 1/k  \big\}\Big)
    \ins \sB[0,\infty) \oti \sB(\O_0) \oti \sB(\OmX)  .
  \eeas
 For any $\fs \ins [0,\infty)$,   since the path-valued random variables $(\oW,\oX)$ on
$\oO$ are $\sB(\O_0) \oti \sB(\OmX)  -$measurable, we can derive from
     the Borel measurability of   $\cI_\vf$ and $ \sT_n $ that the mapping
    \bea
 \oM^{\vf,n}_\fs (t,\oo) \df (\vf \nci \Xi \- \cI_\vf) \big( t, \sT_n \big(t,\oW(\oo),\oX(\oo)\big) \ld n \ld \fs ,\oW(\oo),\oX(\oo) \big)
 \= \big( \oM^t (\vf) \big) \big(\otau^t_n (\oo) \ld (t\+\fs) ,\oo\big) ,
    \label{Sep21_02}
 \eea
 $  (t,\oo) \ins [0,\infty) \ti \oO$  is    $\sB[0,\infty) \oti \sB(\oO) -$measurable,
 where we used the fact $  \otau^t_n (\oo) \= t\+\sT_n \big(t,\oW(\oo),\oX(\oo)\big) \ld  n   $.

 Let    $\th \df \big(  \vf, n, (\fs,\fr) , \{(s_i,\cO_i )\}^k_{i=1} \big) \ins
   \fC(\hR^{d+l}) \ti \hN  \ti \hQ^{2,<}_+ \ti \wh{\sO} (\hR^{d+l}) $.
  Since
  $ \ol{\ff}_\th (t,\oo) \df \big( \oM^{\vf,n}_\fr (t,\oo) \- \oM^{\vf,n}_\fs (t,\oo) \big) \ti   \underset{i=1}{\overset{k}{\prod}} \\  \b1_{\{ \Xi  ( t,s_i \land \fs,\oW(\oo),\oX(\oo)  )   \in \cO_i \}  }   $,
  $  (t,\oo) \ins [0,\infty) \ti   \oO$ is   $\sB[0,\infty) \oti \sB(\oO) -$measurable
   by \eqref{Sep21_02},
  applying Lemma \ref{lem_A1}   yields that the mapping
  $  (t,\oP) \mto \int_{\oo \in \oO} \, \ol{\ff}_\th (t,\oo) \oP(d \, \oo)
    $
  is $ \sB[0,\infty)   \oti \sB\big(\fP\big(\oO\big)\big) 
  - $measurable and the set
  $ \Big\{ (t,\bx,\oP) \ins [0,\infty) \ti \OmX \ti \fP\big(\oO\big) \n :   E_\oP \Big[  \big( \oM^t_{\otau^t_n \land (t+\fr)} (\vf ) \- \oM^t_{\otau^t_n \land (t+\fs)} (\vf ) \big)  \underset{i=1}{\overset{k}{\prod}}   \b1_{ \{(\oW^t_{t+s_i \land \fs }, \oX_{t+s_i \land \fs }) \in \cO_i \}}   \Big]   \= 0 \Big\} $
  is thus Borel-measurable.
  Letting $\th$ run through the countable collection $\fC(\hR^{d+l}) \ti \hN  \ti \hQ^{2,<}_+ \ti \wh{\sO} (\hR^{d+l})$ shows
  $ \big\lan\n\big\lan  \ocP  \big\ran\n\big\ran_2 \ins  \sB[0,\infty) \oti \sB(\OmX)   \oti \sB\big(\fP\big(\oO\big)\big)$.

  \no {\bf 3)}
  \if{0}

   Since $\O_0$ is a Polish space,
 Proposition 7.13 of \cite{Bertsekas_Shreve_1978} shows that   $ \O_0   \ti \hT$  with the product topology
 is   a Borel space. According to Corollary 7.25.1 of \cite{Bertsekas_Shreve_1978},  
   $ \fP \big(\O_0   \ti \hT\big) $ is also a Borel space.

   \fi
 We know from Lemma \ref{lem_082020_15} and   Lemma \ref{lem_082020_17}   that
 the mapping $\Ga \n :   \fS \ni \tau  \mto P_0 \nci  ( W,  \tau    )^{-1} \ins \fP \big(\O_0 \ti \hT\big)$
  is a continuous injection from the Polish space $   \fS $ to 
 $   \fP \big(\O_0   \ti \hT\big)$ and
  the image $\Ga(  \fS)$ is thus  a Lusin subset of $ \fP \big(\O_0   \ti \hT\big) $.
  According  to Theorem A.6 of \cite{Takesaki_1979}, $\Ga(  \fS)$ is even a Borel subset of the Borel space $ \fP \big(\O_0   \ti \hT\big) $.
  Then Lemma \ref{lem_082020_19} implies
  $ \big\lan\n\big\lan  \ocP  \big\ran\n\big\ran_3
  \= \big\{\big(t,\bx,\oP \big)  \ins [0,\infty) \ti \OmX \ti \fP\big(\oO\big) \n : \ol{\Ga} (t,   \oP)  \ins \Ga(  \fS) \big\}
     \ins  \sB[0,\infty) \oti \sB(\OmX)   \oti \sB\big(\fP\big(\oO\big)\big) $.
   Totally, $ \big\lan\n\big\lan  \ocP  \big\ran\n\big\ran  \= \big\lan\n\big\lan  \ocP  \big\ran\n\big\ran_1 \Cp \big\lan\n\big\lan  \ocP  \big\ran\n\big\ran_2 \Cp \big\lan\n\big\lan  \ocP  \big\ran\n\big\ran_3 $ is   a Borel subset of $[0,\infty) \ti \OmX \ti \fP\big(\oO\big)$. \qed

\no {\bf Proof of Corollary  \ref{cor_graph_ocP}: 1)}
  Let $i \ins \hN$.
  By the measurability of functions $g_i$ and $\fl_2$ \big(defined in \eqref{020723_11}\big),
  the mapping  $  \fg_i  (t,s,r,\omX) \df   \b1_{\{  t \land s \le r \le s\}}  g_i \big(r, \fl_2(r,\omX)\big) $
   is Borel-measurable in $   (t,s,r,\omX) \ins [0,\infty) \ti [0,\infty) \ti (0,\infty) \ti \OmX  $.
  It follows that
   \bea \label{Sep21_04}
     \ol{\cI}_{g_i}  (t, \oo) \df   \int_0^\infty \n  \fg_i  \big(t,\oT(\oo),r,\oX(\oo)\big) dr
  \= \int_{ \oT(\oo) \land t }^{\oT(\oo)  }   g_i \big(r,\oX_{r \land \cd} (\oo)\big) dr  \, , \q     (t, \oo) \ins [0,\infty)   \ti \oO
  \eea
  is $\sB[0,\infty) \oti \sB(\oO) -$measurable.
  Lemma \ref{lem_A1}  implies that
 $ \oPhi_{g_i}  (t, \oP) \df \n \int_{\oo \in \oO}    \ol{\cI}_{g_i}   (t, \oo) \oP(d \, \oo)
   \= E_\oP \big[ \n \int_{\oT \land t}^{\oT  } \n   g_i    \big(r,\oX_{r \land \cd}  \big) dr \big]   $,
 $  (t, \oP) \ins [0,\infty)   \ti \fP\big(\oO\big) $
  is $ \sB[0,\infty)    \oti \sB\big(\fP\big(\oO\big)\big)   - $measurable.
  \if{0}

    It follows that
   \bea \label{Sep21_04}
     \ol{\cI}^\pm_{g_i}  (t, \oo) \df   \int_0^\infty \n  \fg^\pm_i  \big(t,\oT(\oo),r,\oX(\oo)\big) dr
  \= \int_{ \oT(\oo) \land t }^{\oT(\oo)  }   g^\pm_i \big(r,\oX_{r \land \cd} (\oo)\big) dr  \, , \q     (t, \oo) \ins [0,\infty)   \ti \oO
  \eea
  is $\sB[0,\infty) \oti \sB(\oO) -$measurable.
  Lemma \ref{lem_A1}  implies that
 $ \oPhi_{g_i}  (t, \oP) \df \n \int_{\oo \in \oO}    \ol{\cI}^+_{g_i}   (t, \oo) \oP(d \, \oo) \- \int_{\oo \in \oO}    \ol{\cI}^-_{g_i}   (t, \oo) \oP(d \, \oo)
   \= E_\oP \big[ \n \int_{\oT \land t}^{\oT  } \n   g_i    \big(r,\oX_{r \land \cd}  \big) dr \big]   $,
 $  (t, \oP) \ins [0,\infty)   \ti \fP\big(\oO\big) $
  is $ \sB[0,\infty)    \oti \sB\big(\fP\big(\oO\big)\big)   - $measurable.

  \fi
  Similarly,
  $  \oPhi_{h_i}  (t, \oP) \df 
   E_\oP \big[ \int_{\oT \land t}^\oT   h_i   \big(r,\oX_{r \land \cd}  \big) dr \big]   $,
  $   (t, \oP) \ins [0,\infty)   \ti \fP\big(\oO\big) $
  is $ \sB[0,\infty)    \oti \sB\big(\fP\big(\oO\big)\big)   - $measurable.  Then the set
  \beas
  \q  \sD  \df \big\{ (t,\bx,y,z,\oP) \ins [0,\infty) \ti \OmX \ti \Re \ti \Re \ti  \fP\big(\oO\big) \n :  \oPhi_{g_i}  (t, \oP)   \ls y_i ,\,
   \oPhi_{h_i}  (t, \oP)     \= z_i,\,\fa i \ins \hN \big\}
  \eeas
  is Borel-measurable.  
  Since   $\big\lan\n\big\lan  \ocP  \big\ran\n\big\ran \ins  \sB[0,\infty) \oti \sB \big(\OmX\big) \oti \sB \big( \fP (\oO) \big)$
  by Proposition \ref{prop_Ptx_char},
  using  the   projection   $\ol{\Pi}_1 (t, \bx,y,z,\oP) \\ \df   \big(t,\bx ,\oP\big)$ yields that
 $\gP   \=   \big\{ (t,\bx, y,z, \oP ) \ins [0,\infty) \ti \OmX \ti \Re \ti \Re \ti \fP\big(\oO\big)  \n :    \oP \ins \ocP_{t,\bx} ; \;   E_\oP \big[ \int_{\oT \land t}^\oT   g_i\big(r,\oX_{r \land \cd}  \big) dr \big]    \ls y_i ,\, E_\oP \big[ \int_{\oT \land t}^\oT   h_i\big(r,   \oX_{r \land \cd}  \big) dr \big]    \= z_i , \, \fa i \ins \hN \big\}  \=     \ol{\Pi}^{-1}_1 \big( \big\lan\n\big\lan  \ocP  \big\ran\n\big\ran \big)     \Cp \sD$
 is a Borel subset of $\oD \ti \fP\big(\oO\big)$.

   \no {\bf 2)}
   \if{0}
   Since the  function $\fl_1 (s,\o_0) \df \o_0(s \ld \cd)  $ is continuous in $ (s,\o_0) \ins  [0,\infty) \ti \O_0  $,   the mapping
 $ \phi_W (t,\bw, \oo) \df  \b1_{\{ \fl_1(t,\oW(\oo) )  - \fl_1(t,\bw )   = 0 \}} $, $ \fa (t,\bw, \oo) \ins [0,\infty) \ti \O_0 \ti  \oO $
 is $\sB[0,\infty) \oti \sB(\O_0) \oti   \sB(\oO)  -$measurable.
 Applying Lemma \ref{lem_A1} again shows that
 \fi
 Similarly  to   $\Psi_{\n X}  (t,\bx,\oP)$,  defined in 
 the proof of Proposition \ref{prop_graph_ocP},
 the mapping $ \Psi_{\n W} (t,\bw, \oP)   \df 
     \oP \big\{  \oW_{\n s}  \=  \bw(s)  $, $\fa s \ins [0,t]  \big\}  $ is Borel-measurable in
 $  (t,\bw, \oP) \ins [0,\infty) \ti \O_0  \ti   \fP\big(\oO\big)$.
  By the projections
   $\ol{\Pi}_2 (t,\bw,\bx,y,z,\oP) \df   \big(t,\bx ,\oP\big)$,
 $\ol{\Pi}_3 (t,\bw,\bx,y,z,\oP) \df   \big(t,\bx,y,z,\oP\big)$ and
 $\ol{\Pi}_4 (t,\bw,\bx,y,z,\oP) \df (t,\bw,  \oP) $,
 we can deduce that
 \if{0}
 \beas
&& \hspace{-1.7cm} \gcP   \= \Big\{ (t, \bw,\bx,y,z, \oP ) \ins [0,\infty) \ti \O_0  \ti \OmX \ti \Re \ti \Re \ti \fP\big(\oO\big)  \n :    \oP \ins \ocP_{t,\bx}  \Big\} \\
&& \hspace{-1cm}  \q \cap \, \Big\{ (t, \bw,\bx,y,z, \oP ) \ins [0,\infty) \ti \O_0  \ti \OmX \ti \Re \ti \Re \ti \fP\big(\oO\big)  \n :   \oP  \big\{  \oW_{\n s}  \=   \bw(s)   ,
    \fa s \ins [0,t]   \big\} \= 1  \Big\} \\
&& \hspace{-1cm}  \q \cap \, \bigg\{ (t, \bw,\bx,y,z, \oP ) \ins [0,\infty) \ti \O_0  \ti \OmX \ti \Re \ti \Re \ti \fP\big(\oO\big)  \n :     E_\oP \Big[ \int_t^{\oT  }   g_i\big(r,\oX_{r \land \cd}  \big) dr \Big]    \ls y_i ,\, E_\oP \Big[ \int_t^{\oT  }   h_i\big(r,\oX_{r \land \cd}  \big) dr \Big]    \= z_i , \, \fa i \ins \hN \bigg\} \\
&& \hspace{-1cm}  \q  \= \ol{\Pi}^{-1}_2 \big( \big\lan\n\big\lan  \ocP  \big\ran\n\big\ran \big)
  \Cp \ol{\Pi}^{-1}_4 \big( \Phi^{-1}_W (1) \big) \Cp \ol{\Pi}^{-1}_3 (\sD)
\eeas
\fi
$  \gcP  
  \= \ol{\Pi}^{-1}_2 \big( \big\lan\n\big\lan  \ocP  \big\ran\n\big\ran \big)
   \Cp \ol{\Pi}^{-1}_3 (\sD) \Cp \ol{\Pi}^{-1}_4 \big( \Psi^{-1}_W (1) \big) $
is a Borel subset of $ \ocD \ti \fP\big(\oO\big) $. \qed

\no {\bf Proof of Theorem \ref{thm_V_usa}:}
 Analogous  to  \eqref{Sep21_04},
 $  \ol{\cI}_f (t, \oo) \df   \int_{\oT(\oo) \land t}^{\oT(\oo)}  f  \big(r,\oX_{r \land \cd} (\oo)\big) dr $
  is 
 Borel-measurable in $   (t,\oo) \ins [0,\infty) \ti   \oO  $.
 Since the measurability of functions $\pi$ and $\fl_2$ \big(defined in \eqref{020723_11}\big)     implies that the mapping
 $ (s,\oo) \mto \pi \big(s,\fl_2(s,\oX(\oo))\big) \= \pi \big(s,\oX_{s \land \cd}(\oo)\big)$  is $\sB(0,\infty) \oti \sB(\oO)-$measurable,
 the random variable
 $ 
    \ol{\phi}_\pi (\oo) \df   \b1_{\{\oT (\oo) < \infty\}}   \pi \big(\oT (\oo), \\ \oX  (\oT (\oo) \ld \cd, \oo) \big)  $, $ \oo \ins \oO $
 is   $\sB(\oO)
       -$measurable.      Lemma \ref{lem_A1} shows that
  $   \ol{\sV}   (t,\oP) \df \int_{\oo \in \oO}   \big( \ol{\cI}_f (t, \oo) \+ \ol{\phi}_\pi (\oo) \big) \oP(d \, \oo)
    \= E_\oP \big[ \, \oR  (t) \big]
  $,
 $ (t,\oP) \ins [0,\infty)  \ti \fP\big(\oO\big)$ is $ \sB[0,\infty)   \oti  \sB\big(\fP\big(\oO\big)\big) 
 - $measurable.
 \if{0}

    Lemma \ref{lem_A1} shows that
  $   \ol{\sV}   (t,\oP) \df \int_{\oo \in \oO}   \big( \ol{\cI}^+_f (t, \oo) \+ \ol{\phi}_\pi (\oo) \big) \oP(d \, \oo)
  \- \int_{\oo \in \oO}    \ol{\cI}^-_f (t, \oo)   \oP(d \, \oo)
   \= E_\oP \Big[ \int_t^{\oT  }  f  \big(r,\oX_{r \land \cd}  \big) dr
   \+   \b1_{\{\oT   < \infty\}}   \pi  \big(\oT  ,\oX_{\oT \land \cd}  \big) \Big] $,
 $ (t,\oP) \ins [0,\infty)  \ti \fP\big(\oO\big)$ is $ \sB[0,\infty)   \oti  \sB\big(\fP\big(\oO\big)\big) 
 - $measurable.

 \fi
 Then  Corollary  \ref{cor_graph_ocP} and  Proposition 7.47 of \cite{Bertsekas_Shreve_1978} yield that
 $  \oV (t,\bx,y,z) \=  \Sup{\oP \in \ocP_{t,\bx}(y,z)}  \ol{\sV}(t,\oP) \= \Sup{(t,\bx,y,z,\oP) \in [[\ocP]]}  \ol{\sV}(t,\oP) $
 is  \usa ~ on $  \oD $ and
   $\oV(t,\bw,\bx,y,z)
   \=  \Sup{(t,\bw,\bx,y,z,\oP) \in \{\n\{\ocP\}\n\}}  \ol{\sV}(t,\oP)  $
   is  \usa ~ on $  \ocD$. \qed

 \no {\bf Proof of Proposition \ref{prop_flow}: } Let us set $t_\oo \df \oga(\oo) \gs t $ for any $\oo \ins \oO$.

 \no {\bf  1)}   We first demonstrate that for   $\oP-$a.s. $ \oo \ins \oO $,  $ \oP^t_{\oga,\oo} $ belongs to $ \ocP^1_{t_\oo,\oX_{\oga \land \cd}  (\oo)} \Cp \ocP^2_{t_\oo}$ and   thus satisfies \(D1\) and \(D2\) in Definition  \ref{def_ocP}  of $\ocP_{t_\oo,\oX_{\oga \land \cd}  (\oo)}$ according to Part \(2a\) of the proof of Proposition \ref{prop_Ptx_char}.

 \no {\bf  1a)}   By (D2) in Definition  \ref{def_ocP} of $ \ocP_{t,\bx} $,
   $\ocN_{\n X} \df \big\{ \oo \ins \oO \n :  \oX_s (\oo)  \nne  \osX^{t,\bx}_s  (\oo)    \hb{ for some } s \ins [0,\infty) \big\}
       \ins \sN_\oP \big(\ocF^t_{\n \infty}\big)$.  
 As $\big\{\osX^{t,\bx}_s \big\}_{s \in [t,\infty)} $ is  an $  \bF^{\oW^t,\oP}  -$adapted continuous   process,
   one can construct  an  $\hR^l-$valued  $  \bF^{\oW^t }  -$predictable  process  $  \big\{ \oK^t_s \big\}_{s \in [t,\infty)}$
   such that $ \ocN_{\n K} \df  \big\{ \oo \ins \oO \n :   \oK^t_s  (\oo)   \nne \osX^{t,\bx}_s  (\oo)   $ for some $s \ins [t,\infty)  \big\}  \ins  \sN_\oP \big(\cF^{\oW^t}_\infty\big)$ \big(see e.g.  Lemma 2.4 of \cite{STZ_2011a}\big).
   Since   $\Ktgo  \df \ccap{r \in  \hQ \cap (t,\infty)}{} \big\{\oo' \ins \oO \n : \oK^t_{\oga \land r} (\oo') \= \oK^t_{\oga \land r} (\oo) \big\} $
    is an $ \cF^{\oW^t}_\oga -$measurable set   including $\oo$,   (R3) in \eqref{R3} shows  that
  $ 
   \oP^t_{\oga,\oo} \big(\Ktgo\big)   \= 1 $, $ \fa \oo \ins \ocN^c_0$. 

   For any $\oo \ins \big(\ocN_{\n X}  \cp \ocN_{\n K}\big)^c$,    we can deduce from \eqref{090520_11} that
   \bea
   \q &    & \hspace{-1.6cm}   \Wtgo  \Cp    \big(\ocN_{\n X} \n \cp \ocN_{\n K}\big)^c  \n  \Cp \Ktgo
   \=   \Wtgo  \Cp   \big(\ocN_{\n X}  \n  \cp \ocN_{\n K}\big)^c  \n  \Cp   \big\{\oo'  \n \ins \oO \n: \oX_{\n s}(\oo') \= \bx(s),   \fa s \ins [0,t]; \,
   \oK^t_{\n \oga(\oo) \land r}   (\oo') \= \oK^t_{\n \oga(\oo) \land r} (\oo),   \fa r \ins  \hQ \Cp (t,\infty)   \big\} \nonumber \\
    & & \=    \Wtgo  \Cp \big(\ocN_{\n X}  \cp \ocN_{\n K}\big)^c \n \Cp   \big\{\oo'  \n \ins \oO \n: \oX_{\n s}(\oo') 
    \= \oX_{\n s}(\oo) ,   \fa s \ins [0,t]; \,
   \oX_{\n \oga(\oo) \land r}   (\oo') \= \oX_{\n \oga(\oo) \land r} (\oo),   \fa r \ins  \hQ \Cp (t,\infty)   \big\} \nonumber \\
  & & \=   \Wtgo  \Cp \big(\ocN_{\n X}  \cp \ocN_{\n K}\big)^c  \n \Cp \big\{\oo'  \n \ins \oO \n:   \oX_r  ( \oo') \= \oX_{\oga \land r} (  \oo), \, \fa r \ins    [0, \oga(\oo) ]   \big\} . \label{091620_21}
\eea
 And     (R2) in \eqref{R2} shows that $\oP^t_{\oga,\oo} \big(\ocN_{\n X} \cp \ocN_{\n K}\big) \=  E_\oP  \big[ \b1_{\ocN_{\n X} \cup \ocN_{\n K}} \big|\cF^{\oW^t}_\oga  \big] (\oo) \= 0$ for all $\oo \ins \oO$ except on a   $ \wh{\cN}_{\n X,K} \ins   \sN_\oP\big(\cF^{\oW^t}_\oga\big)$.

 Set $\ocN_1 \df  \ocN_{\n X} \cp    \ocN_{\n K} \cp \wh{\cN}_{\n X,K} \ins \sN_\oP \big( \ocF^t_{\n \infty} \big) $.
 Given  $\oo \ins \big(\ocN_0  \cp \ocN_1 \big)^c  $,
  taking $\oP(\cd)$ in \eqref{091620_21} and using \eqref{Jan11_03} yield  that
$ \oP^t_{\oga,\oo} \big\{\oo' \ins \oO \n:   \oX_r  ( \oo') \= \oX_{\oga \land r} (  \oo),    \fa r \ins    [0,t_\oo ]   \big\} \= 1 $,
i.e.,  $ \oP^t_{\oga,\oo} \ins \ocP^1_{t_\oo,\oX_{\oga \land \cd}(\oo)} $.

 \ss \no {\bf  1b)}
 For any $\vf \ins \fC(\hR^{d+l})$ and $q \ins \hQ^d$,   define a function $\vf_q (w,x) \df \vf (w\-q,x)$, $(w,x) \ins \hR^{d+l}  $.
We set $\sC \df \{ \vf_q  : \vf \ins \fC(\hR^{d+l}), q \ins \hQ^l\}$, which is a countable sub-collection of $ C^2(\hR^{d+l})$.
 For any $n \ins \hN $, define an  $\obF^t-$stopping time by $ \oz_n (\oo) \df \inf\big\{r \ins [\oga(\oo),\infty) \n : | \oW^t_r (\oo) \- \oW^t_\oga (\oo) |^2  
 \+ |\oX_r (\oo)|^2 \gs n^2  \big\} \ld \big(\oga (\oo)\+n \big)$, $\oo \ins \oO$.
    \if{0}

 In general,  let $\tau$ be an $\bF-$stopping time and let $X$ be an $\bF-$adapted continuous process. Define $\z (\o) \df \inf\big\{r \ins [\tau(\o),\infty) \n: | X_r (\o)|  \gs n \big\}$, $\o \ins \O$. We can show that for any $s \ins [0,\infty)$
 \beas
 \{\z \> s\} \= \ccup{k \in \hN}{} \ccap{r \in \hQ \cap [0,s]}{} \Big(\{r \ls \tau\} \cp \big(\{r \> \tau\} \Cp \{ |X_r| \ls n \- 1/k \}\big)  \Big) \ins \cF_s .
 \eeas
 So $\z $ is also an $\bF-$stopping time.

    \fi

Let  $\th \df \big(  \phi, n,j, (\fs,\fr) , \{(s_i,\cO_i )\}^k_{i=1} \big) \ins
   \sC \ti \hN \ti \hN  \ti \hQ^{2,<}_+ \ti \wh{\sO} (\hR^{d+l})    $.
 Since    $ \big\{ \oM^t_{s \land \otau^t_j } (\phi) \big\}_{ s \in  [t,\infty) } $
 is a bounded $(\obF^t,  \oP)-$martingale by 
applying Proposition \ref{prop_MPF1} with $\big(\O,\cF,P,B,X  \big)\=\big(\oO,\sB(\oO),\oP ,\oW,\oX  \big) $,
 the  optional sampling theorem implies  that
 $  E_\oP  \Big[   \oM^t_{  ( \oga  +\fr) \land \oz_n  \land \otau^t_j} (\phi  )    \Big|\ocF^t_{\n \oga+\fs}  \Big]
  \=  \oM^t_{  ( \oga  +\fs) \land \oz_n \land \otau^t_j } ( \phi  )  $,  $\oP-$a.s.
 Set $ \oxi_\th \df \oM^t_{  ( \oga  +\fr) \land \oz_n  \land \otau^t_j} (\phi  )   \- \oM^t_{  ( \oga  +\fs) \land \oz_n \land \otau^t_j } ( \phi  ) \= \b1_{\{\otau^t_j > \oga \}} \Big( \oM^t_{  ( \oga  +\fr) \land \oz_n  \land \otau^t_j} (\phi  )   \- \oM^t_{  ( \oga  +\fs) \land \oz_n \land \otau^t_j } ( \phi  ) \Big) $
 and set $\oeta_\th \df
 \underset{i=1}{\overset{k}{\prod}}   \b1_{    \{(\oW^t_{ \oga +s_i \land \fs } - \oW^t_\oga,\oX_{ \oga+s_i \land \fs  }) \in \cO_i  \}    } \ins \ocF^t_{\n \oga+\fs} $.
 As $ \cF^{\oW^t}_\oga \sb \ocF^t_{\n \oga} \sb \ocF^t_{\n \oga+\fs} $, the tower property   renders that
 $ E_\oP  \big[  \,  \oxi_\th \oeta_\th \big|\cF^{\oW^t}_\oga  \big]
 \=  E_\oP  \Big[    \oeta_\th  E_\oP  \big[  \,  \oxi_\th   \big|\ocF^t_{\n \oga+\fs}  \big]  \Big|\cF^{\oW^t}_\oga  \Big] \= 0 $, $ \oP-$a.s.
 By (R2) in \eqref{R2} again, there exists an $\ocN_\th \ins \sN_\oP \big(\cF^{\oW^t}_\oga\big)$ such that
 \bea \label{012822_17}
  E_{\oP^t_{\oga,\oo}} \big[ \, \oxi_\th \oeta_\th \big]
   \= E_\oP  \big[  \,  \oxi_\th \oeta_\th \big|\cF^{\oW^t}_\oga  \big] (\oo) \= 0  , \q \fa \oo \ins \ocN^c_\th  .
 \eea

 Define $\ocN_2 \df    \bigcup \big\{\ocN_\th \n :  \th \ins \sC \ti \hN \ti \hN  \ti \hQ^{2,<}_+ \ti \wh{\sO} (\hR^{d+l}) \big\}   \ins \sN_\oP \big(\ocF^t_{\n \infty}\big)$ and fix   $\oo \ins \big( \ocN_0 \cp \ocN_1 \cp \ocN_2 \big)^c $.
 We   let    $\big(\vf,n,(\fs,\fr), \\ \{(s_i,   \cO_i )\}^k_{i=1} \big) \ins \fC(\hR^{d+l}) \ti \hN \ti \hQ^{2,<}_+  \ti  \wh{\sO} (\hR^{d+l}) $
 and let $j \ins \hN$.
 There exists a sequence  $\{q_m\=q_m(\oo)\}_{m \in   \hN} $ of $ \hQ^d$ that converges to $  \oW^t_\oga(\oo)$.

 Let $ m \ins \hN$. We set $\th_m \df \big(  \vf_{q_m}, n,j,   (\fs,   \fr) , \{(s_i,\cO_i )\}^k_{i=1} \big)$ and define $\d^{j,m}_\oo    \df \Sup{|(w,x)| \le j} \Big( \sum^2_{i=0} \big|D^i\vf_{q_m}(w,x) \- D^i \vf \big(w\-\oW^t_\oga(\oo),x\big)\big| \Big)   \= \Sup{|(w,x)| \le j} \Big( \sum^2_{i=0}\big|D^i\vf (w\-q_m,x) \- D^i\vf \big(w\-\oW^t_\oga(\oo),x\big)\big| \Big)$.

  Given $\oo' \ins \Wtgo \Cp \ocN^c_{\n X} \Cp \big\{  \otau^t_j   \> \oga \big\}   $, \eqref{090520_11} implies that
 $\otau^t_j (\oo') \> \oga(\oo') 
 \= t_\oo$ and
 $ 
  \oz_n(\oo') \= \inf\big\{r \ins [t_\oo,\infty) \n : | \oW_r (\oo') \- \oW_{t_\oo} (\oo') |^2
 \+ |\oX_r (\oo')|^2 \gs n^2  \big\} \ld \big(t_\oo\+n \big) \= \otau^{t_\oo}_n (\oo')$. 
 As $ \oW^{t_\oo}_{\n r}(\oo') \= \oW^t_{\n r}(\oo') \- \oW^t_{\n t_\oo}(\oo')  
 \= \oW^t_{\n r}(\oo') \- \oW^t_\oga(\oo) $, $\fa r \ins [t_\oo,\infty)$,  it holds for any   $t_\oo \ls s_1  \ls s_2 \< \infty  $ that
 \beas
\; && \hspace{-0.9cm} \big(\oM^{t_\oo}_{s_2} (\vf) \- \oM^{t_\oo}_{s_1} (\vf)\big)(\oo')   \=   \vf \big(\oW^t_{\n s_2} (\oo') \- \oW^t_{\n \oga} (\oo) , \oX_{\n s_2} (\oo') \big) \- \vf \big( \oW^t_{\n s_1} (\oo') \- \oW^t_{\n \oga} (\oo) , \oX_{\n s_1} (\oo') \big) \\
  &&  \- \n \int_{s_1}^{s_2}  \n  \ol{b}  \big( r, \oX_{r \land \cd} (\oo')  \big) \n \cd \n D \vf \big( \oW^t_{\n r} (\oo') \- \oW^t_\oga (\oo) , \oX_r (\oo') \big) dr
    \-   \frac12 \n \int_{s_1}^{s_2} \n  \ol{\si} \, \ol{\si}^T  \big( r,  \oX_{r \land \cd}(\oo')  \big) \n : \n D^2 \vf  ( \oW^t_{\n r} (\oo') \- \oW^t_\oga (\oo) , \oX_r  (\oo') )   dr .
\eeas
 Since $   \big| \big(\oW^t_r (\oo')  , \oX_r (\oo') \big) \big|  \ls j $ for any $r \ins [t_\oo, \otau^t_j(\oo') ]$,
 an analogy to \eqref{122721_11} shows that for any   $t_\oo \ls s_1  \ls s_2 \ls \otau^t_j(\oo')  $
 \beas
 \hspace{-5mm}
  \Big| \big( \oM^{t_\oo}_{s_2} (\vf) \- \oM^{t_\oo}_{s_1} (\vf) \- \oM^t_{s_2} (\vf_{q_m} )  \+ \oM^t_{s_1} ( \vf_{q_m} ) \big) (\oo') \Big|
 \ls 2 \d^{j,m}_\oo  \+ \d^{j,m}_\oo \n \int_t^{\otau^t_j(\oo')}  \n \Big(   \big| \ol{b}  \big( r, \oX_{r \land \cd} (\oo')  \big) \big|   \+ \frac12   \big| \ol{\si}  \big( r, \oX_{r \land \cd} (\oo')  \big) \big|^2 \Big) dr
  \ls \d^{j,m}_\oo (2   \+ c^j_{t,\bx}) ,
 \eeas
 where $c^j_{t,\bx} \df \big[ d / 2 \+  \k( t\+j)  ( \|\bx\|_t  \+ j    )  \+     \k^2( t\+j)  ( \|\bx\|_t  \+  j     )^2      \big] j \+     \int_t^{t+j}   \big(  |b( r,\bz)| \+ |\si( r,\bz)|^2 \big)  dr \< \infty $.
 Taking $ s_1  \=   \big(( \oga   \+\fs) \ld \oz_n   \ld \otau^t_j\big) (\oo') \=  ( t_\oo  \+\fs) \ld \otau^{t_\oo}_n (\oo') \ld \otau^t_j (\oo')  $ and $ s_2 
 \=   ( t_\oo  \+\fr) \ld \otau^{t_\oo}_n (\oo') \ld \otau^t_j (\oo') $ yields
 $ \Big| \big( \oM^{t_\oo}_{( t_\oo  +\fr) \land \otau^{t_\oo}_n \land \otau^t_j} (\vf)  - \oM^{t_\oo}_{( t_\oo  +\fs) \land \otau^{t_\oo}_n \land \otau^t_j} (\vf)\big) (\oo') \- \oxi_{\th_m}   (\oo') \Big|
 \ls \d^{j,m}_\oo (2   + c^j_{t,\bx}) $.
 As $ \oeta_{\th_m}   (\oo') \=  \underset{i=1}{\overset{k}{\prod}}   \b1_{    \{(\oW^{t_\oo}_{ t_\oo+s_i \land \fs } (\oo'),\oX_{ t_\oo+s_i \land \fs  } (\oo')) \in \cO_i  \}    }    $ by \eqref{090520_11},
  we see from \eqref{Jan11_03}    that
   \beas
   \q E_{\oP^t_{\oga,\oo}}  \Big[ \b1_{\{  \otau^t_j   > \oga \}} \Big| \Big( \oM^{t_\oo}_{  ( t_\oo +\fr) \land \otau^{t_\oo}_n \land \otau^t_j } (\vf )  \- \oM^{t_\oo}_{  ( t_\oo +\fs) \land \otau^{t_\oo}_n \land \otau^t_j } (\vf ) \Big) \underset{i=1}{\overset{k}{\prod}} \,  \b1_{    \{(\oW^{t_\oo}_{ t_\oo+s_i \land \fs },\oX_{ t_\oo+s_i \land \fs  }) \in \cO_i  \}    } \- \oxi_{\th_m} \oeta_{\th_m} \Big|   \Big] \ls \d^{j,m}_\oo (2   \+ c^j_{t,\bx})  .
    \eeas

 The uniform  continuity of $D^i \vf  $'s    over compact sets 
   implies   $ \lmtd{m \to \infty} \d^{j,m}_\oo \= 0$, and one  can  then deduce from  \eqref{012822_17} that
  \bea \label{012922_17}
  \hspace{-0.3cm}
  E_{\oP^t_{\oga,\oo}} \n \Big[ \b1_{\{  \otau^t_j   > \oga \}} \Big( \oM^{t_\oo}_{  ( t_\oo +\fr) \land \otau^{t_\oo}_n \land \otau^t_j } (\vf )  \- \oM^{t_\oo}_{  ( t_\oo +\fs) \land \otau^{t_\oo}_n \land \otau^t_j } (\vf ) \Big) \underset{i=1}{\overset{k}{\prod}} \,  \b1_{    \{(\oW^{t_\oo}_{ t_\oo+s_i \land \fs },\oX_{ t_\oo+s_i \land \fs  }) \in \cO_i  \}    }    \Big]
   \n \=  \n  \lmt{m \to \infty} E_{\oP^t_{\oga,\oo}}  \big[ \,  \oxi_{\th_m} \oeta_{\th_m}  \big]  \n \= 0. \q
  \eea
 Since $  \oP^t_{\oga,\oo} \big\{\oo' \ins \oO \n:   \oX_r  ( \oo') \= \oX_{\oga \land r}  (  \oo),    \fa r \ins    [0,t_\oo ]   \big\} \= 1 $
 by Part (1a),
  applying Proposition \ref{prop_MPF1} with $\big(\O,\cF,P,B,X \big)\=\big(\oO,\sB(\oO),\oP^t_{\oga,\oo},\oW,\oX \big) $
  and   $(t,\bx ) \= \big(t_\oo,\oX_{\oga \land \cd} (\oo) \big)$ renders that
 $  \big\{\oM^{t_\oo}_{  s \land \otau^{t_\oo}_n  } (\vf )\big\}_{s \in [t_\oo,\infty)}  $ is a bounded   process under $\oP^t_{\oga,\oo}$.
 \if{0}
 applying Proposition \ref{prop_MPF1} with $(t,\bx) \= \big(t_\oo,\oX_{\oga \land \cd}(\oo)\big)$ and
 $(\O,\cF,P,B,X) \= \big( \oO,\sB(\oO),\oP^t_{\oga,\oo},\oW,\oX\big)$ renders  that
 $  \big\{\oM^{t_\oo}_{  s \land \otau^{t_\oo}_n  } (\vf )\big\}_{s \in [t_\oo,\infty)}  $ is a bounded $\obF^{t_\oo}-$adapted continuous process under $\oP^t_{\oga,\oo}$.
 \fi
   As $\lmtu{j \to \infty} \otau^t_j (\oo') \= \infty$ for any $\oo' \ins \oO$,
   letting $j \nto \infty$ in \eqref{012922_17}
 and using the  bounded convergence theorem, we obtain  that
 $   E_{\oP^t_{\oga,\oo}} \n \Big[   \Big( \oM^{t_\oo}_{  ( t_\oo +\fr) \land \otau^{t_\oo}_n   } (\vf )  \- \oM^{t_\oo}_{  ( t_\oo +\fs) \land \otau^{t_\oo}_n   } (\vf ) \Big) \underset{i=1}{\overset{k}{\prod}}   \b1_{    \{(\oW^{t_\oo}_{ t_\oo+s_i \land \fs },\oX_{ t_\oo+s_i \land \fs  }) \in \cO_i  \}    }    \Big]    \= 0  $.
 Hence,  $ \oP^t_{\oga,\oo} \ins \ocP^1_{t_\oo,\oX_{\oga \land \cd}    (\oo)} \Cp \ocP^2_{t_\oo}$  for any $\oo \ins \big( \ocN_0 \cp \ocN_1 \cp \ocN_2 \big)^c $.

  \no {\bf  2)}   We next show  that  for $\oP-$a.s.  $ \oo \ins \oO $,  $ \oP^t_{\oga,\oo} $ satisfies \(D3\)
  in  Definition \ref{def_ocP}  of $\ocP_{t_\oo,\oX_{\oga \land \cd}  (\oo)}$.

 By   (D3) in Definition  \ref{def_ocP} of $\ocP_{t,\bx}$,
  there is   a $[t,\infty]-$valued   $ \bF^{W^t,P_0} -$stopping time $\wh{\tau}$    such that
 $ \oP \big\{  \oT  \=   \wh{\tau} (\oW)   \big\} \= 1$.
 \if{0}
 For any $s \ins [t,\infty)$, as $\{\wh{\tau}  \ls s\} \ins \cF^{W^t,P_0}_s $,  applying Lemma \ref{lem_122921_11}   with $t_0 \= t$, $(\O_1, \cF_1, P_1,B^1)   \= \big(\oO ,  \sB(\oO ),  \oP , \oW\big) $, $(\O_2, \cF_2, P_2,B^2) \= \big(\O_0,  \sB(\O_0),  P_0, W\big) $ and $\Phi \= \oW$  implies that  $ \{\wh{\tau} (\oW) \ls s\} \= \oW^{-1} \big( \{\wh{\tau}  \ls s\} \big) \ins \cF^{\oW^t,\oP}_s  $.
 So $ \wh{\tau} (\oW) $ is a $\bF^{\oW^t,\oP}-$stopping time.
 \fi
 Since Lemma \ref{lem_M31_01} (1) implies that $ \wh{\tau} (\oW) $ is a $[t,\infty]-$valued $\bF^{\oW^t,\oP}-$stopping time on $\oO$,
 applying Lemma \ref{lem_012922_11} with $(\fP_t,\tau) \= \big(\{\oP\}, \wh{\tau} \big)$ assures
 that there exists $ \ocA_*  \ins \cF^{\oW^t}_\oga $ satisfying
 \bea \label{020322_11}
 \big\{ \wh{\tau}(\oW) \gs \oga \big\} \D \ocA_* \ins \sN_\oP \big(\cF^{\oW^t}_\infty\big) .
 \eea

     Let $n,i \ins \hN$. Set $s^n_i \df t \+ i 2^{-n}$ and $A^n_i \df \big\{s^n_{i-1} \ls \wh{\tau}  \< s^n_i \big\} \ins \cF^{W^t,P_0}_{s^n_i}$
     with $s^n_0 \df t$.   Using Lemma \ref{lem_013022_11},  we  can find    $   \ocN^n_i \ins \sN_\oP \big(\cF^{\oW^t}_\infty\big) $  such that for any $(s,\oo) \ins [t,s^n_i] \ti \oO$,  there exists   $A^{s,\oo}_{n,i}  \ins \cF^{W^s}_{s^n_i} $   satisfying  $ \b1_{\{ \oW  (\oo')  \in  A^n_i \}}
    \= \b1_{\{ \oW  (\oo')  \in  A^{s,\oo}_{n,i} \}},   \fa  \oo'  \ins  \obW^t_{s,\oo} \Cp \big(\ocN^n_i\big)^c   $.
    For each $ \oo \ins \{\oga \ls s^n_i \} $,      taking     $s 
  \= t_\oo$ yields some $A^\oo_{n,i} \= A^{t_\oo,\oo}_{n,i}  \ins \cF^{W^{t_\oo}}_{s^n_i}$  such that
     \bea \label{052021_21}
      \b1_{\{ \oW  (\oo')  \in  A^n_i \}}    \= \b1_{\{ \oW  (\oo')  \in  A^\oo_{n,i} \}},
      \q \fa  \oo'  \ins  \Wtgo \Cp \big(\ocN^n_i\big)^c    .
     \eea

 Set $\ocN_{\wh{\tau}}  \df \ccup{n,i \in  \hN}{}  \ocN^n_i \ins \sN_\oP \big(\cF^{\oW^t}_\infty\big)   $.
 By (R2) in \eqref{R2}, it holds for any $\oo \ins \oO$ except on an   $  \ocN_3   \ins \sN_\oP \big(\cF^{\oW^t}_\oga\big) $ that
 \bea \label{052021_25}
  \oP^t_{\oga,\oo}  \Big(  \big( \{ \wh{\tau}(\oW) \gs \oga  \} \D \ocA_* \big)   \cp \ocN_{\wh{\tau}} \cp \big\{\oT \nne   \wh{\tau}(\oW )\big\} \Big)
  \= E_\oP \Big[ \b1_{ ( \{ \wh{\tau}(\oW) \ge \oga  \} \D \ocA_*  )   \cup \ocN_{\wh{\tau}} \cup  \{\oT \ne   \wh{\tau}(\oW ) \}} \Big| \cF^{\oW^t}_\oga \Big](\oo)  \=  0   .
 \eea

   Fix $\oo \ins \ocN^c_0 \Cp \ocN^c_3 \Cp \ocA_*  $.    
   For any $n \ins \hN$,   set $\fri_n   (\oo) \df \big\lfloor 2^n ( t_\oo \- t)  \big\rfloor \+ 1 \> 2^n ( t_\oo \- t) $
   and defines a     $(t_\oo,\infty]-$valued $\bF^{W^{t_\oo}}-$stopping time:
 \beas 
 \wh{\tau}^\oo_n (\o_0) \df \sum^\infty_{i=\fri_n \n (\oo)} \b1_{\{\o_0 \in  A^\oo_{n,i}\}}   s^n_i
 \+ \infty \b1_{\big\{\o_0 \in \ccap{i=\fri_n \n (\oo)}{\infty} ( A^\oo_{n,i})^c \big\}}  , \q \fa \o_0 \ins \O_0 .
 \eeas

  As $\bF^{W^{t_\oo},P_0}$ is a right-continuous complete filtration,
 Lemma I.2.11 of \cite{Kara_Shr_BMSC} implies that
  $ \wh{\tau}^\oo   (\o_0) \df \linf{n \to \infty} \wh{\tau}^\oo_n (\o_0) $, $ \fa \o_0 \ins \O_0$
  is a $[t_\oo,\infty]-$valued   $\bF^{W^{t_\oo},P_0}-$stopping time. 

  Let $ \oo' \ins \Wtgo \Cp   \ocN^c_{\wh{\tau}}   \Cp \big\{   \wh{\tau} ( \oW  ) \gs \oga  \big\} $ and   $n \ins \hN$.
   Since \eqref{090520_11} shows that  $ \wh{\tau} \big( \oW  (\oo') \big) \gs \oga (\oo') \= t_\oo
 $,   \eqref{052021_21} renders  that
  \beas
  && \hspace{-1.2cm} \sum_{i \in \hN} \b1_{\{s^n_{i-1} \le \wh{\tau}(\oW   (\oo')) < s^n_i \}}   s^n_i
 \+ \infty \b1_{\{\wh{\tau}(\oW  (\oo'))  = \infty \}}
    \=    \sum^\infty_{i=\fri_n \n (\oo)} \b1_{\{\oW   (\oo') \in  A^n_i\}}   s^n_i
 \+ \infty \b1_{\big\{\oW  (\oo') \in \ccap{i=\fri_n \n (\oo)}{\infty} ( A^n_i )^c \big\}}  \\
 &&  \=   \sum^\infty_{i=\fri_n \n (\oo)} \b1_{\{\oW   (\oo') \in  A^\oo_{n,i}\}}   s^n_i
 \+ \infty \b1_{\big\{\oW  (\oo') \in \ccap{i=\fri_n \n (\oo)}{\infty} ( A^\oo_{n,i})^c \big\}}
 \= \wh{\tau}^\oo_n \big( \oW  (\oo') \big)   .
 \eeas
   Sending $ n \nto \infty$ reaches that
 $ \wh{\tau} \big( \oW (\oo') \big)
 \= \lmtd{n \to \infty}  \wh{\tau}^\oo_n (\oW (\oo'))
 \=   \wh{\tau}^\oo  (\oW (\oo'))    $. So
  $ \Wtgo \Cp   \ocN^c_{\wh{\tau}}  \Cp \big\{   \wh{\tau} ( \oW  ) \gs \oga  \big\} \Cp \big\{ \oT \=   \wh{\tau}^\oo (\oW ) \big\}
 \\ \=   \Wtgo \Cp  \ocN^c_{\wh{\tau}}  \Cp \big\{   \wh{\tau} ( \oW  ) \gs \oga  \big\} \Cp \big\{ \oT \=    \wh{\tau} ( \oW  ) \big\} $.
 Since $\oP^t_{\oga,\oo}   \big\{\oT \nne   \wh{\tau}(\oW )\big\} \= 0$, $   \oP^t_{\oga,\oo}  \big( \{ \wh{\tau}(\oW) \gs \oga  \} \D \ocA_* \big)   \= 0$ by \eqref{052021_25} and since $ \ocA_*  \ins \cF^{\oW^t}_\oga$, we can deduce from \eqref{Jan11_03} and (R3) in \eqref{R3} that
 $ \oP^t_{\oga,\oo} \big( \big\{   \wh{\tau} ( \oW  ) \gs \oga  \big\} \Cp \big\{ \oT \=   \wh{\tau}^\oo (\oW  ) \big\} \big)
 \= \oP^t_{\oga,\oo} \big( \big\{   \wh{\tau} ( \oW  ) \gs \oga  \big\} \Cp \big\{ \oT \=  \wh{\tau} ( \oW  ) \big\} \big)
 \= \oP^t_{\oga,\oo} \big\{   \wh{\tau} ( \oW  ) \gs \oga  \big\}  \= \oP^t_{\oga,\oo} ( \ocA_* )  \=  \b1_{\{\oo \in \ocA_*\}} \= 1  $.
  Hence,   for any $ \oo \ins \ocN^c_0 \Cp \ocN^c_3 \Cp \ocA_* $,
  $\oP^t_{\oga,\oo}   \big\{ \oT \=   \wh{\tau}^\oo (\oW  ) \big\}   \= 1$, i.e.,  $\oP^t_{\oga,\oo} $ satisfies
  (D3) in   Definition \ref{def_ocP}  of $\ocP_{t_\oo,\oX_{\oga \land \cd}  (\oo)}$.

 \no {\bf 3)} Let $i \ins \hN$.
   \if{0}

 Since the $\bF^{\oW^t}-$stopping time $\oga$ is $\cF^{\oW^t}_\oga-$measurable
 and since  $\cF^{\oW^t}_\oga \sb \cF^{\oW^t}_\infty \sb \sB(\oO) $,
  we can deduce from the measurability of function  $ \ol{\cI}_{g_i}  $   defined in \eqref{Sep21_04}  that
  $ \oO \n \ni \n  \oo \mto \ol{\cI}_{g_i} \big(\oga(\oo),\oo\big) \= \int_{\oT(\oo) \land \oga(\oo)}^{\oT(\oo)}  g_i  \big(r,\oX (r \ld \cd,\oo) \big) dr  $
  is a $\sB(\oO) -$measurable random variable on $\oO$.

   \fi
 According to (R2) in \eqref{R2}, it holds for all $\oo \ins \oO$ except on   $\ocN^{\,i}_{\n g,h}   \ins \sN_\oP \big(\cF^{\oW^t}_\oga\big) $   that
 $ E_{\oP^t_{\oga,\oo}} \big[  \int_{\oT \land  \oga }^\oT   g_i  (r, \oX_{r \land \cd}) dr    \big]
 \= \big(\oY^i_{\n \oP}  (\oga) \big) (\oo) $
 \if{0}

 \beas
 E_{\oP^t_{\oga,\oo}} \Big[  \int_{\oT \land  \oga }^\oT   g_i  (r, \oX_{r \land \cd}) dr    \Big]
 & \tn \= & \tn   E_{\oP^t_{\oga,\oo}} \Big[  \int_{\oT \land  \oga }^\oT   g^+_i  (r, \oX_{r \land \cd}) dr    \Big]
 \-  E_{\oP^t_{\oga,\oo}} \Big[  \int_{\oT \land  \oga }^\oT   g^-_i  (r, \oX_{r \land \cd}) dr    \Big] \\
  & \tn \= & \tn   E_\oP \Big[   \int_{\oT \land  \oga }^\oT    g^+_i (r,\oX_{r \land \cd} ) dr \Big| \cF^{\oW^t}_\oga \Big] (\oo)
 \-  E_\oP \Big[   \int_{\oT \land  \oga }^\oT    g^-_i (r,\oX_{r \land \cd} ) dr \Big| \cF^{\oW^t}_\oga \Big] (\oo)
 \= \big(\oY^i_{\n \oP}  (\oga) \big) (\oo)
 \eeas

 \fi
 and $  E_{\oP^t_{\oga,\oo}} \big[  \int_{\oT \land  \oga }^\oT   h_i  (r, \oX_{r \land \cd}) dr    \big]  \= \big(\oZ^i_\oP  (\oga) \big) (\oo) $.
   Given $ \oo \ins   \Big( \ocN_0  \cup \ocN_3 \cup \ocN^{\,i}_{\n g,h}  \Big)^c \Cp \ocA_* $,
 \eqref{090520_11}, \eqref{Jan11_03} and  $\oP^t_{\oga,\oo}   \big\{ \oT \=   \wh{\tau}^\oo (\oW  ) \gs t_\oo \big\}   \= 1$
 from Part (2)  imply that
  $ \big( \oY^i_{\n \oP}  (\oga) \big)  (\oo)
    \=    E_{\oP^t_{\oga,\oo}} \big[  \int_{\oT \land  \oga }^\oT   g_i  (r, \oX_{r \land \cd}) dr    \big]
  \= E_{\oP^t_{\oga,\oo}} \big[ \b1_{\Wtgo} \int_{\oT \land  t_\oo }^\oT   g_i (r, \oX_{r \land \cd}) dr    \big] \\
    \=     E_{\oP^t_{\oga,\oo}} \big[  \int_{  t_\oo }^{\oT  }  g_i (r, \oX_{r \land \cd}) dr    \big]
   $
 and similarly that  $       E_{\oP^t_{\oga,\oo}} \big[  \int_{  t_\oo }^{\oT  } h_i (r, \oX_{r \land \cd}) dr    \big] \= \big( \oZ^i_\oP  (\oga) \big) (\oo) $. Hence,
   \bea   \label{062821_11b}
  \oP^t_{\oga,\oo} \ins \ocP_{ \oga(\oo),\oX_{\oga \land \cd}  (\oo)} \Big(  \big(\oY_{\n \oP}  (\oga )\big)    (\oo), \big(\oZ_\oP  (\oga )\big)    (\oo) \Big), \q \fa \oo \ins \ocA_*   \Cp  \ocN^c_* ,
  \eea
  where $ \ocN_* \df \ocN_0 \cup \ocN_1 \cup \ocN_2 \cup \ocN_3   \cup \Big( \ccup{i \in \hN}{} \ocN^{\,i}_{\n g,h}  \Big) \ins \sN_\oP \big(\ocF^t_{\n \infty}\big) $.
 In particular, \eqref{062821_11}  holds for $\oP-$null set $\ocN \df \ocN_*  \cp   \big\{\oT \nne   \wh{\tau}(\oW )\big\} \cp \big( \{ \wh{\tau}(\oW) \gs \oga  \} \D \ocA_* \big) \ins \sN_\oP  \big(\sB(\oO) \big) $.  \qed

   \no {\bf Proof of Theorem \ref{thm_DPP1}:}
 We first show  the  measurability of $ \b1_{\{\oT \ge \ogaP\}} \oV \big(  \ogaP ,\oX_{ \ogaP  \land \cd} ,    \oY_{\n \oP}  \big( \ogaP \big) , \oZ_\oP  \big( \ogaP \big)   \big) $ for each $\oP \ins \ocP_{t,\bx}(y,z)$ so that the right hand side of \eqref{020422_14} is well-defined.

   Let $\oP \ins \ocP_{t,\bx}(y,z)$ and  simply denote $\ogaP$ by $\oga$.
  Like in Part (1a) of the proof of Proposition \ref{prop_flow},
 we still set   $\ocN_{\n X} \df \big\{ \oo \ins \oO \n :   \oX_s  (\oo)   \nne  \osX^{t,\bx}_s  (\oo)    \hb{ for some } s \ins [0,\infty) \big\}
       \ins \sN_\oP \big(\ocF^t_{\n \infty}\big)$
    and let $  \big\{ \oK^t_s \big\}_{s \in [t,\infty)}$ be  the  $\bF^{\oW^t}-$predictable  process
   such that $ \ocN_{\n K} \df  \big\{ \oo \ins \oO \n :   \oK^t_s  (\oo)   \nne \osX^{t,\bx}_s  (\oo)   $ for some $s \ins [t,\infty)  \big\}  \ins  \sN_\oP \big(\cF^{\oW^t}_\infty\big)$.
  Since $ \oX|_{[0,t)} \= \bx|_{[0,t)}$   and
   $ \oX|_{[t,\infty)} \= \oK^t  $ on $ \big(\ocN_{\n X} \cp \ocN_{\n K}\big)^c $, one can deduce that
    the path-valued random variable $ \oX_{\oga \land \cd} \n : \oO \mto \OmX $ is   $ \si \big( \cF^{\oW^t}_\oga \cp  \sN_\oP  (\ocF^t_{\n \infty} ) \big) \big/   \sB(\OmX)-$measurable.

 Set $\breve{\O} \df [0,\infty) \ti \OmX \times \Re \times \Re \n \supset \n \oD$.
  Let $\wh{\tau}$ be the $[t,\infty]-$valued   $ \bF^{W^t,P_0} -$stopping time with
 $ \oP \big\{  \oT  \=   \wh{\tau} (\oW)   \big\} \= 1$ and let  $ \ocA_*   \ins \cF^{\oW^t}_\oga$, $\ocN_* \ins \sN_\oP \big(\ocF^t_{\n \infty}\big)   $ be  as in \eqref{020322_11} and \eqref{062821_11b}. For any $ \oo \in \ocA_* \cap   \ocN^c_* $, we know from \eqref{062821_11b}   that   $ \Big( \oga (\oo) ,\oX_{\oga \land \cd}  (\oo) ,   \big(\oY_\oP  (\oga )\big)  (\oo),   \big(\oZ_\oP  (\oga )\big)    (\oo) \Big) \ins \oD $. By the measurability of $ \oX_{\oga \land \cd}$,
\bea  \label{020622_23}
\breve{\Psi} (\oo) \df \b1_{\{\oo \in \ocA^c_* \cup   \ocN_*\}} (t,\bx,y,z) \+ \b1_{\{\oo \in \ocA_* \cap   \ocN^c_*\}} \Big( \oga (\oo) ,\oX_{\oga \land \cd}  (\oo) ,   \big(\oY_\oP  (\oga )\big)    (\oo), \big(\oZ_\oP  (\oga )\big)    (\oo) \Big) \ins \oD , \q \fa \oo \ins \oO
\eea
 is  a $ \si \big( \cF^{\oW^t}_\oga \cp   \sN_\oP \big(\ocF^t_{\n \infty}\big) \big) \big/ \sB  (\oD )  -$measurable random variable,
 which induces a probability measure
  $ \breve{P}   \df   \oP  \circ  \breve{\Psi}^{-1}  $ on $ \big( \breve{\O}, \sB  (\breve{\O} ) \big)$.
    Then $\breve{\Psi}$ is further $ \si \big(  \cF^{\oW^t}_\oga \cp   \sN_\oP  (\ocF^t_{\n \infty}   ) \big) \big/  \si \big( \sB  ( \oD  ) \cp \sN_{\breve{P}}  \big( \sB   ( \oD   )  \big)  \big) -$measurable.
  As the universally measurable function $ (t',\bx',   y',   z') \mto \oV (t',\bx',   y',   z')$  is $\si \big( \sB  ( \oD ) \cp \sN_{\breve{P}} \big( \sB  ( \oD ) \big)  \big) \big/ \sB[-\infty,\infty] -$measurable  by Theorem \ref{thm_V_usa},
  \bea  \label{020622_21}
  \breve{V}(\oo) \df \b1_{\{\oo \in \ocA_* \cap   \ocN^c_*\}} \oV \big( \breve{\Psi} (\oo) \big)
  \= \b1_{\{\oo \in \ocA_* \cap   \ocN^c_*\}} \oV \Big(  \oga (\oo) ,\oX_{\oga \land \cd}  (\oo) ,   \big(\oY_\oP  (\oga )\big)    (\oo), \big(\oZ_\oP  (\oga )\big)    (\oo) \Big) 
  , \q \fa \oo \ins \oO
  \eea
   is $ \si \big( \cF^{\oW^t}_\oga \cp   \sN_\oP \big(\ocF^t_{\n \infty}\big) \big) \big/ \sB[-\infty,\infty] -$measurable.
   We see from \eqref{020322_11}   that  $ (\ocA_* \cap   \ocN^c_*) \D \{\oT \gs \oga\} \sb \big(\ocA_* \D \{\oT \gs \oga\}\big) \cp \ocN_*
   \sb \big(\ocA_* \D \{\wh{\tau}(\oW) \gs \oga\}\big) \cup \{\oT \nne \wh{\tau}(\oW)\} \cup \ocN_* \n \ins \n \sN_\oP \big(\sB(\oO)\big)  $,
   where $\wh{\tau}$   is the  $ \bF^{W^t,P_0} -$stopping time with  $ \oP \big\{  \oT  \=   \wh{\tau} (\oW)   \big\} \= 1$.
   It follows that   $  \b1_{\{\oT(\oo) \ge \oga(\oo) \}} \oV \big(  \oga (\oo) ,\oX_{\oga \land \cd}  (\oo) ,   \big(\oY_\oP  (\oga )\big)    (\oo),   \big(\oZ_\oP  (\oga )\big)    (\oo) \big) $, $  \oo \ins \oO$ is $ \si \big(   \cF^{\oW^t}_\oga \cp  \sN_\oP \big(\sB(\oO)\big) \big) \big/ \sB[-\infty,\infty] -$ measurable   and   the right hand side of \eqref{020422_14} is thus well-defined.

 For any $[t,\infty)-$valued $\bF^{\oW^t}-$stopping time $\oz$, we denote  $\oR  (\oz) \df \int_{\oT \land \oz}^\oT   f   (r, \oX_{r \land \cd}  ) dr  \+ \b1_{\{\oT < \infty\}} \pi   \big(\oT, \oX_{\oT \land \cd}\big) $.

  \no {\bf (I) (sub-solution side)} Fix  $\oP \ins \ocP_{t,\bx}(y,z)$   and  simply denote $\ogaP$ by $\oga$.

 \if{0}

 Since the $\bF^{\oW^t}-$stopping time $\oga$ is $\cF^{\oW^t}_\oga-$measurable
 and since  $\cF^{\oW^t}_\oga \sb \cF^{\oW^t}_\infty \sb \sB(\oO) $,
  we can deduce from the measurability of function  $ \ol{\cI}_f  $   defined in \eqref{Sep21_04}  that
  $ \oO \n \ni \n  \oo \mto \ol{\cI}_f \big(\oga(\oo),\oo\big) \= \int_{\oT(\oo) \land \oga(\oo)}^{\oT(\oo)}  f   \big(r,\oX (r \ld \cd,\oo) \big) dr  $
  is a $\sB(\oO) -$measurable random variable on $\oO$.

  \fi

  Let $\wh{\tau}$ be the $[t,\infty]-$valued   $ \bF^{W^t,P_0} -$stopping time with
 $ \oP \big\{  \oT  \=   \wh{\tau} (\oW)   \big\} \= 1$ and let  $ \ocA_*   \ins \cF^{\oW^t}_\oga$, $\ocN_* \ins \sN_\oP \big(\ocF^t_{\n \infty}\big)   $ be  as in \eqref{020322_11} and \eqref{062821_11b}.
  By (R2) in \eqref{R2}, 
  there is a $\ocN_{f,\pi}  \ins \sN_\oP \big(\cF^{\oW^t}_\oga\big)$ such that
 $  E_{\oP^t_{\oga,\oo}} \big[  \oR (\oga)      \big]
  \=  E_\oP \big[   \oR (\oga)    \big| \cF^{\oW^t}_\oga \big] (\oo) $  for any $\oo \ins   \ocN^c_{f,\pi}    $.
  For any $\oo \ins \ocA_* \Cp \big( \ocN_*  \cp \ocN_{f,\pi}   \big)^c$,
   \if{0}

  By  \eqref{062821_11b}, there are  $ \ocA_*   \ins \cF^{\oW^t}_\oga$
   and   $\ocN_* \ins \sN_\oP \big(\ocF^t_{\n \infty}\big)   $ such that

       \bea  \label{Feb01_07b}
  \oP^t_{\oga,\oo} \ins \ocP_{ \oga(\oo),\oX_{\oga \land \cd}  (\oo)} \Big(  \big(\oY_{\n \oP}  (\oga )\big)    (\oo), \big(\oZ_\oP  (\oga )\big)    (\oo) \Big), \q \fa \oo \ins \ocA_*   \Cp  \ocN^c_* .
  \eea

  \fi
 as $\ocN_0 \sb   \ocN_*$,    \eqref{090520_11}, \eqref{Jan11_03} and \eqref{062821_11b} imply  that
 \beas 
 \hspace{-0.5cm}
  E_\oP \big[ \, \oR (\oga)   \big| \cF^{\oW^t}_\oga \big] (\oo)    \=  E_{\oP^t_{\oga,\oo}} \big[ \, \oR (\oga)   \big]
 \= E_{\oP^t_{\oga,\oo}} \big[ \b1_{\Wtzo} \oR (\oga(\oo))  \big]   \=  E_{\oP^t_{\oga,\oo}} \big[ \, \oR (\oga(\oo))   \big]
 \ls \oV \big(  \oga(\oo),\oX_{\oga \land \cd}(\oo) ,  \big(  \oY_{\n \oP}  (\oga ) \big)   (\oo), \big( \oZ_\oP  (\oga ) \big)   (\oo) \big) .
 \eeas

  Since $\b1_{\{\oT    \ge   \oga \}} \= \b1_{\{ \wh{\tau}(\oW)   \ge   \oga \}} \= \b1_{\ocA_*}$, $\oP-$a.s.  by \eqref{020322_11}
 and since $ \ocA_*    \ins \cF^{\oW^t}_\oga$,  
 the tower property renders that
   \beas
  && \hspace{-1.2cm}  E_\oP \Big[  \b1_{\{\oT    \ge   \oga \}} \oV \big( \oga , \oX_{\oga \land \cd} ,     \oY_{\n \oP}  (\oga ) , \oZ_\oP  (\oga )         \big) \Big]
 \= E_\oP \Big[ \b1_{\ocA_*} \oV \big( \oga , \oX_{\oga \land \cd} ,     \oY_{\n \oP}  (\oga ) , \oZ_\oP  (\oga )         \big) \Big]
    \gs      E_\oP \Big[ \b1_{\ocA_*} E_\oP \big[ \, \oR (\oga)   \big| \cF^{\oW^t}_\oga \big] \Big] \\
  &&  \=    E_\oP \Big[  E_\oP \big[ \b1_{\ocA_*} \oR (\oga)   \big| \cF^{\oW^t}_\oga \big] \Big]
     \=    E_\oP \big[  \b1_{\ocA_*} \oR (\oga)   \big]  \= E_\oP \big[   \b1_{\{\oT    \ge   \oga \}} \oR (\oga)  \big]   .
     \eeas
    \if{0}

     This implies that
     $ \b1_{\ocA_*}  \oV^- \big(  \oga ,\oX_{\oga \land \cd}  ,    \oY_{\n \oP}  (\oga )    ,   \oZ_\oP  (\oga )     \big)
     \ls \b1_{\ocA_*}  E_\oP \big[ \, \oR^- (\oga)   \big| \cF^{\oW^t}_\oga \big]  $
      and thus
     \beas
      && \hspace{-1.2cm}  E_\oP \Big[  \b1_{\{\oT    \ge   \oga \}}  \oV^- \big(  \oga ,\oX_{\oga \land \cd}  ,   \oY_{\n \oP}  (\oga )    ,   \oZ_\oP  (\oga )     \big) \Big]
    \=  E_\oP \Big[ \b1_{\ocA_*}  \oV^- \big(  \oga ,\oX_{\oga \land \cd}  ,    \oY_{\n \oP}  (\oga )     ,   \oZ_\oP  (\oga )     \big) \Big] \\
    && \ls   E_\oP \Big[  \b1_{\ocA_*}  E_\oP \big[ \, \oR^- (\oga)   \big| \cF^{\oW^t}_\oga \big]   \Big]
     \ls  E_\oP \big[ \, \b1_{\ocA_*}   \oR^- (\oga)     \big]
     \ls  E_\oP \big[ \, \oR^- (\oga)     \big] \ls E_{P_0} \Big[ \int_0^\infty \n  f^- (r, \, ^o \n X^{t,\bx}_{r \land \cd} )  dr \Big] \+ c_\pi \< \infty  .
     \eeas

    \fi
 It follows that
   $ E_\oP \big[ \, \oR(t)   \big]
    \ls  E_\oP \Big[ \b1_{\{\oT <  \oga \}} \oR(t) \+ \b1_{\{\oT  \ge \oga \}}  \Big( \n \int_t^{\oga } f  (r,\oX_{r \land \cd}   ) dr
   \+     \oV \big( \oga , \oX_{\oga \land \cd} ,     \oY_{\n \oP}  (\oga ) , \oZ_\oP  (\oga )         \big) \Big) \Big]  $.
  Letting $\oP$ vary over $\ocP_{t,\bx}(y,z)$ yields that
  $ \oV(t,\bx,y,z)   \=  \n   \Sup{\oP \in \ocP_{t,\bx}(y,z)} E_\oP \big[ \, \oR(t) \big]  \ls \n  \Sup{\oP \in \ocP_{t,\bx}(y,z)}  E_\oP \Big[
  \b1_{\{\oT < \ogaP \}} \Big( \n \int_t^\oT \n f(r,\oX_{r \land \cd}) dr   \+ \pi \big(\oT, \oX_{\oT \land \cd}\big) \Big)
      \+ \b1_{\{\oT \ge \ogaP \}} \Big( \n \int_t^{\ogaP} \n f(r,\oX_{r \land \cd}) dr  \+ \oV \big(   \ogaP ,\oX_{\ogaP \land \cd} ,    \oY_{\n \oP}   (\ogaP ),  \oZ_\oP   (\ogaP )   \big) \Big) \Big] $.

  \no {\bf (II) (super-solution side)}
  Let   $\oP \ins \ocP_{t,\bx}(y,z)$ and   simply denote $\oga_\oP$ by $\oga$.
   We shall  show that
  \bea \label{081720_15}
   E_\oP \bigg[
  \b1_{\{\oT < \oga  \}} \oR(t) \+ \b1_{\{\oT \ge \oga  \}} \Big( \n \int_t^\oga  \n    f(r,\oX_{r \land \cd}) dr  \+ \oV \big( \oga  ,\oX_{ \oga   \land \cd} ,    \oY_{\n \oP}   (\oga  ) , \oZ_\oP   (\oga  )  \big) \Big) \bigg] \ls  \oV (t,\bx,y,z) .
  \eea

  As $\cF^{\oW^t}_t \= \{\es,\oO\}$,   the $[t,\infty)-$valued  $\bF^{\oW^t} - $stopping time $\oga$ satisfies either $\{\oga \= t\} \= \oO $ or
 $ \{\oga \> t\} \= \oO $.

     Suppose first  that $\{\oga \= t\} \= \oO $:    for any $i \ins \hN$,
 $ 
 \oY^i_{\n \oP} (t)  \= E_\oP \big[ \int_{\oT \land t}^\oT    g_i (r,\oX_{r \land \cd} ) dr \big| \cF^{\oW^t}_t \big] \= E_\oP \big[ \int_t^\oT    g_i(r,\oX_{r \land \cd} ) dr   \big] \ls y_i$ and 
 $ 
 \oZ^i_\oP (t)
 \= E_\oP \big[ \int_t^\oT    h_i(r,\oX_{r \land \cd} ) dr   \big] \= z_i$.
 \if{0}
 For any $(t,\bx,y,z) \ins \oD$ and for any $y' \= \{y'_i\}_{i \in \hN} \ins \Re$ such that $ (t,\bx,y',z) \ins \oD $ and that
 $y_i \ls y'_i$ for any $i \ins \hN$, it is clear that  
 $\oV (t,\bx,y,z) \ls \oV (t,\bx,y',z)$.
 \fi
  Then 
 $  E_\oP \Big[  \b1_{\{\oT < \oga \}} \oR(t)
  \+ \b1_{\{\oT \ge \oga \}} \big(   \int_t^{\oga} \n f(r,\oX_{r \land \cd}) dr  \+  \oV \big(  \oga ,\oX_{\oga \land \cd} ,    \oY_{\n \oP}  (\oga) , \oZ_\oP  (\oga)  \big)  \big) \Big]
  \n \= \n E_\oP \big[  \,   \oV \big( t   ,\oX_{t \land \cd} ,    \oY_{\n \oP}  ( t ) , \oZ_\oP  ( t )   \big)   \big]
    \ls \oV \big( t   , \bx  ,   y ,z  \big)   $.  

 Let us assume  $ \big\{\oga \> t \big\} \= \oO $ in the rest of this proof
 and set $\ocN_{\n X} \df  \big\{ \oo \ins \oO \n :  \oX_s (\oo)  \nne  \osX^{t,\bx}_s (\oo)   \hb{ for some } s \ins [0,\infty) \big\}  \ins \sN_\oP \big(\ocF^t_{\n \infty}\big) $.

 \no {\bf II.a)}
  Define a   random variable  $   \oW^{t,\oga} \n : \oO \mto \O_0   $   by
 $ \oW^{t,\oga}_r   (\oo)  \df   
  \oW^t  \big(  (r \ve t) \ld   \oga(\oo)  , \oo \big) $,
 $ \fa (r,\oo) \ins [0,\infty) \ti \oO  $,
 which is clearly   $   \cF^{\oW^t}_\oga   \big/ \sB(\O_0) -$measurable.

 Set $\ddot{\O} \df [0,\infty) \ti \O_0  \ti \OmX   \ti \Re \ti \Re \n \supset \n \ocD $  and pick up   an arbitrary element $\bw$ from $\O_0$.
  We let $\wh{\tau}$ be the $[t,\infty]-$valued   $ \bF^{W^t,P_0} -$stopping time with
 $ \oP \big\{  \oT  \=   \wh{\tau} (\oW)   \big\} \= 1$ and let  $ \ocA_*   \ins \cF^{\oW^t}_\oga$, $\ocN_* \ins \sN_\oP \big(\ocF^t_{\n \infty}\big)   $ be  as in \eqref{020322_11} and \eqref{062821_11b}.  Since $(t,\bx,y,z) \ins \oD$, Theorem \ref{thm_V=oV} and \eqref{062821_11b} show that
 $(t,\bw,\bx,y,z) \ins \ocD$ and  that  $   \Big(  \oga (\oo) ,  \oW^{t,\oga} (\oo), \oX_{\oga \land \cd}  (\oo) ,  \big(\oY_\oP  (\oga )\big)    (\oo), \\ \big(\oZ_\oP  (\oga )\big)    (\oo) \Big)  \ins \ocD $ for any $\oo \ins \ocA_* \Cp \ocN^c_*  $. Similarly to \eqref{020622_23},
 \bea
 ~ \;  \ddot{\Psi} (\oo) \df \b1_{\{\oo \in \ocA^c_* \cup \ocN_*\}} (t,\bw,\bx,y,z) \+ \b1_{\{\oo \in \ocA_* \cap \ocN^c_*\}} \Big( \oga (\oo) , \oW^{t,\oga} (\oo) ,    \oX_{\oga \land \cd} (\oo) ,
\big(\oY_{\n \oP}  (\oga ) \big)   (\oo), \big(\oZ_\oP  (\oga ) \big)   (\oo) \Big)
 \ins   \ocD  ,   \label{020822_11}
\eea
 $\fa \oo \ins \oO$ is   $ \si \big(  \cF^{\oW^t}_\oga \cp  \sN_\oP  (\ocF^t_{\n \infty} ) \big) \big/ \sB  ( \ocD ) -$measurable,
 which induces a probability measure   $ \ddot{P}    \df   \oP  \circ \ddot{\Psi}^{-1}  $ on $ \big( \ddot{\O},\sB  ( \ddot{\O}  )\big)$.
  Then $\ddot{\Psi} $ is further $ \si \big(  \cF^{\oW^t}_\oga \cp   \sN_\oP  (\ocF^t_{\n \infty}   ) \big) \big/  \si \big( \sB  ( \ocD  ) \cp \sN_{\ddot{P} }  ( \sB   ( \ocD  )  )  \big) -$measurable.

 \no {\bf II.b)} Fix   $\e \ins (0,1)$   through Part (II.e).

   According to Jankov-von Neumann  Theorem (Proposition 7.50 of \cite{Bertsekas_Shreve_1978}),  Corollary  \ref{cor_graph_ocP} and Theorem \ref{thm_V_usa},
  there exists an analytically measurable function $ \ol{\bQ}_\e \n : \ocD \mto \fP\big(\oO\big) $ such that for any  $(\ft,\fw,\fx,\fy,\fz) \ins \ocD $, $  \ol{\bQ}_\e (\ft,\fw,\fx,\fy,\fz) $ belongs to $   \ocP_{\ft,\fw,\fx}(\fy,\fz) $ and satisfies
\bea \label{081620_19}
\hspace{-1cm}
  E_{\ol{\bQ}_\e (\ft,\fw,\fx,\fy,\fz)} \big[ \, \oR(\ft)   \big] \gs
\left\{
\ba{ll}
\oV(\ft,\fw,\fx,\fy,\fz) \- \e , & \hb{ if } \oV(\ft,\fw,\fx,\fy,\fz) \< \infty ; \\
 1/\e ,  & \hb{ if } \oV(\ft,\fw,\fx,\fy,\fz) \= \infty .
 \ea
 \right.
\eea
   As $ \ol{\bQ}_\e$ is    universally measurable,
   it is also  $  \si \big( \sB  ( \ocD  ) \cup \sN_{\ddot{P}}  ( \sB  ( \ocD )  )  \big) \big/  \sB \big( \fP\big(\oO\big) \big) -$measurable and
  $ \oQ^\oo_\e   \df      \b1_{\{\oo \in \ocA^c_*  \cup   \ocN_* \}} \oP  \\ + \n \b1_{\{\oo \in \ocA_* \cap   \ocN^c_* \}}  \ol{\bQ}_\e \big( \ddot{\Psi} (\oo) \big) $, $ \fa  \oo \ins \oO $
    is thus $  \si \big(  \cF^{\oW^t}_\oga \cp   \sN_\oP  (  \ocF^t_{\n \infty}  ) \big) \big/ \sB\big(\fP\big(\oO\big)\big) -$measurable.

   Given a $[0,\infty]-$valued $\sB(\oO)-$measurable random variable $\ol{\phi}$,
   Proposition 7.25 of  \cite{Bertsekas_Shreve_1978} implies that
the mapping
$ \fP\big(\oO\big) \ni \oQ \mto E_\oQ  \big[ \, \ol{\phi} \, \big] $ is $\sB\big(\fP\big(\oO\big)\big)   -$measurable.
 The  
 measurability of $\big\{\oQ^\oo_\e\big\}_{\oo \in \oO} $ renders that
\bea \label{022022_23}
 \hb{the random variable  $ \oO \ni \oo \mto E_{\oQ^\oo_\e} \big[ \, \ol{\phi} \, \big] $
 is $ \si \Big( \cF^{\oW^t}_\oga \cp \sN_\oP \big(\ocF^t_{\n \infty}\big) \Big) -$measurable.}
 \eea

     Let $\oo \ins \ocA_* \Cp \ocN^c_*$ and denote $ t_\oo \df \oga(\oo)$. We know from \eqref{020822_11} that
    \bea
   \oQ^\oo_\e    \= \ol{\bQ}_\e \big( \ddot{\Psi} (\oo) \big)   \ins   \ocP_{\oga(\oo) , \oW^{t,\oga} (\oo) ,  \oX_{\oga \land \cd} (\oo)} \Big( \big(\oY_{\n \oP}  (\oga ) \big)   (\oo), \big(\oZ_\oP  (\oga ) \big)   (\oo) \Big)   .    \label{April07_11}
  \eea
  By  (D3) in Definition \ref{def_ocP} of $   \ocP_{t,\bx} $,
 there is   a $[t_\oo,\infty]-$valued   $ \bF^{W^{t_\oo},P_0} -$stopping time $\wh{\tau}_\loo$ with  
 \bea \label{022122_11}
  \oQ^\oo_\e \big( \big\{  \oT  \=  \wh{\tau}_\loo (\oW  )   \big\} \big) \= 1 .
  \eea

  Set $ \oO^t_{\oga,\oo}   \df    \big\{ \oo' \ins \oO \n: (\oW_{\n s},\oX_s) (\oo') \=  ( \oW^{t,\oga}_s, \oX_s    ) (\oo)   ,    \fa s \ins  [0,\oga(\oo) ]   \big\} $ and $ \oXi^t_{\oga,\oo}
 \df \big\{ \oo' \ins \oO \n:  \oW^t_s  (\oo')  \=  \oW^t_s  (\oo)  , \fa s \ins \big[t, \oga(\oo)\big] ; \,
   \oX_s(\oo') \= \oX_s(\oo), \fa s \ins [0, \oga(\oo)]  \big\} $.
 Since   $   \oO^t_{\oga,\oo}    \sb       \big\{ \oo' \ins \oO \n:  \oW_{\n s}  (\oo') \=  0   ,
    \fa s \ins [0,t] ; \oW_{\n s}  (\oo') \=   \oW^t_{\n s}  (\oo)    ,
    \fa s \ins (t,\oga(\oo)] ; \oX_s(\oo') \= \oX_s(\oo) , \fa s \ins [0,\oga(\oo)] \big\}
      \sb  
     \oXi^t_{\oga,\oo}      \sb \Wtgo  $,
 we see from \eqref{April07_11}   that
    \bea \label{072820_15}
  \oQ^\oo_\e  \big( \oO^t_{\oga,\oo}  \big) \=    1 , \q \hb{and thus} \q  \oQ^\oo_\e \big( \Wtgo \big) \= \oQ^\oo_\e \big( \oXi^t_{\oga,\oo} \big) \=  1 .
   \eea

    Let $\oA \ins \sB (\oO)$.  We claim that
  \bea \label{May06_25}
 \oQ^\oo_\e (\ocA \Cp \oA) \= \b1_{ \{\oo \in \ocA  \}}  \oQ^\oo_\e (\oA) , \q \fa \ocA \ins \ocF^t_{\n \oga} ,
   ~ \fa \oo \ins \ocA_* \Cp \ocN^c_* .
 \eea
 To see this, we  take  $    \ocA  \ins \ocF^t_{\n \oga} $.
 Let $\oo_1 \ins \ocA  \Cp \ocA_* \Cp   \ocN^c_*  $ and   set    $s_1 \df \oga(\oo_1)$.
 Since $  \ocA  \Cp \{\oga \ls s_1\} $ is an $ \ocF^t_{s_1}-$measurable set including $\oo_1 $,
 one can deduce that
 \bea \label{090920_15}
 \hb{$  \oXi^t_{\oga,\oo_1} \= \big\{ \oo' \ins \oO \n:  \oW^t_r  (\oo')  \=  \oW^t_r  (\oo_1)  , \fa r \ins  [t, s_1] ; \,
   \oX_r(\oo') \= \oX_r(\oo_1), \fa r \ins [0, s_1]  \big\}  $ is also contained in $ \ocA  \Cp \{\oga \ls s_1\} $. } \q
 \eea
 By \eqref{072820_15},  $ \oQ^{\oo_1}_\e \big(\ocA \big)   \= 1 $ and thus  $\oQ^{\oo_1}_\e \big( \ocA  \Cp \oA \big) \= \oQ^{\oo_1}_\e \big(   \oA \big) \= \b1_{ \{\oo_1 \in \ocA  \}}  \oQ^{\oo_1}_\e (\oA) $. We next let $\oo_2 \ins \ocA^c \Cp \ocA_* \Cp   \ocN^c_*   $ and   set $s_2  \df \oga(\oo_2)$.
As $  \ocA^c  \Cp \{\oga \ls s_2 \} $ is an $ \ocF^t_{s_2}-$measurable set including $\oo_2 $,
 $ \oXi^t_{\oga,\oo_2} \= \big\{ \oo' \ins \oO \n:  \oW^t_r  (\oo')  \=  \oW^t_r  (\oo_2)  , \fa r \ins  [t, s_2] ; \,
   \oX_r(\oo') \= \oX_r(\oo_2), \fa r \ins [0, s_2]  \big\}  $
 is also included in $ \ocA^c  \Cp \{\oga \ls s_2\} $.
 We correspondingly have  $ \oQ^{\oo_2}_\e (\ocA^c )   \= 1 $
 and thus $\oQ^{\oo_2}_\e \big( \ocA  \Cp \oA \big) \= 0 \= \b1_{ \{\oo_2 \in \ocA  \}}  \oQ^{\oo_2}_\e (\oA) $. 

 Consider a pasted probability measure $\oP_\e \ins \fP\big(\oO\big)$:
    \bea   \label{092520_11}
  \oP_\e (\oA)  \df  \oP \big(\ocA^c_* \Cp \oA \big)   \+
  \int_{\oo \in \ocA_*}   \oQ^\oo_\e (\oA ) \oP(d\oo)    ,
  \q \fa \oA \ins \sB (\oO) .
 \eea
 \if{0}

    Let $\oA \ins \sB \big(\oO\big)$. By Proposition 7.25 of  \cite{Bertsekas_Shreve_1978},
 the function  $\phi_\oA \big(\oQ\big) \df E_\oQ  [\b1_\oA ] \= \oQ (\oA) $, $ \fa  \oQ  \ins   \fP\big(\oO\big) $ is
 $ \sB\big(\fP\big(\oO\big)\big) \big/ \sB[0,1] -$ measurable.
 As  $ \big\{ \oQ^\oo_\e   \big\}_{ \oo \in  \oO } $ is  $  \si \big(  \cF^{\oW^t}_\oga \cp   \sN_\oP  (  \ocF^t_{\n \infty}  ) \big) \big/ \sB\big(\fP\big(\oO\big)\big) -$measurable,     the random variable  $ \oO \ni \oo \mto \oQ^\oo_\e(\oA) \= \phi_\oA \big(\oQ^\oo_\e\big)   $
  is $ \si \big( \cF^{\oW^t}_\oga \cp \sN_\oP  (\ocF^t_{\n \infty} ) \big) \big/ \sB[0,1] -$measurable.

   It is clear that
  $ \oP_\e (\es) \= \oP \big(\ocA^c_* \Cp \es  \big) \+ \int_{\oo \in \ocA_*} \oQ^\oo_\e (\es ) \oP(d\oo)
 \= \oP (\es) \+ E_\oP \big[ \b1_{\ocA_*} 0 \big] \= 0 $ and $ \oP_\e \big(\oO\big) \= \oP \big(\ocA^c_* \Cp \oO  \big) \+
 \int_{\oo \in \ocA_*} \oQ^\oo_\e \big(\oO\big) \oP(d\oo)
 \= \oP \big( \ocA^c_* \big)     \+ E_\oP \big[ \b1_{\ocA_*} \big]
 \= 1  $.
 Given   a disjoint sequence $\{\oA_i\}_{i \in \hN}$ in $ \sB \big(\oO\big) $, the monotone convergence theorem implies that
 \beas
    \oP_\e \Big( \underset{i \in \hN}{\cup}\oA_i \Big)
   & \tn  \=  & \tn   \oP \Big(\ocA^c_* \Cp \Big( \underset{i \in \hN}{\cup}\oA_i \Big)  \Big)
  \+ \int_{\oo \in \ocA_*}
  \oQ^\oo_\e \Big(\underset{i \in \hN}{\cup}\oA_i  \Big) \oP(d\oo)
    \=     \oP \Big(  \underset{i \in \hN}{\cup} \big( \ocA^c_* \Cp \oA_i \big)  \Big)
  \+ \int_{\oo \in \ocA_*} \underset{i \in \hN}{\sum}  \oQ^\oo_\e  \big( \oA_i  \big) \oP(d\oo) \\
    & \tn  \=  & \tn   \underset{i \in \hN}{\sum} \oP   \big( \ocA^c_* \Cp \oA_i \big)
   \+ \underset{i \in \hN}{\sum} \int_{\oo \in \ocA_*} \n \oQ^\oo_\e  \big( \oA_i  \big) \oP(d\oo)
    \=    \underset{i \in \hN}{\sum} \Big( \oP   \big( \ocA^c_* \Cp \oA_i \big)
   \+   \int_{\oo \in \ocA_*} \n \oQ^\oo_\e  \big( \oA_i  \big) \oP(d\oo) \Big)
   \= \underset{i \in \hN}{\sum} \, \oP_\e \big(\oA_i\big) .
 \eeas
 So  $\oP_\e$ is a probability measure  of   $\fP\big(\oO\big)$.

 \fi
 In particular, taking $\oA \= \oO$ in  \eqref{May06_25} renders that
 \bea
  \oP_\e (\ocA)   \=     \oP \big(\ocA^c_* \Cp \ocA \big)   \+
  \int_{\oo \in \ocA_*} \b1_{
  \{\oo \in \ocA\}  }   \oP(d\oo)   \= \oP  (\ocA)  , \q \fa   \ocA \ins \ocF^t_{\n \oga}  .   \label{052420_21}
  \eea

 In the next three parts, we demonstrate that   $\oP_\e $ also belongs to  $\ocP_{t,\bx}(y,z)$,
 i.e.,  the probability class $\ocP_{t,\bx}(y,z)$ is stable under the pasting \eqref{092520_11}.

 \no {\bf II.c)}
 We first show that  $ \oP_\e $ is of $ \ocP^1_{t,\bx} \Cp \ocP^2_t $ and   thus satisfies \(D1\) and \(D2\) in Definition  \ref{def_ocP}   of $\ocP_{t,\bx}$.

 \no {\bf II.c.1)} Set $ \oO_X \df \big\{    \oX_s    \= \bx(s) ,\,\fa s \ins [0,t ] \big\}$.
 By   the proof of Proposition \ref{prop_flow}, $ \oO^c_X \sb \ocN_{\n X} \=   \big\{ \oo \ins \oO \n :  \oX_s (\oo)  \nne  \osX^{t,\bx}_s (\oo)   \hb{ for some } s \ins [0,\infty) \big\} \sb \ocN_1 \sb   \ocN_* $.
 Given $\oo \ins \ocA_* \Cp \ocN^c_* \sb \oO_X$, 
 one has     $  \oX_s (\oo)  \=  \bx(s) $, $\fa s \ins [0,t] $ and thus
    $\oO^t_{\oga,\oo} \sb   \big\{ \oo' \ins \oO \n:   \oX_s(\oo') \= \oX_s(\oo) , \fa s \ins [0,t ] \big\}
\= \oO_X$. As $\oP \big(   \oO_X  \big) \gs \oP (\ocN^c_*) \= 1$,   \eqref{072820_15} implies that
 $  \oP_\e \big( \oO_X \big)  \=  \oP \big(\ocA^c_* \Cp \oO_X \big)   \+
  \int_{\oo \in \ocA_* \cap \ocN^c_*}  1 \n \cd \n \oP(d\oo)  \= \oP \big(\ocA^c_* \big) \+ \oP(\ocA_*) \= 1   $, i.e.,  $ \oP_\e \ins \ocP^1_{t,\bx} $.

  \no {\bf II.c.2)} We need some technical preparation for checking $ \oP_\e \ins \ocP^2_t $:
  Let $\oo \ins \ocA_* \Cp \ocN^c_*$  and set $\fra_{\overset{}{\oo}} \df   \big( \dn - \n \oW^t_\oga( \oo ) ,   \bz \big) \ins \hR^{d+l}$.
  We define an   $\obF^{t_\oo}-$stopping time $ \oz^n_\loo (\oo') \df \inf\big\{s  \ins [t_\oo,\infty) \n : \big|(\oW^{t_\oo}_s,\oX_s)  (\oo') \- \fra_{\overset{}{\oo}}   \big|   \gs n  \big\}  $, $  \oo' \ins \oO$.

  Given $\vf \ins \fC (\hR^{d+l})$,   define a $C^2 (\hR^{d+l})  $ function
 $\vf_\loo(w,x) \df \vf \big(w \+ \oW^t_\oga  ( \oo ), x \big) $, $(w,x) \ins \hR^{d+l}$.
 For $i \= 0, 1,2 $ and $ \oo' \ins  \Wtgo $, since $D^i\vf  \big(  \oW^t_r (\oo'),\oX_r (\oo')\big)   \=   D^i\vf  \big( \oW^t_r(\oo') \- \oW^t  (\oga(\oo),\oo' )
  \+ \oW^t   (\oga(\oo), \oo ) , \oX_r (\oo') \big)
 \= D^i \vf_\loo \big( \oW^{t_\oo}_r ( \oo' ) ,\oX_r (\oo')\big) $, $\fa r \ins  [t_\oo,\infty)$,
 one has
 \bea
 \big( \,\oM^t_{r_2} (\vf)\big)  (\oo'  ) \- \big(\,\oM^t_{r_1} (\vf)\big)  (\oo'  )
 \=  \big( \,\oM^{t_\oo}_{r_2} (\vf_\oo)\big)  (\oo'  ) \- \big( \,\oM^{t_\oo}_{r_1} (\vf_\oo)\big)  (\oo'  ) , \q \fa t_\oo \ls r_1 \ls r_2 \< \infty.   \label{091920_17}
 \eea

 Let   $  \Big\{ \ofX^\oo_s \=  \osX^{t_\oo,\oX_{\oga \land \cd} (\oo)}_s \Big\}_{s \in  [0,\infty)}$ be
 the $\Big\{ \cF^{\oW^{t_\oo},\oQ^\oo_\e}_{s \vee t_\oo } \Big\}_{s \in [0,\infty)} -$adapted continuous process   that
  uniquely  solves  the   following SDE on $\big(\oO , \sB(\oO) , \oQ^\oo_\e\big) \n :  $
 \beas
 \q \osX_s    \=  \oX_\oga (\oo) \+ \n  \int_{t_\oo}^s \n  b ( r, \osX_{r \land \cd})dr
      \+ \n  \int_{t_\oo}^s \n  \si ( r, \osX_{r \land \cd}) d \oW_{\n r}, \q s \ins [t_\oo,\infty) \; \hb{ with initial condition } \osX_s \= \oX_s(\oo)  , \; \fa s \ins [0,t_\oo]    .
 \eeas
   By \eqref{April07_11},
    $ \ocN^{\,\oo}_{\n X} \df   \big\{ \oo' \ins \oO \n :    \oX_s (\oo')  \nne \ofX^\oo_s  ( \oo')   \hb{ for some } s \ins [0,\infty) \big\}   \ins \sN_{\oQ^\oo_\e} \big(\ocF^{ t_\oo }_\infty \big) $.
  And there exists an $\hR^l-$valued  $\bF^{\oW^{t_\oo}}-$ predictable  process  $  \big\{ K^\oo_s \big\}_{s \in [t_\oo,\infty)}$ such that
   $ \ocN^{\,\oo}_K \df  \big\{ \oo' \ins \oO \n : K^\oo_s  (\oo')  \nne  \ofX^\oo_s  ( \oo')   \hb{ for some } s \ins [t_\oo,\infty)  \big\}  \ins  \sN_{\oQ^\oo_\e} \big(\cF^{\oW^{t_\oo}}_\infty\big) $.

  Let $(\vf,n) \ins \fC (\hR^{d+l}) \ti \hN$, $(\fs,\fr) \ins \hQ^{2,<}_+   $ and $   \{(s_i,\cO_i )\}^k_{i=1} \sb \big(\hQ \Cp [0,\fs]\big) \ti \sO (\hR^{d+l}) $.
  Denote $\oM^{t,n}_s (\vf ) \df \oM^t_{s \land \otau^t_n  } (\vf ) $, $s \ins [t,\infty)$ and set $\wh{A} \df \ccap{i=1}{k}  (\oW^t_{ t+s_i  },\oX_{ t+s_i   })^{-1} (\cO_i) \ins \ocF^t_{t+\fs}$.

  \no {\bf (i)} To verify $E_{\oP_\e} \big[   \big( \oM^{t,n}_{t+\fr} (\vf )  \- \oM^{t,n}_{t+\fs} (\vf ) \big)   \b1_{\wh{A}}    \big]   \= 0 $,
  we first show that
   \bea \label{071720_31}
   E_{\oP_\e} \Big[ \b1_{\{\oga >  t+ \fs\}}  \big( \oM^{t,n}_{t+\fr} (\vf )  \- \oM^{t,n}_{t+\fs} (\vf ) \big)   \b1_{\wh{A}}    \Big]   \= 0 .
 \eea

 Since $ \{\oga \> t\+\fs\}  \Cp \wh{A} \= \{\oga \> t\+\fs\}  \Cp \Big( \ccap{i=1}{k}  (\oW^t_{ \oga \land (t+s_i)  },\oX_{ \oga \land (t+s_i)  })^{-1} (\cO_i) \Big) \n \ins \ocF^t_{\n \oga \land (t+\fs)} $ and
$     \oM^{t,n}_{\oga \land (t+\fr)} (\vf ) \- \oM^{t,n}_{\oga \land (t+\fs)} (\vf )   \ins \ocF^t_{ \n \oga \land (t+\fr)}
\sb \ocF^t_{\n \oga}$,  
using \eqref{052420_21} and applying \eqref{020522_17} with $ \big(\fra,\oz_1,\oz_2 \big) \= \big(\bz, \oga \ld (t\+\fs), \oga \ld (t\+\fr) \big)$
 yield that
  \bea
   E_{\oP_\e} \Big[   \big( \oM^{t,n}_{\oga \land (t+\fr)} (\vf ) \- \oM^{t,n}_{t+\fs} (\vf ) \big)   \b1_{\{\oga  > t+\fs\}  \cap \wh{A}}    \Big]
    \=    E_\oP  \Big[   \big( \oM^t_{\otau^t_n \land \oga \land (t+\fr)} (\vf ) \- \oM^t_{\otau^t_n \land \oga \land (t+\fs)} (\vf ) \big)   \b1_{\{\oga  > t+\fs\}  \cap \wh{A}}     \Big]  \= 0 . \qq  \label{071420_21}
  \eea
 And  \eqref{May06_25} implies that
  \bea
 && \hspace{-1.4cm}  E_{\oP_\e} \Big[   \big( \, \oM^{t,n}_{t+ \fr} (\vf ) \- \oM^{t,n}_{ \oga \land (t+\fr)} (\vf ) \big)   \b1_{\{\oga  > t+\fs\}  \cap \wh{A}}     \Big]
\=  E_{\oP_\e} \Big[   \big( \, \oM^{t,n}_{ \oga \vee (t+\fr) } (\vf ) \- \oM^{t,n}_{ \oga} (\vf ) \big)   \b1_{\{\oga  > t+\fs\}  \cap \wh{A}}     \Big]  \nonumber \\
 &&  \hspace{-1cm} \=  E_\oP \Big[ \b1_{\ocA^c_*}   \big(   \oM^{t,n}_{\oga \vee (t+\fr)} (\vf ) \- \n  \oM^{t,n}_{\oga} (\vf ) \big)   \b1_{\{\oga  > t+\fs\}  \cap \wh{A}}    \Big] \n \+ \n \int_{\oo \in \ocA_*} \dn  \b1_{ \{\oga (\oo) > t+\fs \}} \b1_{ \{\oo \in \wh{A}\,\}}  E_{\oQ^\oo_\e} \n \big[ \,   \oM^{t,n}_{\oga \vee (t+\fr)} (\vf ) \- \n \oM^{t,n}_{\oga} (\vf )  \big]  \oP(d\oo) . \qq \label{071720_27}
\eea
 Taking $ \big(\fra,\oz_1,\oz_2 \big) \= \big(\bz, \oga, \oga \ve (t\+\fr) \big)$ in \eqref{020522_17}
 renders   that
\bea
   E_\oP \Big[ \b1_{\ocA^c_*}   \big( \, \oM^{t,n}_{\oga \vee (t+\fr)} (\vf ) \- \oM^{t,n}_{\oga} (\vf ) \big)   \b1_{\{\oga  > t+\fs\}  \cap \wh{A}}     \Big] \=  E_\oP \Big[ \b1_{\ocA^c_*  } \big( \, \oM^t_{\otau^t_n \land (\oga \vee (t+\fr))} (\vf ) \- \oM^t_{\otau^t_n \land \oga} (\vf ) \big)   \b1_{\{\oga  > t+\fs\}  \cap \wh{A}}     \Big]   \= 0 . \q \label{071720_25}
\eea

 Fix $ \oo  \ins \big\{ \otau^t_n   \> \oga   \big\} \Cp \ocA_* \Cp \ocN^c_* $ and set $\fr_\loo \df t_\oo \ve (t\+\fr)$.
 As $ t_\oo \= \oga (\oo) \< \otau^t_n (\oo) \ls t\+n $,    $n_\loo \df \lceil t\+n \- t_\oo \rceil \ins \hN   $.
 Using \eqref{020522_17} with $(t,\bx,\oP,\vf,n,\fra,\oz_1,\oz_2) \= \big(t_\oo ,  \oX_{\oga \land \cd} (\oo), \oQ^\oo_\e ,\vf_\loo,  n_\loo, \fra_{\overset{}{\oo}},  t_\oo, \fr_\loo \ld (t\+n)\big)$ yields that
  \bea \label{071720_23}
  0 \= E_{\oQ^\oo_\e} \Big[   \oMoo_{\oz^n_\loo \land ( t_\oo + n_\loo ) \land \fr_\loo \land (t + n)} (\vf_\loo) \- \oMoo_{\oz^n_\loo \land ( t_\oo + n_\loo )   \land t_\oo} (\vf_\loo) \Big] \=   E_{\oQ^\oo_\e} \Big[   \oMoo_{\oz^n_\loo   \land \fr_\loo \land (t+n) } (\vf_\loo) \- \oM^{t_\oo}_{t_\oo} (\vf_\loo) \Big] .
 \eea

 Because $\oo \ins \big\{ \oo' \ins \oO \n :  \otau^t_n(\oo') \> \oga (\oo) \big\} \ins \ocF^t_{\oga (\oo)}  $, an analogy to \eqref{090920_15} shows that
 $\oXi^t_{\oga,\oo} \sb \big\{ \oo' \ins \oO \n :  \otau^t_n(\oo') \> \oga (\oo) \big\}$.
 Let $\oo' \ins \oXi^t_{\oga,\oo}  $. Since $\inf \n \big\{s \ins [t,\infty) \n : \big| (\oW^t_s,\oX_s) (\oo') \big|   \gs n  \big\} \gs \otau^t_n(\oo') \> \oga(\oo)$, one has $\big| (\oW^t_s,\oX_s) (\oo') \big|\< n $, $\fa s \ins [t,t_\oo] $  and thus
 \beas
   \;\;\;   \inf \n \big\{s \ins [t,\infty) \n : |(\oW^t_s,\oX_s) (\oo')|   \gs n  \big\}
   \=   
    \inf\big\{s \ins [t_\oo,\infty) \n : \big|\oW^t_s (\oo') \- \oW^t (\oga(\oo),\oo'),\oX_s(\oo')\big) \+ \big(\oW^t_\oga (\oo),\bz\big) \big|   \gs n  \big\}
 \= \oz^n_\loo (\oo')   .
 \eeas
  It follows that $  \otau^t_n(\oo')   \=   \oz^n_\loo (\oo')   \ld (t\+n)   $.
  Taking $(r_1,r_2) \= \big(t_\oo, \otau^t_n(\oo') \ld \fr_\loo \big)$ in \eqref{091920_17},
   we can   deduce from \eqref{090520_11}   that
 \beas
&& \hspace{-1cm} \big(\,\oM^t \n (\vf)\big) \big(   \otau^t_n(\oo')    \ld  (\oga (\oo') \ve (t\+\fr) ) ,\oo' \big) \- \big(\,\oM^t \n (\vf ) \big) \big( \oga (\oo'), \oo'  \big)
  \= \big(\,\oM^t \n (\vf)\big) \big(  \otau^t_n(\oo')    \ld  (\oga (\oo) \ve (t\+\fr) ) ,\oo' \big) \- \big(\oM^t \n (\vf ) \big) \big( \oga (\oo), \oo'  \big)  \nonumber \\
 &&  \hspace{-0.5cm}   \= \big( \oM^{t_\oo}  (\vf_\loo) \big) \big( \otau^t_n(\oo') \ld \fr_\loo, \oo' \big) \- \big( \oM^{t_\oo}  (\vf_\loo) \big)  ( t_\oo , \oo'  ) \= \big( \oM^{t_\oo}  (\vf_\loo) \big) \big( \oz^n_\loo (\oo') \ld (t \+ n)  \ld \fr_\loo, \oo' \big) \- \big( \oM^{t_\oo}  (\vf_\loo) \big)  ( t_\oo , \oo'  ) . \q  
\eeas
 As $\{ \otau^t_n   \> \oga   \}  \ins \ocF^t_{\otau^t_n \land \oga} \sb \ocF^t_{\n \oga}$,
   \eqref{May06_25},  \eqref{072820_15}, 
  and \eqref{071720_23} then imply that
 \beas
  && \hspace{-1cm} E_{\oQ^\oo_\e} \big[ \,  \oM^{t,n}_{\oga \vee (t+\fr)} (\vf ) \- \oM^{t,n}_{\oga} (\vf )  \big]
 \= E_{\oQ^\oo_\e} \big[ \b1_{ \{ \otau^t_n   > \oga   \}}  \big( \,   \oM^t_{ \otau^t_n   \land (\oga \vee (t+\fr))} (\vf ) \- \oM^t_\oga (\vf ) \big)   \big]
\= \b1_{ \{ \otau^t_n (\oo)  > \oga (\oo)  \}}  E_{\oQ^\oo_\e} \big[  \,  \oM^t_{ \otau^t_n   \land (\oga \vee (t+\fr))} (\vf ) \- \oM^t_\oga (\vf )   \big] \\
  &&    \=   \b1_{ \{ \otau^t_n(\oo)   > \oga (\oo)  \} }  E_{\oQ^\oo_\e} \Big[   \oMoo_{\oz^n_\loo   \land  (t+n)  \land \fr_\loo }  (\vf_\loo)    \-   \oM^{t_\oo}_{t_\oo}  (\vf_\loo)    \Big] \= 0  , ~ \;  \fa \oo  \ins  \ocA_* \Cp     \ocN^c_* .
 \eeas
 So  $  \int_{\oo \in \ocA_*}  \n  \b1_{ \{\oga (\oo) > t+\fs \}}  \b1_{\{\oo \in \wh{A}\,\} }   E_{\oQ^\oo_\e} \big[  \, \oM^{t,n}_{\oga \vee (t+\fr)} (\vf )  \-   \oM^{t,n}_{\oga} \n (\vf )  \big]  \oP(d\oo)   \= 0 $,
 which together with \eqref{071420_21}$-$\eqref{071720_25}   leads to \eqref{071720_31}.

  \no {\bf (ii)} If $\fs \= 0$,    as $\{\oga \> t\} \= \oO$,    \eqref{071720_31} directly gives
 $  E_{\oP_\e} \big[   \big( \, \oM^{t,n}_{t+\fr} (\vf ) \- \oM^{t,n}_t (\vf ) \big)   \b1_{\wh{A}}  \big]  \= 0  $.

 Next, let $\fs \> 0$. In this case, we can assume with loss of generality that  $ 0 \=   s_1 \< \cds   \< s_k \= \fs $ with $k \gs 2 $.
 As $ \ocA^c_* \ins \cF^{\oW^t}_\oga \sb \ocF^t_{\n \oga}  $, one has
  $ \ocA^c_* \Cp   \{\oga \ls t \+ \fs  \} \ins \ocF^t_{t+\fs} $. Applying \eqref{020522_17} with $(\fra,\oz_1,\oz_2 ) \= (\bz,t\+\fs,t\+\fr) $
   yields that
 \bea
 \hspace{-0.5cm}
  E_\oP \Big[ \b1_{\ocA^c_* \cap \{\oga \le t+\fs\}} \big( \, \oM^{t,n}_{t+\fr} (\vf ) \- \oM^{t,n}_{t+\fs} (\vf ) \big)   \b1_{ \wh{A}}     \Big]
 \= E_\oP \Big[ \b1_{\ocA^c_* \cap \{\oga \le t+\fs\}} \big( \, \oM^t_{\otau^t_n \land (t+\fr)} (\vf ) \- \oM^t_{\otau^t_n \land (t+\fs)} (\vf ) \big)  \b1_{ \wh{A}}    \Big]  \= 0 . \q  \label{071820_27}
\eea

 Fix $i \ins \{ 1, \cds \n , k \- 1\}$   and fix $\oo \ins \big\{ \otau^t_n   \> \oga   \big\} \Cp \big\{t\+s_i \< \oga  \ls t\+ s_{i+1} \big\} \Cp \ocA_* \Cp   \ocN^c_* $.
     Since $ \Wtgo \sb \{\oga \> t\+s_i  \} $ by  \eqref{090520_11},   $\oA_i \df \ccap{j=1}{i}    \big( \oW^t_{ \oga \land (t+s_j) }  , \oX_{ \oga \land (t+s_j) }      \big)^{-1} (\cO_j) \ins 
 \ocF^t_{\n \oga} $ satisfies
\bea
    \Wtgo \Cp \Big(  \ccap{j=1}{i}   \big( \oW^t_{t+ s_j  }, \oX_{t+ s_j  }  \big)^{-1} (\cO_j)  \Big)
 \= \Wtgo \Cp \oA_i .   \label{070120_11}
\eea
 Also,  \eqref{090520_11} shows that $ \Wtgo \sb \{\oga \ls t \+ \fs   \} $ and thus $\Wtgo \Cp \{\otau^t_n \ls \oga\} 
  \sb \{\otau^t_n \ls t \+ \fs\}$. By   \eqref{072820_15},
  \bea  \label{021322_11}
    E_{\oQ^\oo_\e} \big[ \b1_{\{\otau^t_n \le \oga\}}\big( \oM^{t,n}_{t+\fr} (\vf ) \- \oM^{t,n}_{t+\fs} (\vf ) \big)  \b1_{ \wh{A} }  \big]
   \ls  E_{\oQ^\oo_\e} \big[ \b1_{\{\otau^t_n \le t  + \fs\}}\big( \oM^t_{\otau^t_n \land (t+\fr)} (\vf ) \- \oM^t_{\otau^t_n \land (t+\fs)} (\vf ) \big)  \b1_{ \wh{A} }  \big] \= 0  .
   \eea

  Define  $A^\oo_i   \df   \ccap{j=i+1}{k}   \big( \oW^{t_\oo}_{t+s_j   }, K^\oo_{t+s_j } \big)^{-1}  ( \cO_{j,\oo}  )   \ins \cF^{\oW^{t_\oo}}_{t+\fs   } $
 with   $\cO_{j,\oo} \df \big\{\fx \+ \fra_{\overset{}{\oo}} \n : \fx \ins \cO_j \big\} \ins \sB(\hR^{d+l})$.
 Since $ t_\oo \= \oga (\oo) \< \otau^t_n (\oo) \ls t\+n $ and since $ t_\oo    \ls t \+ s_{i+1} \ls t \+ \fs $,
  we   set $n_\loo $ as in Step (i)   and using \eqref{020522_17} with $(t,\bx,\oP,\vf,n,\fra,\oz_1,\oz_2) \= \big(t_\oo ,   \oX_{\oga \land \cd} (\oo), \oQ^\oo_\e ,\vf_\loo, n_\loo, \fra_{\overset{}{\oo}} ,    t\+n \ld \fs   ,    t\+n \ld \fr  \big)$ renders that
\bea
  0 & \tn \= & \tn  E_{\oQ^\oo_\e} \Big[ \Big(    \oMoo_{\oz^n_\loo   \land (t_\oo+n_\loo)   \land (t +n \land \fr)}(\vf_\loo)    \- \oMoo_{\oz^n_\loo   \land (t_\oo+n_\loo)   \land (t +n \land \fs )}(\vf_\loo) \Big) \b1_{A^\oo_i } \Big] \nonumber \\
 & \tn \= & \tn  E_{\oQ^\oo_\e} \Big[ \Big(    \oMoo_{\oz^n_\loo   \land (t +n \land \fr)}(\vf_\loo)
 \- \oMoo_{\oz^n_\loo     \land (t +n \land \fs )}(\vf_\loo) \Big) \b1_{A^\oo_i } \Big] .   \label{071820_25}
\eea
 Given $j \ins \{ i\+1,\cds \n , k \} $ and $\oo' \ins  \Wtgo  \Cp \big( \ocN^{\,\oo}_{\n X} \cp \ocN^{\,\oo}_K \big)^c$,
 we can derive that  $ (\oW^t  , \oX  ) ( t  \+ s_j,\oo')  \ins  \cO_j $ if and only if
 $ \big( \oW^{t_\oo} , \\ K^\oo \big)  ( t  \+ s_j, \oo' )
     \=   \big( \oW^t_{t  + s_j} ( \oo') \- \oW^t \big(\oga(\oo) ,\oo'\big),
\ofX^\oo_{t + s_j}   ( \oo' )  \big)
\=  (\oW^t  , \oX  ) ( t  \+ s_j,\oo') \+ \fra_{\overset{}{\oo}} \ins  \cO_{j,\oo}  $.
 So \eqref{070120_11} implies  that
 \bea
     \hspace{-1.2cm} \wh{A}   \Cp \Wtgo  \Cp \big( \ocN^{\,\oo}_{\n X} \cp  \ocN^{\,\oo}_{\n K} \big)^c
     \=  \oA_i \Cp  A^\oo_i   \Cp  \Wtgo \Cp \big( \ocN^{\,\oo}_{\n X} \cp  \ocN^{\,\oo}_{\n K} \big)^c   . \label{051020_33}
 \eea

 Let $\oo' \ins \oXi^t_{\oga,\oo}  $.
 Like in Step (i), we still have $  \otau^t_n(\oo')   \=   \oz^n_\loo (\oo')   \ld (t\+n)   $ since
  $  \oga(\oo) \< \otau^t_n (\oo) $.
   Taking $(r_1,r_2) \= \big(\otau^t_n(\oo') \ld (t\+\fs), \otau^t_n(\oo') \ld (t\+\fr) \big)$ in \eqref{091920_17} shows that
  $ \big(\oM^t \n (\vf )\big) \big(\otau^t_n(\oo') \ld (t\+\fr),\oo'\big) \- \big(\oM^t \n (\vf )\big) \big(\otau^t_n(\oo') \ld (t\+\fs),\oo'\big) \= \big( \oM^{t_\oo}  (\vf_\loo) \big) \big( \oz^n_\loo (\oo')   \ld (t\+n \ld \fr), \oo' \big) \- \big( \oM^{t_\oo}  (\vf_\loo) \big) \big( \oz^n_\loo (\oo')   \ld (t\+n \ld \fs), \oo' \big)  $.
   Then we can deduce from \eqref{021322_11}, \eqref{051020_33},  \eqref{072820_15},    \eqref{May06_25} and \eqref{071820_25}  that
 \beas
&& \hspace{-1.2cm}  E_{\oQ^\oo_\e} \big[ \big( \oM^{t,n}_{t+\fr} (\vf ) \- \oM^{t,n}_{t+\fs} (\vf ) \big)  \b1_{ \wh{A} }  \big]
\= E_{\oQ^\oo_\e} \big[ \b1_{\{\otau^t_n > \oga\}} \big( \oM^{t,n}_{t+\fr} (\vf ) \- \oM^{t,n}_{t+\fs} (\vf ) \big)  \b1_{ \wh{A} }  \big]
 \= E_{\oQ^\oo_\e} \big[  \b1_{\oA_i \cap  A^\oo_i} \b1_{\{ \otau^t_n  > \oga \}}\big( \oM^t_{\otau^t_n \land (t+\fr)} (\vf )
 \- \oM^t_{\otau^t_n \land ( t+\fs ) } (\vf ) \big)  \big] \\
&&  \hspace{-0.3cm}  \= \b1_{\{\oo \in \oA_i\}} \b1_{\{ \otau^t_n (\oo)  >  \oga (\oo) \}}  E_{\oQ^\oo_\e} \Big[  \Big(     \oMoo_{\oz^n_\loo \land (t +n \land \fr)   }(\vf_\loo)    \- \oMoo_{\oz^n_\loo  \land (t +n \land \fs)    }(\vf_\loo) \Big) \b1_{  A^\oo_i} \Big] \= 0 , \q   \fa \oo \ins \{t\+s_i \< \oga   \ls t\+s_{i+1}\} \Cp \ocA_* \Cp     \ocN^c_* ,
\eeas
 and thus
 $ \int_{\oo \in \ocA_*} \n  \b1_{ \{ t+s_i < \oga(\oo) \le t+s_{i+1}  \}}
       E_{\oQ^\oo_\e} \big[ \big( \oM^{t,n}_{t+\fr} (\vf ) \- \oM^{t,n}_{t+\fs} (\vf ) \big)  \b1_{ \wh{A} } \big]  \oP(d\oo) \= 0 $.
 Taking summation    from $i\=1$ through $i\=k \- 1$, we obtain from \eqref{071820_27}    that
$ E_{\oP_\e} \big[ \b1_{\{  \oga \le t+\fs \}} \big( \oM^{t,n}_{t+\fr} (\vf ) \- \oM^{t,n}_{t+\fs} (\vf ) \big)  \b1_{ \wh{A} }    \big]
    \=  E_\oP \big[ \b1_{\ocA^c_* \cap \{  \oga \le t+\fs\} \cap \wh{A}}  \big( \oM^{t,n}_{t+\fr} (\vf ) \- \oM^{t,n}_{t+\fs} (\vf ) \big)   \big]
\+  \int_{\oo \in \ocA_*}  \n   \b1_{ \{   \oga(\oo) \le t+\fs  \}}
       E_{\oQ^\oo_\e} \big[ \big( \oM^{t,n}_{t+\fr} (\vf ) \- \oM^{t,n}_{t+\fs} (\vf ) \big) \b1_{ \wh{A} } \big]  \oP(d\oo)
\= 0 $.
 Adding it to \eqref{071720_31} yield
 $ E_{\oP_\e} \big[ \big( \oM^{t,n}_{t+\fr} (\vf ) \- \oM^{t,n}_{t+\fs} (\vf ) \big)  \b1_{\wh{A}} \big] \= 0  $.

 Hence, $ \oP_\e \ins \ocP^2_t $.
 According to Part \(2a\) of the proof of Proposition \ref{prop_Ptx_char}, $ \oP_\e $  satisfies \(D1\) and \(D2\) in Definition  \ref{def_ocP}   of $\ocP_{t,\bx}$.

\no {\bf II.d)} In this part, we show that
 $ \oP_\e \big\{  \oT  \=  \wh{\tau}_\e (\oW)   \big\} \= 1$ for some   $[t,\infty]-$valued   $ \bF^{W^t,P_0} -$stopping time $\wh{\tau}_\e$, i.e.,  $ \oP_\e $ satisfies (D3) in Definition  \ref{def_ocP}  of $   \ocP_{t,\bx} $.

\no {\bf II.d.1)} For any $s \ins [t,\infty)$,   there is a $[0,1]-$valued $ \cF^{W^t}_s -$measurable random variable  $\vth^\e_s$  on $\O_0$   such that
\bea    \label{J30_03}
\vth^\e_s \big( \oW  (\oo) \big) \=  E_{\oP_\e}   \Big[ \b1_{\{\oga \le s \}}   \b1_{\{\oT \in [\oga, s ]\}} \big| \cF^{\oW^t}_{\n s} \Big] ( \oo ),
\q \fa   \oo \ins \oO .
\eea
  Since $\oW^t$ is a Brownian motion under $\oP_\e$ by Part (II.c),
  applying  Lemma \ref{lem_122921_11} with $t_0 \= t$, $(\O_1, \cF_1, P_1,B^1)   \= \big(\oO ,  \sB(\oO ),  \oP_\e, \oW \big) $, $(\O_2, \cF_2, P_2,B^2) \= \big(\O_0,  \sB(\O_0),  P_0,  W \big) $  and $\Phi \= \oW$ yields that
  $\big\{ \vth^\e_s (\oW) \big\}_{s \in [t,\infty)}$ is an $\bF^{\oW^t}-$adapted process   and that
    $E_{P_0} [\vth^\e_s] \= E_{\oP_\e} \big[ \vth^\e_s (\oW) \big]
 \= E_{\oP_\e} \big[  \b1_{\{\oga \le s \}}  \b1_{\{\oT \in [\oga, s ]\}}   \big] $ is right-continuous in $s \ins [t,\infty)$.
 As  $\bF^{W^t,P_0}$ is a right-continuous complete filtration,
the process $ \{\vth^\e_s\}_{s \in [t,\infty)}$ admits a $[0,1]-$valued $\bF^{W^t,P_0}-$adapted \cad modification $ \big\{\wh{\vth}^\e_s\big\}_{s \in [t,\infty)} $.  
 Define a  $[t,\infty]-$valued $\bF^{W^t,P_0}-$stopping time by
 \bea \label{022022_17}
  \btau_{\n \e} (\o_0) \df \inf \big\{s \ins [t,\infty) \n :  \wh{\vth}^\e_s (\o_0) \= 1 \big\}  .
 \eea

  As $\oW^t$ is also a Brownian motion under $\oP_\e$ by Part  (II.c),
  using Lemma \ref{lem_122921_11}   with $t_0 \= t$, $(\O_1, \cF_1, P_1,B^1)   \= \big(\oO ,  \sB(\oO ),  \oP_\e, \oW \big) $, $(\O_2, \cF_2, P_2,B^2) \= \big(\O_0,  \sB(\O_0),  P_0,  W \big) $  and $\Phi \= \oW$ implies that
    $ \wh{\tau}(\oW)$ and $ \btau_{\n \e}(\oW)$ are  $[t,\infty]-$valued
 $\bF^{\oW^t,\oP_\e}-$stopping times. Then
 \beas
 \otau_\e \df \wh{\tau}(\oW) \b1_{\{\wh{\tau}(\oW)  < \oga \}} \+ \big(\btau_{\n \e} (\oW) \ve \oga \big)  \b1_{\{\wh{\tau}(\oW) \ge \oga \}}
 \eeas
 is also a   $[t,\infty]-$valued $\bF^{\oW^t,\oP_\e}-$stopping time.
 According to Lemma \ref{lem_012922_11}, there exists two $[t,\infty]-$valued $\sB(\oO)-$measurable random variables $\oxi$ and $\oxi_\e$ such that
 \bea \label{022022_21}
  \big\{\wh{\tau} (\oW) \nne \oxi \big\}   \ins \sN_\oP \big(\cF^{\oW^t}_\infty\big) \Cp \sN_{\oP_\e} \big(\cF^{\oW^t}_\infty\big)
  \aand \big\{\btau_{\n \e} (\oW) \nne \oxi_\e \big\}   \ins \sN_\oP \big(\cF^{\oW^t}_\infty\big) \Cp \sN_{\oP_\e} \big(\cF^{\oW^t}_\infty\big) .
  \eea
 We can also update \eqref{020322_11} to:
 \bea \label{020322_11b}
 \big\{ \wh{\tau}(\oW) \gs \oga \big\} \D \ocA_* \ins \sN_\oP \big(\cF^{\oW^t}_\infty\big) \Cp \sN_{\oP_\e} \big(\cF^{\oW^t}_\infty\big) .
 \eea

 \if{0}

 Note $ \otau_\e \df (\wh{\tau}(\oW) \ld \oga) \b1_{\{\wh{\tau}(\oW)  < \oga \}} \+ \big(\btau_{\n \e} (\oW) \ve \oga \big)  \b1_{\{\wh{\tau}(\oW) \ge \oga \}} $. For $s \ins [t,\infty)$, since $ \{\wh{\tau}(\oW)  \< \oga \} \ins \cF^{\oW^t, \oP_\e}_{\wh{\tau}(\oW) \land \oga}$, $\{\wh{\tau}(\oW)  \< \oga \} \Cp \big\{  \wh{\tau}(\oW) \ld \oga  \ls s \big\} \ins \cF^{\oW^t, \oP_\e}_s $ and $\{\wh{\tau}(\oW)  \gs \oga \}   \Cp \big\{   \oga  \ls s \big\} \ins \cF^{\oW^t, \oP_\e}_s$.
 \beas
 \{\otau_\e \ls s\} & \tn \= & \tn  \Big(\{\wh{\tau}(\oW)  \< \oga \} \Cp \big\{  \wh{\tau}(\oW) \ld \oga  \ls s \big\}\Big) \cp
 \Big(\{\wh{\tau}(\oW)  \gs \oga \} \Cp \big\{  \btau_{\n \e} (\oW) \ve \oga  \ls s \big\}\Big) \\
  & \tn \= & \tn  \Big(\{\wh{\tau}(\oW)  \< \oga \} \Cp \big\{  \wh{\tau}(\oW) \ld \oga  \ls s \big\}\Big) \cp
 \Big(\{\wh{\tau}(\oW)  \gs \oga \} \Cp \big\{  \btau_{\n \e} (\oW)    \ls s \big\} \Cp \big\{   \oga  \ls s \big\}\Big) \ins \cF^{\oW^t, \oP_\e}_s .
 \eeas

 \fi

 Since $  \oQ^\oo_\e \big(\ocA^c_* \Cp\{\oT \= \oxi \}\big) \= \b1_{ \{\oo \in \ocA^c_*  \}}  \oQ^\oo_\e \{\oT \= \oxi \} \= 0 $, $\fa \oo \ins \ocA_* \Cp \ocN^c_*  $ by  \eqref{May06_25},
  one has $ \oP_\e \big( \ocA^c_* \Cp\{\oT \= \oxi \} \big) \= \oP  \big( \ocA^c_* \Cp\{\oT \= \oxi \} \big) $.
 It follows from \eqref{020322_11b}, \eqref{022022_21} and  $\oP \big\{\oT \= \wh{\tau}(\oW) \big\}   \= 1$ that
 \bea
   \oP_\e \big( \ocA^c_* \Cp\{\oT \= \otau_\e\} \big)
   & \tn \= & \tn  \oP_\e \big(  \{ \wh{\tau}(\oW) \< \oga  \} \Cp\{\oT \= \otau_\e\} \big)
   \= \oP_\e \big(  \{ \wh{\tau}(\oW) \< \oga  \} \Cp\{\oT \= \wh{\tau}(\oW)\} \big)
   \= \oP_\e \big(  \ocA^c_* \Cp\{\oT \= \oxi\} \big) \nonumber \\
    & \tn \= & \tn  \oP  \big( \ocA^c_* \Cp\{\oT \= \oxi \} \big)
   \= \oP  \big( \ocA^c_* \Cp\{\oT \= \wh{\tau} (\oW) \} \big) \= \oP  \big( \ocA^c_* \big)   . \label{022022_31}
   \eea

  \no {\bf II.d.2)} We next show that   $ \oP_\e \big( \ocA_* \Cp\{\oT \= \otau_\e\} \big)  \= \oP(\ocA_*) $ and thus $ \oP_\e\{\oT \= \otau_\e\} \= 1 $.

  As $  \big\{\btau_{\n \e} (\oW) \nne \oxi_\e \big\}   \ins   \sN_{\oP_\e} \big(\cF^{\oW^t}_\infty\big) $, 
 there exists   $\oA^\e_\xi \ins \cF^{\oW^t}_\infty \sb \sB(\oO) $ such that $ \big\{\btau_{\n \e} (\oW) \nne \oxi_\e \big\} \sb \oA^\e_\xi$ and $ \oP_\e \big(\oA^\e_\xi\big) \= 0 $.
By \eqref{022022_23}, the random variable  $   \oo \mto \oQ^\oo_\e \big(\oA^\e_\xi\big)$ is   $ \si \Big( \cF^{\oW^t}_\oga \cp \sN_\oP \big(\ocF^t_{\n \infty}\big) \Big)  -$measurable. Since
 $0 \ls  \int_{\oo \in \ocA_*} \oQ^\oo_\e \big(\oA^\e_\xi\big) \oP(d\oo) \ls \oP_\e \big(\oA^\e_\xi\big) \= 0 $,
   we can find   $\ocN^\e_\xi \ins \sN_\oP \big(\ocF^t_{\n \infty}\big)$
 such that
\bea \label{022022_25}
\oQ^\oo_\e \big(\oA^\e_\xi\big) \= 0  \hb{ and thus }   \oQ^\oo_\e \big\{\btau_{\n \e} (\oW) \nne \oxi_\e \big\} \= 0  , \q \fa \oo \ins \ocA_* \Cp \big( \ocN^\e_\xi \big)^c .
\eea

  Let $s \ins \hQ \cap [t,\infty)$ and pick    a countable Pi-system $   \big\{\ocO_j\big\}_{j \in \hN}$ that generates $ \cF^{\oW^t}_s $.
  We also let  $j \ins \hN$ and $\oA  \ins \cF^{\oW^t}_\oga$.
  As $\ocA_* \cap \oA \cap \{\oga \le s\} \ins \cF^{\oW^t}_s$,  it holds $\oP_\e-$a.s. that
   $ \b1_{ \ocA_* \cap \oA \cap \ocO_j} \vth^\e_s (\oW)
   \= \b1_{\ocA_* \cap \oA \cap \ocO_j \cap \{\oga \le s\}}  E_{\oP_\e} \Big[    \b1_{\{\oT \in [\oga,s]\}} \big| \cF^{\oW^t}_{\n s} \Big]
    \=   E_{\oP_\e} \Big[ \b1_{\ocA_* \cap \oA \cap \ocO_j }    \b1_{\{\oga \le s \}}  \b1_{\{\oT \in [\oga,s]\}} \big| \cF^{\oW^t}_{\n s} \Big] $.
 Then   \eqref{May06_25} and \eqref{J30_03} imply  that
 \beas
 && \hspace{-1.2cm} \int_{\oo \in \oO}  \b1_{\{\oo \in \ocA_* \cap \oA   \}}    E_{\oQ^\oo_\e} \big[  \b1_{ \ocO_j} \vth^\e_s (\oW) \big]  \oP(d\oo)
 \= E_{\oP_\e} \big[ \b1_{\ocA_* \cap \oA \cap \ocO_j} \vth^\e_s (\oW) \big]
 \= E_{\oP_\e} \Big[  E_{\oP_\e} \Big[   \b1_{\ocA_* \cap \oA \cap \ocO_j}  \b1_{\{\oga \le s \}}  \b1_{\{\oT \in [\oga,s]\}} \big| \cF^{\oW^t}_{\n s} \Big]  \Big] \nonumber \\
 &&   \= E_{\oP_\e} \big[  \b1_{\ocA_* \cap \oA \cap \ocO_j}  \b1_{\{\oga \le s \}}  \b1_{\{\oT \in [\oga,s]\}} \big]
  \= \int_{\oo \in \oO}  \b1_{\{\oo \in \ocA_* \cap \oA   \}}    E_{\oQ^\oo_\e} \big[  \b1_{ \ocO_j}   \b1_{\{\oga \le s \}}  \b1_{\{\oT \in [\oga,s]\}} \big]  \oP(d\oo)     . 
\eeas
 So $ \cF^{\oW^t}_\oga \cp \sN_\oP\big(\sB\big(\oO\big)\big)$ is contained in the Lambda-system 
$ \ol{\L}^\e_{s,j}   \df   \big\{\oA \ins \sB_\oP(\oO) \n :  \int_{\oo \in \oO}  \b1_{\{\oo \in \ocA_* \cap \oA   \}}    E_{\oQ^\oo_\e} \big[  \b1_{ \ocO_j} \vth^\e_s (\oW) \big]  \oP(d\oo) \=   \int_{\oo \in \oO}  \b1_{\{\oo \in \ocA_* \cap \oA   \}}    E_{\oQ^\oo_\e} \big[  \b1_{ \ocO_j}   \b1_{\{\oga \le s \}}  \b1_{\{\oT \in [\oga,s]\}} \big]  \oP(d\oo)  \big\} $.
As $  \cF^{\oW^t}_\oga \cp \sN_\oP \big(\sB\big(\oO\big)\big)   $ is closed under intersection,
 Dynkin's Pi-Lambda Theorem shows that
$\si \Big( \cF^{\oW^t}_\oga \cp \sN_\oP \big(\sB\big(\oO\big)\big) \Big) \sb \ol{\L}^\e_{s,j}$, i.e.,  for any $\oA \ins \si \Big( \cF^{\oW^t}_\oga \cp \sN_\oP \big(\sB\big(\oO\big)\big) \Big)$
\bea
\hspace{-0.3cm}
  \int_{\oo \in \oO}  \b1_{\{\oo \in \ocA_* \cap \oA   \}}    E_{\oQ^\oo_\e} \big[  \b1_{ \ocO_j} \vth^\e_s (\oW) \big]  \oP(d\oo)   \= \n \int_{\oo \in \oO}  \b1_{\{\oo \in \ocA_* \cap \oA   \}}    E_{\oQ^\oo_\e} \big[  \b1_{ \ocO_j}   \b1_{\{\oga \le s \}}  \b1_{\{\oT \in [\oga,s]\}} \big]  \oP(d\oo)   .  \q  \label{070520_29}
\eea

 Since $\b1_{ \ocO_j}   \vth^\e_s (\oW)  $ and $  \b1_{ \ocO_j}  \b1_{\{\oga \le s \}}  \b1_{\{\oT \in [\oga,s]\}}$ are $\sB(\oO)-$measurable,
 we see from \eqref{022022_23}   that
    the random variables $ \oO \ni \oo \mto E_{\oQ^\oo_\e} \big[  \b1_{ \ocO_j} \vth^\e_s (\oW) \big] $
    \if{0}

  To use \eqref{022022_23}, we can actually take $ \vth^\e_s \big( \oW  (\oo) \big) \=  E_{\oP_\e}   \Big[ \b1_{\{\oga \le s \}}   \b1_{\{\oT \in [\oga, s ]\}} \big| \cF^{\oW^t}_{\n s} \Big] ( \oo ) \ve 0 \gs 0 $, $ \fa   \oo \ins \oO  $

    \fi
 and $  \oO \ni  \oo \mto E_{\oQ^\oo_\e} \big[  \b1_{ \ocO_j}    \b1_{\{\oga \le s \}}  \b1_{\{\oT \in [\oga,s]\}}  \big]  $
 are $ \si \Big( \cF^{\oW^t}_\oga \cp \sN_\oP \big(\ocF^t_{\n \infty}\big) \Big) -$measurable.
  Letting $\oA$ vary over $\si \Big( \cF^{\oW^t}_\oga \cp \sN_\oP \big(\ocF^t_{\n \infty}\big) \Big)$ in \eqref{070520_29}
   yields  that
  $\b1_{ \{\oo \in \ocA_* \}}    E_{\oQ^\oo_\e} \big[  \b1_{ \ocO_j} \vth^\e_s (\oW) \big]
\= \b1_{ \{\oo \in \ocA_* \}}  E_{\oQ^\oo_\e} \big[  \b1_{ \ocO_j}   \b1_{\{\oga \le s \}}  \b1_{\{\oT \in [\oga,s]\}} \big]$
 for all $\oo \ins \oO$ except on some $ \ocN^\e_{s,j} \ins \sN_\oP \big(\ocF^t_{\n \infty}\big)$. It then follows from
   \eqref{022122_11} that
 \bea
   E_{\oQ^\oo_\e} \big[  \b1_{ \ocO_j} \vth^\e_s (\oW) \big]
\=    E_{\oQ^\oo_\e} \big[  \b1_{ \ocO_j}  \b1_{\{\oga \le s \}}   \b1_{\{\wh{\tau}_\loo(\oW) \in [\oga,s]\}} \big] , \q \fa \oo \in \ocA_* \Cp \ocN^c_* \Cp \big(\ocN^\e_{s,j}\big)^c.     \label{070620_23}
\eea
 By Lemma \ref{lem_122921_11}, $  \ocN^\vth_s \df \oW^{-1} \big(   \big\{ \wh{\vth}^\e_s \nne  \vth^\e_s \big\} \big) $ belongs to $ \sN_{\oP_\e} (\cF^{\oW^t}_\infty) $.
 An analogy to \eqref{022022_25} gives   $\ocN^\e_{s,0} \ins \sN_\oP \big(\ocF^t_{\n \infty}\big)$
 such that $ \oQ^\oo_\e \big(\ocN^\vth_s\big) \= 0  $ for any $\oo \ins \ocA_* \Cp \big( \ocN^\e_{s,0} \big)^c $.

 Set $\ocN^\e_* \df \ocN_* \cp \Big(\ccup{s \in \hQ \cap [t,\infty) }{} \ccup{j =0}{\infty}  \ocN^\e_{s,j}\Big)  \cp   \ocN^\e_\xi \ins \sN_\oP \big(\ocF^t_{\n \infty}\big)  $.
 We fix $\oo \ins \ocA_* \Cp \big( \ocN^\e_* \big)^c $
 and let $s \ins \hQ \Cp [t,\infty)$.

 \bul When $s \< t_\oo  $: Since \eqref{070620_23}, \eqref{090520_11} and \eqref{072820_15}  show that
 $ E_{\oQ^\oo_\e} \big[  \b1_{ \ocO_j} \vth^\e_s (\oW) \big]
 \=    E_{\oQ^\oo_\e} \big[ \b1_\Wtgo \b1_{ \ocO_j}  \b1_{\{t_\oo \le s \}}   \b1_{\{\wh{\tau}_\loo(\oW) \in [t_\oo, s ]\}} \big]
 \=  0 $ for any $ j \ins \hN$,
   Dynkin's Pi-Lambda Theorem implies that  $ E_{\oQ^\oo_\e} \big[  \b1_{ \ocE} \vth^\e_s (\oW) \big]
 \=   0 $ for any $ \ocE \ins \cF^{\oW^t}_s $.  Letting $\cE$ vary over $\cF^{\oW^t}_s$ reaches that
 $ \vth^\e_s \big(\oW(\oo')\big) \= 0 $ for all $\oo' \ins \oO$ except on
 some  $ \ofN^\oo_{s,1} \ins \sN_{\oQ^\oo_\e} \big(\cF^{\oW^t}_\infty\big) $.

 \bul When $s \gs t_\oo  $:
 Applying  Lemma \ref{lem_122921_11} with $t_0 \= t_\oo$, $(\O_1, \cF_1, P_1,B^1)   \= \big(\oO ,  \sB(\oO ),  \oQ^\oo_\e, \oW \big) $, $(\O_2, \cF_2, P_2,B^2) \= \big(\O_0,  \sB(\O_0),   P_0,  W \big) $  and $\Phi \= \oW$ yields that
   $\wh{\tau}_\loo(\oW)$ is a $[t_\oo,\infty]-$valued   $ \bF^{\oW^{t_\oo},\oQ^\oo_\e} -$stopping time and thus
   $\big\{\wh{\tau}_\loo(\oW) \ins [\oga \ld s, s ]\big\} \ins   \cF^{\oW^{t_\oo},\oQ^\oo_\e}_s $.
     By    Problem 2.7.3 of \cite{Kara_Shr_BMSC},
  there is   $\ocA^\oo_s   \ins  \cF^{\oW^{t_\oo} }_s $   such that
  $ \ofN^\oo_{s,2} \df  \ocA^\oo_s \, \D \big\{\wh{\tau}_\loo(\oW) \ins [\oga \ld s, s ]\big\}   \ins  \sN_{\oQ^\oo_\e} \big(\cF^{\oW^{t_\oo} }_\infty \big)  $.
  Then  we see from \eqref{070620_23}   that
 $ E_{\oQ^\oo_\e} \big[  \b1_{ \ocO_j} \vth^\e_s (\oW) \big]
 \=    E_{\oQ^\oo_\e} \big[  \b1_{ \ocO_j}  \b1_{\{\oga \le s \}}   \b1_{\{\wh{\tau}_\loo(\oW) \in [\oga \land s, s ]\}} \big]
 \=    E_{\oQ^\oo_\e} \big[  \b1_{ \ocO_j}  \b1_{\{\oga \le s \}}   \b1_{\ocA^\oo_s} \big] $ for any $ j \ins \hN$,
 and we know from Dynkin's Pi-Lambda Theorem   that  $ E_{\oQ^\oo_\e} \big[  \b1_{ \ocE} \vth^\e_s (\oW) \big]
 \=    E_{\oQ^\oo_\e} \big[  \b1_{ \ocE}  \b1_{\{\oga \le s \}} \b1_{\ocA^\oo_s} \big] $ for any $ \ocE \ins \cF^{\oW^t}_s $.
  As   $\cF^{\oW^{t_\oo} }_s \= \si \big( \oW^{t_\oo}_r; r \ins [t_\oo,s] \big)
 \= \si \big( \oW^t_r \- \oW^t_{t_\oo} ; r \ins [t_\oo,s] \big) \sb \si \big( \oW^t_r   ; r \ins [t,s] \big) \= \cF^{\oW^t}_s$,
 letting $\cE$ run over $\cF^{\oW^t}_s$ renders   that
 $ \vth^\e_s \big(\oW(\oo')\big) \= \b1_{\{\oga (\oo') \le s \}} \b1_{\{\oo' \in  \ocA^\oo_s\}} $ for all $\oo' \ins \oO$ except on
 some  $ \ofN^\oo_{s,3} \ins \sN_{\oQ^\oo_\e} \big(\cF^{\oW^t}_\infty\big) $.

   Let $\oo' \ins \Wtgo \Cp \Big( \ccup{s \in \hQ \cap [t,\infty) }{} \ocN^\vth_s   \Big)^c \Cp \Big( \ccup{s \in \hQ \cap [t,t_\oo) }{} \ofN^\oo_{s,1}   \Big)^c \Cp \Big( \ccup{s \in \hQ \cap [t_\oo,\infty) }{}   \ofN^\oo_{s,2} \cp  \ofN^\oo_{s,3} \Big)^c \Cp \big\{\btau_{\n \e} (\oW) \nne \oxi_\e \big\}$.
 The above analysis   and   \eqref{090520_11} show that
  \beas
   \wh{\vth}^\e_s \big(\oW(\oo')\big) \=  \vth^\e_s \big(\oW(\oo')\big) \n   \= \b1_{\{s \ge t_\oo\}}   \b1_{\{\oo' \in  \ocA^\oo_s\}}
   \n \= \b1_{\{s \ge t_\oo\}} \b1_{ \big\{\wh{\tau}_\loo(\oW(\oo')) \in  [\oga (\oo') \land s, s ]  \big\}}
   \n \=  \b1_{\{s \ge t_\oo\}} \b1_{ \big\{\wh{\tau}_\loo(\oW(\oo')) \in [ t_\oo, s ]  \big\}}   , ~ \fa s \ins \hQ \Cp [t,\infty) .
   \eeas
  So the right-continuity of process $\wh{\vth}^\e_\cd$ gives that
  $ \wh{\vth}^\e_s \big(\oW(\oo')\big)
  \=  \b1_{\{s \ge t_\oo\}} \b1_{ \big\{\wh{\tau}_\loo(\oW(\oo')) \in  [ t_\oo, s ]  \big\}}  $, $ \fa s \ins   [t,\infty)  $,
  and we  can deduce from \eqref{022022_17} that
  \beas
   \oxi_\e (\oo') \= \btau_{\n \e} (\oW(\oo')) \= \inf \n \big\{s \ins [t,\infty) \n :  \wh{\vth}^\e_s (\oW(\oo')) \= 1 \big\} \= \wh{\tau}_\loo(\oW(\oo'))
  \gs   t_\oo 
  \=  \oga(\oo') .
  \eeas
  In particular, one has  $ \oxi_\e (\oo') \ve \oga(\oo') \= \wh{\tau}_\loo \big(\oW(\oo')\big) $,
  which together with    \eqref{072820_15}, \eqref{022022_25} and \eqref{022122_11} implies that
 $ 1   \=    \oQ^\oo_\oo \big\{ \oxi_\e \ve \oga \= \wh{\tau}_\loo(\oW)\big\}
 \=    \oQ^\oo_\oo \big\{ \oxi_\e \ve \oga \= \oT\big\} $, $ \fa \oo \ins \ocA_* \Cp \big( \ocN^\e_* \big)^c   $.
 Then  \eqref{020322_11b},  \eqref{022022_21} and \eqref{May06_25} render that
 \beas
  \oP_\e \big( \ocA_* \Cp\{\oT  \=   \otau_\e\} \big)  & \tn \=  & \tn \oP_\e \big(  \{ \wh{\tau}(\oW) \gs \oga  \} \Cp\{\oT \= \otau_\e\} \big)
  \= \oP_\e \big(  \{ \wh{\tau}(\oW) \gs \oga  \} \Cp\{\oT \= \btau_{\n \e} (\oW) \ve \oga\} \big) \\
   & \tn  \=  & \tn  \oP_\e \big(  \ocA_* \Cp \{\oT \= \oxi_\e \ve \oga\} \big)
   \= \int_{\oo \in  \ocA_*} \oQ^\oo_\e \{\oT \= \oxi_\e \ve \oga\} \oP(d \oo)     \= \oP(\ocA_*) .
  \eeas
 Adding it to \eqref{022022_31} yields  $ \oP_\e \{\oT \= \otau_\e\} \= 1 $.
 Moreover,  applying Lemma \ref{lem_M31_01} (2) with $(\O , \cF , P ,B ) \= \big(\oO ,  \sB(\oO ),  \oP_\e, \oW \big) $  and $\Phi \= \oW$,
 we can find   a $[t,\infty]-$valued $\bF^{W^t,P_0}-$stopping time $\wh{\tau}_\e$ on $\O_0$ such that  $\otau_\e \= \wh{\tau}_\e (\oW)$,
 $\oP_\e-$a.s. Hence, $ \oP_\e $ satisfies (D3) in Definition  \ref{def_ocP}  of   $\ocP_{t,\bx}$.

 \no {\bf II.e)} Fix $i \ins \hN$.  Since $ \big\{  \int_t^s g_i  (r, \osX^{t,\bx}_{r \land \cd}  ) dr  \big\}_{s \in [t,\infty)} $
 and $ \big\{  \int_t^s h_i  (r, \osX^{t,\bx}_{r \land \cd}  ) dr  \big\}_{s \in [t,\infty)} $ are two
 $ \bF^{\oW^t,\oP}  - $adapted continuous processes,
  Lemma 2.4 of \cite{STZ_2011a} assures two   $  \bF^{\oW^t }  -$predictable  processes  $  \big\{ \ol{\Phi}^i_s  \big\}_{s \in [t,\infty)}$ and
  $  \big\{ \ol{\Psi}^i_s  \big\}_{s \in [t,\infty)}$
   such that   $ \ocN^{\, i,1}_{g,h} \df \Big\{\oo \ins \oO \n : \ol{\Phi}^i_s (\oo) \nne \int_t^s g_i  \big(r, \osX^{t,\bx}_{r \land \cd} (\oo) \big) dr \hb{ or } \ol{\Psi}^i_s (\oo) \nne \int_t^s h_i  \big(r, \osX^{t,\bx}_{r \land \cd} (\oo) \big) dr \hb{ for some }   s \ins [t,\infty) \Big\} \ins  \sN_\oP \big(\cF^{\oW^t}_\infty \big)$.
  \if{0}

 Since $ \big\{  \int_t^s g^\pm_i  (r, \osX^{t,\bx}_{r \land \cd}  ) dr  \big\}_{s \in [t,\infty)} $
 and $ \big\{  \int_t^s h^\pm_i  (r, \osX^{t,\bx}_{r \land \cd}  ) dr  \big\}_{s \in [t,\infty)} $ are   $[0,\infty]-$valued
 $ \bF^{\oW^t,\oP}  - $adapted continuous processes,  Lemma 2.4 of \cite{STZ_2011a} assures   $[0,\infty]-$valued  $  \bF^{\oW^t }  -$predictable  processes  $  \big\{ \ol{\Phi}^{i,\pm}_s  \big\}_{s \in [t,\infty)}$ and
  $  \big\{ \ol{\Psi}^{i,\pm}_s  \big\}_{s \in [t,\infty)}$
   such that all $\oo \ins \oO$ except on some $ \ocN^\pm_{\n g_i,h_i} \ins  \sN_\oP \big(\cF^{\oW^t}_\infty \big)$,
   $\ol{\Phi}^{i,\pm}_s (\oo) \= \int_t^s g^\pm_i  \big(r, \osX^{t,\bx}_{r \land \cd} (\oo) \big) dr$
   and $\ol{\Psi}^{i,\pm}_s (\oo) \= \int_t^s h^\pm_i  \big(r, \osX^{t,\bx}_{r \land \cd} (\oo) \big) dr$ for any $ s \ins [t,\infty) $.
   So for any $\oo \ins \big( \ocN^+_{\n g_i,h_i} \cp \ocN^-_{\n g_i,h_i}\big)^c$,
    \beas
    \ol{\Phi}^i_s (\oo) \df \ol{\Phi}^{i,+}_s (\oo) \- \ol{\Phi}^{i,-}_s (\oo) \= \int_t^s g_i  \big(r, \osX^{t,\bx}_{r \land \cd} (\oo) \big) dr
    \aand
   \ol{\Psi}^i_s (\oo) \df \ol{\Psi}^{i,+}_s (\oo) \- \ol{\Psi}^{i,-}_s (\oo) \= \int_t^s h_i  \big(r, \osX^{t,\bx}_{r \land \cd} (\oo) \big) dr , \q \fa s \ins [t,\infty) .
    \eeas

     \fi
 By Remark \ref{rem_ocP2} (1),
 $E_\oP \big[ \int_t^\infty   \n g^-_i  ( r,  \oX_{r \land \cd}   ) \ve h^-_i  ( r,  \oX_{r \land \cd}  )   dr \big] \< \infty$.
 So  it holds for any $\oo \ins \oO$ except on some $ \ocN^{\, i,2}_{g,h} \ins \sN_\oP \big(\cF^\oX_\infty\big) $
  that $ \int_t^\infty   \n g^-_i  \big( r,  \oX_{r \land \cd} (\oo)   \big) \ve h^-_i  \big( r,  \oX_{r \land \cd} (\oo)  \big)   dr \< \infty$.

  For any $\oo \ins \ocA_*   \Cp  \ocN^c_* \Cp \ocN^c_{\n X} \Cp \big(  \ocN^{\, i,1}_{g,h} \cp \ocN^{\, i,2}_{g,h} \big)^c   $,
 since $ \oO^t_{\oga,\oo} \sb \big\{ \oo' \ins \oO \n:  \oX_s  (\oo') \=  \oX_s   (\oo)   ,
    \fa s \ins \big[0,\oga(\oo)\big]   \big\} $,
     \eqref{072820_15} and \eqref{April07_11} show that
  \bea
    E_{\oQ^\oo_\e} \Big[ \int_t^\oT   \n g_i \big( r,  \oX_{r \land \cd}   \big)   dr \Big]
     & \tn \=  & \tn     \int_t^{\oga(\oo)}   \n g_i \big(   r,  \oX_{r \land \cd} (\oo)  \big)   dr
    \+ E_{\oQ^\oo_\e} \Big[  \int_{\oga(\oo)}^\oT   \n g_i \big( r,  \oX_{r \land \cd}   \big)   dr   \Big] \label{091020_15} \\
     & \tn  \ls  & \tn   \int_t^{\oga(\oo)}   \n g_i \big(   r,  \osX^{t,\bx}_{r \land \cd} (\oo)  \big)   dr
    \+   \big( \oY^i_{\n \oP}  ( \oga ) \big) (\oo)
    \= \ol{\Phi}^i_\oga (\oo) \+  E_\oP \Big[ \int_{\oT \land  \oga }^\oT    g_i (r,\oX_{r \land \cd} ) dr \Big| \cF^{\oW^t}_\oga \Big] (\oo)
    \nonumber
    \eea
    and similarly that
     $  E_{\oQ^\oo_\e} \Big[ \int_t^\oT   \n h_i \big( r,  \oX_{r \land \cd}   \big)   dr \Big]
     \=     \int_t^{\oga(\oo)}   \n h_i \big(   r,  \oX_{r \land \cd}   (\oo)  \big)   dr
     \+ E_{\oQ^\oo_\e} \Big[  \int_{\oga(\oo)}^\oT   \n h_i \big( r,  \oX_{r \land \cd}   \big)   dr   \Big]
     \=  \int_t^{\oga(\oo)}   \n h_i \big(   r,  \osX^{t,\bx}_{r \land \cd} (\oo)  \big)   dr \\
     \+   \big( \oZ^i_{\n \oP}  ( \oga ) \big) (\oo)
    \= \ol{\Psi}^i_\oga (\oo) \+  E_\oP \Big[ \int_{\oT \land  \oga }^\oT    h_i (r,\oX_{r \land \cd} ) dr \Big| \cF^{\oW^t}_\oga \Big] (\oo) $.
     \if{0}

    \beas
     E_{\oQ^\oo_\e} \Big[ \int_t^\oT   \n g_i \big( r,  \oX_{r \land \cd}   \big)   dr \Big]
     \= \int_{\oo' \in \oO} \b1_{\{\oo' \in \oO^t_{\oga,\oo}\}}\Big(  \int_t^{\oga(\oo)}   \n g_i \big( r,  \oX_{r \land \cd} (\oo')  \big)   dr \+ \int_{\oga(\oo)}^{\oT(\oo')}   \n g_i \big( r,  \oX_{r \land \cd} (\oo')  \big)   dr \Big) \oQ^\oo_\e (d \oo') \\
     \= \int_{\oo' \in \oO} \b1_{\{\oo' \in \oO^t_{\oga,\oo}\}}\Big(  \int_t^{\oga(\oo)}   \n g_i \big( r,  \oX_{r \land \cd} (\oo)  \big)   dr \+ \int_{\oga(\oo)}^{\oT(\oo')}   \n g_i \big( r,  \oX_{r \land \cd} (\oo')  \big)   dr \Big) \oQ^\oo_\e (d \oo') \\
     \= \int_{\oo' \in \oO}  \Big(  \int_t^{\oga(\oo)}   \n g_i \big( r,  \oX_{r \land \cd} (\oo)  \big)   dr \+ \int_{\oga(\oo)}^{\oT(\oo')}   \n g_i \big( r,  \oX_{r \land \cd} (\oo')  \big)   dr \Big) \oQ^\oo_\e (d \oo') \\
    \=    \int_t^{\oga(\oo)}   \n g_i \big(   r,  \oX_{r \land \cd} (\oo)  \big)   dr
    \+ E_{\oQ^\oo_\e} \big[  \int_{\oga(\oo)}^\oT   \n g_i \big( r,  \oX_{r \land \cd}   \big)   dr   \big] ,
    \eeas

     \fi
     Since $ \ocA_* ,  \ol{\Phi}^i_\oga  \ins \cF^{\oW^t}_\oga$ and
  since $\b1_{\ocA_*}  \= \b1_{\{ \wh{\tau}(\oW)   \ge   \oga \}} \=  \b1_{\{\oT    \ge   \oga \}}   $, $\oP-$a.s.  by \eqref{020322_11},
  we can deduce from the tower property that
    \beas
 && \hspace{-1.5cm} \int_{\oo \in \ocA_*}   E_{\oQ^\oo_\e} \Big[  \int_t^\oT   g_i \big(r,\oX_{r \land \cd}  \big) dr \Big]     \oP(d\oo)
     \ls  
     E_\oP \bigg[    E_\oP \Big[ \b1_{\ocA_*} \Big( \ol{\Phi}^i_\oga
    \+ \int_{\oT \land \oga}^\oT    g_i(r,\oX_{r \land \cd} ) dr \Big) \Big| \cF^{\oW^t}_\oga \Big]   \bigg] \\
   &&  \=    E_\oP \Big[ \b1_{\ocA_*} \Big( \int_t^\oga  g_i  \big(r, \osX^{t,\bx}_{r \land \cd}   \big) dr
    \+ \int_\oga^\oT    g_i(r,\oX_{r \land \cd} ) dr \Big)    \Big]
     \=  E_\oP \Big[   \b1_{\ocA_*}  \int_t^\oT   \n g_i \big( r,  \oX_{r \land \cd}    \big)   dr \Big]
 \eeas
 and thus $ E_{\oP_\e} \big[    \int_t^\oT   g_i \big(r,\oX_{r \land \cd}  \big) dr \big]
 \ls   E_\oP  \big[    \int_t^\oT   g_i \big(r,\oX_{r \land \cd}  \big) dr \big] \ls y_i $.
 Analogously, we have  $ E_{\oP_\e} \big[    \int_t^\oT  h_i \big(r,\oX_{r \land \cd}  \big) dr \big]
 \=  E_\oP  \big[    \int_t^\oT  h_i \big(r,\oX_{r \land \cd}  \big) dr \big] \= z_i $.  Hence,   $\oP_\e $ belongs to $  \ocP_{t,\bx}(y,z)$.

 \no {\bf II.f)} Let $\breve{V}$ be the    function    defined in \eqref{020622_21}
 and set $D^V_\infty \df \big\{\oo \ins \oO \n : \breve{V}   (\oo)   \= \infty  \big\}   \= \big\{\oo \ins \ocA_* \cap   \ocN^c_* \n :    \oV \big( \breve{\Psi} (\oo) \big) \= \infty  \big\}
  \ins  \si \Big(  \cF^{\oW^t}_\oga \cp   \sN_\oP \big(\ocF^t_{\n \infty}  \big) \Big)  $.
 By Theorem \ref{thm_V=oV},  $D^V_\infty  $ is also equal to $   \big\{\oo \ins \ocA_* \cap   \ocN^c_* \n :    \oV \big( \ddot{\Psi} (\oo) \big) \= \infty  \big\}$.
  As $E_\oP \big[ \int_t^\infty   \n f^-   ( r,  \oX_{r \land \cd}   )    dr \big] \< \infty$,
 there exists a $ \ocN_{\n f} \ins \sN_\oP \big(\cF^\oX_\infty\big) $ such that
 $ \int_t^\infty   \n f^-   \big( r,  \oX_{r \land \cd} (\oo)   \big)    dr \< \infty$
 for any $\oo \ins \ocN^c_{\n f}$.

  Let $\e \ins (0,1)$.
    \if{0}

    since \eqref{081620_19} implies that
\beas 
    E_{\oQ^\oo_\e} \Big[ \int_{\oga(\oo)}^\oT   \n f  \big( r,  \oX_{r \land \cd}   \big)   dr   \+    \b1_{\{\oT < \infty\}} \pi  \big( \oT,\oX_{\oT \land \cd} \big) \Big] \gs
\left\{
\ba{ll}
1/\e , & \fa \oo \ins   D^V_\infty ;\\
\oV \big( \ddot{\Psi} (\oo) \big)    \- \e , & \fa \oo \ins   \big(D^V_\infty\big)^c,
 \ea
 \right.
 \eeas

   \fi
  For any  $\oo \ins \ocA_* \Cp \ocN^c_* \Cp \ocN^c_{\n f} $,
   an analogy to \eqref{091020_15}, \eqref{081620_19} and Theorem \ref{thm_V=oV}  imply    that
   \beas
    E_{\oQ^\oo_\e} \big[ \, \oR(t) \big]
   & \tn \=  & \tn    \int_t^{\oga(\oo)}   \n f  \big(  r,  \oX_{r \land \cd} (\oo)  \big)   dr \+ E_{\ol{\bQ}_\e  ( \ddot{\Psi} (\oo) )} \big[   \, \oR(\oga(\oo)) \big] \\
    & \tn  \gs  & \tn   \int_t^{\oga(\oo)}   \n f  \big(  r,  \oX_{r \land \cd} (\oo)  \big)   dr \+  \b1_{ \{\oo \in  (D^V_\infty )^c  \}} \Big( \oV \big( \oga(\oo)   ,  \oX_{\oga \land \cd}(\oo)  ,   \big(\oY_{\n \oP}  ( \oga )\big) (\oo),    \big( \oZ_\oP  ( \oga ) \big) (\oo)   \big)   \- \e \Big)
   \+   \frac{1}{\e} \b1_{ \{\oo \in D^V_\infty  \}}    .
 \eeas
  Since $\oP_\e \ins \ocP_{t,\bx}(y,z)$  and since $\b1_{\ocA_*}  \= \b1_{\{ \wh{\tau}(\oW)   \ge   \oga \}} \=  \b1_{\{\oT    \ge   \oga \}}   $, $\oP-$a.s.  by \eqref{020322_11},
   \bea
  && \hspace{-1cm}\oV (t,\bx,y,z)    \gs    E_{\oP_\e} \n  \big[ \, \oR(t) \big]
  \gs    E_\oP   \bigg[ \b1_{\ocA^c_*} \oR(t) \+ \b1_{\ocA_*} \Big(  \int_t^{\oga }   \n f  \big(  r,  \oX_{r \land \cd}   \big)   dr \+
   \b1_{  (D^V_\infty )^c  } \big[ \, \oV \big( \oga   ,  \oX_{\oga \land \cd}  ,   \oY_{\n \oP}  ( \oga ) ,     \oZ_\oP  ( \oga )    \big)   \- \e \big]  \n \+   \frac{1}{\e} \b1_{  D^V_\infty  }  \Big)     \bigg] \nonumber \\
   & & \hspace{-0.5cm} \=    E_\oP   \bigg[ \b1_{\{\oT   <  \oga \}} \oR(t) \+ \b1_{\{\oT    \ge   \oga \}} \Big( \int_t^{\oga }   \n f  \big(  r,  \oX_{r \land \cd}   \big)   dr \+
   \b1_{  (D^V_\infty )^c  } \big[ \, \oV \big( \oga   ,  \oX_{\oga \land \cd}  ,   \oY_{\n \oP}  ( \oga ) ,     \oZ_\oP  ( \oga )    \big)   \- \e \big]    \+   \frac{1}{\e} \b1_{  D^V_\infty  }  \Big)     \bigg] \nonumber \\
    && \hspace{-0.5cm}  \gs    E_\oP   \bigg[ \b1_{\{\oT   <  \oga \}} \oR(t) \+ \b1_{\{\oT    \ge   \oga \}} \Big( \int_t^{\oga }   \n f  \big(  r,  \oX_{r \land \cd}   \big)   dr \+
   \b1_{  (D^V_\infty )^c  }   \oV \big( \oga   ,  \oX_{\oga \land \cd}  ,   \oY_{\n \oP}  ( \oga ) ,     \oZ_\oP  ( \oga )    \big)       \Big)     \bigg] \- \e \+ \frac{1}{\e} \oP \big(   \{\oT    \gs   \oga \} \Cp D^V_\infty \big)
   .  \qq \qq
   \label{081720_17}
 \eea

 To verify \eqref{081720_15}, we set $\ocI^t_\oP \df \b1_{\{\oT < \oga  \}} \oR(t) \+ \b1_{\{\oT \ge \oga  \}} \Big( \n \int_t^\oga  \n    f(r,\oX_{r \land \cd}) dr  \+ \oV \big( \oga  ,\oX_{ \oga   \land \cd} ,    \oY_{\n \oP}   (\oga  ) , \oZ_\oP   (\oga  )  \big) \Big) $.

\bul  If $ \oP \big( \{\oT    \gs   \oga \} \cap D^V_\infty \big) \= 0  $, then
   $ \oV (t,\bx,y,z)   \gs  E_\oP \Big[ \b1_{\{\oT   <  \oga \}} \oR(t) \+ \b1_{\{\oT    \ge   \oga \}} \Big( \int_t^{\oga }   \n f  \big(  r,  \oX_{r \land \cd}   \big)   dr \+
       \oV \big( \oga   ,  \oX_{\oga \land \cd}  ,   \oY_{\n \oP}  ( \oga ) ,     \oZ_\oP  ( \oga )    \big) \Big) \Big] \- \e    $
 holds for any $\e \ins (0,1)$.    Letting $\e \nto 0$ gives   \eqref{081720_15}.

 \bul   If $ \oP \big( \{\oT    \gs   \oga \} \Cp D^V_\infty \big) \> 0  $  and
 $  E_\oP \big[ \big(\ocI^t_\oP\big)^- \big] \= \infty$, then $ E_\oP \big[  \ocI^t_\oP  \big] \= -\infty \ls  \oV (t,\bx,y,z) $,
 so \eqref{081720_15} holds automatically.

 \bul If $ \oP \big( \{\oT    \gs   \oga \} \cap D^V_\infty \big) \> 0  $  and $ E_\oP \big[ \big(\ocI^t_\oP\big)^- \big] \< \infty $,
 since Remark \ref{rem_ocP2} (1) shows that  $E_\oP   \Big[ \-\b1_{\{\oT   <  \oga \}} \oR(t) \- \b1_{\{\oT    \ge   \oga \}} \Big(  \n \int_t^{\oga }   \n f  \big(  r,  \oX_{r \land \cd}   \big)   dr \\ \+
   \b1_{  (D^V_\infty )^c  }   \oV \big( \oga   ,  \oX_{\oga \land \cd}  ,   \oY_{\n \oP}  ( \oga ) ,     \oZ_\oP  ( \oga )      \big) \Big)        \Big]
      \= E_\oP   \Big[ \- \b1_{  (D^V_\infty )^c  }  \ocI^t_\oP
     \- \b1_{  D^V_\infty  }    \int_t^{\oga \land \oT}   \n f  \big(  r,  \oX_{r \land \cd}   \big)   dr
     \- \b1_{\{\oT   <  \oga \} \cap  D^V_\infty   }   \pi   \big(\oT, \oX_{\oT \land \cd}\big)   \Big]
     \ls E_\oP   \big[   \big( \ocI^t_\oP \big)^-
     \+    \int_t^\infty   \n f^-  \big(  r,  \oX_{r \land \cd}   \big)   dr \big] \- c_\pi \< \infty  $,
    we can deduce from   \eqref{081720_17}  that
 \beas
  \oV (t,\bx,y,z)
  \gs - E_\oP   \Big[   \big( \ocI^t_\oP \big)^-
     \+    \int_t^\infty   \n f^-  \big(  r,  \oX_{r \land \cd}   \big)   dr \Big] \+ c_\pi \- \e \+ \frac{1}{\e} \oP \big(   \{\oT    \gs   \oga \} \Cp D^V_\infty \big) , \q  \fa \e \ins (0,1)  .
  \eeas
  Sending $\e \nto 0$   yields $ \oV (t,\bx,y,z) \= \infty $, so \eqref{081720_15} still holds. This completes the proof of Theorem \ref{thm_DPP1}. \qed

\appendix
\renewcommand{\thesection}{A}
\refstepcounter{section}
\makeatletter
\renewcommand{\theequation}{\thesection.\@arabic\c@equation}
\makeatother

\section{Appendix}

\begin{lemm} \label{lem_122921_11}
Let  $t_0 \ins [0,\infty)$. For $i \= 1,2$,  let $(\O_i, \cF_i, P_i)$ be   a    probability space
and  let $B^i\= \{B^i_s\}_{s \in [0,\infty)}$ be an  $\hR^d-$valued continuous process on $\O$ with $B^i_0 \= \bz$ such that
 $\fB^i_s \df B^i_s \- B^i_{t_0}$, $  s \ins [t_0,\infty)$ is a  Brownian motion on $(\O_i, \cF_i, P_i)$.
 Let $\Phi \n: \O_1 \mto \O_2$ be a mapping such that  $ \fB^2_s (\Phi(\o))\= \fB^1_s(\o)  $  for any $(s,\o) \ins [t_0,\infty) \ti \O_1 $, then
 \(i\) $\Phi^{-1}\big(\cF^{\fB^2}_s\big) \= \cF^{\fB^1}_s  $, $\fa s \ins [t_0,\infty]$;
 \(ii\) $\Phi^{-1}\big(\sN_{P_2}(\cF^{\fB^2}_\infty)\big) \sb \sN_{P_1}(\cF^{\fB^1}_\infty) $;
 \(iii\) $\Phi^{-1}\big(\cF^{\fB^2,P_2}_s\big) \sb \cF^{\fB^1,P_1}_s  $, $\fa s \ins [t_0,\infty]$
 and \(iv\) $P_1 \nci \Phi^{-1} (A) \= P_2(A)$ for any $ A \ins \cF^{\fB^2,P_2}_\infty$.

\end{lemm}

  \no {\bf Proof:}   Let $  s \ins [t_0,\infty]$.  Since
  $ (\fB^1_r)^{-1} ( \cE )  \= \{\o \ins \O_1 \n : \fB^1_r(\o) \ins \cE \}  \= \big\{\o \ins \O_1 \n : \fB^2_r \big(\Phi(\o)\big) \ins \cE \big\} \= \Phi^{-1} \big( (\fB^2_r)^{-1} (\cE) \big) $ for any $ r \ins [t_0,s] \Cp \hR  $ and $ \cE \ins \sB(\hR^d) $,
  we see that   all generating sets of $\cF^{\fB^1}_s$ are included in  the  sigma-field $ \Phi^{-1}(\cF^{\fB^2}_s) \= \big\{\Phi^{-1}(A) \n : A \ins \cF^{\fB^2}_s \big\} $ and that all generating sets of $\cF^{\fB^2}_s $ are contained in the sigma-field $\{A \sb \O_2 \n : \Phi^{-1} (A) \ins \cF^{\fB^1}_s \}$. It follows that  $\Phi^{-1} \big(\cF^{\fB^2}_s\big) \= \cF^{\fB^1}_s$.

  Since  the $P_1-$distribution of $\fB^1$ is equal to the $P_2-$distribution of $\fB^2$, 
  it holds for any   $\{(s_i,\cE_i)\}^k_{i=1} \sb [t_0,\infty)\ti\sB(\hR^d)$ that
 $ P_2 \Big( \underset{i=1}{\overset{k}{\cap}}  (\fB^2_{s_i})^{-1} ( \cE_i ) \Big)
 \= P_1 \Big( \underset{i=1}{\overset{k}{\cap}} (\fB^1_{s_i})^{-1} ( \cE_i ) \Big)
  \= P_1 \Big( \underset{i=1}{\overset{k}{\cap}} \Phi^{-1} \big( (\fB^2_{s_i})^{-1} (\cE_i)\big) \Big)
  \= P_1 \Big( \Phi^{-1} \Big( \underset{i=1}{\overset{k}{\cap}}   (\fB^2_{s_i})^{-1} (\cE_i)\Big) \Big)  $.
  So the Lambda-system  $\L \df \big\{ A \ins \cF^{\fB^2}_\infty  \n : P_2(A) \= P_1  ( \Phi^{-1} (A)) \big\}$
  contains the Pi-system  $\Big\{ \underset{i=1}{\overset{k}{\cap}}  (\fB^2_{s_i})^{-1} ( \cE_i ) \n : \{(s_i,\cE_i)\}^k_{i=1} \sb [t_0,\infty)\ti\sB(\hR^d) \Big\}$, which generates $\cF^{\fB^2}_\infty$.
   Dynkin's Pi-Lambda Theorem   yields  that $\cF^{\fB^2}_\infty \= \L$, i.e.,
   \bea \label{122921_11}
    P_2(A) \= P_1  ( \Phi^{-1} (A)) , \q \fa  A \ins \cF^{\fB^2}_\infty  .
    \eea
  Given $ \cN \ins \sN_{P_2}(\cF^{\fB^2}_\infty)$, there exists $A \ins \cF^{\fB^2}_\infty$ such that $\cN  \sb A  $ and  $ P_2(A) \= 0$.
  Since $ \Phi^{-1} (A) \ins 
  \cF^{\fB^1}_\infty $ and since   $ P_1 \big( \Phi^{-1} (A) \big)   \= P_2 (A) \= 0 $,
  we obtain that  $ \Phi^{-1} (\cN) \ins \sN_{P_1} (\cF^{\fB^1}_\infty) $.

 Let  $s \ins [t_0,\infty]$. We have known  from the above that
the sigma-field $  \big\{A \sb \O_2 \n :  \Phi^{-1} (A) \ins \cF^{\fB^1,P}_s\big\}$
  contains both $\cF^{\fB^2}_s$ and $\sN_{P_2}(\cF^{\fB^2}_\infty)$.
  It thus includes $\cF^{\fB^2,P_2}_s$, i.e.,  $\Phi^{-1} \big(\cF^{\fB^2,P_2}_s\big) \sb  \cF^{\fB^1,P_1}_s$.

   Now, let   $ A \ins  \cF^{\fB^2,P_2}_\infty  $. According to Problem 2.7.3 of \cite{Kara_Shr_BMSC},
there exists $A' \ins \cF^{\fB^2}_\infty $ such that $ \cN_2  \df  A \D A' \ins  \sN_{P_2}(\cF^{\fB^2 }_\infty)  $.
So $ \cN_1  \df \Phi^{-1}(A) \D \Phi^{-1}(A') \= \Phi^{-1}(A \D A') \ins \sN_{P_1}(\cF^{\fB^1 }_\infty) $.
 Then we can deduce from \eqref{122921_11} that
   $ P_1 \big(\Phi^{-1}(A)\big) \= P_1 \big( \Phi^{-1}(A') \D \cN_1  \big)
   \= P_1 \big(\Phi^{-1}(A')\big) \= P_2(A') \= P_2(A \D \cN_2 ) \= P_2(A)$. \qed

 \begin{lemm} \label{lem_M31_01}
 Let $(\O, \cF, P)$ be   a    probability space and let $\,t \ins [0,\infty)$.
 Let $B\= \{B_s\}_{s \in [0,\infty)}$ be an $\hR^d-$valued continuous process on $\O$ with $B_0 \= \bz$ such that
 $B^t_s \df B_s \- B_t$, $  s \ins [t,\infty)$ is a  Brownian motion on $(\O, \cF, P)$.

\no \(1\) For any $[t,\infty]-$valued $ \bF^{W^t,P_0}-$stopping time $\wh{\tau}$ on $\O_0$,
$\wh{\tau}(B)$ is an    $  \bF^{B^t,P}-$stopping time on $\O$.

\no \(2\)  Let $\Phi \n: \O  \mto \O_0$ be a mapping such that  $ W^t_s (\Phi(\o))\= B^t_s (\o)  $  for any $(s,\o) \ins [t,\infty) \ti \O  $.
 For any   $[t,\infty]-$valued $  \bF^{B^t,P}-$stopping time $\tau$ on $\O$,  there exists a $[t,\infty]-$valued $ \bF^{W^t,P_0}-$stopping time $\wh{\tau} $ on $\O_0$  such that $\tau \= \wh{\tau}  (\Phi )$, $P-$a.s.

 \end{lemm}

\no {\bf Proof: 1)} Since $ W^t_s (B(\o)) \= B_s(\o) \- B_t(\o) \=  B^t_s(\o)$  for any $(s,\o) \ins [t,\infty) \ti \O$,
 applying Lemma \ref{lem_122921_11} with $t_0 \= t$,  $(\O_1, \cF_1, P_1,B^1)   \= \big(\O, \cF, P , B \big) $, $(\O_2, \cF_2, P_2,B^2) \= \big(\O_0, \sB(\O_0) , P_0 , W\big) $ and $\Phi \= B$ shows that  $B^{-1} \big(\cF^{W^t,P_0}_s\big) \sb \cF^{B^t,P}_s$  for any $s \ins [t,\infty]$.

   Let $\wh{\tau}$ be a $[t,\infty]-$valued $ \bF^{W^t,P_0}-$stopping time on $\O_0$. For any $s \ins [t,\infty)$,
 as $  \{\wh{\tau} \ls s\} \ins \cF^{W^t,P_0}_s$,  one can deduce that     $\{ \wh{\tau}(B) \ls s \} \=  B^{-1} \big(\{\wh{\tau} \ls s\}\big) \ins \cF^{B^t,P}_s $.
  Hence,  $\wh{\tau}(B) $ is a $[t,\infty]-$valued $\bF^{B^t,P}-$stopping time on $\O$.

  \no {\bf 2)}  Let $\Phi \n: \O  \mto \O_0$ be a mapping such that  $ W^t_s (\Phi(\o))\= B^t_s (\o)  $  for any $(s,\o) \ins [t,\infty) \ti \O  $.    Lemma \ref{lem_122921_11} also implies that
        $\Phi^{-1} \big(\cF^{W^t}_s\big) \= \cF^{B^t}_s$, $ \fa s \ins [t,\infty]$.

    Let $ \tau $ be a $[t,\infty]-$valued $  \bF^{B^t,P}-$stopping time on $\O$.
 We   fix $n \ins \hN$ and set $s^n_0 \df t$.
  Given $ i \ins \hN$, we set   $s^n_i \df  t + i 2^{-n}$ and $ \ddot{A}^n_i \df \big\{ s^n_{i-1} \ls \tau \< s^n_i  \big\} \ins \cF^{B^t,P}_{s^n_i} $.
 By Problem 2.7.3 of \cite{Kara_Shr_BMSC}, there is $ A^n_i \ins  \cF^{B^t}_{s^n_i} $ such that  $ \cN^n_i \df \ddot{A}^n_i \D A^n_i \ins \sN_P \big(\cF^{B^t}_\infty\big) $.
 As  $ \cF^{B^t}_{s^n_i} \= \Phi^{-1} \big(\cF^{W^t}_{s^n_i}\big)   $,
 we can find $\cA^n_i \ins \cF^{W^t}_{s^n_i}$ such that $  A^n_i  \=  \Phi^{-1}(\cA^n_i)$.
 Set $\wcA^n_i \df \cA^n_i \Big\backslash \Big( \underset{j < i}{\cup} \cA^n_j \Big) \ins \cF^{W^t}_{s^n_i} $\,.

 Define a $(t,\infty]-$valued $\bF^{W^t}-$stopping time
 $ \wh{\tau}^n \df \sum_{i \in  \hN} s^n_i \b1_{\wcA^n_i} \+ \infty \b1_{\big(\ccup{i \in \hN}{}  \wcA^n_i \big)^c  } $
 and set $\O_n \df \Big(\ccup{i \in \hN}{} \cN^n_i\Big)^c$.
 For any $ i \ins \hN $, since $ \O_n \Cp A^n_i \= \O_n \Cp \ddot{A}^n_i $ and $ \O_n \Cp (A^n_i)^c \= \O_n \Cp \big(\ddot{A}^n_i\big)^c  $, one has
  $ \O_n \Cp  \Phi^{-1} \big(\wcA^n_i\big)   \= \O_n \Cp \Big(\Phi^{-1}(\cA^n_i) \Big\backslash \Big( \underset{j < i}{\cup} \Phi^{-1}(\cA^n_j) \Big) \Big)
 \= \O_n \Cp \Big(A^n_i  \Big\backslash \Big( \underset{j < i}{\cup}  A^n_j  \Big) \Big)
 \= \O_n \Cp \Big(\ddot{A}^n_i  \Big\backslash \Big( \underset{j < i}{\cup}  \ddot{A}^n_j  \Big) \Big)
 \=  \O_n \Cp  \ddot{A}^n_i   $,
 and thus $ \O_n \Cp \Big(\ccup{i \in \hN}{}  \Phi^{-1} \big(\wcA^n_i\big) \Big)^c
 \= \O_n \Cp \Big(\ccup{i \in \hN}{}   \ddot{A}^n_i    \Big)^c
 \=  \O_n \Cp \{\tau \= \infty\} $.
 It follows that
 \bea \label{122820_14}
 \b1_{\O_n} \, \wh{\tau}^n (\Phi) \= \b1_{\O_n}  \Big( \sum_{i \in  \hN} s^n_i \b1_{\Phi^{-1}(\wcA^n_i)} \+ \infty \b1_{\big(\ccup{i \in \hN}{}  \Phi^{-1}(\wcA^n_i) \big)^c  } \Big)
 \= \b1_{\O_n}  \Big( \sum_{i \in  \hN} s^n_i \b1_{ \ddot{A}^n_i } \+ \infty \b1_{\{\tau = \infty\} } \Big)  .
 \eea

 Since   $\bF^{W^t,P_0}$ is a right-continuous complete filtration,
 Lemma 1.2.11 of \cite{Kara_Shr_BMSC} renders  that
 $\wh{\tau}  \df \linf{n \to \infty} \wh{\tau}^n$  is a $[t,\infty]-$valued  $\bF^{W^t,P_0}-$stopping time on $\O_0$.
 Given $\o \ins \ccap{n \in \hN}{} \O_n$, we can deduce from \eqref{122820_14} that
 \beas
  \hspace{1cm} \tau(\o) \= \lmtd{n \to \infty} \Big( \sum_{i \in  \hN} s^n_i \b1_{ \{ s^n_{i-1} \le \tau(\o)  < s^n_i \} } \+ \infty \b1_{\{\tau(\o) = \infty\} } \Big) \= \lmtd{n \to \infty} \wh{\tau}^n ( \Phi (\o))
 \= \linf{n \to \infty} \wh{\tau}^n ( \Phi (\o)) \= \wh{\tau} \big( \Phi (\o)\big) . \hspace{1.5cm} \hb{\qed}
 \eeas

   \if{0}

 \begin{lemm} \label{lem_M31_01_b}
 Let $(\O, \cF, P)$ be   a    probability space and let $\,t \ins [0,\infty)$.
 Let $B\= \{B_s\}_{s \in [0,\infty)}$ be an $\hR^d-$valued continuous process on $\O$ 
 such that $\beta^t_\fs \df   B_{t+\fs} \- B_t $, $  \fs \ins [0,\infty)$ is a  Brownian motion on $(\O, \cF, P)$.
 Set also $B^t_s \df B_s \- B_t $, $\fa s \ins [t,\infty)$.

\no \(1\) For any $[0,\infty]-$valued $ \bF^{W,P_0}-$stopping time $\wh{\tau}$ on $\O_0$,
$ t\+\wh{\tau}(\beta^t)$ is a   $[t,\infty]-$valued  $  \bF^{B^t,P}-$stopping time on $\O$.

\no \(2\)
 For any   $[t,\infty]-$valued $  \bF^{B^t,P}-$stopping time $\tau$ on $\O$,  there exists a $[0,\infty]-$valued $ \bF^{W,P_0}-$stopping time $\wh{\tau} $ on $\O_0$  such that $\tau \=  t \+ \wh{\tau}  (\beta^t )$, $P-$a.s.

 \end{lemm}

\no {\bf Proof: 1)} Since $ W_\fs (\beta^t(\o)) \= \beta^t_\fs(\o)  $  for any $(\fs,\o) \ins [0,\infty) \ti \O$,
 applying Lemma \ref{lem_122921_11} with $t_0 \= 0$,  $(\O_1, \cF_1, P_1,B^1)   \= \big(\O, \cF, P , \beta^t \big) $, $(\O_2, \cF_2, P_2,B^2) \= \big(\O_0, \sB(\O_0) , P_0 , W\big) $ and $\Phi \= \beta^t $ shows that  $(\beta^t)^{-1} \big(\cF^{W,P_0}_\fs\big) \sb \cF^{\beta^t,P}_\fs$  for any $\fs \ins [0,\infty]$.

 Since $\beta^t_\fs \= B^t_{t+\fs}$, $\fa \fs \ins [0,\infty)$, one has
 $\cF^{\beta^t}_\fs \= \si\big(\beta^t_r; r \ins [0,\fs]\big) \= \si\big(B^t_{t+r}; r \ins [0,\fs]\big)
 \= \si\big(B^t_{r'}; r' \ins [t,t\+\fs]\big) \= \cF^{B^t}_{t+\fs} $ for any $\fs \ins [0,\infty)$.
 It follows that $\cF^{\beta^t}_\infty \= \cF^{B^t}_\infty$ and thus
 $ \cF^{\beta^t,P}_\fs \= \cF^{B^t,P}_{t+\fs} $ for any $\fs \ins [0,\infty)$.

   Let $\wh{\tau}$ be a $[0,\infty]-$valued $ \bF^{W,P_0}-$stopping time on $\O_0$. For any $s \ins [t,\infty)$,
 as $  \{\wh{\tau} \ls s\-t\} \ins \cF^{W,P_0}_{s-t}$,  one can deduce that     $\{ t\+\wh{\tau}(\beta^t) \ls s \} \=  (\beta^t)^{-1} \big(\{\wh{\tau} \ls s \- t\}\big) \ins \cF^{\beta^t,P}_{s-t} \= \cF^{B^t,P}_s $.
  So  $t\+\wh{\tau}(\beta^t) $ is a $[t,\infty]-$valued $\bF^{B^t,P}-$stopping time on $\O$.

  \no {\bf 2)} For $\fs \ins [0,\infty)$, one has
  $ (\beta^t)^{-1} \big(\cF^W_\fs\big) \= (\beta^t)^{-1} \big(\si \big\} W^{-1}_r (A) \n: r \ins [0,\fs],A \ins \sB(\hR^d) \big\} \big)
 \=   \si \big\{ (\beta^t)^{-1} \big( W^{-1}_r (A)\big) \n: r \ins [0,\fs],A \ins \sB(\hR^d) \big\}
 \=   \si \big( (B^t_{t+r})^{-1}   (A)  \n: r \ins [0,\fs],A \ins \sB(\hR^d) \big)
 \= \cF^{B^t}_{t+\fs}$.

    Let $ \tau $ be a $[t,\infty]-$valued $  \bF^{B^t,P}-$stopping time on $\O$.
 We   fix $n \ins \hN$ and set $s^n_0 \df 0$.
  Given $ i \ins \hN$, we set   $s^n_i \df    i 2^{-n}$ and $ \ddot{A}^n_i \df \big\{t+ s^n_{i-1} \ls \tau \< t+s^n_i  \big\} \ins \cF^{B^t,P}_{t+s^n_i} $.
 By Problem 2.7.3 of \cite{Kara_Shr_BMSC}, there is $ A^n_i \ins  \cF^{B^t}_{t+s^n_i} $ such that  $ \cN^n_i \df \ddot{A}^n_i \D A^n_i \ins \sN_P \big(\cF^{B^t}_\infty\big) $.
 As  $ \cF^{B^t}_{t+s^n_i} \= (\beta^t)^{-1} \big(\cF^W_{s^n_i}\big)   $,
 we can find $\cA^n_i \ins \cF^W_{ s^n_i}$ such that $  A^n_i  \=  (\beta^t)^{-1}(\cA^n_i)$.
 Set $\wcA^n_i \df \cA^n_i \Big\backslash \Big( \underset{j < i}{\cup} \cA^n_j \Big) \ins \cF^W_{s^n_i} $\,.

 Define a $(0,\infty]-$valued $\bF^W-$stopping time
 $ \wh{\tau}^n \df \sum_{i \in  \hN} s^n_i \b1_{\wcA^n_i} \+ \infty \b1_{\big(\ccup{i \in \hN}{}  \wcA^n_i \big)^c  } $
 and set $\O_n \df \Big(\ccup{i \in \hN}{} \cN^n_i\Big)^c$.
 For any $ i \ins \hN $, since $ \O_n \Cp A^n_i \= \O_n \Cp \ddot{A}^n_i $ and $ \O_n \Cp (A^n_i)^c \= \O_n \Cp \big(\ddot{A}^n_i\big)^c  $, one has
  $ \O_n \Cp  (\beta^t)^{-1} \big(\wcA^n_i\big)   \= \O_n \Cp \Big((\beta^t)^{-1}(\cA^n_i) \Big\backslash \Big( \underset{j < i}{\cup} (\beta^t)^{-1}(\cA^n_j) \Big) \Big)
 \= \O_n \Cp \Big(A^n_i  \Big\backslash \Big( \underset{j < i}{\cup}  A^n_j  \Big) \Big)
 \= \O_n \Cp \Big(\ddot{A}^n_i  \Big\backslash \Big( \underset{j < i}{\cup}  \ddot{A}^n_j  \Big) \Big)
 \=  \O_n \Cp  \ddot{A}^n_i   $,
 and thus $ \O_n \Cp \Big(\ccup{i \in \hN}{}  (\beta^t)^{-1} \big(\wcA^n_i\big) \Big)^c
 \= \O_n \Cp \Big(\ccup{i \in \hN}{}   \ddot{A}^n_i    \Big)^c
 \=  \O_n \Cp \{\tau \= \infty\} $.
 It follows that
 \bea \label{122820_14}
 \b1_{\O_n} \, \wh{\tau}^n (\beta^t) \= \b1_{\O_n}  \Big( \sum_{i \in  \hN} s^n_i \b1_{(\beta^t)^{-1}(\wcA^n_i)} \+ \infty \b1_{\big(\ccup{i \in \hN}{}  (\beta^t)^{-1}(\wcA^n_i) \big)^c  } \Big)
 \= \b1_{\O_n}  \Big( \sum_{i \in  \hN} s^n_i \b1_{ \ddot{A}^n_i } \+ \infty \b1_{\{\tau = \infty\} } \Big)  .
 \eea

 Since   $\bF^{W,P_0}$ is a right-continuous complete filtration,
 Lemma 1.2.11 of \cite{Kara_Shr_BMSC} renders  that
 $\wh{\tau}  \df \linf{n \to \infty} \wh{\tau}^n$  is a $[0,\infty]-$valued  $\bF^{W,P_0}-$stopping time on $\O_0$.
 Given $\o \ins \ccap{n \in \hN}{} \O_n$, we can deduce from \eqref{122820_14} that
 \beas
     \tau(\o) \= \lmtd{n \to \infty} \Big( \sum_{i \in  \hN} (t\+s^n_i) \b1_{ \{ t+s^n_{i-1} \le \tau(\o)  < t+s^n_i \} } \+ \infty \b1_{\{\tau(\o) = \infty\} } \Big) \= t\+ \lmtd{n \to \infty} \wh{\tau}^n ( \beta^t (\o))
 \= t\+  \linf{n \to \infty} \wh{\tau}^n ( \beta^t (\o)) \= t \+ \wh{\tau} \big( \beta^t (\o)\big) . \hspace{0.5cm} \hb{\qed}
 \eeas

\begin{lemm} \label{lem_010922}
The mapping $   (t,\o_0) \mto   \sW^t (\o_0)   $
   is continuous from $[0,\infty) \ti \O_0$ to $ \O_0    $.

\end{lemm}

 \no {\bf Proof:} The metric  of locally uniform convergence on $\O_0$ is
 \bea
  \Rho{\O_0}  (\o_0,\o'_0)  & \tn \df & \tn  \sum_{n \in \hN} \Big( 2^{-n} \ld \Sup{r \in [0,n]}  |\o_0(r)\- \o'_0(r)| \Big) , \q \fa \o_0 , \o'_0 \ins \O_0   \label{062720_23}
  \eea

 Fix $(t,\o_0) \ins [0,\infty) \ti \O_0   $ and $\e \ins (0,1)$. We set $\fn_1 \df \lceil  3 \- \log_2 \e \rceil $
 and $\fn_2 \df \fn_1 \+ \lceil t \rceil \+ 1$.
 As  $ \o_0 (r) $ is uniformly continuous over $r \ins [0,\fn_2]$,
  there exists $  \d   \ins (0,2^{-\fn_2}/\fn_1)$ such that
 \bea \label{010922_23}
  |\o_0(r) \- \o_0(r')|   \< \frac{\e}{8\fn_1}  ,  \q \fa r,r' \ins [0, \fn_2] \hb{ with } |r \- r'| \< \d     .
  \eea
  For any $(t',\o'_0) \ins [0,\infty) \ti \O_0  $ with $ |t'\-t| \ve \Rho{\O_0}  (\o_0,\o'_0) \< \d$, since
 $ 2^{-\fn_2} \ld \Sup{r \in [0,\fn_2]}  |\o_0(r)\- \o'_0(r)| \ls \Rho{\O_0}  (\o_0,\o'_0)
  \< \d $, we have
   \bea \label{010922_21}
   \Sup{r \in [0,\fn_2]}  |\o_0(r)\- \o'_0(r)| \< \d \< 2^{-\fn_2}/\fn_1 \< 2^{-\fn_1}/\fn_1 \ls \frac{\e}{8\fn_1}  .
   \eea

 Let $ \fr \ins [0,\fn_1] $. since $t' \+ \fr \< t \+ \d \+ \fn_1 \< \fn_1 \+ \lceil t \rceil \+ 1 \= \fn_2 $,
  we can deduce from \eqref{010922_21} and \eqref{010922_23} that
 \beas
 \big| \sW^{t'}_\fr (\o'_0) \- \sW^t_\fr (\o_0) \big|
& \tn \= & \tn  \big|\o'_0(t'\+\fr) \- \o'_0(t') \- \o_0(t\+\fr) \+ \o_0(t) \big| \\
& \tn \ls & \tn  \big|\o'_0(t'\+\fr) \- \o_0(t'\+\fr)| \+ \big|\o'_0(t') \- \o_0(t')\big|
\+ \big|\o_0(t'\+\fr) \- \o_0(t\+\fr)|\+\big|  \o_0(t')  \- \o_0(t) \big| \< \frac{\e}{2\fn_1} .
\eeas
  It follows that
 \beas
   \hspace{1cm} \Rho{\O_0} \big(\sW^{t'} (\o'_0), \sW^t (\o_0)\big)
   & \tn \=  & \tn  \sum_{n \in \hN} \Big( 2^{-n} \ld \Sup{\fr \in [0,n]}  | \sW^{t'}_\fr (\o'_0) \- \sW^t_\fr (\o_0) | \Big)
   \ls   \sum^{\fn_1}_{n = 1}    \Sup{\fr \in [0,n]}  | \sW^{t'}_\fr (\o'_0) \- \sW^t_\fr (\o_0) |
   \+ \sum^\infty_{n = \fn_1+1}   2^{-n}   \\
    & \tn \ls & \tn  \e/2 \+ 2^{-\fn_1}\< \frac58 \e .  \hspace{10cm} \hb{\qed}
  \eeas

 Similar to the metric $\Rho{+}$ on $\hT \= [0,\infty]$, we define a metric $\ddot{\rho}$ on $ \ddot{\hR} \df [-\infty,\infty]$ by
 $\ddot{\rho}(a_1,a_2) \df \big|\arctan(a_1) \- \arctan(a_2)\big|$, $\fa a_1,a_2 \ins \ddot{\hR}$.

\begin{lemm} \label{lem040221}
The mapping $(t,\oo) \mto \oT(\oo)\-t$ is continuous from $[0,\infty) \ti \oO$ to $ \ddot{\hR} $.
\end{lemm}

\no {\bf Proof:} Let $ t  \ins [0,\infty)$, $\oo  \= (\o_0,\omX,\ft) \ins  \oO$ and   $\e \ins (0,1)$. We discuss by   two cases:

\no (i) When $\ft \< \infty$, we set $\dis \d \= \d(\e,\ft) \df \frac{\e}{2 (1\+(\ft\+1)^2)} \ld \big(\arctan(\ft\+1) \- \arctan(\ft)\big)   $.
Let $ t' \ins (t \- \d,t\+ \d)\Cp [0,\infty)$ and $\oo'  \= (\o'_0,\omX',\ft') \ins  \oO$ with
 $ \rho^2_{\overset{}{\oO}}   (\oo,\oo') \=   \rho^2_{\overset{}{\O_0}} (\o_0,\o'_0)   \+ \rho^2_{\overset{}{\OmX}} (\omX,\omX') \+ \rho^2_{\overset{}{+}}(\ft,\ft')  \< \d^2 $.
 As $ \arctan(\ft\+1) \- \arctan(\ft) \gs \d \> \Rho{+}(\ft,\ft') \= \big|\arctan(\ft) \- \arctan(\ft')\big|
 \gs \arctan(\ft') \- \arctan(\ft) $, one has that $\ft' \< \ft\+1$.
  Since
 \bea \label{040521_11}
 \frac{b \- a}{1+b^2} \ls \arctan(b) \- \arctan(a) \= \int_a^b \frac{1}{1+x^2}dx \ls \frac{b \- a}{1+a^2} \ls b\-a  , \q \hb{ for } 0 \ls a \ls b \< \infty ,
\eea
we can also deduce that
 $\dis  \frac{\e}{2 (1\+(\ft\+1)^2)} \gs \d \>   \big|\arctan(\ft) \- \arctan(\ft')\big|
 \gs \frac{ |\ft'\-\ft|}{1+(\ft'\ve \ft)^2} \> \frac{ |\ft'\-\ft|}{1+(  \ft\+ 1)^2}$ and thus $  |\ft'\-\ft| \< \e/2$.
 As   $\d \< \e/2$,  \eqref{040521_11} also implies that
$\ddot{\rho} \big( \oT(\oo')\-t' , \oT(\oo)\-t\big)
\= \big|   \arctan \big(\ft'\-t'\big) \- \arctan \big(\ft\-t\big) \big|
   \ls    \big| \ft'\-t' \- \ft \+ t \big|  \ls | t'\- t| \+ |\ft'\-\ft|
\< \d \+ \e/2 \< \e $.

\no (ii) When $\ft \= \infty$, we set $\dis \wt{\d} \= \wt{\d}(\e,t ) \df \frac{\e}{2} \ld \Big( \frac{\pi}{2} \- \arctan \big( t\+1\+ \sqrt{  2(t\+1) \e^{-1}\-1 } \big) \Big)$.  Let $ t' \ins \big(t \- \wt{\d},t\+ \wt{\d} \big)\Cp [0,\infty)$ and $\oo'  \= (\o'_0,\omX',\ft') \ins  \oO$ with
 $ \rho^2_{\overset{}{\oO}}   (\oo,\oo') \=   \rho^2_{\overset{}{\O_0}} (\o_0,\o'_0)   \+ \rho^2_{\overset{}{\OmX}} (\omX,\omX') \+ \rho^2_{\overset{}{+}}(\ft,\ft')  \< \wt{\d}^2 $.
 Since
$\frac{\pi}{2} \- \arctan \big( t\+1\+ \sqrt{  2(t\+1) \e^{-1}\-1 } \big) \gs \wt{\d} \> \Rho{+}(\ft,\ft') \= \big|\arctan(\infty) \- \arctan(\ft')\big| \= \frac{\pi}{2} \- \arctan(\ft') $ and since $t' \< t\+ \wt{\d} \ls t \+ \frac{\e}{2} \< t \+ 1$,
we see that
$\ft' \> t\+1\+ \sqrt{  2(t\+1) \e^{-1}\-1 } \> t' \+ \sqrt{  2(t\+1) \e^{-1}\-1 }$.
 By \eqref{040521_11} again,
 $  \arctan  (\ft' ) \-  \arctan \big(\ft'\-t'\big)
\ls \frac{t'}{ 1+ (\ft'\-t')^2 } \< \frac{t\+1}{1+ ( 2(t\+1) \e^{-1}\-1) } \= \frac{\e}{2} $.
 Adding it to the inequality  $\frac{\pi}{2} \- \arctan(\ft') \= \Rho{+}(\ft,\ft') \< \wt{\d} \ls \e/2 $ yields that
  $ \ddot{\rho} \big( \oT(\oo')\-t' , \oT(\oo)\-t\big)
 \= \big|   \arctan \big(\ft'\-t'\big) \- \arctan \big(\ft\-t\big) \big|
 \= \big|   \arctan \big(\ft'\-t'\big) \- \arctan  (\infty ) \big|
 \= \frac{\pi}{2} \- \arctan \big(\ft'\-t'\big) \< \e $.    \qed

   \fi

 \if{0}

  \begin{lemm} \label{lem_integral_measurable}
 Given a measure space $(\fX,\cF_\fX,\fm)$ and a measurable space $(\fY,\cF_\fY)$,  let   $\ff \n : \fX \ti \fY \mto [-\infty,\infty]$
is a $\cF_\fX \oti \cF_\fY -$measurable function. Then the function
$\phi_\ff (y) \df   \int_{x \in \fX} \ff(x,y) \fm(dx) $, $ \fa y \ins \fY $
 is $\cF_\fY / \sB[-\infty,\infty]-$measurable.
 \end{lemm}

\no {\bf Proof:} Let $\cE \ins \cF_\fX$ and $ A \ins \cF_\fY $.
 Since the function
 $\phi_{\cE,A} (y) \df   \int_{x \in \fX} \b1_\cE (x)  \b1_A (y) \fm(dx)
\=  \b1_A (y) \int_{x \in \cE} \fm(dx) $, $ \fa y \ins \fY $
 is $\cF_\fY / \sB[0,\infty)-$measurable, we see that all measurable rectangles of $ \cF_\fX \oti \cF_\fY $
are included in
\beas
\q \L \df \big\{D \ins \cF_\fX \oti \cF_\fY   \big| \; \phi_D (y) \df   \int_{x \in \fX} \b1_D (x,y) \fm(dx), ~ y \ins \fY   \hb{ is $\cF_\fY / \sB[0,\infty)-$measurable} \big\} ,
\eeas
which is  a Lambda-system  by the monotone convergence theorem.
 Then Dynkin's Pi-Lambda Theorem shows that $ \cF_\fX \oti \cF_\fY \= \L $.

 Next, let $\ff \n : \fX \ti \fY \mto [0,\fk]$ be a $\cF_\fX \oti \cF_\fY -$measuable   function for some $\fk \ins (1,\infty)$.    Given $n \ins \hN$, we set $a^n_i \df   i 2^{-n} \fk $ for $i \= 0, 1, \cds, 1\+2^n$  and define
 $ \ff_n (x,y) \df \underset{i=0}{\overset{2^n}{\sum}} a^n_i \b1_{\{(x,y) \in D^n_i\}} $, $
  \fa (x,y) \ins \fX \ti \fY $,
   where $D^n_i \df \big\{\ff \ins [a^n_i,a^n_{i+1}) \big\} \ins  \cF_\fX \oti \cF_\fY $.
   From the equality $ \cF_\fX \oti \cF_\fY \= \L$, we see that
  $ \phi_{\ff_n} (y)   \= \int_{x \in \fX} \ff_n(x,y) \fm(dx)
 \= \underset{i=0}{\overset{2^n}{\sum}} a^n_i \int_{x \in \fX} \b1_{\{(x,y) \in D^n_i\}} \fm(dx)
 \= \underset{i=0}{\overset{2^n}{\sum}} a^n_i \phi_{D^n_i} (y) $, $ \fa y \ins \fY$
 is $\cF_\fY / \sB[0,\infty)-$measurable. As $ \ff(x,y) \= \lmtu{n \to \infty} \ff_n(x,y)$, $\fa (x,y ) \ins \fX \ti \fY $,
the monotone convergence theorem shows that the mapping
$\phi_\ff (y)   \= \int_{x \in \fX} \ff (x,y) \fm(dx) \= \lmtu{n \to \infty} \int_{x \in \fX} \ff_n(x,y) \fm(dx)
\= \lmtu{n \to \infty}  \phi_{\ff_n} (y) $, $ \fa y \ins \fY$
  is also $\cF_\fY / \sB[0,\infty)-$measurable.

 Now, let $\ff \n : \fX \ti \fY \mto [-\infty,\infty]$ be a  general $\cF_\fX \oti \cF_\fY -$measuable   function, since the functions
 $ \phi^\pm_{\ff,n} (y) \df   \int_{x \in \fX} n \ld \ff^\pm (x,y)   \fm(dx) $, $ \fa y \ins \fY$
  is $\cF_\fY / \sB[0,\infty) -$measurable for each $n \ins \hN$, the monotone convergence theorem  
 implies that the mapping
 $ \phi_{\ff^\pm} (y) \df   \int_{x \in \fX} \ff^\pm (x,y)   \fm(dx) \= \lmtu{n \to \infty}  \int_{x \in \fX} n \ld \ff^\pm (x,y) \fm(dx) \= \lmtu{n \to \infty} \phi^\pm_{\ff,n} (y) \ins [0,\infty] $, $ \fa y \ins \fY$
  is $\cF_\fY / \sB[0,\infty]-$measurable. Eventually, the mapping
 \beas
 \phi_\ff (y)  & \tn \= & \tn   \int_{x \in \fX} \ff  (x,y)   \fm(dx) \=  \int_{x \in \fX} \ff^+  (x,y)   \fm(dx) \-  \int_{x \in \fX} \ff^-  (x,y)   \fm(dx) \= \phi_{\ff^+} (y) \- \phi_{\ff^-} (y) \\
 & \tn \=  & \tn   \b1_{\{\phi_{\ff^-} (y)\}} (-\infty) \+
  \b1_{\{\phi_{\ff^-} (y) < \infty \}}   \lmt{n \to \infty} \big( \phi_{\ff^+} (y) \ld n \- \phi_{\ff^-} (y) \ld n \big) \ins [-\infty,\infty]  , ~\; \fa y \ins \fY
 \eeas
 is $\cF_\fY / \sB[-\infty,\infty] -$measurable.  \qed

 \fi

\begin{lemm}   \label{lem_A1}
Let $ \fX $ be a topological space and let $\fY$ be a Borel space. If $\ff \n : \fX \ti \fY \mto (-\infty,\infty]$
is a $\sB(\fX) \oti \sB(\fY) -$measurable function bounded from below, then
$\phi_\ff (x,P) \df   \int_{y \in \fY} \ff(x,y) P(dy) $, $   (x,P) \ins \fX \ti \fP(\fY)$
is $\sB(\fX) \oti \sB\big(\fP(\fY)\big)  -$measurable.

\end{lemm}

\no {\bf Proof:} Let $\cE \ins \sB(\fX)$ and $ A \ins \sB(\fY) $.
Since Proposition 7.25 of  \cite{Bertsekas_Shreve_1978} shows that
the function  $\phi_A (P) \df E_P  [\b1_A ]$, $    P  \ins   \fP(\fY) $ is
 $ \sB\big(\fP(\fY)\big)   -$measurable, the mapping
 $ \phi_{\cE,A} (x,P) \df   \int_{y \in \fY} \b1_\cE (x)  \b1_A (y,z) P(dy) $, $   (x,P) \ins \fX \ti \fP(\fY) $
 is $\sB(\fX) \oti \sB\big(\fP(\fY)\big)  -$measurable. So all measurable rectangles of $ \sB(\fX) \oti \sB(\fY) $
 are included in the Lambda-system
$ \L \df \big\{D \ins \sB(\fX) \oti \sB(\fY)   \big| \; \phi_D (x,P) \df   \int_{y \in \fY} \b1_D (x,y) P(dy), ~ (x,P) \ins \fX \ti \fP(\fY) \hb{ is $\sB(\fX) \oti \sB\big(\fP(\fY)\big) -$measurable} \big\} $.
  An application of Dynkin's Pi-Lambda Theorem yields that $ \sB(\fX) \oti \sB(\fY) \= \L $.

 Next, let $\ff \n : \fX \ti \fY \mto \hR$ be a $\sB(\fX) \oti \sB(\fY)-$measuable   function
 taking values in a finite  subset $\{a_1 \<  \cds \< a_N\} $ of $\hR$  
 and define $D_i \df  \{\ff \=  a_i \} \ins  \sB(\fX) \oti \sB(\fY)$ for $i \= 1,\cds \n ,N$.
   By the equality $ \sB(\fX) \oti \sB(\fY) \= \L $, we see that
$ \phi_\ff  (x,P)   \= \int_{y \in \fY} \ff (x,y) P(dy)
 \= \underset{i=1}{\overset{N}{\sum}} a_i \int_{y \in \fY} \b1_{\{(x,y) \in D_i\}} P(dy)
 \= \underset{i=1}{\overset{N}{\sum}} a_i \phi_{D_i} (x,P) $, $  (x,P) \ins \fX \ti \fP(\fY)$
 is $\sB(\fX) \oti \sB\big(\fP(\fY)\big)  -$measurable.
 Then   we can use  the standard approximation to obtain such a Borel measurability
 for   general $(-\infty,\infty]-$valued, $\sB(\fX) \oti \sB(\fY)-$measuable   functions bounded from below.  \qed

\begin{lemm} \label{lem_012922_11}
 Given $t \ins [0,\infty)$,   let $\tau $ be a $[t,\infty]-$valued  $\bF^{W^t,P_0}-$stopping time
 and  let $  \fP_t $ be a subset of $\fP(\oO)$ such that
 $ \oW^t $ is a   Brownian motion under each $\oP \ins \fP_t$.
 There is a $[t,\infty]-$valued  $\sB(\oO)-$measurable random variable $\oxi$ such that
 $ \{\tau (\oW) \nne \oxi\} \ins \ccap{\oP \in \fP_t}{} \sN_{\oP}\big(\cF^{\oW^t}_\infty\big) $.
 If $\oga$ is a  $[t,\infty)-$valued  $\bF^{\oW^t} -$stopping time,
 one can find $ \oA \ins \cF^{\oW^t}_\oga$ such that
$  \{\tau (\oW) \gs \oga\} \D  \oA \ins \ccap{\oP \in \fP_t}{} \sN_{\oP}\big(\cF^{\oW^t}_\infty\big)$.
\end{lemm}

 \no {\bf Proof:}   Let $n \ins \hN$. We set $s^n_i \df t \+ i 2^{-n}$, $\fa i \ins \hN \cp \{0\}$ and define
 an $\bF^{W^t,P_0}-$stopping time
   $ \tau_n \df \sum_{i \in \hN} s^n_i \b1_{\{s^n_{i-1} \le  \tau   < s^n_i  \}} \\ \+ \infty \b1_{\{ \tau  = \infty\}} $.
    By similar arguments to those   leading  to \eqref{031221_11},
    we can construct an $\bF^{W^t} -$stopping time $\wt{\tau}_n$ such that $ \cN_n \df \{\tau_n \nne \wt{\tau}_n\} \ins \sN_{P_0}(\cF^{W^t}_\infty)$.
     \if{0}

   Let $n \ins \hN$ and set $s^n_0 \df t$. Given $ i \ins \hN$, we set $s^n_i \df t \+ i 2^{-n}$ and $A^n_i \df \big\{s^n_{i-1} \ls  \tau  \< s^n_i \big\} \ins \cF^{W^t,P_0}_{s^n_i}$.
   By Problem 2.7.3 of \cite{Kara_Shr_BMSC},
 there exists   $\wA^n_i  \ins \cF^{W^t}_{s^n_i}$   such that $ \cN^n_i \df  A^n_i   \D    \wA^n_i  \ins  \sN_{P}  (\cF^{W^t}_\infty )  $.
    Define    $\wcA^n_i  \df  \wt{A}^n_i \Big\backslash \Big(\underset{j < i}{\cup}
     \wt{A}^n_j \Big)  \ins  \cF^{W^t}_{ s^n_i }  $
     and $ \wcA_n  \df  \ccup{i \in \hN}{ }  \wcA^n_i
      \= \ccup{i \in \hN}{ }  \wt{A}^n_i  \ins  \cF^{W^t}_\infty $.
 The $\bF^{W^t,P_0}-$stopping time
   $ \tau_n \df \sum_{i \in \hN} s^n_i \b1_{A^n_i}  \+ \infty \b1_{\{ \tau  = \infty\}} $ coincides with the $\bF^{W^t} -$stopping time
   $ \wt{\tau}_n \df \sum_{i \in \hN} s^n_i \b1_{\wcA^n_i}   \+ \infty \b1_{\wcA^c_n } $  on
   $ \O_n \df \Big(\ccup{i \in \hN}{} \big( A^n_i \Cp \wcA^n_i \big) \Big) \cp \Big( \{ \tau  = \infty\} \Cp \wcA^c_n \Big)$.
   Using similar arguments to those that lead  to \eqref{031221_11}, one has $\O^c_n \ins \sN_{P_0}(\cF^{W^t}_\infty)$.

    \fi

\no {\bf 1)} Define $\oxi (\oo) \df \linf{n \in \hN} \wt{\tau}_n \big(\oW (\oo)\big) \ins [t,\infty] $, $\fa \oo \ins \oO$.
Let $\oP \ins \fP_t$.
 Applying  Lemma \ref{lem_122921_11} with $t_0 \= t$, $(\O_1, \cF_1, P_1,B^1)   \= \big(\oO ,  \sB(\oO ),  \oP, \oW \big) $, $(\O_2, \cF_2,  P_2,B^2) \= \big(\O_0,  \sB(\O_0),   P_0,  W \big) $  and $\Phi \= \oW$ yields that
 \bea \label{022022_11}
 \hb{$\wt{\tau}_n (\oW)  $ is   an $\bF^{\oW^t} -$stopping time and
 $ \oW^{-1}(\cN_n) \= \big\{\tau_n (\oW) \nne \wt{\tau}_n (\oW) \big\} \ins \sN_{\oP}(\cF^{\oW^t}_\infty) , \q \fa n \ins \hN$.}
 \eea
 So $\oxi   \df \linf{n \in \hN} \wt{\tau}_n  (\oW)$ is $ \cF^{\oW^t}_\infty -$measurable and is further $\sB(\oO)-$measurable.
 Since
 \bea \label{022022_14}
 \tau \big(\oW(\oo)\big)   \=     \lmtu{n \to \infty}  \tau_n \big(\oW(\oo)\big)     \= \linf{n \to \infty}   \wt{\tau}_n \big(\oW(\oo)\big) ,
 \q \fa \oo \ins \Big(\ccup{n \in \hN}{} \oW^{-1}(\cN_n)\Big)^c \= \ccap{n \in \hN}{} \oW^{-1}(\cN^c_n) ,
 \eea
 one has  $ \{\tau (\oW) \nne \oxi\} \sb \ccup{n \in \hN}{} \oW^{-1}(\cN_n) $, i.e.,
 $ \{\tau (\oW) \nne \oxi\}  \ins   \sN_{\oP}\big(\cF^{\oW^t}_\infty\big) $.

\no {\bf 2)} Let $\oga$ be a $[t,\infty)-$valued  $\bF^{\oW^t} -$stopping time.
 We set $ \oA   \df     \Big\{ \, \linf{n \in \hN} \wt{\tau}_n (\oW) \gs \oga\Big\}$ and
 let $\oP \ins \fP_t$. As \eqref{022022_11} shows that $ \big\{   \wt{\tau}_k (\oW) \+ a \gs \oga   \big\} \ins \cF^{\oW^t}_{(\wt{\tau}_k (\oW)  + a ) \land \oga} \sb \cF^{\oW^t}_\oga $ for any $k  \ins \hN$ and $a \ins [0,\infty)$, we can deduce that
  $ \oA   \=   \Big\{ \, \Sup{n \in \hN} \Big( \Inf{k \ge n} \wt{\tau}_k (\oW)\Big)  \gs \oga\Big\}
\=  \ccap{m \in \hN}{} \ccup{n \in \hN}{} \Big\{ \, \Inf{k \ge n} \wt{\tau}_k (\oW) \gs \oga  \-  1/m\Big\}
   \=   \ccap{m \in \hN}{} \ccup{n \in \hN}{} \ccap{k \ge n}{}\big\{   \wt{\tau}_k (\oW) \+ 1/m \gs \oga   \big\} \ins \cF^{\oW^t}_\oga $.
  Using \eqref{022022_14} also implies that
 $ \big\{\tau (\oW) \gs \oga \big\} \D \oA  \sb \big\{ \tau (\oW) \nne \linf{n \to \infty}  \wt{\tau}_n (\oW)   \big\} \sb \ccup{n \in \hN}{} \oW^{-1}(\cN_n) $, so $\big\{\tau (\oW) \gs \oga \big\} \D \oA  \ins \sN_{\oP}(\cF^{\oW^t}_\infty) $.   \qed

\begin{lemm} \label{lem_013022_11}
Let $t \ins [0,\infty)$ and $\oP \ins \fP\big(\oO\big)$.
For any $(s,\oo) \ins [t,\infty) \ti \oO$,  set  $\obW^t_{s,\oo}   \df \big\{\oo' \ins \oO \n : \oW^t_a (\oo') \= \oW^t_a (\oo), \;  \fa a \ins [t,s] \big\}$. Then for any $r \ins [t,\infty)$, $ \cF^{W^t,P_0}_r  \sb   \sS_r     \df     \Big\{  A \sb   \O_0   \n: \exists \, \ocN_r \ins \sN_\oP \big(\cF^{\oW^t}_\infty\big) $  such that for any $(s,\oo) \ins [t,r] \ti \oO$,  there exists   $A^{s,\oo}  \ins \cF^{W^s}_r $   satisfying
     $ \b1_{\{ \oW  (\oo')  \in  A \}}
     \= \b1_{\{ \oW   (\oo')  \in  A^{s,\oo} \}},
     \fa  \oo'  \ins  \obW^t_{s,\oo} \Cp \ocN^c_r   \Big\} $.
\end{lemm}

  \no {\bf Proof:} Fix $r \ins [t,\infty)$.
  Clearly,  $\es \ins \sS_r $ with $A^{s,\oo} \= \es$ for any  $(s,\oo) \ins [t,r] \ti \oO$.
  Given $A \ins \sS_r  $, by taking the complement of each $A^{s,\oo}$, we see that  $A^c \ins \sS_r $.
  \if{0}

  As $ \b1_{\{ \oW (\oo')  \in  \es \}} \= 0 \= \b1_{\{ \oW (\oo')  \in  \es\}}  $
   and $ \b1_{\{ \oW (\oo')  \in  \O_0 \}} \= 1 \= \b1_{\{ \oW (\oo')  \in  \O_0\}} $
 for any  $   \oo'   \ins   \oO $,
 it is clear that   both $\es$ and $\O_0  $ belong to $\sS_r $.

 When $A \ins \sS_r  $, one can find an $\ocN_r\ins \sN_\oP \big(\cF^{\oW^t}_\infty\big) $  such that
 for each  $(s,\oo) \ins [t,r] \ti \oO$, there exists   $A^{s,\oo}  \ins \cF^{W^s}_r$   satisfying
    $ \b1_{\{ \oW (\oo')  \in  A \}}
     \= \b1_{\{ \oW (\oo')  \in  A^{s,\oo} \}}$,
     $\fa  \oo'  \ins  \obW^t_{s,\oo} \Cp \ocN^c_r $. Then for any $(s,\oo) \ins [t,r] \ti \oO$,
   $ (A^{s,\oo} )^c \ins \cF^{W^s}_r $  satisfies
   $ \b1_{\{ \oW (\oo')  \in  A^c \}} \= 1\- \b1_{\{ \oW (\oo')  \in  A \}}
     \= 1 \-\b1_{\{ \oW (\oo')  \in A^{s,\oo} \}} \= \b1_{\{ \oW (\oo')  \in  (A^{s,\oo} )^c \}} $, $ \fa \oo'  \ins  \obW^t_{s,\oo} \Cp \ocN^c_r $,
   so $A^c \ins \sS_r $.

 \fi
   Let $\{A_n\}_{n \in \hN} \sb \sS_r $. For any $n \ins \hN$,  we can find
        $\ocN^n_r \ins \sN_\oP \big(\cF^{\oW^t}_\infty\big) $ such that
        for any $(s,\oo) \ins [t,r] \ti \oO$,
        there exists   $A^{s,\oo}_n \ins \cF^{W^s}_r$   satisfying
    $  \b1_{ \{\oW  (\oo')\in A_n \} } \= \b1_{\{ \oW (\oo')  \in A^{s,\oo}_n \} }    $ for any $\oo'  \ins  \obW^t_{s,\oo} \cap \big(\ocN^n_r\big)^c $.
    Then for any $(s,\oo) \ins [t,r] \ti \oO$,
      $ \underset{n \in \hN}{\cap} A^{s,\oo}_n \ins \cF^{W^s}_r $ satisfies
    $  \b1_{\big\{\oW  (\oo')\in \underset{n \in \hN}{\cap}  A_n \big\} }
    \= \underset{n \in \hN}{\prod}  \b1_{ \{\oW  (\oo')\in A_n \} }
    \= \underset{n \in \hN}{\prod} \b1_{\{ \oW (\oo')  \in A^{s,\oo}_n \} }
    \= \b1_{\big\{ \oW  (\oo')  \in \underset{n \in \hN}{\cap} A^{s,\oo}_n \big\} }  $, $ \fa \oo'  \ins  \obW^t_{s,\oo} \Cp \Big( \ccup{n \in \hN}{} \ocN^n_r \Big)^c  $,
  which shows   $ \underset{n \in \hN}{\cap}  A_n \ins \sS_r  $. 
  So $\sS_r $ is a sigma$-$field of $\O_0$.

  Let $(a,\cE) \ins [t,r] \ti \sB(\hR^d)$.
  We verify that  $ (W^t_a)^{-1} (\cE) \ins \sS_r $ by three cases: Let $(s,\oo) \ins [t,r] \ti \oO $.

 \no (i) If $s \gs a$ and $ \oW^t_a   (\oo) \ins \cE  $, then
 $ \b1_{\{ \oW (\oo')  \in (W^t_a)^{-1} (\cE) \}}
  \=  \b1_{\{ \oW^t_a (\oo')  \in  \cE  \}}
  \=  \b1_{\{ \oW^t_a (\oo)  \in  \cE  \}} \= 1
     \= \b1_{\{ \oW  (\oo')  \in  \O_0\}} $,  $  \fa \oo'  \ins  \obW^t_{s,\oo} $.

  \no (ii) If $ s \gs  a   $ but $ \oW^t_a   (\oo) \n  \notin \n \cE  $, then
 $  \b1_{\{ \oW (\oo')  \in (W^t_a)^{-1} (\cE) \}} 
  \=  \b1_{\{ \oW^t_a (\oo)  \in  \cE  \}} \= 0 \= \b1_{\{ \oW  (\oo')  \in  \es\}} $, $ \fa \oo'  \ins  \obW^t_{s,\oo} $.

  \no (iii) If $ s \<  a   $,   set $ \cE_{s,\oo} \df   \big\{\fx \- \oW^t_s   (\oo) \n : \fa \fx \ins \cE \} \ins \sB(\hR^d) $, then $ (W^s_a)^{-1} (\cE_{s,\oo}) \ins 
  \cF^{W^s}_r $ satisfies that
  $ \b1_{\{ \oW (\oo')  \in (W^t_a)^{-1} (\cE) \}} 
  \=  \b1_{\{ \oW^t_a (\oo')  - \oW^t_s (\oo')  \in  \cE_{s,\oo}  \}} 
  \= \b1_{\{ \oW  (\oo')   \in (W^s_a)^{-1}  (\cE_{s,\oo}) \}} $, $ \fa \oo'  \ins  \obW^t_{s,\oo} $.
  So $ (W^t_a)^{-1} (\cE) \ins \sS_r $, 
  and it follows    $ 
   \cF^{W^t}_r \sb \sS_r $.

  Let $\cN \ins \sN_{P_0}(\cF^{W^t}_\infty)$.
  Applying Lemma \ref{lem_122921_11}   with $t_0 \= t$, $(\O_1, \cF_1, P_1,B^1)   \= \big(\oO ,  \sB(\oO ),  \oP , \oW\big) $, $(\O_2, \cF_2, P_2,B^2) \= \big(\O_0,  \sB(\O_0),  P_0, W\big) $ and $\Phi \= \oW$
  yields that $ \oW^{-1} (\cN) \ins  \sN_\oP \big( \cF^{\oW^t}_\infty \big) $.
  For any $(s,\oo) \ins [t,r] \ti \oO$ and   any $ \oo'  \ins  \obW^t_{s,\oo} \Cp  \oW^{-1} (\cN^c) $, one has
 $ \b1_{\{ \oW (\oo')  \in  \cN \}}
 \= \b1_{\{ \oo' \in   \oW^{-1} (\cN) \}} \= 0
     \= \b1_{\{ \oW (\oo')  \in  \es \}} $.
     This shows  $ \cN \ins \sS_r $  
    and thus  $ 
     \cF^{W^t,P_0}_r \sb \sS_r $. \qed

\if{0}

   \begin{lemm}  \label{lemm_FW_rv}
Given $t \ins [0,\infty)$,
 let $\oW^t $ be a $d-$dimensional 
 Brownian motion on $\big(\oO,\sB(\oO)\big)$ under  some $\oP \ins \fP \big(\oO\big)$.
 Let $s \ins [t,\infty]$ and let $\oxi$ be a $\hU-$valued, $\cF^{\oW^t}_s-$measurable random variable on $\oO$.
There exists a $\hU-$valued, $\cF^{W^t}_s-$measurable random variable  $\xi_0$ on $\O_0$ such that
$\xi_0 \big(\oW (\oo)\big)  \= \oxi (\oo)$, $\fa \oo \ins \oO$.
Moreover, it holds for any Borel-measurable function $\psi \n : \hU \mto \hR$ that
 $E_\oP \big[ \psi (\oxi)  \big] \= E_{P_0} \big[ \psi (\xi_0) \big]  $.

\end{lemm}

  \no {\bf Proof:} Let $s \ins [t,\infty]$.

  \ss \no {\bf 1)} Since it holds for any $r \ins [t,s] \Cp \hR $, $\cE \ins \sB(\hR^d)$ and $\oo \ins \oO$ that
  $  \b1_{ (W^t_r)^{-1}(\cE) } \big(\oW (\oo)\big) 
  \= \b1_{\{W^t_r (\oW (\oo)) \in   \cE \}}
  \= \b1_{\{ \oW^t_r (\oo)  \in   \cE \}}
  \=  \b1_{   (\oW^t_r)^{-1}(\cE) } (\oo) $, all generating sets of $\cF^{\oW^t}_s$ is included in
   \beas
    \ol{\L}_s \df \big\{\oA \sb \oO \n :  \b1_{\oA} \= \b1_A  (\oW)
   \hb{ for some } A \ins \cF^{W^t}_s    \big\} .
   \eeas
   Clearly, $\oO \ins \ol{\L}_s$ with $A \= \O_0$. 
   If $\oA \ins \ol{\L}_s$, one can find $\wA  \ins \cF^{W^t}_s$ such that
   $ \b1_{\oA} (\oo) \= \b1_{\wA}  \big(\oW (\oo)\big)$ for any $  \oo \ins \oO$.
   Then
   \beas
   \b1_{\oA^c} (\oo) \= 1 \- \b1_{\oA} (\oo) \= 1 \- \b1_\wA  \big(\oW (\oo)\big) \= \b1_{\wA^c}  \big(\oW (\oo)\big) , \q \fa  \oo \ins \oO .
   \eeas
   So $ \oA^c \ins \ol{\L}_s$ with $A \= \wA^c $.

   Let $\{\oA_n\}_{n \in \hN} \sb \ol{\L}_s$. For any $n \ins \hN$, there exists $A_n \ins \cF^{W^t}_s$ such that
   $ \b1_{\oA_n} (\oo) \= \b1_{A_n}  \big(\oW (\oo)\big)$ for any $  \oo \ins \oO$. We can deduce that
   \beas
   \b1_{\ccap{n \in \hN}{}\oA_n} (\oo) \= \underset{n \in \hN}{\prod} \b1_{ \oA_n} (\oo)
   \= \underset{n \in \hN}{\prod} \b1_{A_n}  \big(\oW (\oo)\big)
   \= \b1_{\ccap{n \in \hN}{} A_n}  \big(\oW (\oo)\big) , \q \fa  \oo \ins \oO ,
   \eeas
   which shows that $\ccap{n \in \hN}{}\oA_n \ins \ol{\L}_s$ with $A \= \ccap{n \in \hN}{} A_n $.
   Hence, $\ol{\L}_s$  is a sigma-field of $\oO$ including $ \cF^{\oW^t}_s$.  

  \ss \no {\bf 2)}  Let $\{u_i\}_{i \in \hN}$ be a countable dense subset of $\big(\hU, \Rho{\hU}\big)$.
   For any $i,n \ins \hN$, we set $O^n_i \df \big\{u \ins \hU: \Rho{\hU} (u,u_i) \< 2^{-n} \big\}$
 as  the open ball centered at $u_i$ with radius $2^{-n}$.

 Let $ \oxi$ be a general $\hU-$valued, $\cF^{\oW^t}_s$ measurable random variable on $\oO$
 and let $n \ins \hN$.
 For any $i \ins \hN$, we define $\oA^n_i \df \{ \oxi   \ins \cO^n_i \} \ins \cF^{\oW^t}_s$ with
 $\cO^n_1 \df O^n_1 $ and $\cO^n_i \df O^n_i \big\backslash \Big( \underset{j < i }{\cup} O^n_j \Big) $ for $i \gs 2$.
 By Part 1), there exists  $A^n_i  \ins \cF^{W^t}_s $ such that $ \b1_{\oA^n_i}(\oo) \= \b1_{A^n_i}  \big(\oW (\oo) \big) $, $\fa \oo \ins \oO$.

 Setting $\cA^n_1 \df A^n_1 $ and $\cA^n_i \df A^n_i \big\backslash \Big( \underset{j < i }{\cup} A^n_j \Big) \ins \cF^{W^t}_s $ for $i \gs 2$.
 We define an $\fE-$valued, $\cF^{W^t}_s-$measurable random variable on $\O_0$ by
 \beas
  \beta_n (\o_0) \df \b1_{\big\{\o_0 \in   \ccap{i \in \hN}{}  (A^n_i)^c \big\}}  \sI   (u_1)
  \+ \sum_{i \in \hN} \b1_{\{  \o_0  \in  \cA^n_i \}}  \sI   (u_i)   ,  \q \fa \o_0  \ins   \O_0 .
 \eeas
 Then  $ \ul{\beta} (\o_0) \df \linf{n \to \infty} \beta_n (\o_0) $
 and  $ \ol{\beta} (\o_0) \df \lsup{n \to \infty} \beta_n (\o_0) $, $ \o_0 \ins \O_0 $,
 are two $[0,1]-$valued, $\cF^{W^t}_s-$measurable random variables on $\O_0$. It follows that
\beas
\beta (\o_0) \df \ul{\beta} (\o_0) \b1_{\{ \ul{\beta} (\o_0) =  \ol{\beta} (\o_0) \} \cap \{\ul{\beta} (\o_0) \in \fE\}}
\+  \sI (u_1)  \b1_{\{ \ul{\beta} (\o_0) \ne  \ol{\beta} (\o_0) \} \cup \{\ul{\beta} (\o_0) \notin \fE\}}, \q \o_0 \ins \O_0
\eeas
 is an $\fE-$valued, $\cF^{W^t}_s-$measurable random variable on $\O_0$. Moreover,
 $\xi_0  (\o_0) \df \sI^{-1} (\beta (\o_0))$, $\o_0 \ins \O_0$ is
 a $\hU-$valued, $\cF^{W^t}_s-$measurable random variable on $\O_0$.

 Let  $ \oo  \ins  \oO$ and $ n \ins \hN $.
 Since $ \oW^{-1}(\cA^n_i) \= \oW^{-1}(A^n_i) \Big\backslash \Big( \underset{j < i }{\cup} \oW^{-1}(A^n_j) \Big)
\= \oA^n_i \big\backslash \Big( \underset{j < i }{\cup} \oA^n_j \Big) \= \oA^n_i $ for any  $i \ins \hN$,
 we can deduce that  $   \oW^{-1}\Big(\ccup{i \in \hN}{}\cA^n_i\Big)
 \= \ccup{i \in \hN}{} \oW^{-1} (\cA^n_i ) \= \ccup{i \in \hN}{} \oA^n_i \= \oO $ and thus
 \beas
 \sI^{-1} \big(  \beta_n (\oW (\oo)) \big) \=   \sum_{i \in \hN} \b1_{\{  \oW (\oo)  \in  \cA^n_i \}}  u_i
 \= \sum_{i \in \hN} \b1_{\{   \oo   \in \oW^{-1} (\cA^n_i) \}}  u_i \= \sum_{i \in \hN} \b1_{\{   \oo   \in  \oA^n_i \}}  u_i .
\eeas
 It follows that  $ \Rho{\hU} \Big(   \oxi(\oo) ,   \sI^{-1}  \big(  \beta_n (\oW (\oo)) \big) \Big)
 \< 2^{-n} $.
 Letting $n \to \infty$ yields that $\lmt{n \to \infty}   \sI^{-1}  \big(  \beta_n (\oW (\oo)) \big) \=  \oxi(\oo)  $ in $\hU$
 and then $\sI(\oxi(\oo)) \= \lmt{n \to \infty} \beta_n (\oW (\oo)) $,
 the latter of which implies that $\ul{\beta} (\oW (\oo)) \= \ol{\beta} (\oW (\oo)) \= \sI(\oxi(\oo)) \ins \fE$.
  So $\beta (\oW (\oo)) \= \ul{\beta} (\oW (\oo)) \= \sI(\oxi(\oo))$ or
  $\oxi(\oo) \= \sI^{-1} \big( \beta (\oW (\oo)) \big) \= \xi_0  \big(\oW (\oo)\big) $.

    \if{0}

  \ss \no {\bf Another Method:} When $\hU \= \hR$, we can use   the monotone class theorem.

  Fix $t,s \ins [0,\infty]$ and let
  $\fH_s $ denote the collection of all bounded real-valued random variables $\oxi$ on $\oO$ such that
  $\xi(\oW) \= \oxi $ for some  $\cF^{W^t}_s-$measurable random variable  $\xi$ on $\O_0$.
  Clearly, $\fH_s$ is a vector space that contains all constants 
  and that is closed under multiplication between random variables.

  We show that $\fH_s$ is also closed under increasing monotone convergence:
  Let $\big\{\oxi_n \big\}_{n \in \hN}$ be a uniformly bounded increasing sequence of positive random variables in $\fH_s$.
  For any $n \ins \hN$, let $\xi_n $ be the $\cF^{W^t}_s-$measurable random variable   on $\O_0$
  satisfying $\xi_n(\oW) \= \oxi_n $. For any $\o_0 \ins \O_0$, the increasing monotonicity of
  sequence $ \big\{\oxi_n \big\}_{n \in \hN} $ shows that
  $ \xi_n (\o_0) \= \oxi_n \big(\{\o_0\} \ti [0,1] \ti \hJ  \ti \OmX  \ti [0,\infty]\big)$ is increasing in $n \ins \hN$,
   and the uniform  boundedness of $ \big\{\oxi_n \big\}_{n \in \hN} $ implies that
   $\xi_* (\o_0) \df \lmtu{n \to \infty} \xi_n (\o_0)
   \= \lmtu{n \to \infty} \oxi_n \big(\{\o_0\} \ti [0,1] \ti \hJ  \ti \OmX  \ti [0,\infty]\big) \ls \big\| \oxi \big\|
   \df \lmtu{n \to \infty} \underset{\oo \in \oO}{\sup} \oxi_n (\oo) \< \infty $.
     As $\xi_* \= \lmtu{n \to \infty} \xi_n$ is an $\cF^{W^t}_s-$measurable random variable   on $\O_0$
     satisfying   $ \xi_* (\oW ) \= \lmtu{n \to \infty} \xi_n (\oW ) \= \lmtu{n \to \infty} \oxi_n $,
     we see that $ \lmtu{n \to \infty} \oxi_n \ins \fH_s $.

  For any $r \ins [0,s] \Cp \hR $ and $\cE \ins \sB(\hR^d)$, since
  $  \b1_{ W^{-1}_r(\cE) } \big(\oW (\oo)\big) 
  \= \b1_{\{W_r (\oW (\oo)) \in   \cE \}}
  \= \b1_{\{ \oW_r (\oo)  \in   \cE \}}
  \=  \b1_{ \oW_r^{-1}(\cE) } (\oo) $, $\fa \oo \ins \oO$,
  one has  $  \b1_{ \oW_r^{-1}(\cE) } \ins \fH_s$. As $\fH_s$ is closed under multiplication,
  \beas
  \sC_s \df \Big\{  \b1_{\big(\underset{i=1}{\overset{n}{\cap}}  \oW_{s_i}^{-1}(\cE_i)\big)} \n :
   s_i \ins [0,s] \Cp \hR \hb{ and } \cE_i \ins \sB(\hR^d) \hb{ for } i \= 1,\cds, n \Big\} \sb \fH_s .
  \eeas
  Then the monotone class theorem yields that
  $\fH_s$ includes all bounded random variables measurable with respect to  $\si(\sC_s) \= \si (\oW_r, r \ins [0,s] \Cp \hR)
  \= \cF^{\oW^t}_s $.

  Now, let $\oxi$ be a general real-valued, $\cF^{\oW^t}_s$ measurable random variable on $\oO$.
  For any $n \in \hN$, as $\oxi^\pm_n \df \oxi^\pm \ld n \ins \fH_s $, we can find
  some  $\cF^{W^t}_s-$measurable random variable  $\xi^\pm_n$ on $\O_0$ such that
  $\xi^\pm_n (\oW) \= \oxi^\pm_n $.
  For any $\o_0 \ins \O_0$, the increasing monotonicity of
  sequence $ \big\{\oxi^\pm_n \big\}_{n \in \hN} $ shows that
  $ \xi^\pm_n (\o_0) \= \oxi^\pm_n \big(\{\o_0\} \ti [0,1] \ti \hJ  \ti \OmX  \ti [0,\infty]\big)$
  is increasing in $n \ins \hN$. So
   \beas
   \xi^\pm (\o_0) \df \lmtu{n \to \infty} \xi^\pm_n (\o_0)
   \= \lmtu{n \to \infty} \oxi^\pm_n \big(\{\o_0\} \ti [0,1] \ti \hJ  \ti \OmX  \ti [0,\infty]\big)
   \= \oxi^\pm  \big(\{\o_0\} \ti [0,1] \ti \hJ  \ti \OmX  \ti [0,\infty]\big)   \< \infty  .
   \eeas
   It follows that $\xi \df \xi^+ \- \xi^- \= \lmt{n \to \infty} \xi^+_n \- \lmt{n \to \infty} \xi^-_n $ is an $\cF^{W^t}_s-$measurable random variable  on $\O_0$ satisfying
   $ \xi (\oW ) \= \lmt{n \to \infty} \xi^+_n (\oW ) \- \lmt{n \to \infty} \xi^-_n  (\oW )
   \= \lmt{n \to \infty} \oxi^+_n   \- \lmt{n \to \infty} \oxi^-_n
   \= \oxi^+ \- \oxi^- \= \oxi $.

    \fi

  Moreover, let  $\psi \n : \hU \mto \hR$ be a $\sB(\hU)/\sB(\hR)-$measurable function.
  Since  an application of \ref{lem_122921_11} with $t_0 \= t$, $(\O_1, \cF_1, P_1,B^1)   \= \big(\oO ,  \sB(\oO ),  \oP, \oW \big) $, $(\O_2, \cF_2, P_2,B^2) \= \big(\O_0,  \sB(\O_0),  P_0,  W \big) $  and $\Phi \= \oW $ yields that
  $ E_\oP \big[ \psi^\pm( \oxi ) \ld n \big] \= E_\oP \big[ \psi^\pm \big( \xi_0(\oW)  \big) \ld n \big]
  \= E_{P_0} \big[ \psi^\pm  ( \xi_0   ) \ld n \big]$, $\fa n \ins \hN $,
    the monotone convergence theorem renders  that
   $ E_\oP \big[ \psi^\pm( \oxi )   \big] \= \lmtu{n \to \infty} E_\oP \big[ \psi^\pm( \oxi ) \ld n \big]
  \= \lmtu{n \to \infty} E_{P_0} \big[ \psi^\pm  ( \xi_0   ) \ld n \big] \= E_{P_0} \big[ \psi^\pm  ( \xi_0   )   \big]$.

   If $ E_\oP \big[ \psi^-( \oxi )   \big] \= E_{P_0} \big[ \psi^-  ( \xi_0   )   \big] \= \infty $,
   then $E_\oP \big[ \psi ( \oxi )   \big] \= E_{P_0} \big[ \psi   ( \xi_0   )   \big] \= - \infty $.
   If $ E_\oP \big[ \psi^-( \oxi )   \big] \= E_{P_0} \big[ \psi^-  ( \xi_0   )   \big] \< \infty $,
   then $E_\oP \big[ \psi ( \oxi )   \big]
   \= E_\oP \big[ \psi^+ ( \oxi )   \big] \- E_\oP \big[ \psi^-( \oxi )   \big]
   \= E_{P_0} \big[ \psi^+  ( \xi_0   )   \big] \- E_{P_0} \big[ \psi^-  ( \xi_0   )   \big]
   \=  E_{P_0} \big[ \psi   ( \xi_0   )   \big]$.  \qed

   \fi

 \ss \no {\bf Acknowledgments } We are grateful to  Xiaolu Tan for helpful comments.

\bibliographystyle{siam}
\bibliography{OSMC_bib}

  \end{document}